\documentclass[11pt,a4paper]{amsart}
\usepackage{amssymb,amsxtra}
\usepackage[english,french]{babel}
\newtheorem{thm}{Theorem}[section]
\newtheorem{lem}[thm]{Lemma}
\newtheorem{prop}[thm]{Proposition}

\theoremstyle{definition}
\newtheorem{defn}[thm]{Definition}

\theoremstyle{remark}

\theoremstyle{plain}

\theoremstyle{remark}

\newtheorem*{example}{Example}
\numberwithin{equation}{section}

\begin{document}
\title{\textbf{ EXPONENTIAL SUMS AND RANK OF DOUBLE  PERSYMMETRIC MATRICES OVER $\mathbb{F}_{2} $}}
\author{Jorgen~Cherly}
\address{D\'epartement de Math\'ematiques, Universit\'e de    
    Brest, 29238 Brest cedex~3, France}
\email{Jorgen.Cherly@univ-brest.fr}
\thanks{}
 \maketitle 
 
\begin{abstract}
 Soit $\mathbb{K}^{2} $  le  $\mathbb{K}$ - espace vectoriel de dimension 2  o\`{u} $\mathbb{K}$  d\'{e}note 
    le corps des s\'eries  de Laurent formelles  $ \mathbb{F}_{2}((T^{-1})). $ Nous calculons en particulier des sommes exponentielles 
    (dans $\mathbb{K}^{2}$) de la forme \\
     $  \sum_{\deg Y \leq  k-1}\sum_{\deg Z \leq s-1}E(tYZ)\sum_{\deg U \leq s+m-1}E(\eta YU) $
   o\`{u} $ (t,\eta ) $ est dans la boule  unit\'{e} de $\mathbb{K}^{2}.$\\
     Nous d\'{e}montrons qu'elles d\'{e}pendent uniquement
   du rang de matrices doubles persym\'{e}triques  avec des entr\'{e}es dans $\mathbb{F}_{2},$
    c'est-\`{a}-dire des matrices de la forme  $\left[A\over B\right]  $ o\`{u}  A  est une matrice  $ s \times k $
    persym\'{e}trique et B une matrice  $ (s+m) \times k $  persym\'{e}trique  (une matrice \;$ [\alpha _{i,j}]  $ est 
      persym\'{e}trique  si  $ \alpha _{i,j} = \alpha _{r,s} $ \; pour  \; i+j = r+s). En outre, nous \'{e}tablissons plusieurs 
      formules concernant des propri\'{e}t\'{e}s de rang de partitions de matrices doubles persym\'{e}triques, ce qui nous 
      conduit \`{a} une formule r\'{e}currente du nombre   $  \Gamma _{i}^{\Big[\substack{s \\ s+m }\Big] \times k} $ 
   des matrices de rang i de la forme  $\left[A\over B\right]. $ Nous d\'{e}duisons de cette formule r\'{e}currente 
   que si  $  0\leq i\leq\inf(s-1, k-1),$ le nombre    $  \Gamma _{i}^{\Big[\substack{s \\ s+m }\Big] \times k} $ d\'{e}pend uniquement de i.
   D'autre part, si   $ i\geq s+1, k\geq i,\; \Gamma _{i}^{\Big[\substack{s \\ s+m }\Big] \times k}  $ peut \^{e}tre  calcul\'{e}
   \`{a} partir du nombre  $  \Gamma _{s' +1}^{\Big[\substack{s' \\ s'+m' }\Big] \times k'} $de matrices  de rang (s'+1) de la forme  
       $\left[A'\over B'\right] $ o\`{u}  A' est une matrice   $ s' \times k' $ persym\'{e}trique  et  B' une matrice  $ (s'+m') \times k' $ persym\'{e}trique,
        o\`{u}  s', m' et k' d\'{e}pendent  de  i, s, m et k. La preuve de ce r\'{e}sultat est bas\'{e}e sur une formule (donn\'{e}e dans  [4])
        du nombre de matrices de rang i de la forme     $\left[A\over b_{-}\right] $   o\`{u}  A est  persym\'{e}trique  et  $ b_{-} $
 une matrice ligne avec entr\'{e}es dans $ \mathbb{F}_{2}.$ Nous montrons \'{e}galement que le nombre R de repr\'{e}sentations dans 
  $ \mathbb{F}_{2}[T] $ 
 des \'{e}quations polynomiales
   \[\left\{\begin{array}{cc}
YZ + Y_{1}Z_{1} + \ldots + Y_{q-1}Z_{q-1}= 0 \\
YU + Y_{1}U_{1} + \ldots + Y_{q-1}U_{q-1}= 0
    \end{array}\right.\]\\ 
  associ\'{e}es aux sommes exponentielles \\
   $\sum_{degY\leq  k-1}\sum_{degZ\leq s-1}E(tYZ)\sum_{\deg U\leq s+m-1}E(\eta YU)  $
  est donn\'{e} par une int\'{e}grale sur la boule unit\'{e} de $\mathbb{K}^{2}$ et est une combinaison lin\'{e}aire de
  $ \Gamma _{i}^{\Big[\substack{s \\ s+m }\Big] \times k}$ pour $ i\geq 0. $ Nous pouvons alors calculer explicitement le nombre R.
 \end{abstract}
\allowdisplaybreaks

 \selectlanguage{english}
 \begin{abstract}
Let $ \mathbb{K}^{2}  $ be the 2-dimensional vectorspace over $ \mathbb{K}$  where $ \mathbb{K}$ denotes the field of
 Laurent Series  $ \mathbb{F}_{2}((T^{-1})). $ We compute in particular   exponential sums, (in $\mathbb{K}^{2} $)  of the form \\
  $  \sum_{\deg Y \leq  k-1}\sum_{\deg Z \leq s-1}E(tYZ)\sum_{\deg U \leq s+m-1}E(\eta YU) $ where 
$ (t,\eta ) $ is in the unit interval of  $\mathbb{K}^{2}.$  We  show that they only depend  on the rank  of some associated  double  persymmetric matrices  
with entries in  $\mathbb{F}_{2}, $ that is matrices of the form $\left[A\over B\right] $ 
 where A is  a  $ s \times k $ persymmetric matrix and B a  $ (s+m) \times k $ persymmetric matrix. 
  (A matrix $\; [\alpha _{i,j}]  $ is  persymmetric if $ \alpha _{i,j} = \alpha _{r,s}  \; for  \; i+j = r+s ).$
   Besides, we establish several formulas concerning  rank properties of  partitions of  double  persymmetric matrices, which leads to a 
    recurrent formula for the number $  \Gamma _{i}^{\Big[\substack{s \\ s+m }\Big] \times k} $ of rank i  matrices of the form $ \left[{A\over B}\right].$  
    We deduce from the recurrent formula that  if  $  0\leq i\leq\inf(s-1, k-1) $ then   $  \Gamma _{i}^{\Big[\substack{s \\ s+m }\Big] \times k} $ 
     depends only on  i. On the other hand,  if $ i\geq s+1, k\geq i,\; \Gamma _{i}^{\Big[\substack{s \\ s+m }\Big] \times k}  $
      can be computed from the number  $  \Gamma _{s' +1}^{\Big[\substack{s' \\ s'+m' }\Big] \times k'} $ of rank (s'+1) matrices of the form 
       $\left[A'\over B'\right] $ 
 where A' is  a  $ s' \times k' $ persymmetric matrix and B' a  $ (s'+m') \times k' $ persymmetric matrix, where s', m' and k' depend on i, s, m and k.
 The proof of this result is based on a formula (given in [4]) of the number of rank i matrices of the form  $\left[A\over b_{-}\right] $ 
 where A is persymmetric and $ b_{-} $
 a one-row matrix with entries in  $ \mathbb{F}_{2}.$
   We also  prove that the number R of representations in  $ \mathbb{F}_{2}[T] $ 
 of the polynomial equations
   \[\left\{\begin{array}{cc}
YZ + Y_{1}Z_{1} + \ldots + Y_{q-1}Z_{q-1}= 0 \\
YU + Y_{1}U_{1} + \ldots + Y_{q-1}U_{q-1}= 0
    \end{array}\right.\]\\ 
  associated to the exponential sums \\
   $\sum_{degY\leq  k-1}\sum_{degZ\leq s-1}E(tYZ)\sum_{\deg U\leq s+m-1}E(\eta YU)  $
   is given by an integral over the unit interval of $\mathbb{K}^{2}$, and is a linear combination of the
    $ \Gamma _{i}^{\Big[\substack{s \\ s+m }\Big] \times k}\;for\;  i\geq 0 . $ We can  then compute 
 explicitly the number R.
   \end{abstract}

\newpage
\tableofcontents
\newpage

\section{\textbf{NOTATION}}
\label{sec 1}
\subsection{\textbf{ANALYSIS ON  $\mathbb{K} $} }
\label{subsec 1.1}

We denote by $ \mathbb{F}_{2}\big(\big({\frac{1}{T}}\big) \big)
 = \mathbb{K} $ the completion
 of the field $\mathbb{F}_{2}(T), $ the field of  rational fonctions over the
 finite field\; $\mathbb{F}_{2}$, for the  infinity  valuation \;
 $ \mathfrak{v}=\mathfrak{v}_{\infty }$ \;defined by \;
 $ \mathfrak{v}\big(\frac{A}{B}\big) = degB -degA $ \;
 for each pair (A,B) of non-zero polynomials.
 Then every element non-zero t in  $\mathbb{F}_{2}\big(\big({\frac{1}{T}}\big) \big) $
 can be expanded in a unique way in a convergent Laurent series
                              $  t = \sum_{j= -\infty }^{-\mathfrak{v}(t)}t_{j}T^j
                                 \; where\; t_{j}\in \mathbb{F}_{2}. $\\
  We associate to the infinity valuation\; $\mathfrak{v}= \mathfrak{v}_{\infty }$
   the absolute value \; $\vert \cdot \vert_{\infty} $\; defined by \;
  \begin{equation*}
  \vert t \vert_{\infty} =  \vert t \vert = 2^{-\mathfrak{v}(t)}. \\
\end{equation*}
    We denote  E the  Character of the additive locally compact group
$  \mathbb{F}_{2}\big(\big({\frac{1}{T}}\big) \big) $ defined by \\
\begin{equation*}
 E\big( \sum_{j= -\infty }^{-\mathfrak{v}(t)}t_{j}T^j\big)= \begin{cases}
 1 & \text{if      }   t_{-1}= 0, \\
  -1 & \text{if      }   t_{-1}= 1.
    \end{cases}
\end{equation*}
  We denote $\mathbb{P}$ the valuation ideal in $ \mathbb{K},$ also denoted the unit interval of  $\mathbb{K},$ i.e.
  the open ball of radius 1 about 0 or, alternatively, the set of all Laurent series 
   $$ \sum_{i\geq 1}\alpha _{i}T^{-i}\quad (\alpha _{i}\in  \mathbb{F}_{2} ) $$ and, for every rational
    integer j,  we denote by $\mathbb{P}_{j} $
     the  ideal $\left\{t \in \mathbb{K}|\; \mathfrak{v}(t) > j \right\}. $
     The sets\; $ \mathbb{P}_{j}$\; are compact subgroups  of the additive
     locally compact group  $ \mathbb{K}. $\\
      All $ t \in \mathbb{F}_{2}\Big(\Big(\frac{1}{T}\Big)\Big) $ may be written in a unique way as
$ t = [t] + \left\{t\right\}, $   $  [t] \in \mathbb{F}_{2}[T] ,
 \; \left\{t\right\}\in \mathbb{P}  ( =\mathbb{P}_{0}). $\\
 We denote by dt the Haar measure on  $ \mathbb{K} $\; chosen so that \\
  $$ \int_{\mathbb{P}}dt = 1. $$\\
   \begin{defn}
\label{defn 1.1}We introduce the following definitions in  $ \mathbb{K}: $
\begin{itemize}
\item Let s, m and k denote rational integers such that $ s\geq 2,\; m\geq 0 \; and \; k\geq 1. $\\

\item  A matrix\;$ D = [\alpha _{i,j}]  $ is said to be persymmetric if $ \alpha _{i,j} = \alpha _{r,s} $ \; whenever \; i+j = r+s. \vspace{0.5 cm}
\item Set $ t = \sum_{i\geq 1}\alpha _{i}T^{-i}\in \mathbb{P}, $ 
we denote by $ D_{  s \times k}(t) $   the following $ s \times k $ persymmetric matrix 
 $$   \left ( \begin{array} {cccccc}
\alpha _{1} & \alpha _{2} & \alpha _{3} &  \ldots & \alpha _{k-1}  &  \alpha _{k} \\
\alpha _{2 } & \alpha _{3} & \alpha _{4}&  \ldots  &  \alpha _{k} &  \alpha _{k+1} \\
\vdots & \vdots & \vdots    &  \vdots & \vdots  &  \vdots \\
\alpha _{s-1} & \alpha _{s} & \alpha _{s +1} & \ldots  &  \alpha _{s+k-3} &  \alpha _{s+k-2}  \\
\alpha _{s} & \alpha _{s+1} & \alpha _{s +2} & \ldots  &  \alpha _{s+k-2} &  \alpha _{s+k-1} 
 \end{array}  \right). $$ \vspace{0.5 cm}
 \item  Set $ t = \sum_{i\geq 1}\alpha _{i}T^{-i}\in \mathbb{P}, $ 
we denote by $ D_{  s \times (k-j+1)}^{j}(t) $   the following $ s \times (k-j+1) $ persymmetric matrix 
    $$   \left ( \begin{array} {cccccc}
\alpha _{j} & \alpha _{j+1} & \alpha _{j+2} &  \ldots & \alpha _{k-1}  &  \alpha _{k} \\
\alpha _{j+1 } & \alpha _{j+2} & \alpha _{j+3}&  \ldots  &  \alpha _{k} &  \alpha _{k+1} \\
\vdots & \vdots & \vdots    &  \vdots & \vdots  &  \vdots \\
\alpha _{j+s-2} & \alpha _{j+s-1} & \alpha _{j+s} & \ldots  &  \alpha _{k+s-3} &  \alpha _{k+s-2}  \\
\alpha _{j+s-1} & \alpha _{j+s} & \alpha _{j+s+1} & \ldots  &  \alpha _{k+s-2} &  \alpha _{k+s-1}\\
\end{array}  \right). $$ 
 \item  Set $ \eta  = \sum_{i\geq 1}\beta _{i}T^{-i}\in \mathbb{P}, $ 
we denote by $ D_{ (s +m )\times k}(\eta ) $   the following $ (s+m)\times k $ persymmetric matrix 
 $$   \left ( \begin{array} {cccccc}
\beta  _{1} & \beta  _{2} & \beta  _{3} & \ldots  &  \beta_{k-1} &  \beta _{k}  \\
\beta  _{2} & \beta  _{3} & \beta  _{4} & \ldots  &  \beta_{k} &  \beta _{k+1}  \\
\vdots & \vdots & \vdots    &  \vdots & \vdots  &  \vdots \\
\beta  _{m+1} & \beta  _{m+2} & \beta  _{m+3} & \ldots  &  \beta_{k+m-1} &  \beta _{k+m}  \\
\vdots & \vdots & \vdots    &  \vdots & \vdots  &  \vdots \\
\beta  _{s+m-1} & \beta  _{s+m} & \beta  _{s+m+1} & \ldots  &  \beta_{s+m+k-3} &  \beta _{s+m+k-2}  \\
\beta  _{s+m} & \beta  _{s+m+1} & \beta  _{s+m+2} & \ldots  &  \beta_{s+m+k-2} &  \beta _{s+m+k-1} 
\end{array}  \right). $$ \vspace{0.5 cm}
\end{itemize}
\end{defn}

\subsection{\textbf{ANALYSIS ON  THE  TWO-DIMENSIONAL  $\mathbb{K}$-VECTORSPACE}}
\label{subsec 1.2}
  Let $\mathbb{K}\times \mathbb{K}  = \mathbb{K}^2 $
 be the 2-dimensional vector space 
  over $ \mathbb{K}. $ Let $ (t,\eta )\in \mathbb{K}^2  $  and  $ \vert (t,\eta )\vert =
  sup \left\{\vert t \vert ,\vert \eta \vert \right\} =  2^{-inf(\mathfrak{v}(t),\mathfrak{v}(\eta ))}. $\\
 It is easy to see that $ (t,\eta )\longrightarrow   \vert (t,\eta )\vert $
  is an ultrametric valuation on $ \mathbb{K}^2, $ 
  that is,   $ (t,\eta )\longrightarrow   \vert (t,\eta )\vert $ is a norm and
  $  \vert((t,\eta ) + (t',\eta '))\vert \leq max\left\{\vert(t,\eta ) \vert ,\vert (t',\eta ')\vert \right\}. $ \\
   We denote by  $ d (t,\eta )=dtd\eta $ the  Haar measure on   $ \mathbb{K}^2 $
  chosen so that the measure on the unit interval of $ \mathbb{K}^2  $ is equal to one, thus 
 $$ \iint_{\mathbb{P}\times \mathbb{P}} d (t,\eta )=
   \int_{\mathbb{P}}dt\int_{\mathbb{P}}d\eta = 1\cdot 1 =1. $$\\
   
 Let  $(t,\eta ) = \big( \sum_{i= -\infty }^{-\mathfrak{v}(t)}t_{i}T^i ,
\sum_{i = -\infty }^{-\mathfrak{v}(\eta )}\eta _{i}T^i   \big) \in  \mathbb{K}^2, $ 
 we denote $\chi $  the  Character on  $(\mathbb{K}^2, +) $ defined by \\
\begin{equation*}
 \chi\big( \sum_{i= -\infty }^{-\mathfrak{v}(t)}t_{i}T^i ,
\sum_{ i = -\infty }^{-\mathfrak{v}(\eta )}\eta _{i}T^i   \big) =
  E\big( \sum_{i = -\infty }^{-\mathfrak{v}(t)}t_{i}T^i\big)\cdot
  E\big( \sum_{i = -\infty }^{-\mathfrak{v}(\eta )}\eta _{i}T^i\big) =
\begin{cases}
 1 & \text{if      }   t_{-1} +  \eta _{-1}= 0, \\
  -1 & \text{if      }   t_{-1} + \eta _{-1}=1.
    \end{cases}
\end{equation*}

\begin{defn}
\label{defn 1.2}We introduce the following definitions in the two-dimensional  $ \mathbb{K}$ -vectorspace.
\begin{itemize}
\item Let k, s and m denote rational integers such that $ k\geq 1,\; s\geq 2 \; and \; m\geq 0. $\\

\item We denote by $\mathbb{P}/\mathbb{P}_{i}\times \mathbb{P}/\mathbb{P}_{j} $
a complete set of coset representatives of $\mathbb{P}_{i}\times\mathbb{P}_{j} $
in  $\mathbb{P}\times\mathbb{P},$ for instance $\mathbb{P}/\mathbb{P}_{s+k-1}\times \mathbb{P}/\mathbb{P}_{s+m +k-1} $
denotes  a complete set of coset representatives of $\mathbb{P}_{s+k-1}\times\mathbb{P}_{s+m +k-1} $ in  $\mathbb{P}\times\mathbb{P}.$\\

\item $Set \;(t,\eta )= (\sum_{i\geq 1}\alpha _{i}T^{-i},\sum_{i\geq 1}\beta  _{i}T^{-i})
\in \mathbb{P}\times\mathbb{P}. $\\
  We denote by  $  D^{\left[\stackrel{s}{s+m}\right] \times k }(t,\eta  ) $ any  $(2s+m)\times k $   matrix, 
such that  after a rearrangement of the rows, if necessary,  we can  obtain  the following double persymmetric matrix 
$ \left[{D_{s  \times k}(t)\over D_{(s+m )\times k}(\eta )}\right] $ \\

 $$   \left ( \begin{array} {cccccc}
\alpha _{1} & \alpha _{2} & \alpha _{3} &  \ldots & \alpha _{k-1}  &  \alpha _{k} \\
\alpha _{2 } & \alpha _{3} & \alpha _{4}&  \ldots  &  \alpha _{k} &  \alpha _{k+1} \\
\vdots & \vdots & \vdots    &  \vdots & \vdots  &  \vdots \\
\alpha _{s-1} & \alpha _{s} & \alpha _{s +1} & \ldots  &  \alpha _{s+k-3} &  \alpha _{s+k-2}  \\
\alpha _{s} & \alpha _{s+1} & \alpha _{s +2} & \ldots  &  \alpha _{s+k-2} &  \alpha _{s+k-1}  \\
\hline \\
\beta  _{1} & \beta  _{2} & \beta  _{3} & \ldots  &  \beta_{k-1} &  \beta _{k}  \\
\beta  _{2} & \beta  _{3} & \beta  _{4} & \ldots  &  \beta_{k} &  \beta _{k+1}  \\
\vdots & \vdots & \vdots    &  \vdots & \vdots  &  \vdots \\
\beta  _{m+1} & \beta  _{m+2} & \beta  _{m+3} & \ldots  &  \beta_{k+m-1} &  \beta _{k+m}  \\
\vdots & \vdots & \vdots    &  \vdots & \vdots  &  \vdots \\
\beta  _{s+m-1} & \beta  _{s+m} & \beta  _{s+m+1} & \ldots  &  \beta_{s+m+k-3} &  \beta _{s+m+k-2}  \\
\beta  _{s+m} & \beta  _{s+m+1} & \beta  _{s+m+2} & \ldots  &  \beta_{s+m+k-2} &  \beta _{s+m+k-1} 
\end{array}  \right). $$ \vspace{0.5 cm}
We recall that the rank of a matrix does not change under elementary row operations.
\item   Let j be a rational integer such that $1\leq j\leq k-1,$
 $set \;(t,\eta )= (\sum_{i\geq 1}\alpha _{i}T^{-i},\sum_{i\geq 1}\beta  _{i}T^{-i})
\in \mathbb{P}\times\mathbb{P}. $\\
We denote by  $  D_{j}^{\left[\stackrel{s}{s+m}\right] \times (k-j+1) }(t,\eta  ) $ any 
$ (2s +m)\times (k-j+1) $ matrix such  that after a rearrangement of the rows, if necessary,  we  can obtain the following double persymmetric matrix
$ \left[{D_{s  \times (k-j+1)}^{j}(t)\over D_{(s+m )\times (k-j+1)}^{j}(\eta )}\right] $ \\

  $$   \left ( \begin{array} {cccccc}
\alpha _{j} & \alpha _{j+1} & \alpha _{j+2} &  \ldots & \alpha _{k-1}  &  \alpha _{k} \\
\alpha _{j+1 } & \alpha _{j+2} & \alpha _{j+3}&  \ldots  &  \alpha _{k} &  \alpha _{k+1} \\
\vdots & \vdots & \vdots    &  \vdots & \vdots  &  \vdots \\
\alpha _{j+s-2} & \alpha _{j+s-1} & \alpha _{j+s} & \ldots  &  \alpha _{k+s-3} &  \alpha _{k+s-2}  \\
\alpha _{j+s-1} & \alpha _{j+s} & \alpha _{j+s+1} & \ldots  &  \alpha _{k+s-2} &  \alpha _{k+s-1}\\
\hline \\
\beta  _{j} & \beta  _{j+1} & \beta  _{j+2} & \ldots  &  \beta_{k-1} &  \beta _{k}  \\
\beta  _{j+1} & \beta  _{j+2} & \beta  _{j+3} & \ldots  &  \beta_{k} &  \beta _{k+1}  \\
\vdots & \vdots & \vdots    &  \vdots & \vdots  &  \vdots \\
\beta  _{j+s-2} & \beta  _{j+s-1} & \beta  _{j+s} & \ldots  &  \beta_{k+s-3} &  \beta _{k+s-2}  \\
\beta  _{j+s-1} & \beta  _{j+s} & \beta  _{j+s+1} & \ldots  &  \beta_{k+s-2} &  \beta _{k+s-1}  \\
\beta  _{j+s} & \beta  _{j+s+1} & \beta  _{j+s+2} & \ldots  &  \beta_{k+s-1} &  \beta _{k+s}  \\
\vdots & \vdots & \vdots    &  \vdots & \vdots  &  \vdots \\
\beta  _{j+s+m-3} & \beta  _{j+s+m-2} & \beta  _{j+s+m-1} & \ldots  &  \beta_{k+s+m -4} &  \beta _{k+s+m-3}  \\
\beta  _{j+s+m-2} & \beta  _{j+s+m-1} & \beta  _{j+s+m} & \ldots  &  \beta_{k+s+m-3} &  \beta _{k+s+m-2}  \\
\beta  _{j+s+m-1} & \beta  _{j+s+m} & \beta  _{j+s+m+1} & \ldots  &  \beta_{k+s+m-2} &  \beta _{k+s+m-1}
\end{array}  \right). $$ 
For instance we can denote by  $  D_{j}^{\left[\stackrel{s}{s+m}\right] \times (k-j+1) }(t,\eta  ) $ the following matrix \\
 $$   \left ( \begin{array} {cccccc}
\alpha _{j} & \alpha _{j+1} & \alpha _{j+2} &  \ldots & \alpha _{k-1}  &  \alpha _{k} \\
 \beta  _{j} & \beta  _{j+1} & \beta  _{j+2} & \ldots  &  \beta_{k-1} &  \beta _{k}  \\
\alpha _{j+1 } & \alpha _{j+2} & \alpha _{j+3}&  \ldots  &  \alpha _{k} &  \alpha _{k+1} \\
\beta  _{j+1} & \beta  _{j+2} & \beta  _{j+3} & \ldots  &  \beta_{k} &  \beta _{k+1}  \\
\vdots & \vdots & \vdots    &  \vdots & \vdots  &  \vdots \\
\alpha _{j+s-2} & \alpha _{j+s-1} & \alpha _{j+s} & \ldots  &  \alpha _{k+s-3} &  \alpha _{k+s-2}  \\
\beta  _{j+s-2} & \beta  _{j+s-1} & \beta  _{j+s} & \ldots  &  \beta_{k+s-3} &  \beta _{k+s-2}  \\
\beta  _{j+s-1} & \beta  _{j+s} & \beta  _{j+s+1} & \ldots  &  \beta_{k+s-2} &  \beta _{k+s-1}  \\
\beta  _{j+s} & \beta  _{j+s+1} & \beta  _{j+s+2} & \ldots  &  \beta_{k+s-1} &  \beta _{k+s}  \\
\vdots & \vdots & \vdots    &  \vdots & \vdots  &  \vdots \\
\beta  _{j+s+m-3} & \beta  _{j+s+m-2} & \beta  _{j+s+m-1} & \ldots  &  \beta_{k+s+m -4} &  \beta _{k+s+m-3}  \\
\beta  _{j+s+m-2} & \beta  _{j+s+m-1} & \beta  _{j+s+m} & \ldots  &  \beta_{k+s+m-3} &  \beta _{s+m+k-2}  \\
\alpha _{j+s-1} & \alpha _{j+s} & \alpha _{j+s+1} & \ldots  &  \alpha _{k+s-2} &  \alpha _{k+s-1} \\
\beta  _{j+s+m-1} & \beta  _{j+s+m} & \beta  _{j+s+m+1} & \ldots  &  \beta_{k+s+m-2} &  \beta _{k+s+m-1}
\end{array}  \right). $$ \\

If j = 1 we denote  $  D_{1}^{\left[\stackrel{s}{s+m}\right] \times k }(t,\eta  ) $   by 
 $  D^{\left[\stackrel{s}{s+m}\right] \times k }(t,\eta  ). $ \\
 
\item We denote by  
 $  D^{\left[\stackrel{s-1}{\stackrel{s+m-1}{\alpha_{s -} + \beta_{s+m -} }}\right] \times k }(t,\eta  ) $
the following $(2s+m-1)\times k $ matrix, where the submatrix formed by  the first (2s +m -2) rows is equal to the matrix 
  $  D^{\left[\stackrel{s-1}{s+m-1}\right] \times k }(t,\eta  ), $ and the last row form
a  $ 1\times k $ persymmetric matrix of the form 
 $ (\alpha_{s } + \beta_{s+m }, \alpha_{s +1 } + \beta_{s+m+1 },\ldots, 
 \alpha_{s +k-1} + \beta_{s+k+m-1 })  $ \vspace{0.1 cm}\\
  $$   \left ( \begin{array} {cccccc}
\alpha _{1} & \alpha _{2} & \alpha _{3} &  \ldots & \alpha _{k-1}  &  \alpha _{k} \\
\alpha _{2 } & \alpha _{3} & \alpha _{4}&  \ldots  &  \alpha _{k} &  \alpha _{k+1} \\
\vdots & \vdots & \vdots    &  \vdots & \vdots  &  \vdots \\
\alpha _{s-1} & \alpha _{s} & \alpha _{s +1} & \ldots  &  \alpha _{s+k-3} &  \alpha _{s+k-2}  \\
\beta  _{1} & \beta  _{2} & \beta  _{3} & \ldots  &  \beta_{k-1} &  \beta _{k}  \\
\beta  _{2} & \beta  _{3} & \beta  _{4} & \ldots  &  \beta_{k} &  \beta _{k+1}  \\
\vdots & \vdots & \vdots    &  \vdots & \vdots  &  \vdots \\
\beta  _{m+1} & \beta  _{m+2} & \beta  _{m+3} & \ldots  &  \beta_{k+m-1} &  \beta _{k+m}  \\
\vdots & \vdots & \vdots    &  \vdots & \vdots  &  \vdots \\
\beta  _{s+m-1} & \beta  _{s+m} & \beta  _{s+m+1} & \ldots  &  \beta_{s+m+k-3} &  \beta _{s+m+k-2}  \\
\hline
\alpha _{s} + \beta  _{s+m} & \alpha _{s+1} + \beta  _{s+m+1} & \alpha _{s+2} + \beta  _{s+m+2} & \ldots 
& \alpha _{k+s-2} + \beta  _{k+s+m-2} &  \alpha _{k+s-1} + \beta  _{k+s+m-1}
\end{array}  \right). $$ \vspace{0.5 cm}
\item We denote by $\ker $  D the nullspace of the matrix D and r(D) the rank of the matrix D.\\

\item To simplify the notations concerning the exponential sums used in the proofs, we introduce the following definitions.\\

\item  Let $  g_{s,k,m}(t,\eta ) $ be the quadratic  exponential sum in $\mathbb{P}\times\mathbb{P}$ defined by
$$ (t,\eta ) \in  \mathbb{P}\times \mathbb{P}\longmapsto  
  \sum_{deg Y\leq k-1}\sum_{deg Z \leq  s-1}E(tYZ)\sum_{deg U \leq s+m-1}E(\eta YU) \in \mathbb{Z}.  $$\vspace{0.1 cm}\\
  \item  Let $  g(t,\eta ) $ be the quadratic  exponential sum in $\mathbb{P}\times\mathbb{P}$ defined by
$$ (t,\eta ) \in  \mathbb{P}\times \mathbb{P}\longmapsto  
  \sum_{deg Y\leq k-1}\sum_{deg Z \leq  s-1}E(tYZ)\sum_{deg U \leq s+m-1}E(\eta YU) \in \mathbb{Z}.  $$\vspace{0.1 cm}\\
  \item  Let $  g_{1}(t,\eta ) $ be the quadratic  exponential sum in $\mathbb{P}\times\mathbb{P}$ defined by
$$ (t,\eta ) \in  \mathbb{P}\times \mathbb{P}\longmapsto  
  \sum_{deg Y\leq k-1}\sum_{deg Z =  s-1}E(tYZ)\sum_{deg U \leq s+m-1}E(\eta YU) \in \mathbb{Z}.  $$\vspace{0.1 cm}\\
  \item  Let $  g_{2}(t,\eta ) $ be the quadratic  exponential sum in $\mathbb{P}\times\mathbb{P}$ defined by
$$ (t,\eta ) \in  \mathbb{P}\times \mathbb{P}\longmapsto  
  \sum_{deg Y\leq k-1}\sum_{deg Z \leq  s-1}E(tYZ)\sum_{deg U = s+m-1}E(\eta YU) \in \mathbb{Z}.  $$\vspace{0.1 cm}\\
  \item  Let $  f_{1}(t,\eta ) $ be the quadratic  exponential sum in $\mathbb{P}\times\mathbb{P}$ defined by
$$ (t,\eta ) \in  \mathbb{P}\times \mathbb{P}\longmapsto  
  \sum_{deg Y\leq k-1}\sum_{deg Z = s-1}E(tYZ)\sum_{deg U \leq s+m- 2}E(\eta YU) \in \mathbb{Z}.  $$\vspace{0.1 cm}\\
  \item  Let $  f_{2}(t,\eta ) $ be the quadratic  exponential sum in $\mathbb{P}\times\mathbb{P}$ defined by
$$ (t,\eta ) \in  \mathbb{P}\times \mathbb{P}\longmapsto  
  \sum_{deg Y\leq k-1}\sum_{deg Z \leq  s-2}E(tYZ)\sum_{deg U =  s+m-1}E(\eta YU) \in \mathbb{Z}.  $$\vspace{0.1 cm}\\
  \item  Let $  h(t,\eta ) $ be the quadratic  exponential sum in $\mathbb{P}\times\mathbb{P}$ defined by
$$ (t,\eta ) \in  \mathbb{P}\times \mathbb{P}\longmapsto  
  \sum_{deg Y\leq k-1}\sum_{deg Z =  s-1}E(tYZ)\sum_{deg U = s+m-1}E(\eta YU) \in \mathbb{Z}.  $$\vspace{0.1 cm}\\
  \item  Let $  v(t,\eta ) $ be the quadratic  exponential sum in $\mathbb{P}\times\mathbb{P}$ defined by
$$ (t,\eta ) \in  \mathbb{P}\times \mathbb{P}\longmapsto  
  \sum_{deg Y\leq k-1}\sum_{deg Z \leq  s-2}E(tYZ)\sum_{deg U \leq s+m-2}E(\eta YU) \in \mathbb{Z}.  $$\vspace{0.1 cm}\\
  \item  Let $  \psi (t,\eta ) $ be the quadratic  exponential sum in $\mathbb{P}\times\mathbb{P}$ defined by
$$ (t,\eta ) \in  \mathbb{P}\times \mathbb{P}\longmapsto  
  \sum_{deg Y =  k-1}\sum_{deg Z \leq  s-1}E(tYZ)\sum_{deg U \leq s+m-1}E(\eta YU) \in \mathbb{Z}.  $$\vspace{0.1 cm}\\
  \item   Let $  \phi (t,\eta ) $ be the quadratic  exponential sum in $\mathbb{P}\times\mathbb{P}$ defined by
$$ (t,\eta ) \in  \mathbb{P}\times \mathbb{P}\longmapsto  
  \sum_{deg Y =  k-1}\sum_{deg Z \leq  s-1}E(tYZ)\sum_{deg U = s+m-2}E(\eta YU) \in \mathbb{Z}.  $$\vspace{0.1 cm}\\
    \item  Let $  \phi _{1}(t,\eta ) $ be the quadratic  exponential sum in $\mathbb{P}\times\mathbb{P}$ defined by
$$ (t,\eta ) \in  \mathbb{P}\times \mathbb{P}\longmapsto  
  \sum_{deg Y =  k-1}\sum_{deg Z \leq  s- 1}E(tYZ)\sum_{deg U \leq s+m-2}E(\eta YU) \in \mathbb{Z}.  $$\vspace{0.1 cm}\\
    \item  Let $  \phi _{2}(t,\eta ) $ be the quadratic  exponential sum in $\mathbb{P}\times\mathbb{P}$ defined by
$$ (t,\eta ) \in  \mathbb{P}\times \mathbb{P}\longmapsto  
  \sum_{deg Y\leq k- 2}\sum_{deg Z \leq  s- 1}E(tYZ)\sum_{deg U = s+m- 1}E(\eta YU) \in \mathbb{Z}.  $$\vspace{0.1 cm}\\
  \item  Let $  \theta _{1}(t,\eta ) $ be the quadratic  exponential sum in $\mathbb{P}\times\mathbb{P}$ defined by
$$ (t,\eta ) \in  \mathbb{P}\times \mathbb{P}\longmapsto  
  \sum_{deg Y = k-1}\sum_{deg Z = s-1}E(tYZ)\sum_{deg U \leq s+m-2}E(\eta YU) \in \mathbb{Z}.  $$\vspace{0.1 cm}\\
    \item  Let $  \theta _{2}(t,\eta ) $ be the quadratic  exponential sum in $\mathbb{P}\times\mathbb{P}$ defined by
$$ (t,\eta ) \in  \mathbb{P}\times \mathbb{P}\longmapsto  
  \sum_{deg Y\leq k-2}\sum_{deg Z \leq  s-2}E(tYZ)\sum_{deg U = s+m- 1}E(\eta YU) \in \mathbb{Z}.  $$\vspace{0.1 cm}\\
    \item  Let $  \theta _{3}(t,\eta ) $ be the quadratic  exponential sum in $\mathbb{P}\times\mathbb{P}$ defined by
$$ (t,\eta ) \in  \mathbb{P}\times \mathbb{P}\longmapsto  
  \sum_{deg Y\leq k-2}\sum_{deg Z =  s - 1}E(tYZ)\sum_{deg U \leq s+m-2}E(\eta YU) \in \mathbb{Z}.  $$\vspace{0.1 cm}\\
    \item Consider the following partition of the  matrix $  D^{\left[\stackrel{s}{s+m}\right] \times k }(t,\eta  ) $\\
    
 $$   \left ( \begin{array} {cccccc}
\alpha _{1} & \alpha _{2} & \alpha _{3} &  \ldots & \alpha _{k-1}  &  \alpha _{k} \\
\alpha _{2 } & \alpha _{3} & \alpha _{4}&  \ldots  &  \alpha _{k} &  \alpha _{k+1} \\
\vdots & \vdots & \vdots    &  \vdots & \vdots  &  \vdots \\
\alpha _{s-1} & \alpha _{s} & \alpha _{s +1} & \ldots  &  \alpha _{s+k-3} &  \alpha _{s+k-2}  \\
\beta  _{1} & \beta  _{2} & \beta  _{3} & \ldots  &  \beta_{k-1} &  \beta _{k}  \\
\beta  _{2} & \beta  _{3} & \beta  _{4} & \ldots  &  \beta_{k} &  \beta _{k+1}  \\
\vdots & \vdots & \vdots    &  \vdots & \vdots  &  \vdots \\
\beta  _{m+1} & \beta  _{m+2} & \beta  _{m+3} & \ldots  &  \beta_{k+m-1} &  \beta _{k+m}  \\
\vdots & \vdots & \vdots    &  \vdots & \vdots  &  \vdots \\
\beta  _{s+m-1} & \beta  _{s+m} & \beta  _{s+m+1} & \ldots  &  \beta_{s+m+k-3} &  \beta _{s+m+k-2}  \\
\beta  _{s+m} & \beta  _{s+m+1} & \beta  _{s+m+2} & \ldots  &  \beta_{s+m+k-2} &  \beta _{s+m+k-1}\\
\hline \\
 \alpha _{s} & \alpha _{s+1} & \alpha _{s +2} & \ldots  &  \alpha _{s+k-2} &  \alpha _{s+k-1} 
\end{array}  \right). $$ \vspace{0.5 cm}\\

 We define $$  \sigma _{i,i}^{\left[\stackrel{s-1}{\stackrel{s+m }{\overline{\alpha_{s -}}}}\right] \times k } $$\\
 to be the cardinality of the following set \\
$$ \left\{(t,\eta )\in \mathbb{P}/\mathbb{P}_{k+s-1}\times \mathbb{P}/\mathbb{P}_{k+s+m-1}
\mid  r(D^{\big[\stackrel{s-1}{s+m}\big] \times k }(t,\eta ))  = r(D^{\big[\stackrel{s}{s+m}\big] \times k }(t,\eta )) = i
 \right\}.$$\vspace{0.1 cm}
 \item Consider the following partition of the  matrix $  D^{\left[\stackrel{s}{s+m -1}\right] \times k }(t,\eta  ) $\\
 
 $$   \left ( \begin{array} {cccccc}
\alpha _{1} & \alpha _{2} & \alpha _{3} &  \ldots & \alpha _{k-1}  &  \alpha _{k} \\
\alpha _{2 } & \alpha _{3} & \alpha _{4}&  \ldots  &  \alpha _{k} &  \alpha _{k+1} \\
\vdots & \vdots & \vdots    &  \vdots & \vdots  &  \vdots \\
\alpha _{s-1} & \alpha _{s} & \alpha _{s +1} & \ldots  &  \alpha _{s+k-3} &  \alpha _{s+k-2}  \\
\beta  _{1} & \beta  _{2} & \beta  _{3} & \ldots  &  \beta_{k-1} &  \beta _{k}  \\
\beta  _{2} & \beta  _{3} & \beta  _{4} & \ldots  &  \beta_{k} &  \beta _{k+1}  \\
\vdots & \vdots & \vdots    &  \vdots & \vdots  &  \vdots \\
\beta  _{m+1} & \beta  _{m+2} & \beta  _{m+3} & \ldots  &  \beta_{k+m-1} &  \beta _{k+m}  \\
\vdots & \vdots & \vdots    &  \vdots & \vdots  &  \vdots \\
\beta  _{s+m-1} & \beta  _{s+m} & \beta  _{s+m+1} & \ldots  &  \beta_{s+m+k-3} &  \beta _{s+m+k-2}  \\
\hline \\
 \alpha _{s} & \alpha _{s+1} & \alpha _{s +2} & \ldots  &  \alpha _{s+k-2} &  \alpha _{s+k-1} 
\end{array}  \right). $$ \vspace{0.5 cm}\\

   We define $$  \sigma _{i,i}^{\left[\stackrel{s-1}{\stackrel{s+m -1}{\overline{\alpha_{s -}}}}\right] \times k } $$ \\
 to be the cardinality of the following set \\
$$ \left\{(t,\eta )\in \mathbb{P}/\mathbb{P}_{k+s-1}\times \mathbb{P}/\mathbb{P}_{k+s+m-2}
\mid  r(D^{\big[\stackrel{s-1}{s-1+m}\big] \times k }(t,\eta ))  = r(D^{\big[\stackrel{s}{s+m-1}\big] \times k }(t,\eta )) = i
 \right\}.$$ \\
 
  \item We consider the following  partition of the  matrix $  D^{\left[\stackrel{s-1}{\stackrel{s+m-1}{\alpha_{s -} + \beta_{s+m -} }}\right] \times k }(t,\eta  ) $\vspace{0.1 cm}\\
   $$   \left ( \begin{array} {cccccc}
\alpha _{1} & \alpha _{2} & \alpha _{3} &  \ldots & \alpha _{k-1}  &  \alpha _{k} \\
\alpha _{2 } & \alpha _{3} & \alpha _{4}&  \ldots  &  \alpha _{k} &  \alpha _{k+1} \\
\vdots & \vdots & \vdots    &  \vdots & \vdots  &  \vdots \\
\alpha _{s-1} & \alpha _{s} & \alpha _{s +1} & \ldots  &  \alpha _{s+k-3} &  \alpha _{s+k-2}  \\
\beta  _{1} & \beta  _{2} & \beta  _{3} & \ldots  &  \beta_{k-1} &  \beta _{k}  \\
\beta  _{2} & \beta  _{3} & \beta  _{4} & \ldots  &  \beta_{k} &  \beta _{k+1}  \\
\vdots & \vdots & \vdots    &  \vdots & \vdots  &  \vdots \\
\beta  _{m+1} & \beta  _{m+2} & \beta  _{m+3} & \ldots  &  \beta_{k+m-1} &  \beta _{k+m}  \\
\vdots & \vdots & \vdots    &  \vdots & \vdots  &  \vdots \\
\beta  _{s+m-1} & \beta  _{s+m} & \beta  _{s+m+1} & \ldots  &  \beta_{s+m+k-3} &  \beta _{s+m+k-2}  \\
\hline
\alpha _{s} + \beta  _{s+m} & \alpha _{s+1} + \beta  _{s+m+1} & \alpha _{s+2} + \beta  _{s+m+2} & \ldots 
& \alpha _{k+s-2} + \beta  _{k+s+m-2} &  \alpha _{k+s-1} + \beta  _{k+s+m-1}
\end{array}  \right). $$ \vspace{0.5 cm}\\

 We define $$  \sigma _{i,i}^{\left[\stackrel{s-1}{\stackrel{s+m -1}{\overline{\alpha_{s -}+ \beta _{s+m-}}}}\right] \times k } $$ \\
 to be the cardinality of the following set \\
$$ \left\{(t,\eta )\in \mathbb{P}/\mathbb{P}_{k+s-1}\times \mathbb{P}/\mathbb{P}_{k+s+m-1}
\mid  r(D^{\big[\stackrel{s-1}{s-1+m}\big] \times k }(t,\eta ))  =  r( D^{\left[\stackrel{s-1}{\stackrel{s+m-1}{\alpha_{s -} + \beta_{s+m -} }}\right] \times k }(t,\eta  )) = i
  \right\}.$$
   \item   Consider the following partition of the  matrix $  D^{\left[\stackrel{s}{s+m}\right] \times k }(t,\eta  ) $\\
   
 $$   \left ( \begin{array} {cccccc}
\alpha _{1} & \alpha _{2} & \alpha _{3} &  \ldots & \alpha _{k-1}  &  \alpha _{k} \\
\alpha _{2 } & \alpha _{3} & \alpha _{4}&  \ldots  &  \alpha _{k} &  \alpha _{k+1} \\
\vdots & \vdots & \vdots    &  \vdots & \vdots  &  \vdots \\
\alpha _{s-1} & \alpha _{s} & \alpha _{s +1} & \ldots  &  \alpha _{s+k-3} &  \alpha _{s+k-2}  \\
\beta  _{1} & \beta  _{2} & \beta  _{3} & \ldots  &  \beta_{k-1} &  \beta _{k}  \\
\beta  _{2} & \beta  _{3} & \beta  _{4} & \ldots  &  \beta_{k} &  \beta _{k+1}  \\
\vdots & \vdots & \vdots    &  \vdots & \vdots  &  \vdots \\
\beta  _{m+1} & \beta  _{m+2} & \beta  _{m+3} & \ldots  &  \beta_{k+m-1} &  \beta _{k+m}  \\
\vdots & \vdots & \vdots    &  \vdots & \vdots  &  \vdots \\
\beta  _{s+m-1} & \beta  _{s+m} & \beta  _{s+m+1} & \ldots  &  \beta_{s+m+k-3} &  \beta _{s+m+k-2}  \\
\hline \\
 \alpha _{s} & \alpha _{s+1} & \alpha _{s +2} & \ldots  &  \alpha _{s+k-2} &  \alpha _{s+k-1} \\
 \hline \\
  \beta  _{s+m} & \beta  _{s+m+1} & \beta  _{s+m+2} & \ldots  &  \beta_{s+m+k-2} &  \beta _{s+m+k-1}
\end{array}  \right). $$ \vspace{0.5 cm}\\

   We define   $$  \sigma _{i,i,i}^{\left[\stackrel{s-1}{\stackrel{s+m-1 }
{\overline {\stackrel{\alpha_{s -}}{\beta_{s+m-} }}}}\right] \times k } $$\\ to be the cardinality of the following set \\
 \small
 $$ \left\{(t,\eta )\in \mathbb{P}/\mathbb{P}_{k+s-1}\times \mathbb{P}/\mathbb{P}_{k+s+m-1}
\mid  r(D^{\big[\stackrel{s-1}{s-1+m}\big] \times k }(t,\eta ))  = r(D^{\big[\stackrel{s}{s+m-1}\big] \times k }(t,\eta ))
= r(D^{\big[\stackrel{s-1}{s+m}\big] \times k }(t,\eta )) = i
 \right\}$$\\
  $$ =   \left\{(t,\eta )\in \mathbb{P}/\mathbb{P}_{k+s-1}\times \mathbb{P}/\mathbb{P}_{k+s+m-1}
\mid  r(D^{\big[\stackrel{s-1}{s-1+m}\big] \times k }(t,\eta ))  = r(D^{\big[\stackrel{s}{s+m-1}\big] \times k }(t,\eta ))
= r(D^{\big[\stackrel{s}{s+m}\big] \times k }(t,\eta )) = i
 \right\}.$$\\
   \item  Let  $ ( j_{1}, j_{2},   j_{3})  \in \mathbb{N}^{3}$ we define $$  \sigma _{j_{1},j_{2},j_{3}}^{\left[\stackrel{s-1}{\stackrel{s+m-1 }
{\overline {\stackrel{\alpha_{s -}}{\beta_{s+m-} }}}}\right] \times k } $$\\ to be the cardinality of the following set \\
 \small
 $$\begin{array}{l}\Big\{ (t,\eta ) \in \mathbb{P}/\mathbb{P}_{k+s -1}\times
           \mathbb{P}/\mathbb{P}_{k+s+m-1} 
\mid r(  D^{\left[\stackrel{s-1}{s+m-1}\right] \times k }(t,\eta  ) ) = j_{1}, \quad 
r( D^{\left[\stackrel{s}{s+m-1}\right] \times k }(t,\eta  ) ) = j_{2},  \\
 r(  D^{\left[\stackrel{s}{s+m}\right] \times k }(t,\eta  )  ) = j_{3}  \Big\}.
    \end{array}$$\\
   \item  Consider the following partition of the  matrix $  D^{\left[\stackrel{s}{s+m}\right] \times k }(t,\eta  ) $\\
   
 $$   \left ( \begin{array} {cccccc}
\alpha _{1} & \alpha _{2} & \alpha _{3} &  \ldots & \alpha _{k-1}  &  \alpha _{k} \\
\alpha _{2 } & \alpha _{3} & \alpha _{4}&  \ldots  &  \alpha _{k} &  \alpha _{k+1} \\
\vdots & \vdots & \vdots    &  \vdots & \vdots  &  \vdots \\
\alpha _{s-1} & \alpha _{s} & \alpha _{s +1} & \ldots  &  \alpha _{s+k-3} &  \alpha _{s+k-2}  \\
\beta  _{1} & \beta  _{2} & \beta  _{3} & \ldots  &  \beta_{k-1} &  \beta _{k}  \\
\beta  _{2} & \beta  _{3} & \beta  _{4} & \ldots  &  \beta_{k} &  \beta _{k+1}  \\
\vdots & \vdots & \vdots    &  \vdots & \vdots  &  \vdots \\
\beta  _{m+1} & \beta  _{m+2} & \beta  _{m+3} & \ldots  &  \beta_{k+m-1} &  \beta _{k+m}  \\
\vdots & \vdots & \vdots    &  \vdots & \vdots  &  \vdots \\
\beta  _{s+m-1} & \beta  _{s+m} & \beta  _{s+m+1} & \ldots  &  \beta_{s+m+k-3} &  \beta _{s+m+k-2}  \\
\hline \\
  \beta  _{s+m} & \beta  _{s+m+1} & \beta  _{s+m+2} & \ldots  &  \beta_{s+m+k-2} &  \beta _{s+m+k-1}\\
  \hline \\
   \alpha _{s} & \alpha _{s+1} & \alpha _{s +2} & \ldots  &  \alpha _{s+k-2} &  \alpha _{s+k-1}  
\end{array}  \right). $$ \vspace{0.5 cm}\\

     Let  $ ( j_{1}, j_{2},   j_{3})  \in \mathbb{N}^{3}$ we define $$  \sigma _{j_{1},j_{2},j_{3}}^{\left[\stackrel{s-1}{\stackrel{s+m-1 }
{\overline {\stackrel{\beta _{s +m -}}{\alpha _{s -} }}}}\right] \times k } $$\\ to be the cardinality of the following set \\
 \small
 $$\begin{array}{l}\Big\{ (t,\eta ) \in \mathbb{P}/\mathbb{P}_{k+s -1}\times
           \mathbb{P}/\mathbb{P}_{k+s+m-1} 
\mid r(  D^{\left[\stackrel{s-1}{s+m-1}\right] \times k }(t,\eta  ) ) = j_{1}, \quad 
r( D^{\left[\stackrel{s-1}{s+m}\right] \times k }(t,\eta  ) ) = j_{2},  \\
 r(  D^{\left[\stackrel{s}{s+m}\right] \times k }(t,\eta  )  ) = j_{3}  \Big\}.
    \end{array}$$\\ 
 \item  Consider the following partition of the  matrix $  D^{\left[\stackrel{s}{s+m}\right] \times k }(t,\eta  ) $\\
 
 $$   \left ( \begin{array} {ccccc|c}
\alpha _{1} & \alpha _{2} & \alpha _{3} &  \ldots & \alpha _{k-1}  &  \alpha _{k} \\
\alpha _{2 } & \alpha _{3} & \alpha _{4}&  \ldots  &  \alpha _{k} &  \alpha _{k+1} \\
\vdots & \vdots & \vdots    &  \vdots & \vdots  &  \vdots \\
\alpha _{s-1} & \alpha _{s} & \alpha _{s +1} & \ldots  &  \alpha _{s+k-3} &  \alpha _{s+k-2}  \\
\beta  _{1} & \beta  _{2} & \beta  _{3} & \ldots  &  \beta_{k-1} &  \beta _{k}  \\
\beta  _{2} & \beta  _{3} & \beta  _{4} & \ldots  &  \beta_{k} &  \beta _{k+1}  \\
\vdots & \vdots & \vdots    &  \vdots & \vdots  &  \vdots \\
\beta  _{m+1} & \beta  _{m+2} & \beta  _{m+3} & \ldots  &  \beta_{k+m-1} &  \beta _{k+m}  \\
\vdots & \vdots & \vdots    &  \vdots & \vdots  &  \vdots \\
\beta  _{s+m-1} & \beta  _{s+m} & \beta  _{s+m+1} & \ldots  &  \beta_{s+m+k-3} &  \beta _{s+m+k-2}  \\
\hline \\
  \alpha _{s} & \alpha _{s+1} & \alpha _{s +2} & \ldots  &  \alpha _{s+k-2} &  \alpha _{s+k-1}  \\
  \hline \\
  \beta  _{s+m} & \beta  _{s+m+1} & \beta  _{s+m+2} & \ldots  &  \beta_{s+m+k-2} &  \beta _{s+m+k-1}
\end{array}  \right). $$ \vspace{0.5 cm}
      Let  $ ( j_{1}, j_{2},   j_{3}, j_{4},j_{5}, j_{6}) \in \mathbb{N}^{6}, $ we define \\
\begin{align*}
{}^{\#}\left(\begin{array}{c | c}
           j_{1} & j_{2} \\
           \hline
           j_{3} & j_{4}\\
           \hline
            j_{5} & j_{6} 
           \end{array} \right)_{\mathbb{P}/\mathbb{P}_{k+s -1}\times
           \mathbb{P}/\mathbb{P}_{k+s+m-1} }^{{\alpha  \over \beta }}
            \end{align*}
    to be the cardinality of the following set
    $$\begin{array}{l}\Big\{ (t,\eta ) \in \mathbb{P}/\mathbb{P}_{k+s -1}\times
           \mathbb{P}/\mathbb{P}_{k+s+m-1} 
\mid r(  D^{\left[\stackrel{s-1}{s+m-1}\right] \times (k-1) }(t,\eta  ) ) = j_{1}, \quad 
r( D^{\left[\stackrel{s-1}{s+m-1}\right] \times k }(t,\eta  ) ) = j_{2},  \\
 r(  D^{\left[\stackrel{s}{s+m-1}\right] \times (k-1) }(t,\eta  )  ) = j_{3},\quad  
  r( D^{\left[\stackrel{s}{s+m-1}\right] \times k }(t,\eta  ) ) = j_{4}, \\
   r(  D^{\left[\stackrel{s}{s+m}\right] \times (k-1) }(t,\eta  )  ) = j_{5}, \quad
  r( D^{\left[\stackrel{s}{s+m}\right] \times k }(t,\eta  ) ) = j_{6} \Big\}.
    \end{array}$$\\
    \item  Consider the following partition of the  matrix $  D^{\left[\stackrel{s}{s+m}\right] \times k }(t,\eta  ) $
 $$   \left ( \begin{array} {ccccc|c}
\alpha _{1} & \alpha _{2} & \alpha _{3} &  \ldots & \alpha _{k-1}  &  \alpha _{k} \\
\alpha _{2 } & \alpha _{3} & \alpha _{4}&  \ldots  &  \alpha _{k} &  \alpha _{k+1} \\
\vdots & \vdots & \vdots    &  \vdots & \vdots  &  \vdots \\
\alpha _{s-1} & \alpha _{s} & \alpha _{s +1} & \ldots  &  \alpha _{s+k-3} &  \alpha _{s+k-2}  \\
\beta  _{1} & \beta  _{2} & \beta  _{3} & \ldots  &  \beta_{k-1} &  \beta _{k}  \\
\beta  _{2} & \beta  _{3} & \beta  _{4} & \ldots  &  \beta_{k} &  \beta _{k+1}  \\
\vdots & \vdots & \vdots    &  \vdots & \vdots  &  \vdots \\
\beta  _{m+1} & \beta  _{m+2} & \beta  _{m+3} & \ldots  &  \beta_{k+m-1} &  \beta _{k+m}  \\
\vdots & \vdots & \vdots    &  \vdots & \vdots  &  \vdots \\
\beta  _{s+m-1} & \beta  _{s+m} & \beta  _{s+m+1} & \ldots  &  \beta_{s+m+k-3} &  \beta _{s+m+k-2}  \\
\hline \\
  \beta  _{s+m} & \beta  _{s+m+1} & \beta  _{s+m+2} & \ldots  &  \beta_{s+m+k-2} &  \beta _{s+m+k-1}\\
  \hline \\
   \alpha _{s} & \alpha _{s+1} & \alpha _{s +2} & \ldots  &  \alpha _{s+k-2} &  \alpha _{s+k-1}  
\end{array}  \right). $$ \vspace{0.5 cm}
      Let  $( j_{1}, j_{2},   j_{3}, j_{4},j_{5}, j_{6}) \in \mathbb{N}^{6}, $ we define \\
\begin{align*}
{}^{\#}\left(\begin{array}{c | c}
           j_{1} & j_{2} \\
           \hline
           j_{3} & j_{4}\\
           \hline
            j_{5} & j_{6} 
           \end{array} \right)_{\mathbb{P}/\mathbb{P}_{k+s -1}\times
           \mathbb{P}/\mathbb{P}_{k+s+m-1} }^{{\beta  \over \alpha }}
            \end{align*}
    to be the cardinality of the following set
    $$\begin{array}{l}\Big\{ (t,\eta ) \in \mathbb{P}/\mathbb{P}_{k+s -1}\times
           \mathbb{P}/\mathbb{P}_{k+s+m-1} 
\mid r(  D^{\left[\stackrel{s-1}{s+m-1}\right] \times (k-1) }(t,\eta  ) ) = j_{1}, \quad 
r( D^{\left[\stackrel{s-1}{s+m-1}\right] \times k }(t,\eta  ) ) = j_{2},  \\
 r(  D^{\left[\stackrel{s -1}{s+m}\right] \times (k-1) }(t,\eta  )  ) = j_{3},\quad  
  r( D^{\left[\stackrel{s -1}{s+m}\right] \times k }(t,\eta  ) ) = j_{4}, \\
   r(  D^{\left[\stackrel{s}{s+m}\right] \times (k-1) }(t,\eta  )  ) = j_{5}, \quad
  r( D^{\left[\stackrel{s}{s+m}\right] \times k }(t,\eta  ) ) = j_{6} \Big\}.
    \end{array}$$\\ 
   \item   Consider the following partition of the  matrix $  D^{\left[\stackrel{s}{s+m}\right] \times k }(t,\eta  ) $\\
   
 $$   \left ( \begin{array} {ccccc|c}
\alpha _{1} & \alpha _{2} & \alpha _{3} &  \ldots & \alpha _{k-1}  &  \alpha _{k} \\
\alpha _{2 } & \alpha _{3} & \alpha _{4}&  \ldots  &  \alpha _{k} &  \alpha _{k+1} \\
\vdots & \vdots & \vdots    &  \vdots & \vdots  &  \vdots \\
\alpha _{s-1} & \alpha _{s} & \alpha _{s +1} & \ldots  &  \alpha _{s+k-3} &  \alpha _{s+k-2}  \\
\beta  _{1} & \beta  _{2} & \beta  _{3} & \ldots  &  \beta_{k-1} &  \beta _{k}  \\
\beta  _{2} & \beta  _{3} & \beta  _{4} & \ldots  &  \beta_{k} &  \beta _{k+1}  \\
\vdots & \vdots & \vdots    &  \vdots & \vdots  &  \vdots \\
\beta  _{m+1} & \beta  _{m+2} & \beta  _{m+3} & \ldots  &  \beta_{k+m-1} &  \beta _{k+m}  \\
\vdots & \vdots & \vdots    &  \vdots & \vdots  &  \vdots \\
\beta  _{s+m-1} & \beta  _{s+m} & \beta  _{s+m+1} & \ldots  &  \beta_{s+m+k-3} &  \beta _{s+m+k-2}  \\
  \alpha _{s} & \alpha _{s+1} & \alpha _{s +2} & \ldots  &  \alpha _{s+k-2} &  \alpha _{s+k-1}  \\
  \hline \\
  \beta  _{s+m} & \beta  _{s+m+1} & \beta  _{s+m+2} & \ldots  &  \beta_{s+m+k-2} &  \beta _{s+m+k-1}   
\end{array}  \right). $$ \vspace{0.5 cm}
   Let $  ( j_{1}, j_{2},   j_{3}, j_{4}) \in \mathbb{N}^{4}, $ we define \\
\begin{align*}
{}^{\#}\left(\begin{array}{c | c}
           j_{1} & j_{2} \\
           \hline
           j_{3} & j_{4}
           \end{array} \right)_{\mathbb{P}/\mathbb{P}_{k+s -1}\times
           \mathbb{P}/\mathbb{P}_{k+s+m-1} }^{{\alpha \over \beta }}
            \end{align*}
    to be the cardinality of the following set
    $$\begin{array}{l}\Big\{ (t,\eta ) \in \mathbb{P}/\mathbb{P}_{k+s -1}\times
           \mathbb{P}/\mathbb{P}_{k+s+m-1} 
\mid r(  D^{\left[\stackrel{s}{s+m-1}\right] \times (k-1) }(t,\eta  ) ) = j_{1}, \quad 
r( D^{\left[\stackrel{s}{s+m-1}\right] \times k }(t,\eta  ) ) = j_{2},  \\
 r(  D^{\left[\stackrel{s}{s+m}\right] \times (k-1) }(t,\eta  )  ) = j_{3},\quad  
  r( D^{\left[\stackrel{s}{s+m}\right] \times k }(t,\eta  ) ) = j_{4} \Big\}.
    \end{array}$$\\
\item $Set \;(t,\eta )= (\sum_{i\geq 1}\alpha _{i}T^{-i},\sum_{i\geq 1}\beta  _{i}T^{-i})
\in \mathbb{P}\times\mathbb{P}. $\\
Let $  \Gamma _{i}^{\Big[\substack{s \\ s+m }\Big] \times k} $ denote the number of 
double persymmetric $ (2s+m)\times k $ rank i  matrices  of the form $ \left[{D_{s  \times k}(t)\over D_{(s+m )\times k}(\eta )}\right], $that is \vspace{0.1 cm}\\
 $  \Gamma _{i}^{\Big[\substack{s \\ s+m }\Big] \times k}  =   
Card  \left\{(t,\eta )\in \mathbb{P}/\mathbb{P}_{k+s-1}\times \mathbb{P}/\mathbb{P}_{k+s+m-1}
\mid   r(D^{\big[\stackrel{s}{s+m}\big] \times k }(t,\eta )) = i
 \right\}. $ \\
  \item  Set for $ 0\leq i\leq \inf(2s+m,k) $ $$ \Delta _{i}^{\Big[\substack{s \\ s+m }\Big] \times k} = 
     \sigma _{i,i,i}^{\left[\stackrel{s-1}{\stackrel{s+m-1 }{\overline {\stackrel{\alpha_{s -}}{\beta_{s+m-} }}}}\right] \times k }
  - 3\cdot \sigma _{i-1,i-1,i-1}^{\left[\stackrel{s-1}{\stackrel{s+m-1 }{\overline {\stackrel{\alpha_{s -}}{\beta_{s+m-} }}}}\right] \times k } 
  + 2\cdot \sigma _{i-2,i-2,i-2}^{\left[\stackrel{s-1}{\stackrel{s+m-1 }{\overline {\stackrel{\alpha_{s -}}{\beta_{s+m-} }}}}\right] \times k }.  $$ \vspace{0.1 cm}\\
  \item Set  for $ 0\leq j\leq s+m,\; k\geq s+j+1 $  $$  \Omega _{j+1}(1,m+s-1,k) = 
 \Gamma _{j+1}^{\Big[\substack{1 \\ 1+m+(s-1)}\Big] \times k} - 4\cdot \Gamma _{j}^{\Big[\substack{ 1 \\ 1+m+(s-2) }\Big] \times k}. $$ \vspace{0.1 cm} \\
 \item We define    
     \begin{equation*}
\begin{vmatrix}
k-1 & \vline & k  & \vline \\
\hline
\cdot &\vline & \cdot & \vline & D^{\left[s-1\atop s+m-1 \right]\times \cdot}  \\
\hline
j &\vline & j &  \vline  & \alpha _{s}-\\
\hline 
j &\vline & j+1 & \vline  & \beta _{s+m}- 
\end{vmatrix}
\end{equation*}
    to be  the following subset of $\mathbb{P}^{2}$
    $$\begin{array}{l}\Big\{ (t,\eta  ) = (\sum_{i\geq 1}\alpha _{i}T^{-i},\sum_{i\geq 1}\beta  _{i}T^{-i}) \in \mathbb{P}^2
   \mid r( D^{ \left[s\atop s+m-1\right]\times (k-1)}(t,\eta  ) )  = j, \\
r( D^{ \left[s\atop s+m-1\right]\times k}(t,\eta ) )  = j,\quad
r( D^{ \left[s\atop s+m \right]\times (k-1)}(t,\eta  ) ) = j, \quad 
 r(  D^{ \left[s\atop s+m \right]\times k}(t,\eta  ) ) = j+1 \Big\}  
   \end{array}$$\\ 
\item   
 We define    
     \begin{equation*}
\begin{vmatrix}
k-1 & \vline & k  & \vline \\
\hline
j &\vline & j & \vline & D^{\left[s-1\atop s+m-1 \right]\times \cdot}  \\
\hline
j &\vline & j+1 &  \vline  & \alpha _{s}-\\
\hline 
j &\vline & j+1 & \vline  & \beta _{s+m}- 
\end{vmatrix}
\end{equation*}
    to be  the following subset of $\mathbb{P}^{2}$
    $$\begin{array}{l}\Big\{ (t,\eta  ) = (\sum_{i\geq 1}\alpha _{i}T^{-i},\sum_{i\geq 1}\beta  _{i}T^{-i}) \in \mathbb{P}^2
   \mid r( D^{ \left[s-1\atop s+m-1\right]\times (k-1)}(t,\eta  ) )  = j, \\
r( D^{ \left[s-1\atop s+m-1\right]\times k}(t,\eta ) )  = j,\quad
r( D^{ \left[s\atop s+m-1\right]\times (k-1)}(t,\eta ) )  = j,\quad
r( D^{ \left[s\atop s+m-1\right]\times k}(t,\eta ) )  = j+1,\\
r( D^{ \left[s\atop s+m \right]\times (k-1)}(t,\eta  ) ) = j, \quad 
 r(  D^{ \left[s\atop s+m \right]\times k}(t,\eta  ) ) = j+1 \Big\}  
   \end{array}$$\\ 
\item Similar expressions are defined in a similar way.   
  \item  We denote by  $ R_{q}(k,s,m) $ the number of solutions \\
 $(Y_1,Z_1,U_{1}, \ldots,Y_q,Z_q,U_{q}) $  of the polynomial equations
   \[\left\{\begin{array}{c}
 Y_{1}Z_{1} +Y_{2}Z_{2}+ \ldots + Y_{q}Z_{q} = 0  \\
   Y_{1}U_{1} + Y_{2}U_{2} + \ldots  + Y_{q}U_{q} = 0
 \end{array}\right.\]
  satisfying the degree conditions \\
                   $$  degY_i \leq k-1 , \quad degZ_i \leq s-1 ,\quad degU_{i}\leq s+m-1 \quad for \quad 1\leq i \leq q. $$ \\                           
     \end{itemize}
\end{defn}

\section{\textbf{INTRODUCTION (REFER  TO  SECTION  1 AND SECTION 3)}}
\label{sec 2}

The rational function field  $\mathbb{F}_{2}(T) $ is completed with respect to an appropriate valuation 
to a field $ \mathbb{K} $ (i.e. the field of Laurent Series). The unit interval of $\mathbb{K},$ that is, the open ball
of radius 1 about 0, is a compact additive group. We shall use the Haar integral on this group.
Let $ \mathbb{K}^{2} $ be the 2-dimensional vectorspace over  $ \mathbb{K}. $\vspace{0.2 cm}\\
The main result of this paper is to obtain a formula for the number  $  \Gamma _{i}^{\Big[\substack{s \\ s+m }\Big] \times k} $
 of double persymmetric $ (2s+m)\times k $ rank i  matrices.\vspace{0.1 cm}\\
 In particular we shall show that if $0\leq i\leq \inf(s-1,k-1)$  the number of double persymmetric  $(2s+m)\times k$ rank i matrices 
depends only on i and if $i\geq s+1,\;k\geq i$ the number of rank i matrices can be obtained  from the formula for the number of rank s+1 matrices.
 The computation of the  $  \Gamma _{i}^{\Big[\substack{s \\ s+m }\Big] \times k} $ is based on a formula (given in [4]) of the number of rank i matrices of the form  $\left[A\over b_{-}\right] $ 
 where A is persymmetric and $ b_{-} $ a one-row matrix with entries in  $ \mathbb{F}_{2}.$
We observe that the  $  \Gamma _{i}^{\Big[\substack{s \\ s+m }\Big] \times k}$ where  $ 0\leq i\leq \inf(2s+m,k) $
are solutions to the system 
 \begin{equation}
  \label{eq 2.1}
 \begin{cases} 
 \displaystyle \sum_{ i = 0}^{\inf(2s+m,k)} \Gamma _{i}^{\Big[\substack{s  \\ s+m }\Big] \times k}  = 2^{2k+2s+m-2}, \\
 \displaystyle  \sum_{i = 0}^{\inf(2s+m,k)} \Gamma _{i}^{\Big[\substack{s  \\ s+m }\Big] \times k}\cdot2^{-i}  =  2^{k+2s+m-2} + 2^{2k-2} - 2^{k-2}. 
\end{cases}
    \end{equation}\vspace{0.1 cm}\\
   The first equation is obvious.  \vspace{0.1 cm}\\  
The second equation is a consequence of the identity \vspace{0.1 cm}\\  
    \begin{align*}
 \displaystyle
   \int_{\mathbb{P}\times \mathbb{P}} g(t,\eta ) dt d \eta &  =
 Card \left\{(Y,Z,U), deg Y \leq k-1, deg Z \leq s-1,  deg U \leq s+m-1 \mid Y \cdot Z = Y \cdot U = 0 \right\}\\
 &  = 2^{2s+m} + 2^{k} -1. 
\end{align*}

\subsection{\textbf{SKETCH OF THE PROOF}}
\label{subsec 2.1}

\begin{itemize}
\item  Consider  the following $\left[s\atop s+m\right]\times k $ double persymmetric matrix \\

 $$   \left ( \begin{array} {cccccc}
\alpha _{1} & \alpha _{2} & \alpha _{3} &  \ldots & \alpha _{k-1}  &  \alpha _{k} \\
\alpha _{2 } & \alpha _{3} & \alpha _{4}&  \ldots  &  \alpha _{k} &  \alpha _{k+1} \\
\vdots & \vdots & \vdots    &  \vdots & \vdots  &  \vdots \\
\alpha _{s-1} & \alpha _{s} & \alpha _{s +1} & \ldots  &  \alpha _{s+k-3} &  \alpha _{s+k-2}  \\
\alpha _{s} & \alpha _{s+1} & \alpha _{s +2} & \ldots  &  \alpha _{s+k-2} &  \alpha _{s+k-1}  \\
\hline \\
\beta  _{1} & \beta  _{2} & \beta  _{3} & \ldots  &  \beta_{k-1} &  \beta _{k}  \\
\beta  _{2} & \beta  _{3} & \beta  _{4} & \ldots  &  \beta_{k} &  \beta _{k+1}  \\
\vdots & \vdots & \vdots    &  \vdots & \vdots  &  \vdots \\
\beta  _{m+1} & \beta  _{m+2} & \beta  _{m+3} & \ldots  &  \beta_{k+m-1} &  \beta _{k+m}  \\
\vdots & \vdots & \vdots    &  \vdots & \vdots  &  \vdots \\
\beta  _{s+m-1} & \beta  _{s+m} & \beta  _{s+m+1} & \ldots  &  \beta_{s+m+k-3} &  \beta _{s+m+k-2}  \\
\beta  _{s+m} & \beta  _{s+m+1} & \beta  _{s+m+2} & \ldots  &  \beta_{s+m+k-2} &  \beta _{s+m+k-1} 
\end{array}  \right). $$ \vspace{0.5 cm}
\item  We establish  the following recurrent formula for the number $ \Gamma _{i}^{\Big[\substack{s \\ s+m }\Big] \times k} $ of $ \left[s\atop s+m\right]\times k $ double persymmetric matrices of rank i: \vspace{0.2 cm}\\
   \begin{align*}
(1)\quad  \Gamma _{i}^{\Big[\substack{s \\ s+m }\Big] \times k}  = 2\cdot \Gamma _{i-1}^{\Big[\substack{s -1\\ s-1+(m+1) }\Big] \times k}
+ 4\cdot \Gamma _{i-1}^{\Big[\substack{s \\ s+(m-1) }\Big] \times k} - 8\cdot \Gamma _{i-2}^{\Big[\substack{s-1 \\ s-1+m }\Big] \times k}
 + \Delta _{i}^{\Big[\substack{s \\ s+m }\Big] \times k}.
 \end{align*} 
 We observe that the  $  \Gamma _{i}^{\Big[\substack{s \\ s+m }\Big] \times k},\; for \; 0\leq i\leq \inf(2s+m,k) $
are solutions to the system 
 \begin{equation*}
(2)\quad   \begin{cases} 
 \sum_{ i = 0}^{\inf(2s+m,k)} \Gamma _{i}^{\Big[\substack{s  \\ s+m }\Big] \times k}  = 2^{2k+2s+m-2}, \\
  \sum_{i = 0}^{\inf(2s+m,k)} \Gamma _{i}^{\Big[\substack{s  \\ s+m }\Big] \times k}\cdot2^{-i}
    =  2^{k+2s+m-2} + 2^{2k-2} - 2^{k-2}. 
\end{cases}
    \end{equation*}\vspace{0.1 cm}\\
     \item    We establish that the remainder in (1) satisfies the following equality  $$ (3) \quad \Delta _{i}^{\Big[\substack{s \\ s+m }\Big] \times k} 
  =  \Delta _{i}^{\Big[\substack{s \\ s+m }\Big] \times (i+1)} \quad  for \; all \quad 0\leq i\leq 2s+m,\; k\geq i+1. $$\\
  
  \item     We establish from (1)  and  (3) a recurrent formula for the difference \vspace{0.2 cm}\\
     $$ \Gamma _{i}^{\Big[\substack{s \\ s +m }\Big] \times (k+1)} -
    \Gamma _{i}^{\Big[\substack{s \\ s +m }\Big] \times k} \quad for\; 2\leq i \leq 2s+m,\quad k\geq i+1.$$
    In fact we prove that \vspace{0.2 cm}\\
    \Small
       \begin{align*}
(4) \quad  &  \Gamma _{i}^{\Big[\substack{s \\ s +m }\Big] \times (k+1)} -  \Gamma _{i}^{\Big[\substack{s \\ s +m }\Big] \times k} =  4\cdot \left[\Gamma _{i-1}^{\Big[\substack{s \\ s +(m-1) }\Big] \times (k+1)} - \Gamma _{i-1}^{\Big[\substack{s \\ s +(m-1) }\Big] \times  k}\right]  
 \quad \text{if    }\quad  1 \leq i\leq s-1,\quad k\geq i+1,  \\
(5) \quad  &  \Gamma _{s+j}^{\Big[\substack{s \\ s +m }\Big] \times (k+1)} -  \Gamma _{s+j}^{\Big[\substack{s \\ s +m }\Big] \times k} =  4\cdot \left[\Gamma _{s+j -1}^{\Big[\substack{s \\ s +(m-1) }\Big] \times (k+1)} - \Gamma _{s+j -1}^{\Big[\substack{s \\ s +(m-1) }\Big] \times  k}\right]  
+ R(j,s,m,k) \quad \text{if  }  0 \leq j\leq \inf(s+m, k-s+1) 
       \end{align*}
  where R(j,s,m,k)  is equal to \\
  \begin{equation*}
    2^{s-1}\left[ \Gamma _{j+1}^{\Big[\substack{1 \\ 1+ (s+m-1)}\Big] \times (k+1)} -  \Gamma _{j+1}^{\Big[\substack{1 \\ 1+(s+m-1) }\Big] \times k}\right] 
   -   2^{s+1}\left[ \Gamma _{j}^{\Big[\substack{1 \\ 1+ (s+m-2)}\Big] \times (k+1)} -  \Gamma _{j}^{\Big[\substack{1 \\ 1+(s+m-2) }\Big] \times k}\right]. 
  \end{equation*}
\normalsize
\item  From a formula of $ \Gamma _{i}^{\Big[\substack{1 \\ 1 +m }\Big] \times k} $ obtained in [4] we compute 
the remainder    R(j,s,m,k) in  (5). \vspace{0.1 cm}\\ 
    In fact we prove  \vspace{0.1 cm}\\ 
    \begin{equation*}
         (6)\quad       R(j,s,m,k)   = \begin{cases}
2^{k+s-1} & \text{if  } j = 0 ,\quad k\geq s+1, \\
- 2^{k+s-1} &  \text{if  }    j =1, \quad k\geq s+2, \\
 0  & \text{if   } 2\leq j\leq s+m-1, \quad k \geq s+j+1, \\
  -3\cdot 2^{2k+2s+m-2}         & \text{if   } j = s+m, \quad k \geq 2s+m +1. 
\end{cases}
\end{equation*}\vspace{0.1 cm}\\ 
\item   From (4) applying (2) we deduce by induction on i  
 \begin{equation*}
(7) \quad  \Gamma _{i}^{\Big[\substack{s \\ s +m }\Big] \times k}  = \begin{cases}
  \Gamma _{i}^{\Big[\substack{i \\ i }\Big] \times (i+1)} = 
21\cdot2^{3i-4} - 3\cdot2^{2i-3}        & \text{if  } \quad 1\leq i\leq s-1,\;k\geq i+1, \\
2^{2s+2i+m-2} -3\cdot2^{3i-4} + 2^{2i-3} &  \text{if  }\quad 1\leq i = k\leq  s+1.  
\end{cases}
\end{equation*}\vspace{0.1 cm}\\ 
We have now computed $\Gamma _{i}^{\Big[\substack{s \\ s +m }\Big] \times k}\; for \;1\leq i\leq s-1,\;k\geq i. $\vspace{0.1 cm}\\ 
\item    From (5) applying  (6) we obtain \vspace{0.1 cm}\\ 
\begin{align*}
(8)\quad   \Gamma _{s+j}^{\Big[\substack{s \\ s +m }\Big] \times (k+1)} -  \Gamma _{s+j}^{\Big[\substack{s \\ s +m }\Big] \times k} =  4\cdot \left[\Gamma _{s+j-1}^{\Big[\substack{s \\ s +(m-1) }\Big] \times (k+1)} - \Gamma _{s+j-1}^{\Big[\substack{s \\ s +(m-1) }\Big] \times  k}\right]  \\
+   \begin{cases}
2^{k+s-1} & \text{if  } j = 0 ,\quad k\geq s+1, \\
- 2^{k+s-1} &  \text{if  }    j =1, \quad k\geq s+2, \\
 0  & \text{if   } 2\leq j\leq s+m-1, \quad k \geq s+j+1, \\
  -3\cdot 2^{2k+2s+m-2}         & \text{if   } j = s+m, \quad k \geq 2s+m +1. 
\end{cases}
\end{align*}
\item  Using successively the recurrent formula in (8) we deduce  \vspace{0.1 cm}\\ 
 In the case $ m \in \left\{0,1\right\} $ \vspace{0.01 cm}\\
 \begin{equation*}
 \Gamma _{s+j}^{\Big[\substack{s \\ s}\Big] \times (k+1)}
  - \Gamma _{s+j}^{\Big[\substack{s \\ s}\Big] \times k}=
  \begin{cases}
 3\cdot 2^{k+s-1}  & \text{if  } j = 0,\; k >s, \\
 21\cdot 2^{k+s+3j-4}   &  \text{if  }    1\leq j\leq s-1,\; k>s+j, \\
3\cdot2^{2k+2s-2}- 3\cdot2^{k+4s - 4}       & \text{if   } j  = s,  \; k>2s,
\end{cases}
\end{equation*}
 \begin{equation*}
 \Gamma _{s+j}^{\Big[\substack{s \\ s+1}\Big] \times (k+1)}
  - \Gamma _{s+j}^{\Big[\substack{s \\ s+1}\Big] \times k}=
  \begin{cases}
  2^{k+s-1}  & \text{if  } j = 0,\; k >s, \\
  11\cdot2^{k+s-1}& \text{if  } j = 1,\; k > s +1, \\
21 \cdot 2^{k+s+3j-5}   &  \text{if  }    2\leq j\leq s ,\; k>s+j, \\
3\cdot2^{2k+2s-1}- 3\cdot2^{k+4s -2}       & \text{if   } j  = s+1 ,  \; k>2s +1.
\end{cases}
\end{equation*}
 In the case $ m\geq 2  $\vspace{0.1 cm}\\
 \begin{equation*}
(9)\quad \Gamma _{s+j}^{\Big[\substack{s \\ s+m}\Big] \times (k+1)}
  - \Gamma _{s+j}^{\Big[\substack{s \\ s+m}\Big] \times k}=
  \begin{cases}
  2^{k+s-1}  & \text{if  } j = 0,\; k >s, \\
  3\cdot2^{k+s +2j -3}& \text{if  } 1\leq j\leq m-1,\; k > s+j,  \\
11\cdot2^{k+ s + 2m -3}              &  \text{if  }  j = m, \; k > s+m,
\end{cases}
\end{equation*}
\begin{equation*}
(10)\quad  \Gamma _{s+m+j}^{\Big[\substack{s \\ s+m}\Big] \times (k+1)}
  - \Gamma _{s+m+j}^{\Big[\substack{s \\ s+m}\Big] \times k}=
  \begin{cases}
  21\cdot2^{k+s +2m +3j -4}& \text{if  } 1\leq j\leq s-1,\; k > s+m+ j,  \\
3\cdot2^{2k+2s+m-2} -  3\cdot2^{k+4s+2m-4}            &  \text{if  }  j = s, \; k > 2s+m.
\end{cases}
\end{equation*}\\

\item    Consider the following array $ \left(\Gamma_{s+j}^{\left[s\atop s+m\right]\times k}\right)_{0\leq j\leq s+m,\; k\geq s+j}: $
  \[ 
     \begin{array} {cccccc}
     \Gamma _{s}^{\Big[\substack{s \\ s+m}\Big] \times s}      & \Gamma _{s}^{\Big[\substack{s \\ s+m}\Big] \times (s+1)} & \Gamma _{s}^{\Big[\substack{s \\ s+m}\Big] \times (s+2)} &  \ldots & \Gamma _{s}^{\Big[\substack{s \\ s+m}\Big] \times k} &  \ldots \quad \text{$ k\geq s $} \\
 \Gamma _{s+1}^{\Big[\substack{s \\ s+m}\Big] \times (s+1)}      & \Gamma _{s+1}^{\Big[\substack{s \\ s+m}\Big] \times (s+2)} & \Gamma _{s +1}^{\Big[\substack{s \\ s+m}\Big] \times (s+3)} &  \ldots & \Gamma _{s +1}^{\Big[\substack{s \\ s+m}\Big] \times k} &  \ldots \quad \text{$ k\geq s +1 $} \\
\vdots & \vdots & \vdots    &  \vdots & \vdots  &  \vdots \\
 \Gamma _{s+j}^{\Big[\substack{s \\ s+m}\Big] \times (s+j)}      & \Gamma _{s+j}^{\Big[\substack{s \\ s+m}\Big] \times (s+j+1)} & \Gamma _{s +j}^{\Big[\substack{s \\ s+m}\Big] \times (s+ j+2)} &  \ldots & \Gamma _{s + j}^{\Big[\substack{s \\ s+m}\Big] \times k} &  \ldots \quad \text{$ k\geq s +j $} \\ 
\vdots & \vdots & \vdots    &  \vdots & \vdots  &  \vdots \\
 \Gamma _{s+m}^{\Big[\substack{s \\ s+m}\Big] \times (s+m)}      & \Gamma _{s+m}^{\Big[\substack{s \\ s+m}\Big] \times (s+m+1)} & \Gamma _{s +m}^{\Big[\substack{s \\ s+m}\Big] \times (s+ m+2)} &  \ldots & \Gamma _{s + m}^{\Big[\substack{s \\ s+m}\Big] \times k} &  \ldots \quad \text{$ k\geq s +m $} \\ 
  \Gamma _{s+m +1}^{\Big[\substack{s \\ s+m}\Big] \times (s+m+1)}      & \Gamma _{s+m +1}^{\Big[\substack{s \\ s+m}\Big] \times (s+m+2)} & \Gamma _{s +m +1}^{\Big[\substack{s \\ s+m}\Big] \times (s+ m+3)} &  \ldots & \Gamma _{s + m+1}^{\Big[\substack{s \\ s+m}\Big] \times k} &  \ldots \quad \text{$ k\geq s +m+1 $} \\ 
\vdots & \vdots & \vdots    &  \vdots & \vdots  &  \vdots \\
  \Gamma _{s+m +1+j}^{\Big[\substack{s \\ s+m}\Big] \times (s+m+ 1+j)}      & \Gamma _{s+m +1+j}^{\Big[\substack{s \\ s+m}\Big] \times (s+m+ 2+j)} & \Gamma _{s +m +1+j}^{\Big[\substack{s \\ s+m}\Big] \times (s+ m+3+j)} &  \ldots & \Gamma _{s + m+1+j}^{\Big[\substack{s \\ s+m}\Big] \times k} &  \ldots \quad \text{$ k\geq s +m+1+j $} \\  
   \Gamma _{2s+m -1}^{\Big[\substack{s \\ s+m}\Big] \times (2s+m-1)}      & \Gamma _{2s+m -1}^{\Big[\substack{s \\ s+m}\Big] \times (2s+m)} & \Gamma _{2s +m -1}^{\Big[\substack{s \\ s+m}\Big] \times (2s+ m+1)} &  \ldots & \Gamma _{2s + m -1}^{\Big[\substack{s \\ s+m}\Big] \times k} &  \ldots \quad \text{$ k\geq 2s+m -1 $} \\ 
      \Gamma _{2s+m }^{\Big[\substack{s \\ s+m}\Big] \times (2s+m)}      & \Gamma _{2s+m }^{\Big[\substack{s \\ s+m}\Big] \times (2s+m+1)} & \Gamma _{2s +m }^{\Big[\substack{s \\ s+m}\Big] \times (2s+ m+2)} &  \ldots & \Gamma _{2s + m }^{\Big[\substack{s \\ s+m}\Big] \times k} &  \ldots  \quad \text{$ k\geq 2s+m  $}.
\end{array}
   \] \\
   
   \item To compute the elements in the array  $ \left(\Gamma_{s+j}^{\left[s\atop s+m\right]\times k}\right)_{0\leq j\leq s+m,\; k\geq s+j}, $
   we proceed by induction on j as follows:\vspace{0.1 cm}\\
   Let $l$ be a rational integer such that $  0\leq l\leq s+m-1 $.\vspace{0.1 cm}\\
   Assume that we  have computed the elements in the following subarray
    $ \left(\Gamma_{s+j}^{\left[s\atop s+m\right]\times k}\right)_{0\leq j\leq l,\; k\geq s+j}. $\vspace{0.1 cm}\\
Recall that from  (7) $ \Gamma _{i }^{\Big[\substack{s \\ s+m}\Big] \times k}$ are known for $ 0\leq i\leq s-1,\;k\geq i . $\vspace{0.1 cm}\\
     From the first equation in (2) with k = s+l+1,  $ \; \sum_{ i = 0}^{s+l+1} \Gamma _{i}^{\Big[\substack{s  \\ s+m }\Big] \times (s+l+1)}  = 2^{4s+2l+m} $ \vspace{0.1 cm}\\
  we deduce  $  \Gamma _{s+l+1}^{\Big[\substack{s \\ s+m}\Big] \times (s+l+1)} $
   since the terms $ \Gamma _{i}^{\Big[\substack{s  \\ s+m }\Big] \times (s+l+1)} $ for $ i\leq s +l $ are known.  \vspace{0.05 cm}\\
   From the  equations in (2) with $ k = s+ l+2 $ we have
   \begin{equation*}
  \begin{cases} 
 \sum_{ i = 0}^{s+ l+2} \Gamma _{i}^{\Big[\substack{s  \\ s+m }\Big] \times (s+l+2)}  = 2^{4s+ 2l+m+2}, \\
  \sum_{i = 0}^{s+ l+2} \Gamma _{i}^{\Big[\substack{s  \\ s+m }\Big] \times (s+l+2)}\cdot2^{-i}  =  2^{3s+m +l} + 2^{2s +2l +2} - 2^{s +l}. 
\end{cases}
    \end{equation*}\vspace{0.05 cm}\\ 
  We then  deduce  $  \Gamma _{s+l+1}^{\Big[\substack{s \\ s+m}\Big] \times (s+1+2)} $ and
    $  \Gamma _{s+l+2}^{\Big[\substack{s \\ s+m}\Big] \times (s+ l+2)} $ since the terms $ \Gamma _{i}^{\Big[\substack{s  \\ s+m }\Big] \times (s+l+2)} $ for $ i\leq s +l $ are known.  \vspace{0.1 cm}\\
  Then from (9) or (10)  we compute $\Gamma _{s+l+1}^{\Big[\substack{s \\ s+m}\Big] \times k}$ knowing $  \Gamma _{s+l+1}^{\Big[\substack{s \\ s+m}\Big] \times (s+1+2)} $ for all $k\geq s+l+2$.\vspace{0.1 cm}\\ 
   We have now computed any element in the $l+1$-th row and the first element in the $l+2$-th row  in the array  $ \left(\Gamma_{s+j}^{\left[s\atop s+m\right]\times k}\right)_{0\leq j\leq s+m,\; k\geq s+j}. $\vspace{0.05 cm}\\
   
   \item Using the procedure above we prove successively by induction on j the following reduction formulas \vspace{0.1 cm}\\
    $ \Gamma _{s+j}^{\Big[\substack{s \\ s+m}\Big] \times k} = 8^{j-1}\cdot \Gamma _{s+1}^{\Big[\substack{s \\ s+ (m-(j-1))}\Big] \times (k-(j-1))}\quad for\quad 1\leq j\leq m,\;k\geq s+j, $\vspace{0.1 cm}\\ 
    $ \Gamma _{s+m+1+j}^{\Big[\substack{s \\ s+m}\Big] \times k} = 8^{2j+m}\cdot \Gamma _{s -j+1}^{\Big[\substack{s -j\\ s -j}\Big] \times (k -m -2j)}\quad for\quad 0\leq j\leq s-1 ,\;k\geq s+m+1+j $.\vspace{0.1 cm}\\
    \end{itemize}

 \subsection{\textbf{ORGANIZATION OF THE PAPER}}
 \label{subsec 2.2}
 We proceed as follows:\\
 
  \textbf{ In Section 1} are  introduced main notations and definitions. \\
  
  \textbf{ In Section 3,} we state the main theorems  in the  two  dimensional $ \mathbb{K} $-vector space.\\
  
   \textbf{ In Section 4,}    we establish  results  on  exponential sums, (in $ \mathbb{P}\times \mathbb{P}$) of the form $ g(t,\eta ), h(t,\eta ) $ or similar sums and we show
    that they only depend on rank properties of some corresponding double  persymmetric matrices with entries in $ \mathbb{F}_{2}.$ \\
Consider in particular the quadratic  exponential sum in $\mathbb{P}\times\mathbb{P}$ defined by
$$ (t,\eta ) \in  \mathbb{P}\times \mathbb{P}\longmapsto  
   g(t,\eta ) =  \sum_{deg Y\leq k-1}\sum_{deg Z \leq  s-1}E(tYZ)\sum_{deg U \leq s+m-1}E(\eta YU) \in \mathbb{Z}.  $$\vspace{0.1 cm}\\
  We associate to the exponential sum $ g(t,\eta ) $ the following double persymmetric matrix  
$ \left[{D_{s  \times k}(t)\over D_{(s+m )\times k}(\eta )}\right], $ 
 $$   \left ( \begin{array} {cccccc}
\alpha _{1} & \alpha _{2} & \alpha _{3} &  \ldots & \alpha _{k-1}  &  \alpha _{k} \\
\alpha _{2 } & \alpha _{3} & \alpha _{4}&  \ldots  &  \alpha _{k} &  \alpha _{k+1} \\
\vdots & \vdots & \vdots    &  \vdots & \vdots  &  \vdots \\
\alpha _{s-1} & \alpha _{s} & \alpha _{s +1} & \ldots  &  \alpha _{s+k-3} &  \alpha _{s+k-2}  \\
\alpha _{s} & \alpha _{s+1} & \alpha _{s +2} & \ldots  &  \alpha _{s+k-2} &  \alpha _{s+k-1}  \\
\hline \\
\beta  _{1} & \beta  _{2} & \beta  _{3} & \ldots  &  \beta_{k-1} &  \beta _{k}  \\
\beta  _{2} & \beta  _{3} & \beta  _{4} & \ldots  &  \beta_{k} &  \beta _{k+1}  \\
\vdots & \vdots & \vdots    &  \vdots & \vdots  &  \vdots \\
\beta  _{m+1} & \beta  _{m+2} & \beta  _{m+3} & \ldots  &  \beta_{k+m-1} &  \beta _{k+m}  \\
\vdots & \vdots & \vdots    &  \vdots & \vdots  &  \vdots \\
\beta  _{s+m-1} & \beta  _{s+m} & \beta  _{s+m+1} & \ldots  &  \beta_{s+m+k-3} &  \beta _{s+m+k-2}  \\
\beta  _{s+m} & \beta  _{s+m+1} & \beta  _{s+m+2} & \ldots  &  \beta_{s+m+k-2} &  \beta _{s+m+k-1} 
\end{array}  \right). $$ \vspace{0.1 cm}\\
We get \vspace{0.1 cm}\\
\begin{equation}
\label{eq 2.2}
  g(t,\eta ) =  \sum_{deg Y\leq k-1}\sum_{deg Z \leq  s-1}E(tYZ)\sum_{deg U \leq s+m-1}E(\eta YU) =
   2^{2s+m+k- r( D^{\left[\stackrel{s}{s+m}\right] \times k }(t,\eta  ) )}. 
\end{equation}
 The proof of \eqref{eq 2.2} is based on the following two similar identities \vspace{0.1 cm}\\
  \begin{align*}
  \sum_{deg Z \leq  s-1}E(tYZ) & = 2^{s}\quad  if\; and\; only\; if \quad  Y\in \ker D_{s\times k}(t),  \\
   \sum_{deg U \leq  s+m-1}E(\eta YU) &  = 2^{s+m}\quad  if\; and\; only\; if \quad  Y\in \ker D_{(s+m)\times k}(\eta ). 
\end{align*}
Furthermore we establish similar results concerning similar exponential sums.\vspace{0.1 cm}\\
For instance we have \vspace{0.1 cm}\\
 Let $ (t,\eta ) \in\mathbb{P}\times \mathbb{P} $ and set \\
$$  g_{2}(t,\eta ) =  \sum_{deg Y\leq k-1}\sum_{deg Z \leq  s-1}E(tYZ)\sum_{deg U = s-1+m}E(\eta YU). $$ \\
 Then
 \begin{equation}
 \label{eq 2.3}
g_{2}(t,\eta )  =   \begin{cases}
 2^{2s +m +k-1 - r( D^{\left[\stackrel{s}{s+m}\right] \times k }(t,\eta  ) ) }  & \text{if }
 r( D^{\left[\stackrel{s}{s+m}\right] \times k }(t,\eta  ) )    =  r( D^{\left[\stackrel{s}{s+(m-1)}\right] \times k }(t,\eta  ) ), \\
  0  & \text{otherwise }.
    \end{cases}
\end{equation}

In the end of the section we obtain by \eqref{eq 2.2}, observing that $ g(t,\eta ) $ is constant on cosets of $ \mathbb{P}_{k+s-1}\times \mathbb{P}_{k+s+m-1} $ \\ 
 \begin{align}
  R(q,k-1,s-1,s+m-1,0,0) & =  \int_{\mathbb{P}\times \mathbb{P}} g^{q}(t,\eta )dtd\eta   \label{eq 2.4} \\
&   = 2^{(2s+m+k)(q-1)}\cdot 2^{-k+2}\cdot \sum_{i = 0}^{\inf(2s+m,k)} \Gamma _{i}^{\Big[\substack{s  \\ s+m }\Big] \times k}\cdot2^{- qi}.  \nonumber
\end{align}\vspace{0.1 cm}\\

  \textbf{ In Section 5,} we establish  a recurrent formula for the number $ \Gamma_{i}^{\left[s\atop s+m\right]\times k} $ of rank i  matrices of the form $ \left[{A\over B}\right], $ 
 where A is  a  $ s \times k $ persymmetric matrix and B a  $ (s+m) \times k $ persymmetric matrix with entries in $ \mathbb{F}_{2}. $ \vspace{0.2 cm}\\
 Consider the  following two partitions of the matrix  $  D^{\left[\stackrel{s}{s+m}\right] \times k }(t,\eta  ) $ \vspace{0.1 cm}\\
 
$ \left( \begin{smallmatrix}
 \alpha _{1} & \alpha _{2} & \alpha _{3} &  \ldots & \alpha _{k-1}  &  \alpha _{k} \\
\alpha _{2 } & \alpha _{3} & \alpha _{4}&  \ldots  &  \alpha _{k} &  \alpha _{k+1} \\
\vdots & \vdots & \vdots    &  \vdots & \vdots  &  \vdots \\
\alpha _{s-1} & \alpha _{s} & \alpha _{s +1} & \ldots  &  \alpha _{s+k-3} &  \alpha _{s+k-2}  \\
\beta  _{1} & \beta  _{2} & \beta  _{3} & \ldots  &  \beta_{k-1} &  \beta _{k}  \\
\beta  _{2} & \beta  _{3} & \beta  _{4} & \ldots  &  \beta_{k} &  \beta _{k+1}  \\
\vdots & \vdots & \vdots    &  \vdots & \vdots  &  \vdots \\
\beta  _{m+1} & \beta  _{m+2} & \beta  _{m+3} & \ldots  &  \beta_{k+m-1} &  \beta _{k+m}  \\
\vdots & \vdots & \vdots    &  \vdots & \vdots  &  \vdots \\
\beta  _{s+m-1} & \beta  _{s+m} & \beta  _{s+m+1} & \ldots  &  \beta_{s+m+k-3} &  \beta _{s+m+k-2}  \\
\hline
 \alpha _{s} & \alpha _{s+1} & \alpha _{s +2} & \ldots  &  \alpha _{s+k-2} &  \alpha _{s+k-1}  \\
 \beta  _{s+m} & \beta  _{s+m+1} & \beta  _{s+m+2} & \ldots  &  \beta_{s+m+k-2} &  \beta _{s+m+k-1}
  \end{smallmatrix}\right) $  and
 $ \left( \begin{smallmatrix}
 \alpha _{1} & \alpha _{2} & \alpha _{3} &  \ldots & \alpha _{k-1}  &  \alpha _{k} \\
\alpha _{2 } & \alpha _{3} & \alpha _{4}&  \ldots  &  \alpha _{k} &  \alpha _{k+1} \\
\vdots & \vdots & \vdots    &  \vdots & \vdots  &  \vdots \\
\alpha _{s-1} & \alpha _{s} & \alpha _{s +1} & \ldots  &  \alpha _{s+k-3} &  \alpha _{s+k-2}  \\
\beta  _{1} & \beta  _{2} & \beta  _{3} & \ldots  &  \beta_{k-1} &  \beta _{k}  \\
\beta  _{2} & \beta  _{3} & \beta  _{4} & \ldots  &  \beta_{k} &  \beta _{k+1}  \\
\vdots & \vdots & \vdots    &  \vdots & \vdots  &  \vdots \\
\beta  _{m+1} & \beta  _{m+2} & \beta  _{m+3} & \ldots  &  \beta_{k+m-1} &  \beta _{k+m}  \\
\vdots & \vdots & \vdots    &  \vdots & \vdots  &  \vdots \\
\beta  _{s+m-1} & \beta  _{s+m} & \beta  _{s+m+1} & \ldots  &  \beta_{s+m+k-3} &  \beta _{s+m+k-2}  \\
\hline
 \beta  _{s+m} & \beta  _{s+m+1} & \beta  _{s+m+2} & \ldots  &  \beta_{s+m+k-2} &  \beta _{s+m+k-1}\\
  \alpha _{s} & \alpha _{s+1} & \alpha _{s +2} & \ldots  &  \alpha _{s+k-2} &  \alpha _{s+k-1}  
   \end{smallmatrix}\right). $ \vspace{0.2 cm}\\
 Obviously $ \Gamma _{i}^{\Big[\substack{s \\ s+m }\Big] \times k}$ can be represented in the following two ways   \\
  \begin{align}
   \Gamma _{i}^{\Big[\substack{s \\ s+m }\Big] \times k}
& = \sigma _{i,i,i}^{\left[\stackrel{s-1}{\stackrel{s+m-1 }{\overline {\stackrel{\alpha_{s -}}{\beta_{s+m-} }}}}\right] \times k } 
 + \sigma _{i-1,i,i}^{\left[\stackrel{s-1}{\stackrel{s+m-1 }{\overline {\stackrel{\alpha_{s -}}{\beta_{s+m-} }}}}\right] \times k }
  +\sigma _{i-1,i-1,i}^{\left[\stackrel{s-1}{\stackrel{s+m-1 }{\overline {\stackrel{\alpha_{s -}}{\beta_{s+m-} }}}}\right] \times k } 
  + \sigma _{i-2,i-1,i}^{\left[\stackrel{s-1}{\stackrel{s+m-1 }{\overline {\stackrel{\alpha_{s -}}{\beta_{s+m-} }}}}\right] \times k },\label{eq 2.5} \\
  & \nonumber \\
  \Gamma _{i}^{\Big[\substack{s \\ s+m }\Big] \times k}
   &   = \sigma _{i,i,i}^{\left[\stackrel{s-1}{\stackrel{s+m-1 }{\overline {\stackrel{\beta_{s+m -}}{\alpha_{s -} }}}}\right] \times k } 
 + \sigma _{i-1,i,i}^{\left[\stackrel{s-1}{\stackrel{s+m-1 }{\overline {\stackrel{\beta_{s+m -}}{\alpha_{s -} }}}}\right] \times k }
  +\sigma _{i-1,i-1,i}^{\left[\stackrel{s-1}{\stackrel{s+m-1 }{\overline {\stackrel{\beta_{s+m -}}{\alpha_{s -} }}}}\right] \times k } 
  + \sigma _{i-2,i-1,i}^{\left[\stackrel{s-1}{\stackrel{s+m-1 }{\overline {\stackrel{\beta_{s+m -}}{\alpha_{s -} }}}}\right] \times k }. \label{eq 2.6}\\
  & \nonumber
    \end{align}
We establish the recurrent formula by obtaining relations between the terms of the  two expressions of $ \Gamma _{i}^{\Big[\substack{s \\ s+m }\Big] \times k}.$
 By considering the following partition of the matrix $ D^{\left[\stackrel{s-1}{\stackrel{s+m-1}{\alpha_{s -} + \beta_{s+m -} }}\right] \times k }(t,\eta  ), $ \vspace{0.1 cm}\\
 
\footnotesize
  $  \bordermatrix{%
                    &                   &                      &                      &                            &    \cr
   r_{1} &    \alpha _{1} & \alpha _{2} & \alpha _{3} &  \ldots & \alpha _{k-1}  &  \alpha _{k} \cr
 r_{2} &   \alpha _{2 } & \alpha _{3} & \alpha _{4}&  \ldots  &  \alpha _{k} &  \alpha _{k+1} \cr
  \vdots &   \vdots & \vdots & \vdots    &  \vdots & \vdots  &  \vdots  \cr
 r_{s-1}   &  \alpha _{s-1} & \alpha _{s} & \alpha _{s +1} & \ldots  &  \alpha _{s+k-3} &  \alpha _{s+k-2}  \cr
r_{s}   & \beta  _{1} & \beta  _{2} & \beta  _{3} & \ldots  &  \beta_{k-1} &  \beta _{k}  \cr
 r_{s+1} & \beta  _{2} & \beta  _{3} & \beta  _{4} & \ldots  &  \beta_{k} &  \beta _{k+1}  \cr
\vdots  &  \vdots & \vdots & \vdots    &  \vdots & \vdots  &  \vdots \cr
r_{s+m} &      \beta  _{m+1} & \beta  _{m+2} & \beta  _{m+3} & \ldots  &  \beta_{k+m-1} &  \beta _{k+m}  \cr
\vdots  & \vdots & \vdots & \vdots    &  \vdots & \vdots  &  \vdots \cr
   r_{2s+m-2} & \beta  _{s+m-1} & \beta  _{s+m} & \beta  _{s+m+1} & \ldots  &  \beta_{s+m+k-3} &  \beta _{s+m+k-2}  \cr
  \hline
r_{2s+m-1} & \alpha _{s} + \beta  _{s+m} & \alpha _{s+1} + \beta  _{s+m+1} & \alpha _{s+2} + \beta  _{s+m+2} & \ldots 
& \alpha _{s+k-2} + \beta  _{s+k+m-2} &  \alpha _{s+k-1} + \beta  _{s+k+m-1} \cr
} $  \vspace{0.2 cm}\\

\normalsize
 and performing the following row operations, adding the s+m+j-th row to the j-th row for $ 0\leq j\leq s-2, $
we obtain by elementary rank considerations  the following identity 
\begin{equation}
\label{eq 2.7}
 \sigma _{i,i}^{\left[\stackrel{s-1}{\stackrel{s+m -1}{\overline{\alpha_{s -}+ \beta _{s+m-}}}}\right] \times k } =
  2\cdot \sigma _{i,i}^{\left[\stackrel{s-1}{\stackrel{s+m -1}{\overline{\alpha_{s -}}}}\right] \times k }\quad  for \;  0\leq i\leq \inf(2s+m-2,k).
  \end{equation}
  Using the binomial formula, we can write \vspace{0.2 cm}\\

\begin{align}
g_{1}^{q} & = (h + f_{1})^{q} = h^{q} + f_{1}^{q} + \sum_{i = 1}^{q -1} \binom{q}{i}h^{i}\cdot f_{1}^{q-i}, \label{eq 2.8} \\
& \nonumber \\
g_{2}^{q} & = (h + f_{2})^{q} = h^{q} + f_{2}^{q} + \sum_{i = 1}^{q -1} \binom{q}{i}h^{i}\cdot f_{2}^{q-i}. \label{eq 2.9}
\end{align}
 Integrating \eqref{eq 2.8} on the unit interval of $\mathbb{K}^{2}$ we get by using results similar  to \eqref{eq 2.3} \vspace{0.2 cm}\\
 \begin{equation}
  \label{eq 2.10}
\sum_{i=1}^{\inf(2s+m-1, k)}\left[  \sigma _{i-1,i,i}^{\left[\stackrel{s-1}{\stackrel{s+m-1 }{\overline {\stackrel{\beta _{s +m -}}{\alpha _{s -} }}}}\right] \times k } 
-2\cdot \sigma _{i-1,i-1,i}^{\left[\stackrel{s-1}{\stackrel{s+m-1 }{\overline {\stackrel{\alpha_{s -}}{\beta_{s+m-} }}}}\right] \times k }  \right]\cdot2^{-qi} =  0  \; for \;
all \; q\geq 2.
\end{equation} \vspace{0.2 cm}\\ 
From \eqref{eq 2.10} we obtain  \vspace{0.2 cm}\\
\begin{equation}
\label{eq 2.11}
  \sigma _{i-1,i,i}^{\left[\stackrel{s-1}{\stackrel{s+m-1 }{\overline {\stackrel{\beta _{s +m -}}{\alpha _{s -} }}}}\right] \times k } 
= 2\cdot \sigma _{i-1,i-1,i}^{\left[\stackrel{s-1}{\stackrel{s+m-1 }{\overline {\stackrel{\alpha_{s -}}{\beta_{s+m-} }}}}\right] \times k } 
\quad for \;  0\leq i\leq \inf(2s+m-1, k).    
  \end{equation}\\
  
Integrating \eqref{eq 2.9} on the unit interval of $\mathbb{K}^{2}$ we get in the same way \vspace{0.2 cm}\\
\begin{equation}
\label{eq 2.12}
  \sum_{i = 1}^{\inf(2s+m-1,k)}\left[  \sigma _{i-1,i,i}^{\left[\stackrel{s-1}{\stackrel{s+m-1 }{\overline {\stackrel{\beta _{s+m -}}{\alpha _{s -} }}}}\right] \times k }
- \sigma _{i-1,i,i}^{\left[\stackrel{s-1}{\stackrel{s+m-1 }{\overline {\stackrel{\alpha_{s -}}{\beta_{s+m-} }}}}\right] \times k } - \sigma _{i-1,i-1,i}^{\left[\stackrel{s-1}{\stackrel{s+m-1 }{\overline {\stackrel{\alpha_{s -}}{\beta_{s+m-} }}}}\right] \times k }+ 
 \sigma _{i-1,i-1,i}^{\left[\stackrel{s-1}{\stackrel{s+m-1 }{\overline {\stackrel{\beta _{s+m -}}{\alpha _{s -} }}}}\right] \times k }       \right]\cdot 2^{-qi} = 0 \quad for \;  q\geq 2.
  \end{equation}
From \eqref{eq 2.12} we get \vspace{0.2 cm}\\
\begin{equation}
\label{eq 2.13}
 \sigma _{i-1,i,i}^{\left[\stackrel{s-1}{\stackrel{s+m-1 }{\overline {\stackrel{\beta _{s+m -}}{\alpha _{s -} }}}}\right] \times k } +  
   \sigma _{i-1,i-1,i}^{\left[\stackrel{s-1}{\stackrel{s+m-1 }{\overline {\stackrel{\beta _{s+m -}}{\alpha _{s -} }}}}\right] \times k }=
    \sigma _{i-1,i,i}^{\left[\stackrel{s-1}{\stackrel{s+m-1 }{\overline {\stackrel{\alpha_{s -}}{\beta_{s+m-} }}}}\right] \times k } +
     \sigma _{i-1,i-1,i}^{\left[\stackrel{s-1}{\stackrel{s+m-1 }{\overline {\stackrel{\alpha_{s -}}{\beta_{s+m-} }}}}\right] \times k } \quad for \; 1\leq i\leq \inf(2s+m-1,k).
 \end{equation}\vspace{0.2 cm}\\
 Combining \eqref{eq 2.5}, \eqref{eq 2.6}, \eqref{eq 2.7}, \eqref{eq 2.11} and  \eqref{eq 2.13} we obtain
 the following recurrent formula for the number of $ (2s+m)\times k $ double persymmetric matrices of rank i  \\
   \begin{align}
  \Gamma _{i}^{\Big[\substack{s \\ s+m }\Big] \times k} & = 2\cdot \Gamma _{i-1}^{\Big[\substack{s -1\\ s-1+(m+1) }\Big] \times k}
+ 4\cdot \Gamma _{i-1}^{\Big[\substack{s \\ s+(m-1) }\Big] \times k} - 8\cdot \Gamma _{i-2}^{\Big[\substack{s-1 \\ s-1+m }\Big] \times k}
 + \Delta _{i}^{\Big[\substack{s \\ s+m }\Big] \times k}\; for \; 0\leq i\leq \inf(2s+m,k) \label{eq 2.14} \\
 & \text{where $\Delta _{i}^{\Big[\substack{s \\ s+m }\Big] \times k}$ is equal to} \nonumber \\
  & \sigma _{i,i,i}^{\left[\stackrel{s-1}{\stackrel{s+m-1 }{\overline {\stackrel{\alpha_{s -}}{\beta_{s+m-} }}}}\right] \times k }
  - 3\cdot \sigma _{i-1,i-1,i-1}^{\left[\stackrel{s-1}{\stackrel{s+m-1 }{\overline {\stackrel{\alpha_{s -}}{\beta_{s+m-} }}}}\right] \times k } 
  + 2\cdot \sigma _{i-2,i-2,i-2}^{\left[\stackrel{s-1}{\stackrel{s+m-1 }{\overline {\stackrel{\alpha_{s -}}{\beta_{s+m-} }}}}\right] \times k }.\label{eq 2.15} 
  \end{align}\vspace{0.2 cm}\\
We observe that we have for instance \vspace{0.2 cm}\\
\begin{align*}
 2\cdot \Gamma _{i-1}^{\Big[\substack{s -1\\ s-1+(m+1) }\Big] \times k}
  =    \sigma _{i-1,i-1,i-1}^{\left[\stackrel{s-1}{\stackrel{s+m-1 }{\overline {\stackrel{\beta_{s+m -}}{\alpha_{s -} }}}}\right] \times k } 
 + \sigma _{i-1,i-1,i}^{\left[\stackrel{s-1}{\stackrel{s+m-1 }{\overline {\stackrel{\beta_{s+m -}}{\alpha_{s -} }}}}\right] \times k }
  +\sigma _{i-2,i-1,i-1}^{\left[\stackrel{s-1}{\stackrel{s+m-1 }{\overline {\stackrel{\beta_{s+m -}}{\alpha_{s -} }}}}\right] \times k } 
  + \sigma _{i-2,i-1,i}^{\left[\stackrel{s-1}{\stackrel{s+m-1 }{\overline {\stackrel{\beta_{s+m -}}{\alpha_{s -} }}}}\right] \times k }. \\
\end{align*}
 \textbf{ In Section 6,} we study rank properties of a partition of double persymmetric matrices by integrating some
  appropriate  exponential sums on the unit interval  of $\mathbb{K}^2 $ with integral equal to zero.\vspace{0.2 cm}\\
 Consider the following partition of the matrix  $  D^{\left[\stackrel{s-1}{\stackrel{s+m-1 }
{\overline {\stackrel{\alpha  _{s -}}{\beta  _{s +m -} }}}}\right] \times k }(t,\eta ),  $ \vspace{0.2 cm}\\
      $$   \left ( \begin{array} {ccccc|c}
\alpha _{1} & \alpha _{2} & \alpha _{3} &  \ldots & \alpha _{k-1}  &  \alpha _{k} \\
\alpha _{2 } & \alpha _{3} & \alpha _{4}&  \ldots  &  \alpha _{k} &  \alpha _{k+1} \\
\vdots & \vdots & \vdots    &  \vdots & \vdots  &  \vdots \\
\alpha _{s-1} & \alpha _{s} & \alpha _{s +1} & \ldots  &  \alpha _{s+k-3} &  \alpha _{s+k-2}  \\
 \beta  _{1} & \beta  _{2} & \beta  _{3} & \ldots  &  \beta_{k-1} &  \beta _{k}  \\
\beta  _{2} & \beta  _{3} & \beta  _{4} & \ldots  &  \beta_{k} &  \beta _{k+1}  \\
\vdots & \vdots & \vdots    &  \vdots & \vdots  &  \vdots \\
\beta  _{m+1} & \beta  _{m+2} & \beta  _{m+3} & \ldots  &  \beta_{k+m-1} &  \beta _{k+m}  \\
\vdots & \vdots & \vdots    &  \vdots & \vdots  &  \vdots \\
\beta  _{s+m-1} & \beta  _{s+m} & \beta  _{s+m+1} & \ldots  &  \beta_{s+m+k-3} &  \beta _{s+m+k-2}  \\
\hline
\alpha _{s} & \alpha _{s+1} & \alpha _{s +2} & \ldots  &  \alpha _{s+k-2} &  \alpha _{s+k-1} \\
\hline
\beta  _{s+m} & \beta  _{s+m+1} & \beta  _{s+m+2} & \ldots  &  \beta_{s+m+k-2} &  \beta _{s+m+k-1}
\end{array}  \right). $$ \vspace{0.5 cm} 
In fact we prove that for all $ j\in [0,\inf(2s+m-2,k-1)] $ \vspace{0.2 cm} we have \\
\small
    \begin{align}
& {}^{\#}\left(\begin{array}{c | c}
           j & j \\
           \hline
           j & j\\
           \hline
            j & j +1 
           \end{array} \right)_{\mathbb{P}/\mathbb{P}_{k+s -1}\times
           \mathbb{P}/\mathbb{P}_{k+s+m-1} }^{{\alpha  \over \beta }} = 
{}^{\#}\left(\begin{array}{c | c}
           j & j \\
           \hline
           j & j\\
           \hline
            j & j 
           \end{array} \right)_{\mathbb{P}/\mathbb{P}_{k+s -1}\times
           \mathbb{P}/\mathbb{P}_{k+s+m-1} }^{{\alpha  \over \beta }},\label{eq 2.16} \\
             &   {}^{\#}\left(\begin{array}{c | c}
           j & j \\
           \hline
           j & j +1 \\
           \hline
            j & j +1 
           \end{array} \right)_{\mathbb{P}/\mathbb{P}_{k+s -1}\times
           \mathbb{P}/\mathbb{P}_{k+s+m-1} }^{{\alpha  \over \beta }} = 
{}^{\#}\left(\begin{array}{c | c}
           j & j \\
           \hline
           j & j\\
           \hline
            j & j 
           \end{array} \right)_{\mathbb{P}/\mathbb{P}_{k+s -1}\times
           \mathbb{P}/\mathbb{P}_{k+s+m-1} }^{{\alpha  \over \beta }} + 
           {}^{\#}\left(\begin{array}{c | c}
           j & j \\
           \hline
           j & j\\
           \hline
            j & j +1
           \end{array} \right)_{\mathbb{P}/\mathbb{P}_{k+s -1}\times
           \mathbb{P}/\mathbb{P}_{k+s+m-1} }^{{\alpha  \over \beta }}.\label{eq 2.17} 
          \end{align} 
We deduce  \eqref{eq 2.16} from the fact that for all $ q\geq 2 $
the integral   $\displaystyle \int_{\mathbb{P}\times \mathbb{P}}   \phi (t,\eta )\cdot  \theta _{3}^{q-1} (t,\eta )  dtd\eta $ is equal to  zero, that is  \vspace{0.1 cm}\\
  \begin{align*}
  \sum_{j=0}^{\inf(2s+m-2,k-1)}a(s,m,k,q)\cdot2^{-j q}\cdot\Bigg[{}^{\#}\left(\begin{array}{c | c}
           j & j \\
           \hline
           j & j\\
           \hline
            j & j 
           \end{array} \right)_{\mathbb{P}/\mathbb{P}_{k+s -1}\times
           \mathbb{P}/\mathbb{P}_{k+s+m-1} }^{{\alpha  \over \beta }}
            -   
{}^{\#}\left(\begin{array}{c | c}
           j & j \\
           \hline
           j & j\\
           \hline
            j & j +1
           \end{array} \right)_{\mathbb{P}/\mathbb{P}_{k+s -1}\times
           \mathbb{P}/\mathbb{P}_{k+s+m-1} }^{{\alpha  \over \beta }}\Bigg] = 0. 
           \end{align*}\vspace{0.2 cm}\\
We deduce  \eqref{eq 2.17}  in a similar way from the fact that for all $ q\geq 2 $ the integral \\
 $ \displaystyle \int_{\mathbb{P}\times \mathbb{P}}\theta _{1} (t,\eta )\cdot  \theta _{2}^{q-1} (t,\eta )  dtd\eta $
is equal to zero, that is \vspace{0.2 cm}\\
  \begin{align*}
  \sum_{j=0}^{\inf(2s+m-2,k-1)}b(s,m,k,q)\cdot2^{-j q}\cdot\Bigg[{}^{\#}\left(\begin{array}{c | c}
           j & j \\
           \hline
           j & j\\
           \hline
            j & j 
           \end{array} \right)_{\mathbb{P}/\mathbb{P}_{k+s -1}\times
           \mathbb{P}/\mathbb{P}_{k+s+m-1} }^{{\alpha  \over \beta }}
            +    
{}^{\#}\left(\begin{array}{c | c}
           j & j \\
           \hline
           j & j\\
           \hline
            j & j +1
           \end{array} \right)_{\mathbb{P}/\mathbb{P}_{k+s -1}\times
           \mathbb{P}/\mathbb{P}_{k+s+m-1} }^{{\alpha  \over \beta }} \\
           - {}^{\#}\left(\begin{array}{c | c}
           j & j \\
           \hline
           j & j+1 \\
           \hline
            j & j +1
           \end{array} \right)_{\mathbb{P}/\mathbb{P}_{k+s -1}\times
           \mathbb{P}/\mathbb{P}_{k+s+m-1} }^{{\alpha  \over \beta }}\Bigg] = 0.
            \end{align*}\vspace{0.1 cm}\\

     \textbf{ In Section 7,} we study some  rank properties of submatrices  of double persymmetric matrices. Consider the following partition of the matrix $  D^{\left[\stackrel{s-1}{\stackrel{s+m-1 }
{\overline {\stackrel{\alpha  _{s -}}{\beta  _{s +m -} }}}}\right] \times k }(t,\eta ),  $ \\

         $$   \left ( \begin{array} {ccccc|c}
\alpha _{1} & \alpha _{2} & \alpha _{3} &  \ldots & \alpha _{k-1}  &  \alpha _{k} \\
\alpha _{2 } & \alpha _{3} & \alpha _{4}&  \ldots  &  \alpha _{k} &  \alpha _{k+1} \\
\vdots & \vdots & \vdots    &  \vdots & \vdots  &  \vdots \\
\alpha _{s-1} & \alpha _{s} & \alpha _{s +1} & \ldots  &  \alpha _{s+k-3} &  \alpha _{s+k-2}  \\
 \beta  _{1} & \beta  _{2} & \beta  _{3} & \ldots  &  \beta_{k-1} &  \beta _{k}  \\
\beta  _{2} & \beta  _{3} & \beta  _{4} & \ldots  &  \beta_{k} &  \beta _{k+1}  \\
\vdots & \vdots & \vdots    &  \vdots & \vdots  &  \vdots \\
\beta  _{m+1} & \beta  _{m+2} & \beta  _{m+3} & \ldots  &  \beta_{k+m-1} &  \beta _{k+m}  \\
\vdots & \vdots & \vdots    &  \vdots & \vdots  &  \vdots \\
\beta  _{s+m-1} & \beta  _{s+m} & \beta  _{s+m+1} & \ldots  &  \beta_{s+m+k-3} &  \beta _{s+m+k-2}  \\
\hline
\alpha _{s} & \alpha _{s+1} & \alpha _{s +2} & \ldots  &  \alpha _{s+k-2} &  \alpha _{s+k-1} \\
\hline
\beta  _{s+m} & \beta  _{s+m+1} & \beta  _{s+m+2} & \ldots  &  \beta_{s+m+k-2} &  \beta _{s+m+k-1}
\end{array}  \right). $$ \vspace{0.5 cm} \\

In fact we prove that for all $ j\in [0,\inf(2s+m-3,k-2)] $ \vspace{0.2 cm} we have \\
  \begin{align}
 {}^{\#}\left(\begin{array}{c | c}
           j & j +1 \\
           \hline
           j & j +1 \\
           \hline
            j & j +1
           \end{array} \right)_{\mathbb{P}/\mathbb{P}_{k+s -1}\times
           \mathbb{P}/\mathbb{P}_{k+s+m-1} }^{{\alpha  \over \beta }} = 0. \label{eq 2.18}
\end{align} \\
We prove \eqref{eq 2.18} by contradiction.\vspace{0.1 cm}\\
 Assume on the contrary that there exist $ j_{0} \in [0,\inf(2s+m-3,k-2)] $  such that \\  
   \begin{align*}
 {}^{\#}\left(\begin{array}{c | c}
           j_{0}  & j_{0}  +1 \\
           \hline
          j_{0}   &  j_{0} +1 \\
           \hline
            j_{0}  &  j_{0} +1
           \end{array} \right)_{\mathbb{P}/\mathbb{P}_{k+s -1}\times
           \mathbb{P}/\mathbb{P}_{k+s+m-1} }^{{\alpha  \over \beta }} > 0. 
\end{align*} \\
   We  show that \\
 \begin{align*}
& {}^{\#}\left(\begin{array}{c | c}
           j_{0} & j_{0}+1 \\
           \hline
           j_{0} & j_{0}+1\\
           \hline
            j_{0} & j_{0}+1 
           \end{array} \right)_{\mathbb{P}/\mathbb{P}_{k+s -1}\times
           \mathbb{P}/\mathbb{P}_{k+s+m-1} }^{{\alpha  \over \beta }} > 0
 \Longrightarrow 
 {}^{\#}\left(\begin{array}{c | c}
           j_{0}-1 & j_{0} \\
           \hline
           j_{0}-1 & j_{0}\\
           \hline
            j_{0}-1 & j_{0}
           \end{array} \right)_{\mathbb{P}_{1}/\mathbb{P}_{k+s -1}\times
           \mathbb{P}_{1}/\mathbb{P}_{k+s+m-1} }^{{\alpha  \over \beta }} > 0   \\
&  \Longrightarrow  \ldots  \Longrightarrow
  {}^{\#}\left(\begin{array}{c | c}
           0 & 1 \\
           \hline
           0 & 1 \\
           \hline
            0 & 1 
           \end{array} \right)_{\mathbb{P}_{j_{0}}/\mathbb{P}_{k+s -1}\times
           \mathbb{P}_{j_{0}}/\mathbb{P}_{k+s+m-1} }^{{\alpha  \over \beta }} > 0
\end{align*}
  which obviously contradicts
 \begin{equation*}
{}^{\#}\left(\begin{array}{c | c}
           0 & 1 \\
           \hline
           0 & 1 \\
           \hline
            0 & 1 
           \end{array} \right)_{\mathbb{P}_{j_{0}}/\mathbb{P}_{k+s -1}\times
           \mathbb{P}_{j_{0}}/\mathbb{P}_{k+s+m-1} }^{{\alpha  \over \beta }} = 0.
 \end{equation*} \\
  These inequalities are a consequence of rank properties of submatrices of double persymmetric matrices.\vspace{0.2 cm}\\
  
    \textbf{ In Section 8,} we  study the remainder  $ \Delta _{j}^{\Big[{s \atop s+m}\Big] \times k}$  in the recurrent formula  and we prove in particular
    that \vspace{0.2 cm}\\
 \begin{equation}
\label{eq 2.19}
  \Delta _{j}^{\Big[{s \atop s+m}\Big] \times (k+1)} -\Delta _{j}^{\Big[{s \atop s+m}\Big] \times k}= 0 \text{ if $ k\geq j+1$}. 
\end{equation} \\

 Consider still  the following partition of the matrix $ D^{\left[\stackrel{s-1}{\stackrel{s+m-1 }
{\overline {\stackrel{\alpha  _{s -}}{\beta  _{s +m -} }}}}\right] \times k }(t,\eta ) $
   $$   \left ( \begin{array} {ccccc|c}
\alpha _{1} & \alpha _{2} & \alpha _{3} &  \ldots & \alpha _{k-1}  &  \alpha _{k} \\
\alpha _{2 } & \alpha _{3} & \alpha _{4}&  \ldots  &  \alpha _{k} &  \alpha _{k+1} \\
\vdots & \vdots & \vdots    &  \vdots & \vdots  &  \vdots \\
\alpha _{s-1} & \alpha _{s} & \alpha _{s +1} & \ldots  &  \alpha _{s+k-3} &  \alpha _{s+k-2}  \\
\beta  _{1} & \beta  _{2} & \beta  _{3} & \ldots  &  \beta_{k-1} &  \beta _{k}  \\
\beta  _{2} & \beta  _{3} & \beta  _{4} & \ldots  &  \beta_{k} &  \beta _{k+1}  \\
\vdots & \vdots & \vdots    &  \vdots & \vdots  &  \vdots \\
\beta  _{m+1} & \beta  _{m+2} & \beta  _{m+3} & \ldots  &  \beta_{k+m-1} &  \beta _{k+m}  \\
\vdots & \vdots & \vdots    &  \vdots & \vdots  &  \vdots \\
\beta  _{s+m-1} & \beta  _{s+m} & \beta  _{s+m+1} & \ldots  &  \beta_{s+m+k-3} &  \beta _{s+m+k-2}  \\
\hline
\alpha _{s} & \alpha _{s+1} & \alpha _{s +2} & \ldots  &  \alpha _{s+k-2} &  \alpha _{s+k-1}\\
\hline
\beta  _{s+m} & \beta  _{s+m+1} & \beta  _{s+m+2} & \ldots  &  \beta_{s+m+k-2} &  \beta _{s+m+k-1}
  \end{array}  \right). $$ \vspace{0.5 cm}\\
  
Let $j\leq \inf(k-2,2s+m - 3). $ Then by elementary rank considerations and using \eqref{eq 2.18} we obtain \\
\begin{align}
\label{eq 2.20}
  \sigma _{j,j,j}^{\left[\stackrel{s-1}{\stackrel{s+m-1 }
{\overline {\stackrel{\alpha_{s -}}{\beta_{s+m-} }}}}\right] \times k } 
&  = 
 {}^{\#}\left(\begin{array}{c | c}
           j & j \\
           \hline
           j & j \\
           \hline
            j & j 
           \end{array} \right)_{\mathbb{P}/\mathbb{P}_{k+s -1}\times
           \mathbb{P}/\mathbb{P}_{k+s+m-1} }^{{\alpha \over \beta }} +
        {}^{\#}\left(\begin{array}{c | c}
             j -1 & j \\
           \hline
           j & j \\
           \hline
            j & j 
           \end{array} \right)_{\mathbb{P}/\mathbb{P}_{k+s -1}\times
           \mathbb{P}/\mathbb{P}_{k+s+m-1} }^{{\alpha  \over \beta }}\\
           &  +
           {}^{\#}\left(\begin{array}{c | c}
           j-1 & j \\
           \hline
           j-1 & j \\
           \hline
            j & j 
           \end{array} \right)_{\mathbb{P}/\mathbb{P}_{k+s -1}\times
           \mathbb{P}/\mathbb{P}_{k+s+m-1} }^{{\alpha  \over \beta }}  + 
             {}^{\#}\left(\begin{array}{c | c}
            j -1 & j \\
           \hline
           j-1 & j \\
           \hline
            j-1 & j 
           \end{array} \right)_{\mathbb{P}/\mathbb{P}_{k+s -1}\times
           \mathbb{P}/\mathbb{P}_{k+s+m-1} }^{{\alpha \over \beta }}  \nonumber  \\
           & =  {}^{\#}\left(\begin{array}{c | c}
           j & j \\
           \hline
           j & j \\
           \hline
            j & j 
           \end{array} \right)_{\mathbb{P}/\mathbb{P}_{k+s -1}\times
           \mathbb{P}/\mathbb{P}_{k+s+m-1} }^{{\alpha \over \beta }}. \nonumber
            \end{align}\\
            
      Further  by  \eqref{eq 2.16}, \eqref{eq 2.17} and    \eqref{eq 2.18}  we have  \\
         
         \begin{align}
\label{eq 2.21}
4\cdot \sigma _{j,j,j}^{\left[\stackrel{s-1}{\stackrel{s+m-1 }
{\overline {\stackrel{\alpha_{s -}}{\beta_{s+m-} }}}}\right] \times (k-1) } 
&  = 
 {}^{\#}\left(\begin{array}{c | c}
           j & j \\
           \hline
           j & j \\
           \hline
            j & j 
           \end{array} \right)_{\mathbb{P}/\mathbb{P}_{k+s -1}\times
           \mathbb{P}/\mathbb{P}_{k+s+m-1} }^{{\alpha  \over \beta  }} +
        {}^{\#}\left(\begin{array}{c | c}
             j  & j \\
           \hline
           j & j+1 \\
           \hline
            j & j +1
           \end{array} \right)_{\mathbb{P}/\mathbb{P}_{k+s -1}\times
           \mathbb{P}/\mathbb{P}_{k+s+m-1} }^{{\alpha \over \beta  }}\\
           &  +
           {}^{\#}\left(\begin{array}{c | c}
           j & j \\
           \hline
           j & j \\
           \hline
            j & j +1
           \end{array} \right)_{\mathbb{P}/\mathbb{P}_{k+s -1}\times
           \mathbb{P}/\mathbb{P}_{k+s+m-1} }^{{\alpha \over \beta }}  + 
             {}^{\#}\left(\begin{array}{c | c}
            j  & j +1\\
           \hline
           j & j +1 \\
           \hline
            j & j +1
           \end{array} \right)_{\mathbb{P}/\mathbb{P}_{k+s -1}\times
           \mathbb{P}/\mathbb{P}_{k+s+m-1} }^{{\alpha \over \beta }}  \nonumber  \\
           & = 4\cdot {}^{\#}\left(\begin{array}{c | c}
           j & j \\
           \hline
           j & j \\
           \hline
            j & j 
           \end{array} \right)_{\mathbb{P}/\mathbb{P}_{k+s -1}\times
           \mathbb{P}/\mathbb{P}_{k+s+m-1} }^{{\alpha  \over \beta }}. \nonumber
            \end{align}
      From  \eqref{eq 2.20},  \eqref{eq 2.21} we deduce \vspace{0.2 cm} \\
      \begin{equation}
      \label{eq 2.22}
     \sigma _{j,j,j}^{\left[\stackrel{s-1}{\stackrel{s+m-1 }
{\overline {\stackrel{\alpha_{s -}}{\beta_{s+m-} }}}}\right] \times k } = \sigma _{j,j,j}^{\left[\stackrel{s-1}{\stackrel{s+m-1 }
{\overline {\stackrel{\alpha_{s -}}{\beta_{s+m-} }}}}\right] \times (j+1) } \text{if  }\; j\leq \inf(k-2 , 2s+m-3).
      \end{equation} \vspace{0.2 cm} \\
  Using \eqref{eq 2.16}, \eqref{eq 2.17} we obtain easily  \vspace{0.2 cm} \\
  \begin{align}
  \sigma _{j,j,j}^{\left[\stackrel{s-1}{\stackrel{s+m-1 }
{\overline {\stackrel{\alpha_{s -}}{\beta_{s+m-} }}}}\right] \times (j+1) } 
&  = 
 {}^{\#}\left(\begin{array}{c | c}
           j & j \\
           \hline
           j & j \\
            \hline
            j & j 
           \end{array} \right)_{\mathbb{P}/\mathbb{P}_{j+s }\times
           \mathbb{P}/\mathbb{P}_{j+s+m} }^{{\alpha \over \beta }},\label{eq 2.23} \\
           & \nonumber \\
    4\cdot \sigma _{j,j,j}^{\left[\stackrel{s-1}{\stackrel{s+m-1 }
{\overline {\stackrel{\alpha_{s -}}{\beta_{s+m-} }}}}\right] \times j } 
&  =  4 \cdot 4\cdot \Gamma _{j}^{\Big[\substack{s -1\\ s-1 +m }\Big] \times j}\label{eq 2.24} \\
  & =   4\cdot {}^{\#}\left(\begin{array}{c | c}
           j & j \\
           \hline
           j & j \\
           \hline
            j & j
           \end{array} \right)_{\mathbb{P}/\mathbb{P}_{j+s }\times
           \mathbb{P}/\mathbb{P}_{j+s+m} }^{{\alpha  \over \beta }} +  
              4 \cdot \Gamma _{j+1}^{\Big[\substack{s -1\\ s-1 +m }\Big] \times (j+1)}.\nonumber        
            \end{align}
From \eqref{eq 2.23}, \eqref{eq 2.24} we deduce Theorem \ref{thm 3.3}.\\

 Recalling that $ \Delta _{i}^{\Big[\substack{s \\ s+m }\Big] \times k} $, is equal to
$$ \sigma _{i,i,i}^{\left[\stackrel{s-1}{\stackrel{s+m-1 }{\overline {\stackrel{\alpha_{s -}}{\beta_{s+m-} }}}}\right] \times k }
  - 3\cdot \sigma _{i-1,i-1,i-1}^{\left[\stackrel{s-1}{\stackrel{s+m-1 }{\overline {\stackrel{\alpha_{s -}}{\beta_{s+m-} }}}}\right] \times k } 
  + 2\cdot \sigma _{i-2,i-2,i-2}^{\left[\stackrel{s-1}{\stackrel{s+m-1 }{\overline {\stackrel{\alpha_{s -}}{\beta_{s+m-} }}}}\right] \times k }  $$  we deduce easily Theorem \ref{thm 3.4}.\vspace{0.2 cm}\\
We get  \eqref{eq 2.19} from Theorem \ref{thm 3.5}.\vspace{0.2 cm}\\
   
    \textbf{ In Section 9,} we establish  a recurrent formula for the difference \vspace{0.2 cm}\\
     $$ \Gamma _{i}^{\Big[\substack{s \\ s +m }\Big] \times (k+1)} -
    \Gamma _{i}^{\Big[\substack{s \\ s +m }\Big] \times k} \quad for\; 2\leq i \leq 2s+m,\quad k\geq i+1.$$
    In fact we prove that \vspace{0.2 cm}\\
       \begin{align}
  &  \Gamma _{i}^{\Big[\substack{s \\ s +m }\Big] \times (k+1)} -  \Gamma _{i}^{\Big[\substack{s \\ s +m }\Big] \times k} =  4\cdot \left[\Gamma _{i-1}^{\Big[\substack{s \\ s +(m-1) }\Big] \times (k+1)} - \Gamma _{i-1}^{\Big[\substack{s \\ s +(m-1) }\Big] \times  k}\right]  
 \quad \text{if    }\quad  2 \leq i\leq s-1,\quad k\geq i+1, \label{eq 2.25} \\
 & \nonumber \\
  &  \Gamma _{s+j}^{\Big[\substack{s \\ s +m }\Big] \times (k+1)} -  \Gamma _{s+j}^{\Big[\substack{s \\ s +m }\Big] \times k} =  4\cdot \left[\Gamma _{s+j-1}^{\Big[\substack{s \\ s +(m-1) }\Big] \times (k+1)} - \Gamma _{s+j-1}^{\Big[\substack{s \\ s +(m-1) }\Big] \times  k}\right]  
  +   2^{s-1}\left[ \Gamma _{j+1}^{\Big[\substack{1 \\ 1+ (s+m-1)}\Big] \times (k+1)} -  \Gamma _{j+1}^{\Big[\substack{1 \\ 1+(s+m-1) }\Big] \times k}\right]\label{eq 2.26} \\
&   -   2^{s+1}\left[ \Gamma _{j}^{\Big[\substack{1 \\ 1+ (s+m-2)}\Big] \times (k+1)} -  \Gamma _{j}^{\Big[\substack{1 \\ 1+(s+m-2) }\Big] \times k}\right] 
     \text{if  }  0 \leq j\leq s+m ,\quad k\geq s+j+1.  \nonumber  
       \end{align}\\
       
 To prove \eqref{eq 2.25}  we proceed as follows   \vspace{0.2 cm}\\   
 Set for $ 1\leq i\leq \inf(2s+m,k) $  $$\Omega _{i}(s,m,k)= \Gamma _{i}^{\Big[\substack{s \\ s +m }\Big] \times k} - 4\cdot \Gamma _{i-1}^{\Big[\substack{s \\ s +(m-1) }\Big] \times k}. $$\\
 
 Then, from the recurrent formula \eqref{eq 2.14} we get  \vspace{0.2 cm}\\  
 \begin{equation}
\label{eq 2.27}
  \Omega _{i}(s,m,k)= 2\cdot \Omega _{i-1}(s-1,m+1,k) +  \Delta _{i}^{\Big[\substack{s \\ s+m }\Big] \times k}. 
\end{equation}\\

Applying succesively  \eqref{eq 2.27} we obtain  \vspace{0.3 cm}\\         
   \begin{align}
\displaystyle   \sum_{j = 0}^{i-2}2^{j}\cdot \Omega _{i-j}(s-j,m+j,k) &  =   \sum_{j = 0}^{i-2}2^{j+1}\cdot \Omega _{i-(j+1)}(s-(j+1),m+(j+1),k) 
  +  \sum_{j = 0}^{i-2}2^{j}\cdot \Delta _{i-j}^{\Big[\substack{s-j \\ s-j+(m+j) }\Big] \times k}  \label{eq 2.28} \\
&   =  \sum_{j = 1}^{i-1}2^{j}\cdot \Omega _{i-j}(s-j,m+ j,k) 
  +  \sum_{j = 0}^{i-2}2^{j}\cdot \Delta _{i-j}^{\Big[\substack{s-j \\ s-j+(m+j) }\Big] \times k}. \nonumber
\end{align}\vspace{0.3 cm}\\     
  Observing that \vspace{0.1 cm}\\
  $$  \Omega _{1}(s-(i-1),m+(i-1),k) = 
 \Gamma _{1}^{\Big[\substack{s -i +1\\ s +m }\Big] \times k} - 4\cdot \Gamma _{0}^{\Big[\substack{s - i+1\\ s +(m-1) }\Big] \times k} = 9 - 4 = 5,  $$\vspace{0.1 cm}\\
   we get by \eqref{eq 2.28} after some simplifications
     \begin{equation}
\label{eq 2.29}
  \Omega _{i}(s,m,k) = 2^{i-1}\cdot \Omega _{1}(s-(i-1),m+(i-1),k) +  \sum_{j = 0}^{i-2}2^{j}\cdot \Delta _{i-j}^{\Big[\substack{s-j \\ s-j+(m+j) }\Big] \times k} 
  =  5\cdot2^{i-1} +   \sum_{j = 0}^{i-2}2^{j}\Delta_{i-j}^{\Big[\substack{s -j\\ s-j +(m+j) }\Big] \times k}.
   \end{equation} \vspace{0.1 cm}\\
    
   Applying  \eqref{eq 2.19} and \eqref{eq 2.29} we obtain   \eqref{eq 2.25}.\vspace{0.1 cm}\\
    
  We prove \eqref{eq 2.26} in a similar way. \vspace{0.1 cm}\\
      
 \textbf{ In Section 10,} we compute explicitly\vspace{0.1 cm}\\ 
  
    $ 2^{s-1}\left[ \Gamma _{j+1}^{\Big[\substack{1 \\ 1+ (s+m-1)}\Big] \times (k+1)} -  \Gamma _{j+1}^{\Big[\substack{1 \\ 1+(s+m-1) }\Big] \times k}\right]
   -   2^{s+1}\left[ \Gamma _{j}^{\Big[\substack{1 \\ 1+ (s+m-2)}\Big] \times (k+1)} -  \Gamma _{j}^{\Big[\substack{1 \\ 1+(s+m-2) }\Big] \times k}\right]  \overset{ \text{def}}{=} R(j,s,m,k) \vspace{0.1 cm}\\ 
     \text{if  }\quad   0 \leq j\leq s+m ,\quad k\geq s+j+1.   $ \vspace{0.1 cm}\\  
   In fact we prove  \vspace{0.1 cm}\\ 
    \begin{equation}
     \label{eq 2.30}
               R(j,s,m,k)   = \begin{cases}
2^{k+s-1} & \text{if  } j = 0 ,\quad k\geq s+1, \\
- 2^{k+s-1} &  \text{if  }    j =1, \quad k\geq s+2, \\
 0  & \text{if   } 2\leq j\leq s+m-1, \quad k \geq s+j+1, \\
  -3\cdot 2^{2k+2s+m-2}         & \text{if   } j = s+m, \quad k \geq 2s+m +1. 
\end{cases}
\end{equation}\vspace{0.1 cm}\\ 
From \eqref{eq 2.25}, \eqref{eq 2.26} and \eqref{eq 2.30} we get  for $ s\geq 2, m\geq 0  $\vspace{0.1 cm}\\ 
\begin{align}
 \Gamma _{i}^{\Big[\substack{s \\ s +m }\Big] \times (k+1)} -  \Gamma _{i}^{\Big[\substack{s \\ s +m }\Big] \times k}  =  4\cdot \left[\Gamma _{i-1}^{\Big[\substack{s \\ s +(m-1) }\Big] \times (k+1)} - \Gamma _{i-1}^{\Big[\substack{s \\ s +(m-1) }\Big] \times  k}\right]  
 & \text{if  }  1 \leq i\leq s-1,\quad k\geq i+1, \label{eq 2.31}  \\
  \nonumber \\
   \Gamma _{s+j}^{\Big[\substack{s \\ s +m }\Big] \times (k+1)} -  \Gamma _{s+j}^{\Big[\substack{s \\ s +m }\Big] \times k}  =  4\cdot \left[\Gamma _{s+j-1}^{\Big[\substack{s \\ s +(m-1) }\Big] \times (k+1)} - \Gamma _{s+j-1}^{\Big[\substack{s \\ s +(m-1) }\Big] \times  k}\right] \label{eq 2.32} \\
+   \begin{cases}
2^{k+s-1} & \text{if  } j = 0 ,\quad k\geq s+1, \\
- 2^{k+s-1} &  \text{if  }    j =1, \quad k\geq s+2, \\
 0  & \text{if   } 2\leq j\leq s+m-1, \quad k \geq s+j+1, \\
  -3\cdot 2^{2k+2s+m-2}         & \text{if   } j = s+m, \quad k \geq 2s+m +1. \nonumber
\end{cases}
\end{align}
We establish \eqref{eq 2.30} by applying the formula for $ \Gamma _{i}^{\Big[\substack{1 \\ 1 + m }\Big] \times k} $  (obtained in 
 [4], see Theorem 3.9, with $ m\rightarrow s+m-1, m\rightarrow s+m-2 $).\vspace{0.1 cm}\\

   \textbf{ In Section 11,}  we establish a formula for   $ \Gamma _{i}^{\Big[\substack{s \\ s + m }\Big] \times k} $ in the case
 $ 1\leq i\leq \inf(s-1,k-1),\quad m\geq 0. $\vspace{0.1 cm}\\ 
 We show precisely   \vspace{0.1 cm}\\ 
  \begin{equation}
\label{eq 2.33}
\Gamma _{i}^{\Big[\substack{s \\ s + m }\Big] \times k}  = \begin{cases}
 1  & \text{if  } i = 0, \\
 \Gamma _{i}^{\Big[\substack{i \\ i }\Big] \times (i+1)}  \overset{ \text{def}}{=} \gamma _{i}   &  \text{if  }  1\leq i\leq \inf(s-1,k-1),  m\geq 0 
\end{cases}
\end{equation}\vspace{0.1 cm}  \\ 
and \vspace{0.1 cm}  \\ 
  \begin{equation}
\label{eq 2.34}
 \Gamma _{i}^{\Big[\substack{i \\ i }\Big] \times (i+1)}  \overset{ \text{def}}{=} \gamma _{i} = 21\cdot2^{3i-4} - 3\cdot2^{2i-3}\;  \text{if  }  1\leq i\leq \inf(s-1,k-1) \; m\geq 0. 
\end{equation}\vspace{0.1 cm}  \\ 
From \eqref{eq 2.33}, \eqref{eq 2.34} we  deduce the following results  \vspace{0.1 cm}  \\  
 \begin{align}
 \Gamma _{k}^{\Big[\substack{s  \\ s+m }\Big] \times k} &  = 2^{2k+2s+m-2} -3\cdot2^{3k-4} + 2^{2k-3}\quad for \; 1\leq k\leq s, \label{eq 2.35}\\
 & \nonumber \\
  \sum_{j = 0}^{s -1} \Gamma _{j}^{\Big[\substack{s  \\ s+m }\Big] \times k} & = 3\cdot2^{3s-4} - 2^{2s-3}\quad for \; k\geq s, \label{eq 2.36}\\
  & \nonumber \\
   \sum_{j = 0}^{s -1} \Gamma _{j}^{\Big[\substack{s  \\ s+m }\Big] \times k}\cdot 2^{-j}  & = 7\cdot2^{2s - 4} - 3\cdot2^{s-3} \quad for \; k\geq s. \label{eq 2.37}
\end{align} 
We  prove \eqref{eq 2.33} by induction on j. \vspace{0.1 cm}  \\ 
Set  \vspace{0.1 cm}  
 \begin{equation}
 \label{eq 2.38}
(H_{j})\quad  \Gamma _{j}^{\Big[\substack{s  \\ s+m }\Big] \times k} =  \Gamma _{j}^{\Big[\substack{j  \\ j}\Big] \times (j+1)} =\gamma _{j}
    \quad   \text{if  \quad   $ s\geq j+1,\; k\geq j+1, \; m\geq 0   $}.\vspace{0.1 cm}  \\
 \end{equation}
 Applying \eqref{eq 2.1} we get, using successively \eqref{eq 2.31}  for  i = 1, 2, 3  \vspace{0.1 cm}  \\
 \begin{align}
  \Gamma _{1}^{\Big[\substack{s  \\ s+m }\Big] \times k}& =  \Gamma _{1}^{\Big[\substack{1  \\ 1}\Big] \times 2} =\gamma _{1}
    \quad   \text{if $ \quad s\geq 2   ,\; k\geq 2, \; m\geq 0 $  }. \label{eq 2.39}\\
    & \nonumber \\
  \text{ Using \eqref{eq 2.39} we get } \quad    \Gamma _{2}^{\Big[\substack{s  \\ s+m }\Big] \times k}& =  \Gamma _{2}^{\Big[\substack{2  \\ 2}\Big] \times 3} =\gamma _{2}
    \quad   \text{ if $ \quad    s\geq 3,\; k\geq 3, \; m\geq 0 $  }. \label{eq 2.40}\\ 
    & \nonumber \\
     \text{ Using  \eqref{eq 2.39},\eqref{eq 2.40} we get } \quad       \Gamma _{3}^{\Big[\substack{s  \\ s+m }\Big] \times k} & =  \Gamma _{3}^{\Big[\substack{3  \\ 3}\Big] \times 4} =\gamma _{3}
    \quad   \text{if $ \quad    s\geq 4,\; k\geq 4, \; m\geq 0 $  }. \label{eq 2.41}
    \end{align}\vspace{0.1 cm}  \\
Thus $ (H_{j}) $  [see \eqref{eq 2.38}] holds for $ j\in\left\{1,2,3\right\}. $ \vspace{0.1 cm}  \\
Assume now that $ (H_{j}) $ holds for $ j\in\left\{1,2,3,\ldots,i-1\right\},$ that is  \vspace{0.1 cm}  \\
\begin{equation}
 \label{eq 2.42}
\bigwedge_{j=1}^{j = i-1} (H_{j}) \hspace{2 cm}  \Gamma _{j}^{\Big[\substack{s  \\ s+m }\Big] \times k} =  \Gamma _{j}^{\Big[\substack{j  \\ j}\Big] \times (j+1)} =\gamma _{j}
    \quad   \text{for all   \quad   $ 1\leq j\leq i-1,\; k\geq j+1, \;s \geq j +1 \; m\geq 0   $}. \vspace{0.1 cm}  \\
 \end{equation}\vspace{0.1 cm}  \\ 
 
 From  $ \bigwedge_{j=1}^{j = i-1} (H_{j}) $ and \eqref{eq 2.31} it follows \vspace{0.1 cm}  \\ 
$ \Gamma _{i}^{\Big[\substack{s \\ s +m }\Big] \times (k+1)} -  \Gamma _{i}^{\Big[\substack{s \\ s +m }\Big] \times k}
 =  4\cdot \left[\Gamma _{i-1}^{\Big[\substack{s \\ s +(m-1) }\Big] \times (k+1)} - \Gamma _{i-1}^{\Big[\substack{s \\ s +(m-1) }\Big] \times  k}\right] = 0 $
  for all $ s \geq i+1, \; k\geq i+1. $ \vspace{0.05 cm}\\
  Hence $ \Gamma _{i}^{\Big[\substack{s \\ s +m }\Big] \times k} = \Gamma _{i}^{\Big[\substack{s \\ s +m }\Big] \times (i+1)} \;\text{for all} \; s\geq i+1, \; k\geq i+1. $\vspace{0.05 cm}\\
To get  \eqref{eq 2.33} we need  only to show that $ \Gamma _{i}^{\Big[\substack{s \\ s +m }\Big] \times (i+1)} \;\text{is equal to} \quad  \Gamma _{i}^{\Big[\substack{i \\ i }\Big] \times (i+1)}
\;\text{for all }\;  s\geq i+1,\; m\geq 0. $\vspace{0.1 cm}\\
To do so we consider the following matrix \vspace{0.1 cm}\\
   $$   \left ( \begin{array} {cccccc}
\alpha _{1} & \alpha _{2} & \alpha _{3} &  \ldots & \alpha _{i}  &  \alpha _{i+1} \\
\alpha _{2 } & \alpha _{3} & \alpha _{4}&  \ldots  &  \alpha _{i+1} &  \alpha _{i+2} \\
\vdots & \vdots & \vdots    &  \vdots & \vdots  &  \vdots \\
\alpha _{s-1} & \alpha _{s} & \alpha _{s +1} & \ldots  &  \alpha _{s+ i-2} &  \alpha _{s+ i-1}  \\
\alpha _{s} & \alpha _{s+1} & \alpha _{s +2} & \ldots  &  \alpha _{s+ i-1} &  \alpha _{s+ i} \\
 \beta  _{1} & \beta  _{2} & \beta  _{3} & \ldots  &  \beta_{i} &  \beta _{i+1}  \\
\beta  _{2} & \beta  _{3} & \beta  _{4} & \ldots  &  \beta_{i+1} &  \beta _{i+2}  \\
\vdots & \vdots & \vdots    &  \vdots & \vdots  &  \vdots \\
\beta  _{m+1} & \beta  _{m+2} & \beta  _{m+3} & \ldots  &  \beta_{i+m} &  \beta _{i+m+1}  \\
\vdots & \vdots & \vdots    &  \vdots & \vdots  &  \vdots \\
\beta  _{s+m-1} & \beta  _{s+m} & \beta  _{s+m+1} & \ldots  &  \beta_{s+m+i-2} &  \beta _{s+m+i-1}  \\
\beta  _{s+m} & \beta  _{s+m+1} & \beta  _{s+m+2} & \ldots  &  \beta_{s+m+ i-1} &  \beta _{s+m+ i}
\end{array}  \right). $$ \vspace{0.5 cm}
  
Using the equations \eqref{eq 2.1} with k = i+1 we obtain \\
    \begin{align}
  \sum_{ j = 0}^{i+1} \Gamma _{j}^{\Big[\substack{s  \\ s+m }\Big] \times (i+1)} & = 2^{2i+2s+m +2} \label{eq 2.43}\\
 &  \text{and}\nonumber \\
  \sum_{j = 0}^{i+1} \Gamma _{j}^{\Big[\substack{s  \\ s+m }\Big] \times (i+1)}\cdot2^{-j} & =  2^{i+2s+m -1} + 2^{2i} - 2^{i-1}. \label{eq 2.44}
\end{align}
  Setting s = i, m = 0   in  \eqref{eq 2.43},  \eqref{eq 2.44} we have \vspace{0.5 cm}  \\
     \begin{align}
  \sum_{ j = 0}^{i+1} \Gamma _{j}^{\Big[\substack{i \\ i}\Big] \times (i+1)} & = 2^{4i} \label{eq 2.45}\\
 &  \text{and}\nonumber \\
  \sum_{j = 0}^{i+1} \Gamma _{j}^{\Big[\substack{i \\ i}\Big] \times (i+1)}\cdot2^{- j} & =  2^{3i -1} + 2^{2i} - 2^{i-1}. \label{eq 2.46}
\end{align} \vspace{0.1 cm}\\
From \eqref{eq 2.43}, \eqref{eq 2.44}, \eqref{eq 2.45}, \eqref{eq 2.46} and   $ \bigwedge_{j=1}^{j = i-1} (H_{j}) $ we deduce \eqref{eq 2.34}. \vspace{0.5 cm}  \\

 \textbf{ In Section 12,} using successively the recurrent formula  \eqref{eq 2.32} for the difference 
$\Gamma _{s+j}^{\Big[\substack{s \\ s+j}\Big] \times (k+1)} - \Gamma _{s+j}^{\Big[\substack{s \\ s+j}\Big] \times k},\\
\text{with}\;  0\leq j\leq s+m,\; k\geq s+j+1 $ we deduce  \vspace{0.1 cm}\\ 
 In the case $ m \in \left\{0,1\right\} $ \vspace{0.01 cm}\\
 \begin{equation}
 \label{eq 2.47}
 \Gamma _{s+j}^{\Big[\substack{s \\ s}\Big] \times (k+1)}
  - \Gamma _{s+j}^{\Big[\substack{s \\ s}\Big] \times k}=
  \begin{cases}
 3\cdot 2^{k+s-1}  & \text{if  } j = 0,\; k >s, \\
 21\cdot 2^{k+s+3j-4}   &  \text{if  }    1\leq j\leq s-1,\; k>s+j, \\
3\cdot2^{2k+2s-2}- 3\cdot2^{k+4s - 4}       & \text{if   } j  = s,  \; k>2s,
\end{cases}
\end{equation}
 \begin{equation}
 \label{eq 2.48}
 \Gamma _{s+j}^{\Big[\substack{s \\ s+1}\Big] \times (k+1)}
  - \Gamma _{s+j}^{\Big[\substack{s \\ s+1}\Big] \times k}=
  \begin{cases}
  2^{k+s-1}  & \text{if  } j = 0,\; k >s, \\
  11\cdot2^{k+s-1}& \text{if  } j = 1,\; k > s +1, \\
21 \cdot 2^{k+s+3j-5}   &  \text{if  }    2\leq j\leq s ,\; k>s+j, \\
3\cdot2^{2k+2s-1}- 3\cdot2^{k+4s -2}       & \text{if   } j  = s+1 ,  \; k>2s +1.
\end{cases}
\end{equation}
 In the case $ m\geq 2  $\vspace{0.1 cm}\\
 \begin{equation}
 \label{eq 2.49}
\Gamma _{s+j}^{\Big[\substack{s \\ s+m}\Big] \times (k+1)}
  - \Gamma _{s+j}^{\Big[\substack{s \\ s+m}\Big] \times k}=
  \begin{cases}
  2^{k+s-1}  & \text{if  } j = 0,\; k >s, \\
  3\cdot2^{k+s +2j -3}& \text{if  } 1\leq j\leq m-1,\; k > s+j,  \\
11\cdot2^{k+ s + 2m -3}              &  \text{if  }  j = m, \; k > s+m,
\end{cases}
\end{equation}
\begin{equation}
\label{eq 2.50}
 \Gamma _{s+m+j}^{\Big[\substack{s \\ s+m}\Big] \times (k+1)}
  - \Gamma _{s+m+j}^{\Big[\substack{s \\ s+m}\Big] \times k}=
  \begin{cases}
  21\cdot2^{k+s +2m +3j -4}& \text{if  } 1\leq j\leq s-1,\; k > s+m+ j,  \\
3\cdot2^{2k+2s+m-2} -  3\cdot2^{k+4s+2m-4}            &  \text{if  }  j = s, \; k > 2s+m.
\end{cases}
\end{equation}\vspace{0.1 cm}\\
For instance to compute $ \Gamma _{s+j}^{\Big[\substack{s \\ s}\Big] \times (k+1)}
  - \Gamma _{s+j}^{\Big[\substack{s \\ s}\Big] \times k} $ for $0\leq j\leq s-1,\; k\geq s+j+1$ we proceed as follows : \vspace{0.1 cm}\\
  From \eqref{eq 2.32} with m = 0, j = 0, $k\geq s+1$ we get, using \eqref{eq 2.33} \vspace{0.1 cm}\\
\begin{align}
& \Gamma _{s}^{\Big[\substack{s \\ s  }\Big] \times (k+1)} -  \Gamma _{s}^{\Big[\substack{s \\ s  }\Big] \times k}  = 
 4\cdot \left[\Gamma _{s-1}^{\Big[\substack{s -1 \\ s -1 +1 }\Big] \times (k+1)} - \Gamma _{s-1}^{\Big[\substack{s -1  \\ s -1+1  }\Big] \times  k}\right]  + 2^{k+s-1}, \label{eq 2.51}\\
 & \nonumber \\
 &  \Gamma _{s -1}^{\Big[\substack{s -1 \\ s-1+1  }\Big] \times (k+1)} -  \Gamma _{s-1}^{\Big[\substack{s-1 \\ s -1+1 }\Big] \times k}  = 
 4\cdot \left[\Gamma _{s-2}^{\Big[\substack{s -1 \\ s -1  }\Big] \times (k+1)} - \Gamma _{s-2 }^{\Big[\substack{s -1  \\ s -1  }\Big] \times  k}\right]  + 2^{k+s-2} =  2^{k+s-2}.\label{eq 2.52}
 \end{align}\vspace{0.1 cm}\\
From \eqref{eq 2.51},  \eqref{eq 2.52} we obtain \vspace{0.1 cm}\\
\begin{equation}
\label{eq 2.53}
 \Gamma _{s}^{\Big[\substack{s \\ s  }\Big] \times (k+1)} -  \Gamma _{s}^{\Big[\substack{s \\ s  }\Big] \times k}  = 
 4\cdot2^{k+s-2} + 2^{k+s-1} = 3\cdot2^{k+s-1}.
\end{equation}\\

 We get in a similar way, using \eqref{eq 2.53} with $s\rightarrow s-1$ \vspace{0.1 cm}\\
\begin{equation}
\label{eq 2.54}
 \Gamma _{s+1}^{\Big[\substack{s \\ s  }\Big] \times (k+1)} -  \Gamma _{s+1}^{\Big[\substack{s \\ s  }\Big] \times k}  = 
  21\cdot2^{k+s-1}. 
\end{equation}\vspace{0.1 cm}\\
 Let  $ 2 \leq j\leq s-1, \; k > s+j $. \vspace{0.1 cm}\\
 By \eqref{eq 2.32} we obtain  \vspace{0.1 cm}\\
\begin{align*}
 \Gamma _{s +j}^{\Big[\substack{s \\ s  }\Big] \times (k+1)} -  \Gamma _{s +j}^{\Big[\substack{s \\ s  }\Big] \times k} &  = 
 4\cdot \left[\Gamma _{s +j-1}^{\Big[\substack{s -1 \\ (s-1) +1 }\Big] \times (k+1)} - 
 \Gamma _{s -1+j}^{\Big[\substack{s -1  \\ (s-1)+1  }\Big] \times  k}\right]  && \text{if $ 2\leq j\leq s-1 $}, \\
 \\
   \Gamma _{s -1+j}^{\Big[\substack{s -1\\ (s-1) +1}\Big] \times (k+1)} -  \Gamma _{s -1+j}^{\Big[\substack{s-1 \\ (s-1)+1 }\Big] \times k} &  = 
 4\cdot \left[\Gamma _{s -1+(j-1)}^{\Big[\substack{s -1 \\  s-1 }\Big] \times (k+1)} - 
 \Gamma _{s -1+(j-1)}^{\Big[\substack{s -1  \\ s-1  }\Big] \times  k}\right] && \text{if $ 2\leq j\leq s-1+1-1= s-1 $}.  \\
\end{align*}\vspace{0.1 cm}\\
From the above equations we get \vspace{0.1 cm}\\
\begin{equation}
\label{eq 2.55}
 \Gamma _{s +j}^{\Big[\substack{s \\ s  }\Big] \times (k+1)} -  \Gamma _{s +j}^{\Big[\substack{s \\ s  }\Big] \times k}   = 
 4^{2}\cdot \left[\Gamma _{s -1 +(j-1)}^{\Big[\substack{s -1 \\ s-1 }\Big] \times (k+1)} - 
 \Gamma _{s -1+(j-1)}^{\Big[\substack{s -1  \\ s-1  }\Big] \times  k}\right]  \quad  if\;  2\leq j\leq s-1. 
\end{equation}\vspace{0.1 cm} \\
Using successively \eqref{eq 2.55} we get from \eqref{eq 2.54} with $ s\rightarrow s-j+1  $ \vspace{0.1 cm} \\
\begin{align}
\label{eq 2.56}
\Gamma _{s +j}^{\Big[\substack{s \\ s  }\Big] \times (k+1)} -  \Gamma _{s +j}^{\Big[\substack{s \\ s  }\Big] \times k} & =
 4^{2(j-1)}\cdot \left[\Gamma _{s -j +2}^{\Big[\substack{s -j+1) \\ s -j+1) }\Big] \times (k+1)} - 
 \Gamma _{s -j+2}^{\Big[\substack{s -j+1  \\ s-j+1  }\Big] \times  k}\right] \quad if \; 2\leq j\leq s-1 \\
 & =  4^{2(j-1)}\cdot21\cdot2^{k+s-j+1-1}= 21\cdot2^{k+s+3j-4}.\nonumber\\
 & \nonumber
\end{align}
 \textbf{ In Sections 13, 14, 15, 16 and 17} we compute  $\Gamma_{s+j}^{\left[s\atop s+m\right]\times k}
\;\text{for}\; 0\leq j\leq s+m,\; k \geq s+j,\;m\geq 0. $ \vspace{0.1 cm} \\
To do so, consider the following array  $ \left(\Gamma_{s+j}^{\left[s\atop s+m\right]\times k}\right)_{0\leq j\leq s+m,\; k\geq s+j}: $
\begin{equation}
\label{eq 2.57}
    \begin{array} {cccccc}
     \Gamma _{s}^{\Big[\substack{s \\ s+m}\Big] \times s}      & \Gamma _{s}^{\Big[\substack{s \\ s+m}\Big] \times (s+1)} & \Gamma _{s}^{\Big[\substack{s \\ s+m}\Big] \times (s+2)} &  \ldots & \Gamma _{s}^{\Big[\substack{s \\ s+m}\Big] \times k} &  \ldots \quad \text{$ k\geq s $} \\
 \Gamma _{s+1}^{\Big[\substack{s \\ s+m}\Big] \times (s+1)}      & \Gamma _{s+1}^{\Big[\substack{s \\ s+m}\Big] \times (s+2)} & \Gamma _{s +1}^{\Big[\substack{s \\ s+m}\Big] \times (s+3)} &  \ldots & \Gamma _{s +1}^{\Big[\substack{s \\ s+m}\Big] \times k} &  \ldots \quad \text{$ k\geq s +1 $} \\
\vdots & \vdots & \vdots    &  \vdots & \vdots  &  \vdots \\
 \Gamma _{s+j}^{\Big[\substack{s \\ s+m}\Big] \times (s+j)}      & \Gamma _{s+j}^{\Big[\substack{s \\ s+m}\Big] \times (s+j+1)} & \Gamma _{s +j}^{\Big[\substack{s \\ s+m}\Big] \times (s+ j+2)} &  \ldots & \Gamma _{s + j}^{\Big[\substack{s \\ s+m}\Big] \times k} &  \ldots \quad \text{$ k\geq s +j $} \\ 
\vdots & \vdots & \vdots    &  \vdots & \vdots  &  \vdots \\
 \Gamma _{s+m}^{\Big[\substack{s \\ s+m}\Big] \times (s+m)}      & \Gamma _{s+m}^{\Big[\substack{s \\ s+m}\Big] \times (s+m+1)} & \Gamma _{s +m}^{\Big[\substack{s \\ s+m}\Big] \times (s+ m+2)} &  \ldots & \Gamma _{s + m}^{\Big[\substack{s \\ s+m}\Big] \times k} &  \ldots \quad \text{$ k\geq s +m $} \\ 
  \Gamma _{s+m +1}^{\Big[\substack{s \\ s+m}\Big] \times (s+m+1)}      & \Gamma _{s+m +1}^{\Big[\substack{s \\ s+m}\Big] \times (s+m+2)} & \Gamma _{s +m +1}^{\Big[\substack{s \\ s+m}\Big] \times (s+ m+3)} &  \ldots & \Gamma _{s + m+1}^{\Big[\substack{s \\ s+m}\Big] \times k} &  \ldots \quad \text{$ k\geq s +m+1 $} \\ 
\vdots & \vdots & \vdots    &  \vdots & \vdots  &  \vdots \\
  \Gamma _{s+m +1+j}^{\Big[\substack{s \\ s+m}\Big] \times (s+m+ 1+j)}      & \Gamma _{s+m +1+j}^{\Big[\substack{s \\ s+m}\Big] \times (s+m+ 2+j)} & \Gamma _{s +m +1+j}^{\Big[\substack{s \\ s+m}\Big] \times (s+ m+3+j)} &  \ldots & \Gamma _{s + m+1+j}^{\Big[\substack{s \\ s+m}\Big] \times k} &  \ldots \quad \text{$ k\geq s +m+1+j $} \\  
   \Gamma _{2s+m -1}^{\Big[\substack{s \\ s+m}\Big] \times (2s+m-1)}      & \Gamma _{2s+m -1}^{\Big[\substack{s \\ s+m}\Big] \times (2s+m)} & \Gamma _{2s +m -1}^{\Big[\substack{s \\ s+m}\Big] \times (2s+ m+1)} &  \ldots & \Gamma _{2s + m -1}^{\Big[\substack{s \\ s+m}\Big] \times k} &  \ldots \quad \text{$ k\geq 2s+m -1 $} \\ 
      \Gamma _{2s+m }^{\Big[\substack{s \\ s+m}\Big] \times (2s+m)}      & \Gamma _{2s+m }^{\Big[\substack{s \\ s+m}\Big] \times (2s+m+1)} & \Gamma _{2s +m }^{\Big[\substack{s \\ s+m}\Big] \times (2s+ m+2)} &  \ldots & \Gamma _{2s + m }^{\Big[\substack{s \\ s+m}\Big] \times k} &  \ldots  \quad \text{$ k\geq 2s+m  $}.
\end{array}
     \end{equation}
     \begin{itemize}
     \item
   To compute the elements in the array  $ \left(\Gamma_{s+j}^{\left[s\atop s+m\right]\times k}\right)_{0\leq j\leq s+m,\; k\geq s+j}, $
   we proceed by induction on j as follows:\vspace{0.1 cm}\\
   Let $l$ be a rational integer such that $  0\leq l\leq s+m-1. $\vspace{0.1 cm}\\
   Assume that we  have computed the elements in the following subarray
    $ \left(\Gamma_{s+j}^{\left[s\atop s+m\right]\times k}\right)_{0\leq j\leq l,\; k\geq s+j}. $\vspace{0.1 cm}\\
Recall that from \eqref{eq 2.33}, \eqref{eq 2.34} $ \Gamma _{i }^{\Big[\substack{s \\ s+m}\Big] \times k}$ are known for $ 0\leq i\leq s-1,\;k\geq i . $\vspace{0.1 cm}\\
     From the first equation in \eqref{eq 2.1} with k = s+l+1,    $ \quad \sum_{ i = 0}^{s+l+1} \Gamma _{i}^{\Big[\substack{s  \\ s+m }\Big] \times (s+l+1)}  = 2^{4s+2l+m} $ \vspace{0.1 cm}\\
  we deduce  $  \Gamma _{s+l+1}^{\Big[\substack{s \\ s+m}\Big] \times (s+l+1)} $ since the terms $ \Gamma _{i}^{\Big[\substack{s  \\ s+m }\Big] \times (s+l+1)} $ for $ i\leq s +l $ are known.  \vspace{0.1 cm}\\
   From the  equations in \eqref{eq 2.1} with $ k = s+ l+2 $ we get 
   \begin{equation*}
  \begin{cases} 
 \sum_{ i = 0}^{s+ l+2} \Gamma _{i}^{\Big[\substack{s  \\ s+m }\Big] \times (s+l+2)}  = 2^{4s+ 2l+m+2}, \\
  \sum_{i = 0}^{s+ l+2} \Gamma _{i}^{\Big[\substack{s  \\ s+m }\Big] \times (s+l+2)}\cdot2^{-i}  =  2^{3s+m +l} + 2^{2s +2l +2} - 2^{s +l}. 
\end{cases}
    \end{equation*}\vspace{0.1 cm}\\ 
  We deduce  $  \Gamma _{s+l+1}^{\Big[\substack{s \\ s+m}\Big] \times (s+1+2)} $ and  $  \Gamma _{s+l+2}^{\Big[\substack{s \\ s+m}\Big] \times (s+ l+2)} $ since the terms $ \Gamma _{i}^{\Big[\substack{s  \\ s+m }\Big] \times (s+l+2)} $ for $ i\leq s +l $ are known. \vspace{0.1 cm}\\
  Then, from \eqref{eq 2.49} or \eqref{eq 2.50}, in the case  $ m\geq 2, $  we compute $\Gamma _{s+l+1}^{\Big[\substack{s \\ s+m}\Big] \times k}$ knowing $  \Gamma _{s+l+1}^{\Big[\substack{s \\ s+m}\Big] \times (s+1+2)} $ for all $k\geq s+l+2. $\vspace{0.1 cm}\\ 
   We have now computed any element in the $l+1$-th row and the first element in the $l+2$-th row  in the array  $ \left(\Gamma_{s+j}^{\left[s\atop s+m\right]\times k}\right)_{0\leq j\leq s+m,\; k\geq s+j}. $\\

   \item  From the first equation in \eqref{eq 2.1}  with k = s,    $ \quad \sum_{ i = 0}^{s} \Gamma _{i}^{\Big[\substack{s  \\ s+m }\Big] \times s}  = 2^{4s+m-2} $ \vspace{0.1 cm}\\
  we deduce  $  \Gamma _{s}^{\Big[\substack{s \\ s+m}\Big] \times s} $ since the terms $ \Gamma _{i}^{\Big[\substack{s  \\ s+m }\Big] \times s} $ for $ i\leq s-1 $ are known by \eqref{eq 2.33}, \eqref{eq 2.34}. \vspace{0.1 cm}\\
   From the  equations in \eqref{eq 2.1} with k = s+1 we obtain
    \begin{equation*}
  \begin{cases} 
 \sum_{ i = 0}^{s+1} \Gamma _{i}^{\Big[\substack{s  \\ s+m }\Big] \times (s+1)}  = 2^{4s+m}, \\
  \sum_{i = 0}^{s+1} \Gamma _{i}^{\Big[\substack{s  \\ s+m }\Big] \times (s+1)}\cdot2^{-i}  =  2^{3s+m- 1} + 2^{2s} - 2^{s-1}. 
\end{cases}
    \end{equation*}\vspace{0.1 cm}\\ 
  We deduce  $  \Gamma _{s}^{\Big[\substack{s \\ s+m}\Big] \times (s+1)} $ and  $  \Gamma _{s+1}^{\Big[\substack{s \\ s+m}\Big] \times (s+1)} $ since the terms $ \Gamma _{i}^{\Big[\substack{s  \\ s+m }\Big] \times (s+1)} $ for $ i\leq s-1 $ are known.  \vspace{0.1 cm}\\
  Then, from \eqref{eq 2.49} in the case $m\geq 2 $ with j = 0 we compute $\Gamma _{s}^{\Big[\substack{s \\ s+m}\Big] \times k}$ for all $k\geq s+1.$\vspace{0.1 cm}\\ 
   We have now computed any element in the subarray  $ \left(\Gamma_{s+j}^{\left[s\atop s+m\right]\times k}\right)_{0\leq j\leq l,\; k\geq s+j}\quad \text{with}\; l = 0.  $\\
   
     \item 
   From the  equations in \eqref{eq 2.1} with k = s+ 2, we have 
    \begin{equation*}
  \begin{cases} 
 \sum_{ i = 0}^{s+2} \Gamma _{i}^{\Big[\substack{s  \\ s+m }\Big] \times (s+2)}  = 2^{4s+m+2}, \\
  \sum_{i = 0}^{s+2} \Gamma _{i}^{\Big[\substack{s  \\ s+m }\Big] \times (s+2)}\cdot2^{-i}  =  2^{3s+m} + 2^{2s+2} - 2^{s}. 
\end{cases}
    \end{equation*}\vspace{0.1 cm}\\ 
  We deduce  $  \Gamma _{s+1}^{\Big[\substack{s \\ s+m}\Big] \times (s+2)} $ and  $  \Gamma _{s+2}^{\Big[\substack{s \\ s+m}\Big] \times (s+2)} $ since the terms $ \Gamma _{i}^{\Big[\substack{s  \\ s+m }\Big] \times (s+2)} $ for $ i\leq s $ are known.  \vspace{0.1 cm}\\
   Then from \eqref{eq 2.49}   in the case  $m\geq 2$  with j = 1  we compute $\Gamma _{s+1}^{\Big[\substack{s \\ s+m}\Big] \times k}$  knowing   $  \Gamma _{s+1}^{\Big[\substack{s \\ s+m}\Big] \times (s+2)} $  for all $k\geq s+2.$\vspace{0.1 cm}\\ 
   We have now computed any element in the subarray  $ \left(\Gamma_{s+j}^{\left[s\atop s+m\right]\times k}\right)_{0\leq j\leq l,\; k\geq s+j}\quad \text{with}\; l=1. $\vspace{0.1 cm}\\ 
   
\item    From the  equations in \eqref{eq 2.1} with k = s+ 3, we get 
 \begin{equation*}
  \begin{cases} 
 \sum_{ i = 0}^{s+3} \Gamma _{i}^{\Big[\substack{s  \\ s+m }\Big] \times (s+3)}  = 2^{4s+m+4}, \\
  \sum_{i = 0}^{s+3} \Gamma _{i}^{\Big[\substack{s  \\ s+m }\Big] \times (s+3)}\cdot2^{-i}  =  2^{3s+m +1} + 2^{2s+4} - 2^{s +1}. 
\end{cases}
    \end{equation*}\vspace{0.1 cm}\\ 
  We deduce  $  \Gamma _{s+2}^{\Big[\substack{s \\ s+m}\Big] \times (s+3)} $ and  $  \Gamma _{s+3}^{\Big[\substack{s \\ s+m}\Big] \times (s+3)} $ since the terms $ \Gamma _{i}^{\Big[\substack{s  \\ s+m }\Big] \times (s+3)} $ for $ i\leq s +1$ are known.  \vspace{0.1 cm}\\
   Then from \eqref{eq 2.46} with j = 2 we compute $\Gamma _{s+2}^{\Big[\substack{s \\ s+m}\Big] \times k}$ knowing   $  \Gamma _{s+2}^{\Big[\substack{s \\ s+m}\Big] \times (s+3)} $  for all $k\geq s+3. $\vspace{0.1 cm}\\ 
   We have now computed any element in the 
     subarray  $ \left(\Gamma_{s+j}^{\left[s\atop s+m\right]\times k}\right)_{0\leq j\leq 2,\; k\geq s+j} $\quad with $l = 2.$\vspace{0.1 cm}\\ 
   Further we deduce that  \begin{equation}
   \label{eq 2.58}
 \Gamma _{s+2}^{\Big[\substack{s \\ s+m}\Big] \times k} = 8\cdot \Gamma _{s+1}^{\Big[\substack{s \\ s+m -1}\Big] \times (k-1)}\quad \text{for}\; k\geq s+2,\; m\geq 2. 
\end{equation}\\

\item To compute the elements in the subarray  $ \left(\Gamma_{s+j}^{\left[s\atop s+m\right]\times k}\right)_{2\leq j\leq  m,\; k\geq s+j} $
  we proceed by induction on j. \vspace{0.1 cm}\\
To do so, let $l$ be a rational integer such that $ 2\leq l\leq m.$\vspace{0.1 cm}\\
Assume $ (H_{l-1}) \quad  \Gamma _{s+j}^{\Big[\substack{s \\ s+m}\Big] \times k} = 8^{j-1}\cdot \Gamma _{s+1}^{\Big[\substack{s \\ s+ (m-(j-1))}\Big] \times (k-(j-1))}\quad \text{for}\quad 1\leq j\leq l-1,\;k\geq s+j. $\vspace{0.1 cm}\\ 
$ (H_{l-1})$ holds for $l = 2$ and for $l = 3$ by \eqref{eq 2.58}. \vspace{0.1 cm}\\
To establish $(H_{l})$ we need only to show that \vspace{0.1 cm}\\
 \begin{equation}
 \label{eq 2.59}
 \Gamma _{s+l}^{\Big[\substack{s \\ s+m}\Big] \times k} = 8^{l-1}\cdot \Gamma _{s+1}^{\Big[\substack{s \\ s+ (m-(l-1))}\Big] \times (k-(l-1))}\quad \text{for}\quad k\geq s+ l. 
\end{equation}
In the case $k>s+l$ we proceed as follows:\vspace{0.1 cm}\\
  From the  equations in \eqref{eq 2.1} with $ k = s+ l+1,$ we have 
   \begin{equation*}
  \begin{cases} 
 \sum_{ i = 0}^{s+l+1} \Gamma _{i}^{\Big[\substack{s  \\ s+m }\Big] \times (s+l+1)}  = 2^{4s+m+2l}, \\
  \sum_{i = 0}^{s+l+1} \Gamma _{i}^{\Big[\substack{s  \\ s+m }\Big] \times (s+l+1)}\cdot2^{-i}  =  2^{3s+m + l-1} + 2^{2s+2l} - 2^{s +l-1}. 
\end{cases}
    \end{equation*}\vspace{0.1 cm}\\ 
     We deduce  $  \Gamma _{s+l}^{\Big[\substack{s \\ s+m}\Big] \times (s+l+1)} $ and  $  \Gamma _{s+l+1}^{\Big[\substack{s \\ s+m}\Big] \times (s+l+1)} $ since by $(H_{l-1})$ the terms $ \Gamma _{i}^{\Big[\substack{s  \\ s+m }\Big] \times (s+l+1)} $ for $ i\leq s + l-1 $ are known.  \vspace{0.1 cm}\\ 
 Then from \eqref{eq 2.49} with $j = l $ we compute    $  \Gamma _{s+l}^{\Big[\substack{s \\ s+m}\Big] \times k} $ for all $ k\geq s+l+1 $ and \eqref{eq 2.59} follows.\vspace{0.1 cm}\\ 
 The case $k = s+l $ is obtained in a similar way.\vspace{0.1 cm}\\ 
 
 \item To compute the elements in the remaining  subarray  $ \left(\Gamma_{s +m+1+j}^{\left[s\atop s+m\right]\times k}\right)_{0\leq j\leq  s-1,\; k\geq s+ m+1+j} $
  we proceed still  by induction on j.\vspace{0.1 cm}\\
  Let j = 0, a similar proof of \eqref{eq 2.58} gives the following identity:\vspace{0.1 cm}\\
  \begin{equation}
  \label{eq 2.60}
 \Gamma _{s+m+1}^{\Big[\substack{s \\ s+m}\Big] \times k} = 8^{m}\cdot \Gamma _{s +1}^{\Big[\substack{s \\ s }\Big] \times (k -m)}\quad \text{for}\quad k\geq s+m+1. 
\end{equation}
 Let $l$ be a rational integer such that $ 1\leq l\leq s-1.$\vspace{0.1 cm}\\
Assume $ (H_{l-1}) \quad  \Gamma _{s+m+1+j}^{\Big[\substack{s \\ s+m}\Big] \times k} = 8^{2j+m}\cdot \Gamma _{s -j+1}^{\Big[\substack{s -j\\ s -j}\Big] \times (k -m -2j)}\quad \text{for}\quad 0\leq j\leq l-1,\;k\geq s+m+1+j. $\vspace{0.1 cm}\\
From \eqref{eq 2.60} $ (H_{l-1})$  holds for  $l = 1.$ \vspace{0.1 cm}\\ 
To establish $(H_{l})$ we need only to show that  \vspace{0.1 cm}\\
\begin{equation}
\label{eq 2.61}
\Gamma _{s+m+1+l}^{\Big[\substack{s \\ s+m}\Big] \times k} = 8^{2l+m}\cdot \Gamma _{s - l+1}^{\Big[\substack{s - l\\ s -l}\Big] \times (k-(l-1))}\quad \text{for}\quad k\geq s+ m+1+l.
\end{equation}
In the case $k>s+m+1+l$ we proceed as follows :\vspace{0.1 cm}\\
  From the  equations in \eqref{eq 2.1} with $ k = s+m+l+2,$ we get 
   \begin{equation*}
  \begin{cases} 
 \sum_{ i = 0}^{s+m+l+2} \Gamma _{i}^{\Big[\substack{s  \\ s+m }\Big] \times (s+ m+l+2)}  = 2^{4s+3m+2l+2}, \\
  \sum_{i = 0}^{s+m+l+2} \Gamma _{i}^{\Big[\substack{s  \\ s+m }\Big] \times (s+m+l+2)}\cdot2^{-i}  =  2^{3s+ 2m +l} + 2^{2s+2m+2l+2} - 2^{s + m+l}. 
\end{cases}
    \end{equation*}\vspace{0.1 cm}\\ 
     We deduce  $  \Gamma _{s+m+l+1}^{\Big[\substack{s \\ s+m}\Big] \times (s+m+l+2)} $ and  $  \Gamma _{s+ m+l+2}^{\Big[\substack{s \\ s+m}\Big] \times (s+ m+l+2)} $ since by $(H_{l-1})$ the terms $ \Gamma _{i}^{\Big[\substack{s  \\ s+m }\Big] \times (s+m+l+2)} $ for $ i\leq s + m+l $ are known.  \vspace{0.1 cm}\\ 
 Then from \eqref{eq 2.50} with $j = l $ knowing  $  \Gamma _{s+m+l+1}^{\Big[\substack{s \\ s+m}\Big] \times (s+m+l+2)} $ we compute    $  \Gamma _{s+ m+l+1}^{\Big[\substack{s \\ s+m}\Big] \times k} $ for all $ k\geq s+m+l+1 $ and \eqref{eq 2.61} follows.\vspace{0.1 cm}\\ 
 The case $k = s+m+l+1 $ is obtained in a similar way.\vspace{0.1 cm}\\ 
 \end{itemize}

\section{\textbf{STATEMENT OF RESULTS}}
\label{sec 3}
\begin{thm}
\label{thm 3.1}Let q be a rational integer $ \geq 1,$ then \\
 \begin{align}
    g_{k,s,m}(t,\eta )  =  g(t,\eta ) & = \sum_{deg Y\leq k-1}\sum_{deg Z \leq  s-1}E(tYZ)\sum_{deg U \leq s+m-1}E(\eta YU)  =
   2^{2s+m+k- r( D^{\left[\stackrel{s}{s+m}\right] \times k }(t,\eta  ) )}, \label{eq 3.1} \\
   & \nonumber \\
    \int_{\mathbb{P}\times \mathbb{P}} g_{k,s,m}^{q}(t,\eta )dtd\eta 
 & =2^{(2s+m+k)(q-1)}\cdot 2^{-k+2}\cdot \sum_{i = 0}^{\inf(2s+m,k)} \Gamma _{i}^{\Big[\substack{s  \\ s+m }\Big] \times k}\cdot2^{- qi}. \label{eq 3.2} \\
  & \nonumber 
 \end{align}
\end{thm}

  \begin{thm}
 \label{thm 3.2} Let $ s\geq 2, \; m\geq 0, \; k\geq 1 $ and  $ 0\leq i\leq \inf{(2s+m,k)}. $ Then we have the following recurrent formula for the number $ \Gamma _{i}^{\Big[\substack{s \\ s+m }\Big] \times k} $ of rank i matrices 
 of the form $ \left[{A\over B}\right],$ such that A is a $ s\times k $ persymmetric matrix and B a  $ (s+m)\times k $ persymmetric matrix  with entries in  $ \mathbb{F}_{2}: $\\
   \begin{align}
  \Gamma _{i}^{\Big[\substack{s \\ s+m }\Big] \times k} & = 2\cdot \Gamma _{i-1}^{\Big[\substack{s -1\\ s-1+(m+1) }\Big] \times k}
+ 4\cdot \Gamma _{i-1}^{\Big[\substack{s \\ s+(m-1) }\Big] \times k} - 8\cdot \Gamma _{i-2}^{\Big[\substack{s-1 \\ s-1+m }\Big] \times k}
 + \Delta _{i}^{\Big[\substack{s \\ s+m }\Big] \times k} \label{eq 3.3} \\
 & \text{where the remainder  $\Delta _{i}^{\Big[\substack{s \\ s+m }\Big] \times k}$ is equal to} \nonumber \\
  & \sigma _{i,i,i}^{\left[\stackrel{s-1}{\stackrel{s+m-1 }{\overline {\stackrel{\alpha_{s -}}{\beta_{s+m-} }}}}\right] \times k }
  - 3\cdot \sigma _{i-1,i-1,i-1}^{\left[\stackrel{s-1}{\stackrel{s+m-1 }{\overline {\stackrel{\alpha_{s -}}{\beta_{s+m-} }}}}\right] \times k } 
  + 2\cdot \sigma _{i-2,i-2,i-2}^{\left[\stackrel{s-1}{\stackrel{s+m-1 }{\overline {\stackrel{\alpha_{s -}}{\beta_{s+m-} }}}}\right] \times k }.\label{eq 3.4} 
  \end{align}
    Recall that   $$  \sigma _{i,i,i}^{\left[\stackrel{s-1}{\stackrel{s+m-1 }
{\overline {\stackrel{\alpha_{s -}}{\beta_{s+m-} }}}}\right] \times k } $$\\  is equal to  the cardinality of the following set \\
 \small
 $$   \left\{(t,\eta )\in \mathbb{P}/\mathbb{P}_{k+s-1}\times \mathbb{P}/\mathbb{P}_{k+s+m-1}
\mid  r(D^{\big[\stackrel{s-1}{s-1+m}\big] \times k }(t,\eta ))  = r(D^{\big[\stackrel{s}{s+m-1}\big] \times k }(t,\eta ))
= r(D^{\big[\stackrel{s}{s+m}\big] \times k }(t,\eta )) = i
 \right\}. $$\\ 
   \end{thm}
\begin{thm}
\label{thm 3.3}Let $ s\geq 2 $ and $ m\geq 0, $ we have in the following two cases : \vspace{0.1 cm}\\
\underline {The case $1 \leq  k \leq 2s+m-2 $}\\
\begin{equation}
\label{eq 3.5}
 \sigma _{i,i,i}^{\left[\stackrel{s}{\stackrel{s+m }{\overline {\stackrel{\alpha_{s -}}{\beta_{s+m-} }}}}\right] \times k }= 
  \begin{cases}
1 & \text{if  } i = 0,\quad k\geq 1, \\
   4 \cdot \Gamma _{i}^{\Big[\substack{s -1\\ s-1 +m }\Big] \times i} -  \Gamma _{i+1}^{\Big[\substack{s-1 \\ s-1+m }\Big] \times (i+1)} 
      & \text{if   } 1 \leq i \leq k-1, \\
   4 \cdot \Gamma _{k}^{\Big[\substack{s -1\\ s-1 +m }\Big] \times k}    & \text{if   } i = k.
\end{cases}
\end{equation}

\underline {The case $ k\geq  2s+m-2 $ }\\
\begin{equation}
\label{eq 3.6}
 \sigma _{i,i,i}^{\left[\stackrel{s}{\stackrel{s+m }{\overline {\stackrel{\alpha_{s -}}{\beta_{s+m-} }}}}\right] \times k }= 
  \begin{cases}
1 & \text{if  } i = 0,          \\
   4 \cdot \Gamma _{i}^{\Big[\substack{s -1\\ s-1 +m }\Big] \times i} -  \Gamma _{i+1}^{\Big[\substack{s-1 \\ s-1+m }\Big] \times (i+1)} 
      & \text{if   } 1 \leq i \leq 2s+m-3,  \\
   4 \cdot \Gamma _{2s+m-2}^{\Big[\substack{s -1\\ s-1 +m }\Big] \times (2s+m-2)}    & \text{if   } i = 2s+m-2. 
\end{cases}
\end{equation}
\end{thm}

\begin{thm}The remainder $\Delta _{i}^{\Big[\substack{s \\ s+m }\Big] \times k} $in the recurrent formula is equal to \vspace{0.1 cm}\\
\label{thm 3.4}
   \begin{equation}
\label{eq 3.7}
  \begin{cases}
1 & \text{if  } i = 0,\; k \geq 1,         \\
   4 \cdot \Gamma _{1}^{\Big[\substack{s -1\\ s-1 +m }\Big] \times 1} -  \Gamma _{2}^{\Big[\substack{s-1 \\ s-1+m }\Big] \times 2} 
      & \text{if   } i = 1,\; k\geq 2, \\
    4 \cdot \Gamma _{1}^{\Big[\substack{s -1\\ s-1 +m }\Big] \times 1} - 3   & \text{if   } i = 1,\; k =1, \\
     7 \cdot \Gamma _{2}^{\Big[\substack{s -1\\ s-1 +m }\Big] \times 2} -  12\cdot \Gamma _{1}^{\Big[\substack{s-1 \\ s-1+m }\Big] \times 1}
     -  \Gamma _{3}^{\Big[\substack{s -1\\ s-1 +m }\Big] \times 3} + 2  & \text{if   } i = 2,\; k\geq 3, \\
        7 \cdot \Gamma _{2}^{\Big[\substack{s -1\\ s-1 +m }\Big] \times 2} -  12\cdot \Gamma _{1}^{\Big[\substack{s-1 \\ s-1+m }\Big] \times 1} +2
         & \text{if   } i = 2,\; k = 2, \\
         7 \cdot \Gamma _{i}^{\Big[\substack{s -1\\ s-1 +m }\Big] \times i} -  14\cdot \Gamma _{i-1}^{\Big[\substack{s-1 \\ s-1+m }\Big] \times (i-1)}
     +8\cdot \Gamma _{i-2}^{\Big[\substack{s -1\\ s-1 +m }\Big] \times (i-2)}  -  \Gamma _{i+1}^{\Big[\substack{s -1\\ s-1 +m }\Big] \times (i+1)}   & \text{if } 3\leq i\leq 2s+m-3,  k\geq i+1, \\
        7 \cdot \Gamma _{i}^{\Big[\substack{s -1\\ s-1 +m }\Big] \times i} -  14\cdot \Gamma _{i-1}^{\Big[\substack{s-1 \\ s-1+m }\Big] \times (i-1)}
     +8\cdot \Gamma _{i-2}^{\Big[\substack{s -1\\ s-1 +m }\Big] \times (i-2)}    & \text{if   }  3\leq i\leq 2s+m-3, \; k = i, \\
       7 \cdot \Gamma _{2s+m-2}^{\Big[\substack{s -1\\ s-1 +m }\Big] \times (2s+m-2)} -  14\cdot \Gamma _{2s+m-3}^{\Big[\substack{s-1 \\ s-1+m }\Big] \times (2s+m-3)}
     +8\cdot \Gamma _{2s+m-4}^{\Big[\substack{s -1\\ s-1 +m }\Big] \times (2s+m-4)}    & \text{if   }  i = 2s+m-2, \; k \geq  i, \\
     - 14\cdot \Gamma _{2s+m-2}^{\Big[\substack{s-1 \\ s-1+m }\Big] \times (2s+m-2)}
     +8\cdot \Gamma _{2s+m-3}^{\Big[\substack{s -1\\ s-1 +m }\Big] \times (2s+m-3)}    & \text{if   }  i = 2s+m-1, \; k \geq  i, \\
      8\cdot \Gamma _{2s+m-2}^{\Big[\substack{s -1\\ s-1 +m }\Big] \times (2s+m-2)}    & \text{if   }  i = 2s+m , \; k \geq  i. 
\end{cases}
\end{equation}
  \end{thm}

     \begin{thm}
\label{thm 3.5}
We have
  \begin{align}
\Delta _{i}^{\Big[\substack{s \\ s+m }\Big] \times k}&  =  \Delta _{i}^{\Big[\substack{s \\ s+m }\Big] \times (i+1)} & \;  for \; i\in [0,2s+m-3],  \; k\geq i +1, \label{eq 3.8}\\
 \Delta _{i}^{\Big[\substack{s \\ s+m }\Big] \times k} & =  \Delta _{i}^{\Big[\substack{s \\ s+m }\Big] \times i} & \;  for \; i\in \left\{2s+m-2, 2s+m-1,2s+m \right\}, \; k\geq i. \label{eq 3.9}
\end{align}
\end{thm} 
 
  \begin{thm}
 \label{thm 3.6}We have for all $ m\geq 0 $\vspace{0.01 cm}\\
  \begin{equation}
\label{eq 3.10}
 \Gamma _{j}^{\Big[\substack{s \\ s +m }\Big] \times (k+1)}
  - \Gamma _{j}^{\Big[\substack{s \\ s +m}\Big] \times k}= 0  \quad \text{if $ \quad 0\leq j\leq s-1,\;k >j. $}
 \end{equation}
  We have in the cases  $ m \in \left\{0,1\right\} $ \vspace{0.01 cm}\\
 \begin{equation}
\label{eq 3.11}
 \Gamma _{s+j}^{\Big[\substack{s \\ s}\Big] \times (k+1)}
  - \Gamma _{s+j}^{\Big[\substack{s \\ s}\Big] \times k}=
  \begin{cases}
 3\cdot 2^{k+s-1}  & \text{if  } j = 0,\; k >s, \\
 21\cdot 2^{k+s+3j-4}   &  \text{if  }    1\leq j\leq s-1,\; k>s+j, \\
3\cdot2^{2k+2s-2}- 3\cdot2^{k+4s - 4}       & \text{if   } j  = s,  \; k>2s,
\end{cases}
\end{equation}
 \begin{equation}
\label{eq 3.12}
 \Gamma _{s+j}^{\Big[\substack{s \\ s+1}\Big] \times (k+1)}
  - \Gamma _{s+j}^{\Big[\substack{s \\ s+1}\Big] \times k}=
  \begin{cases}
  2^{k+s-1}  & \text{if  } j = 0,\; k >s, \\
  11\cdot2^{k+s-1}& \text{if  } j = 1,\; k > s +1, \\
21 \cdot 2^{k+s+3j-5}   &  \text{if  }    2\leq j\leq s ,\; k>s+j, \\
3\cdot2^{2k+2s-1}- 3\cdot2^{k+4s -2}       & \text{if   } j  = s+1 ,  \; k>2s +1.
\end{cases}
\end{equation}
 In the case $ m\geq 2  $\vspace{0.1 cm}\\
 \begin{equation}
\label{eq 3.13}
 \Gamma _{s+j}^{\Big[\substack{s \\ s+m}\Big] \times (k+1)}
  - \Gamma _{s+j}^{\Big[\substack{s \\ s+m}\Big] \times k}=
  \begin{cases}
  2^{k+s-1}  & \text{if  } j = 0,\; k >s, \\
  3\cdot2^{k+s +2j -3}& \text{if  } 1\leq j\leq m-1,\; k > s+j,  \\
11\cdot2^{k+ s + 2m -3}              &  \text{if  }  j = m, \; k > s+m,
\end{cases}
\end{equation}
\begin{equation}
\label{eq 3.14}
 \Gamma _{s+m+j}^{\Big[\substack{s \\ s+m}\Big] \times (k+1)}
  - \Gamma _{s+m+j}^{\Big[\substack{s \\ s+m}\Big] \times k}=
  \begin{cases}
  21\cdot2^{k+s +2m +3j -4}& \text{if  } 1\leq j\leq s-1,\; k > s+m+ j,  \\
3\cdot2^{2k+2s+m-2} -  3\cdot2^{k+4s+2m-4}            &  \text{if  }  j = s, \; k > 2s+m.
\end{cases}
\end{equation}
 \end{thm}
 
 \begin{thm}
 \label{thm 3.7}We have for $ m\geq 1 $ \vspace{0.01 cm}\\
   \begin{align}
  \Gamma _{s +j}^{\Big[\substack{s \\ s +m }\Big] \times k} &  =  8^{j-1}\cdot \Gamma _{s+1}^{\Big[\substack{s \\ s +(m-(j-1)) }\Big] \times (k-(j-1))}  && \text{if\;$1\leq j\leq m,\;k\geq s+j, $} \label{eq 3.15}\\
  & \nonumber \\
    \Gamma _{s+1}^{\Big[\substack{s \\ s +(m-(j-1)) }\Big] \times (s+1)} & =  2^{4s+(m-(j-1))} - 3\cdot2^{3s -1} + 2^{2s -1} &&\text{if\; $ 1\leq j\leq m,\; k = s+j, $} \label{eq 3.16}\\ 
    & \nonumber \\
     \Gamma _{s+1}^{\Big[\substack{s \\ s +(m-(j-1)) }\Big] \times (k-(j-1))}& = 3\cdot2^{k - j +s} +21\cdot[2^{3s-1}-2^{2s-1}] && \text{if\;$1\leq j\leq m-1,\;k>s+j, $} \label{eq 3.17} \\
      & \nonumber \\
  \Gamma _{s +1}^{\Big[\substack{s \\ s + 1 }\Big] \times (k-(m-1))}& = 11\cdot2^{k- m +s} + 21\cdot2^{3s- 1} - 11\cdot2^{2s- 1} && \text{if\;$ j = m,\;k>s+m. $} \label{eq 3.18} \\
   & \nonumber 
\end{align}
 \end{thm}
  \begin{thm}
\label{thm  3.8}We have for $ m\geq 0 $ \vspace{0.1 cm}
\begin{align}
 \Gamma _{s +m +1 +j}^{\Big[\substack{s \\ s +m }\Big] \times k} &  =  8^{2j + m}\cdot \Gamma _{s-j+1}^{\Big[\substack{s -j\\ s -j }\Big] \times (k- m-2j)}  &&\text{if   \; $ 0\leq j\leq s-1,\;  k\geq s + m+1+j $}, \label{eq 3.19}\vspace{0.1 cm} \\
  & \nonumber \\
 \Gamma _{s - j+1}^{\Big[\substack{s -j\\ s -j }\Big] \times ( s - j+1)} & = 2^{4s-4j} -3\cdot2^{3s-3j-1} + 2^{2s-2j-1} &&\text{if   \; $ 0\leq j\leq s-1,\;  k =  s + m+1+j $}, \label{eq 3.20} \vspace{0.1 cm} \\
  & \nonumber \\
 \Gamma _{s-j+1}^{\Big[\substack{s -j\\ s -j }\Big] \times (k- m-2j)} 
   & = 21\cdot[2^{k-m -3j+s-1} + 2^{3s-3j-1} -5\cdot2^{2s-2j-1}] &&\text{ if $ \; 0\leq j\leq s-2,\; k > s+m+ 1+j $}, \label{eq 3.21}\vspace{0.1 cm} \\
    & \nonumber \\
    \Gamma _{2}^{\Big[\substack{1 \\ 1  }\Big] \times (k-m -2s +2)}  & =  2^{2(k-m) -4s +4} -3\cdot2^{k-m -2s +2} +2  &&\text{if   \; $ j =  s-1,\;  k > 2s + m $}. \label{eq 3.22} \\
     & \nonumber 
  \end{align}
\end{thm}  
 \begin{thm}
\label{thm  3.9}  We have \vspace{0.1 cm} \\
  \begin{equation}
\label{eq 3.23} \Gamma _{i}^{\Big[\substack{s \\ s  }\Big] \times k} = 
\begin{cases}
1 & \text{if  } i = 0,\; k \geq 1,         \\
 21\cdot2^{3i-4}  - 3\cdot2^{2i-3}  & \text{if   } 1\leq i\leq s-1,\; k > i, \\
 3\cdot2^{k+s-1} + 21\cdot2^{3s-4} - 27\cdot2^{2s-3}  & \text{if   }  i = s ,\; k > s, \\
 21\cdot[ 2^{k - 2s +3i - 4}  + 2^{3i -4} - 5\cdot2^{4i-2s -5} ] & \text{if   } s+1\leq i\leq 2s-1,\; k > i, \\
 2^{2k + 2s  -2}  -3\cdot 2^{k+4s - 4}  + 2^{6s  -5} & \text{if   }  i = 2s ,\; k > 2s. 
  \end{cases}    
   \end{equation}
     \end{thm} 
     
 \begin{thm}
\label{thm  3.10}  We have \vspace{0.1 cm} \\
  \begin{equation}
\label{eq 3.24} \Gamma _{i}^{\Big[\substack{s \\ s  }\Big] \times i} = 
\begin{cases}
2^{2s +2i -2}  - 3\cdot2^{3i-4} + 2^{2i-3} & \text{if   } 1\leq i\leq s,  \\
2^{2s +2i -2}  - 3\cdot2^{3i-4} + 2^{4i -2s -5} & \text{if   } s+1\leq i\leq 2s. 
  \end{cases}    
   \end{equation} 
   \end{thm} 
    
 \begin{thm}
\label{thm 3.11}  We have \vspace{0.1 cm} \\
  \begin{equation}
\label{eq 3.25} \Gamma _{i}^{\Big[\substack{s \\ s +1 }\Big] \times k} = 
\begin{cases}
1 & \text{if  } i = 0,\; k \geq 1,         \\
 21\cdot2^{3i-4}  - 3\cdot2^{2i-3}  & \text{if   } 1\leq i\leq s-1,\; k > i, \\
 2^{k+s-1} + 21\cdot2^{3s-4} - 11\cdot2^{2s-3}  & \text{if   }  i = s ,\; k > s, \\
 11\cdot 2^{k+s-1} + 21\cdot2^{3s- 1} - 53 \cdot2^{2s - 1}  & \text{if   }  i = s +1 ,\; k > s +1, \\
 21\cdot[ 2^{k - 2s +3i -  5}  + 2^{3i -4} - 5\cdot2^{4i-2s - 6} ] & \text{if   } s+ 2\leq i\leq 2s ,\; k > i, \\
 2^{2k + 2s  - 1}  -3\cdot 2^{k+4s - 2}  + 2^{6s  - 2} & \text{if   }  i = 2s +1 ,\; k > 2s +1.
  \end{cases}    
   \end{equation}
     \end{thm} 
     
 \begin{thm}
\label{thm  3.12}  We have  \vspace{0.1 cm} \\
  \begin{equation}
\label{eq 3.26} \Gamma _{i}^{\Big[\substack{s \\ s +1 }\Big] \times i} = 
\begin{cases}
2^{2s +2i - 1}  - 3\cdot2^{3i-4} + 2^{2i-3} & \text{if   } 1\leq i\leq s +1, \\
2^{2s +2i -2}  - 3\cdot2^{3i-4} + 2^{4i -2s - 6} & \text{if   } s+2\leq i\leq 2 s +1.
  \end{cases}    
   \end{equation} 
   \end{thm} 
    
 \begin{thm}
\label{thm  3.13}  We have for $ m\geq 2 $ \vspace{0.1 cm} \\
  \begin{equation}
\label{eq 3.27} \Gamma _{i}^{\Big[\substack{s \\ s + m }\Big] \times k} = 
\begin{cases}
1 & \text{if  } i = 0,\; k \geq 1,         \\
 21\cdot2^{3i-4}  - 3\cdot2^{2i-3}  & \text{if   } 1\leq i\leq s-1,\; k > i, \\
 2^{k+s-1} + 21\cdot2^{3s-4} - 11\cdot2^{2s-3}  & \text{if   }  i = s ,\; k > s, \\
 3\cdot2^{k-s+2i-3} + 21\cdot[2^{3i-4} - 2^{3i -s-4}] & \text{if   }  s+1 \leq i\leq s+m-1 ,\; k > i,\\
 11\cdot 2^{k+s+2m-3} + 21\cdot2^{3s +3m-4} - 53 \cdot2^{2s  +3m- 4}  & \text{if   }  i = s +m ,\; k > s +m, \\
 21\cdot[ 2^{k - 2s +3i - m -4}  + 2^{3i -4} - 5\cdot2^{4i-2s - m-5} ] & \text{if   } s +m +1\leq i\leq 2s +m -1 ,\; k > i, \\
 2^{2k + 2s  +m -2}  -3\cdot 2^{k+4s +2m-4}  + 2^{6s  + 3m -5} & \text{if   }  i = 2s + m ,\; k > 2s + m.
  \end{cases}    
   \end{equation}
     \end{thm} 
     
 \begin{thm}
\label{thm  3.14}  We have for $ m\geq 2  $ \vspace{0.1 cm} \\
  \begin{equation}
\label{eq 3.28} \Gamma _{i}^{\Big[\substack{s \\ s +m }\Big] \times i} = 
\begin{cases}
2^{2s +2i +m -2}  - 3\cdot2^{3i-4} + 2^{2i-3} & \text{if   } 1\leq i\leq s +1, \\
2^{2s +2i  +m -2 }  - 3\cdot2^{3i-4} + 2^{3i - s - 4} & \text{if   } s+2\leq i\leq s+m+1,\\
2^{2s +2i  +m -2 }  - 3\cdot2^{3i-4} + 2^{4i - 2s - m -5 } & \text{if   } s+m+2\leq i\leq 2s+m. 
  \end{cases}    
   \end{equation} 
   \end{thm}

 \begin{thm}
\label{thm 3.15}
 We denote by  $ R_{q}(k,s,m) $ the number of solutions \\
 $(Y_1,Z_1,U_{1}, \ldots,Y_q,Z_q,U_{q}) $  of the polynomial equations
   \[\left\{\begin{array}{c}
 Y_{1}Z_{1} +Y_{2}Z_{2}+ \ldots + Y_{q}Z_{q} = 0,  \\
   Y_{1}U_{1} + Y_{2}U_{2} + \ldots  + Y_{q}U_{q} = 0,
 \end{array}\right.\]
  satisfying the degree conditions \\
                   $$  degY_i \leq k-1 , \quad degZ_i \leq s-1 ,\quad degU_{i}\leq s+m-1 \quad for \quad 1\leq i \leq q. $$ \\                           
Then \\
\begin{align}
  R_{q}(q,k,s,m) & =  \int_{\mathbb{P}\times \mathbb{P}} g_{k,s,m}^{q}(t,\eta )dtd\eta  \label{eq 3.29}\\
&   = 2^{(2s+m+k)(q-1)}\cdot 2^{-k+2}\cdot \sum_{i = 0}^{\inf(2s+m,k)} \Gamma _{i}^{\Big[\substack{s  \\ s+m }\Big] \times k}\cdot2^{- qi}. \nonumber 
\end{align}
  \end{thm}
\begin{example} Let q = 3, k = 4, s = 3, m = 2. Then 
   \begin{equation*}
   \Gamma _{i}^{\Big[\substack{3 \\  3 + 2 }\Big] \times 4} = 
\begin{cases}
1 & \text{if  } i = 0,        \\
 9  & \text{if   } i = 1, \\
 78   & \text{if   }  i = 2, \\
 648   & \text{if   }  i = 3,   \\
 15648 &  \text{if   }  i = 4. 
  \end{cases}    
   \end{equation*}\vspace{0.1 cm}\\
Hence  the number $ R_{3}(4,3,2) $ of solutions \\
 $(Y_1,Z_1,U_{1}, Y_2,Z_2,U_{2},Y_3,Z_3,U_{3}) $  of the polynomial equations
   \[\left\{\begin{array}{c}
 Y_{1}Z_{1} +Y_{2}Z_{2} + Y_{3}Z_{3} = 0,  \\
   Y_{1}U_{1} + Y_{2}U_{2}  + Y_{3}U_{3} = 0,
 \end{array}\right.\]
  satisfying the degree conditions \\
                   $$  degY_i \leq 3 , \quad degZ_i \leq 2 ,\quad degU_{i}\leq 4 \quad for \quad 1\leq i \leq 3 $$ \\                           
  is equal to \vspace{0.1 cm}\\
      \begin{align*}
   &  \int_{\mathbb{P}\times \mathbb{P}} g_{4,3,2}^{3}(t,\eta )dtd\eta 
 =  2^{22}\cdot \sum_{i = 0}^{4} \Gamma _{i}^{\Big[\substack{3 \\ 3 +2 }\Big] \times 4}\cdot2^{-3i}  \\
 & =  2^{22}\cdot[1 + 9\cdot2^{-3} + 78\cdot2^{-6} + 648\cdot2^{-9} + 15648\cdot2^{-12}] \\
&  = 35356672.
\end{align*}
\end{example}
\begin{example} Let  q = 4, k = 6, s = 5, m = 0. Then  
   \begin{equation*}
   \Gamma _{i}^{\Big[\substack{5 \\  5}\Big] \times 6} = 
\begin{cases}
1 & \text{if  } i = 0,        \\
 9  & \text{if   } i = 1, \\
 78   & \text{if   }  i = 2, \\
 648   & \text{if   }  i = 3,   \\
 5280 &  \text{if  }  i = 4, \\
 42624 & \text{if  } i  = 5, \\
  999936  &  \text{if   }  i = 6.
  \end{cases}    
   \end{equation*}\vspace{0.1 cm}\\
   Hence  the number  $  R_{4}(6,5,0) $ of solutions  \vspace{0.1 cm}\\
 $(Y_1,Z_1,U_{1}, Y_2,Z_2,U_{2},Y_3,Z_3,U_{3}, Y_4, Z_4,U_{4} ) $  of the polynomial equations
   \[\left\{\begin{array}{c}
 Y_{1}Z_{1} +Y_{2}Z_{2} + Y_{3}Z_{3} + Y_{4}Z_{4} = 0,  \\
   Y_{1}U_{1} + Y_{2}U_{2}  + Y_{3}U_{3} + Y_{4}U_{4}  = 0,
 \end{array}\right.\]
  satisfying the degree conditions \\
                   $$  degY_i \leq 5 , \quad degZ_i \leq 4 ,\quad degU_{i}\leq 4 \quad for \quad 1\leq i \leq 4 $$ \\                           
  is equal to \vspace{0.1 cm}\\
     \begin{align*}
&   \int_{\mathbb{P}\times \mathbb{P}} g_{4,5,0}^{4}(t,\eta )dtd\eta \\
&  =  2^{44}\cdot \sum_{i = 0}^{6} \Gamma _{i}^{\Big[\substack{5 \\  5 }\Big] \times 6}\cdot2^{-4i}  \\
 & =  2^{44}\cdot[1 + 9\cdot2^{-4} + 78\cdot2^{-8} + 648\cdot2^{-12} + 5280 \cdot2^{-16} + 42624\cdot 2^{-20} + 999936\cdot 2^{-24}] \\
&  = 37014016\cdot 2^{20}.
\end{align*}
\end{example}  
 \begin{example} The fraction of square  double persymmetric  $\left[s\atop s+m\right]\times (2s+m) $ matrices which are 
 invertible is equal to $ \displaystyle \frac{ \Gamma _{2s+m}^{\left[s\atop s+m\right] \times (2s+m)}}{\sum_{i = 0}^{2s+m}\Gamma _{i}^{\left[s\atop s+m\right] \times (2s+m)} } = \frac{3}{8}. $
   \end{example}

    \section{\textbf{EXPONENTIAL SUMS FORMULAS ON  $ \mathbb{K}^2 $}}
   \label{sec 4}
 In this section we compute exponential quadratic sums in  $\mathbb{P}^2 $ and show that they only depend on rank properties of some 
 associated double persymmetric matrices.  
 The following propositions are proved in [2]. \\
 
\begin{prop}
\label{prop 4.1}
The following holds :
\begin{itemize}
\item  For every rational integer j, the measure of  $\mathbb{P}_{j} \; is \; 2^{-j}.$
\item   For every \; $ A \in \mathbb{F}_{2}[T]  $,  \;E(A) = 1. 
\item   For $ u \in \mathbb{K} \quad \nu(u)\geq  2 \Rightarrow  $ \; E(u) = 1.
\end{itemize}
\end{prop}

\begin{prop}
\label{prop 4.2}
Let j be a rational integer and \; $ u \in \mathbb{K} $. Then
  \[ \int_{\mathbb{P}_{j}}E(ut)dt   =  \left\{ \begin{array}{ccc}
                         2^{-j}    &  if  &  \nu(u) > - j,    \\
                          0  &       &  otherwise.
                          \end{array}\right.\]
\end{prop}

 \begin{prop}
\label{prop 4.3}
 Let j be a rational integer and let \; $ u \in \mathbb{K}.$  Then
    \[ \sum_{deg B\leq  j}E(Bu)  =  \left\{ \begin{array}{ccc}
                         2^{j+1}    &  if  &   \nu\left(\left\{u\right\}\right) > j +1,  \\
                          0  &       &  otherwise.
                          \end{array}\right.\]
\end{prop}

\begin{lem}
\label{lem 4.4}
Let $Y \in\mathbb{F}_{2}[T] $ and a a rational integer $ \geq 1 $. Then \\
  \[ \sum_{deg Z\leq a -1} E(tYZ)   =  \left\{ \begin{array}{ccc}
                         2^{a}    &  if  &   \nu\left(\left\{tY\right\}\right) > a, \\
                          0  &       &  otherwise.
                          \end{array}\right.\]
\end{lem}
 \begin{proof}
The Lemma follows from the Proposition \ref{prop 4.3} with u = tY and j = a -1.
\end{proof}

\begin{lem}
\label{lem 4.5}
Let $ t \in\mathbb{P}\; and \; Y \in \mathbb{F}_{2}[T]  ,\; deg Y\leq b-1 $. Then \\
$$\nu\left(\left\{tY\right\}\right) > a
 \Longleftrightarrow  Y\in \ker D_{a\times b}(t). $$.
\end{lem}
 \begin{proof}
 Let  $ Y=\sum_{j=0}^{b -1}\gamma _{j}T^j     \; , \gamma _{j}\in \mathbb{F}_{2} $. Then\\
 $$ tY =(\sum_{i\geq 1}\alpha _{i}T^{-i})(\sum_{j=0}^{b-1} \gamma  _{j}  T^{j})
  = \sum_{i\geq 1}\sum_{j=0}^{b-1}\alpha_{i}\gamma_{j}T^{-(i-j)} $$ and  \\
 $$\left\{tY\right\} = \left(\sum_{i=1}^{b}
 \alpha_{i}\gamma_{i-1}\right)T^{-1}
+\left(\sum_{i=2}^{b+1}\alpha_{i}\gamma_{i-2}\right)T^{-2} +
\ldots + \left(\sum_{i= s}^{b+a -1}\alpha_{i}\gamma_{i-a}\right)T^{- a} + \ldots $$
Therefore $  \mathcal{\nu}(\left\{tY\right\}) \geq a  $  if and only if 
  \def\BA{ \begin{array}{ccccccc}}
  \def\BB{ \begin{array}{c}}
  \def\EA{\end{array}}
  \begin{displaymath}
\left[\BA
\alpha  _{1} & \alpha_{2} &\alpha _{3} & \alpha_{4} &\alpha _{5}  & \ldots  & \alpha_{b}\\
\alpha _{2} &\alpha _{3} &\alpha _{4} & \alpha_{5} & \alpha_{6}  & \ldots  & \alpha_{b+1} \\
 \vdots & \vdots & \vdots &\vdots & \vdots & \ldots                   & \vdots \\
 \vdots & \vdots & \vdots & \vdots & \vdots  & \ldots                   &\vdots \\
 \alpha _{a} &\alpha _{a+1} & \alpha_{a+2} &\alpha _{a+3} & \alpha_{a+4} &\ldots  &\alpha _{b+a -1} \\  \EA\right ]
 \left[ \BB \gamma _{0}\\ \gamma _{1}\\ \gamma _{2}\\ \gamma _{3}\\ \gamma _{4}
  \\ \vdots  \\\gamma _{b-1} \EA\right] =
  \left[ \BB 0 \\ 0 \\ 0 \\  0 \\ 0  \\ \vdots \\  0
    \EA\right]
 \end{displaymath} \\
\begin{displaymath} \Longleftrightarrow
D_{a\times b }(t) \left(\begin{array}{c}
\gamma _{0}\\
\gamma _{1}\\
\vdots \\
\gamma _{b -1}\\
\end{array}\right) =\left(\begin{array}{c}
0\\
0\\
\vdots\\
0\\
\end{array}\right)\Longleftrightarrow
Y= \left(\begin{array}{c}
\gamma _{0}\\
\gamma _{1}\\
\vdots \\
\gamma _{b -1}\\
\end{array}\right)   \in Ker \, D_{a\times b}(t).
\end{displaymath}
\end{proof}

\begin{lem}
\label{lem 4.6}
Let $ (t,\eta ) \in\mathbb{P}\times \mathbb{P}. $ Then \\
\begin{equation}
\label{eq 4.1}
  g(t,\eta ) =  \sum_{deg Y\leq k-1}\sum_{deg Z \leq  s-1}E(tYZ)\sum_{deg U \leq s+m-1}E(\eta YU) =
   2^{2s+m+k- r( D^{\left[\stackrel{s}{s+m}\right] \times k }(t,\eta  ) )}. 
\end{equation}
  \end{lem}
 \begin{proof}
 By lemma \ref{lem 4.4} and lemma \ref{lem 4.5} we obtain with $ a\rightarrow s \; and \; b\rightarrow k $
$$ \sum_{deg Z \leq  s-1}E(tYZ) = 2^{s}\quad  if\; and\; only\; if \quad  Y\in \ker D_{s\times k}(t). $$ \\
Similarly with $ t\rightarrow\eta, \;  a\rightarrow s+m \; and \; b\rightarrow k $ we get \\
$$ \sum_{deg U \leq  s+m-1}E(\eta YU) = 2^{s+m}\quad  if\; and\; only\; if \quad  Y\in \ker D_{(s+m)\times k}(\eta ). $$ \\
We then deduce \\
$$\sum_{deg Z \leq  s-1}E(tYZ)\sum_{deg U \leq s+m-1}E(\eta YU) = 2^{2s+m}\quad  if\; and\; only\; if \quad
Y\in \ker D^{\left[\stackrel{s}{s+m}\right] \times k }(t,\eta  )  $$\\
and we obtain \\
 \begin{equation*}
  \sum_{deg Y\leq k-1}\sum_{deg Z \leq  s-1}E(tYZ)\sum_{deg U \leq s+m-1}E(\eta YU) 
   =2^{2s + m }\sum_{{deg Y\leq k-1\atop Y\in \ker D^{\left[\stackrel{s}{s+m}\right] \times k }(t,\eta  ) }}1
 = 2^{2s + m }\cdot2^{k - r( D^{\left[\stackrel{s}{s+m}\right] \times k }(t,\eta  ) ) }.
\end{equation*}
 \end{proof}
\begin{lem}
\label{lem 4.7}
Let $ (t,\eta ) \in\mathbb{P}\times \mathbb{P} $ and set \\
$$  g_{1}(t,\eta ) =  \sum_{deg Y\leq k-1}\sum_{deg Z = s-1}E(tYZ)\sum_{deg U \leq s+m-1}E(\eta YU). $$
  Then
 \begin{equation}
 \label{eq 4.2}
g_{1}(t,\eta )  =   \begin{cases}
 2^{2s +m +k-1 - r( D^{\left[\stackrel{s}{s+m}\right] \times k }(t,\eta  ) ) }  & \text{if }
 r( D^{\left[\stackrel{s-1}{s+m}\right] \times k }(t,\eta  ) )    =  r( D^{\left[\stackrel{s}{s+m}\right] \times k }(t,\eta  ) ), \\
  0  & \text{otherwise }.
    \end{cases}
\end{equation} 
\end{lem}
\begin{proof}We have 
$$ g_{1}(t,\eta ) =  \sum_{deg Y\leq k-1}\sum_{deg Z \leq  s-1}E(tYZ)\sum_{deg U \leq s+m-1}E(\eta YU) -
 \sum_{deg Y\leq k-1}\sum_{deg Z \leq  s - 2}E(tYZ)\sum_{deg U \leq s+m-1}E(\eta YU). $$
 By  \eqref{eq 4.1} we obtain \\
 $$  g_{1}(t,\eta ) = 2^{2s+m+k- r( D^{\left[\stackrel{s}{s+m}\right] \times k }(t,\eta  ) )} -
   2^{2s+m+k-1 - r( D^{\left[\stackrel{s-1}{s+m}\right] \times k }(t,\eta  ) )}. $$
Now either $$ r( D^{\left[\stackrel{s}{s+m}\right] \times k }(t,\eta  ) ) = 
r( D^{\left[\stackrel{s-1}{s+m}\right] \times k }(t,\eta  ) ) $$
or  $$ r( D^{\left[\stackrel{s}{s+m}\right] \times k }(t,\eta  ) ) = 
r( D^{\left[\stackrel{s-1}{s+m}\right] \times k }(t,\eta  ) ) +1. $$\vspace{0.1 cm}\\
Lemma \ref{lem 4.7} follows.
\end{proof}

  \begin{lem}
\label{lem 4.8}
Let $ (t,\eta ) \in\mathbb{P}\times \mathbb{P} $ and set \\
$$  g_{2}(t,\eta ) =  \sum_{deg Y\leq k-1}\sum_{deg Z \leq  s-1}E(tYZ)\sum_{deg U = s-1+m}E(\eta YU), $$ \\
$$  f_{1}(t,\eta ) =  \sum_{deg Y\leq k-1}\sum_{deg Z = s-1}E(tYZ)\sum_{deg U \leq  s-2 +m}E(\eta YU), $$ \\
$$  f_{2}(t,\eta ) =  \sum_{deg Y\leq k-1}\sum_{deg Z \leq  s-2}E(tYZ)\sum_{deg U =  s -2 +(m+1)}E(\eta YU). $$ \\
  Then
 \begin{equation}
 \label{eq 4.3}
g_{2}(t,\eta )  =   \begin{cases}
 2^{2s +m +k-1 - r( D^{\left[\stackrel{s}{s+m}\right] \times k }(t,\eta  ) ) }  & \text{if }
 r( D^{\left[\stackrel{s}{s+m}\right] \times k }(t,\eta  ) )    =  r( D^{\left[\stackrel{s}{s+(m-1)}\right] \times k }(t,\eta  ) ), \\
  0  & \text{otherwise},
    \end{cases}
\end{equation}
 \begin{equation}
 \label{eq 4.4}
f_{1}(t,\eta )  =   \begin{cases}
 2^{2s +m +k - 2 - r( D^{\left[\stackrel{s-1}{s -1+m}\right] \times k }(t,\eta  ) ) }  & \text{if }
 r( D^{\left[\stackrel{s-1}{s -1+m}\right] \times k }(t,\eta  ) )    =  r( D^{\left[\stackrel{s}{s+(m-1)}\right] \times k }(t,\eta  ) ), \\
  0  & \text{otherwise},
    \end{cases}
\end{equation}
 \begin{equation}
 \label{eq 4.5}
f_{2}(t,\eta )  =   \begin{cases}
 2^{2s +m +k-2 - r( D^{\left[\stackrel{s-1}{s-1+m}\right] \times k }(t,\eta  ) ) }  & \text{if }
 r( D^{\left[\stackrel{s-1}{s-1+m}\right] \times k }(t,\eta  ) )    =  r( D^{\left[\stackrel{s-1}{s-1+(m+1)}\right] \times k }(t,\eta  ) ), \\
  0  & \text{otherwise}.
    \end{cases}
\end{equation}
\end{lem}
\begin{proof}
Similarly to the proof of Lemma \ref{lem 4.7}
\end{proof}

 \begin{lem}
\label{lem 4.9}
Set   $(t, \eta ) = (\sum_{j\geq 1}\alpha  _{j}T^{-j},\sum_{j\geq 1}\beta _{j}T^{-j}) \in \mathbb{P}\times\mathbb{P}
\;  and \;
Y = \sum_{j=1}^{k-1}\delta _{j}T^j\in \mathbb{F}_{2}[T] ,\; degY\leq k-1.  $ \\
  Then
  \begin{equation}
  \label{eq  4.6}
E(Y(tT^{s-1}+\eta T^{s+m-1})) = \begin{cases}
1 & \text{if           }                      \sum_{j=1}^{k-1}(\alpha _{s+j}+\beta _{s+m+j})\delta _{j} = 0, \\
-1 & \text{if          }                       \sum_{j=1}^{k-1}(\alpha _{s+j}+\beta _{s+m+j})\delta _{j} = 1.
\end{cases}
\end{equation}
  \end{lem}
\begin{proof} We have \\
\begin{align*}
Y(tT^{s-1} + \eta T^{s+m-1}) & = \sum_{j = 0}^{k-1}\delta _{j}T^{j}\cdot
\left[\sum_{i \geq 1}\alpha _{i}T^{-i}\cdot T^{s-1} +
 \sum_{l\geq 1}\beta _{l}T^{-l}\cdot T^{s +m -1}\right] \\
 & = \sum_{j = 0}^{k-1}\left[\sum_{i\geq 1}\alpha _{i}\delta _{j}T^{s-1-i+j} +
  \sum_{l\geq 1}\beta _{l}\delta _{j}T^{s+m-1-l+j} \right]\\
  & = \ldots +\left[\sum_{j = 0}^{k-1}(\alpha _{s+j}+ \beta _{s+m+j})\delta _{j}\right]\cdot T^{-1} + \ldots
\end{align*}\vspace{0.1 cm}
Now \begin{align*}
\nu\left(\left\{Y(tT^{s-1} + \eta T^{s+m-1})\right\}\right) > 1
  & \Longleftrightarrow \sum_{j = 0}^{k -1} ( \alpha _{s+j}+ \beta _{s+m+j})\delta _{j}  = 0,  \\
\nu\left(\left\{Y(tT^{s-1} + \eta T^{s+m-1})\right\}\right) = 1
& \Longleftrightarrow  \sum_{j = 0}^{k -1} ( \alpha _{s+j}+ \beta _{s+m+j})\delta _{j}  = 1, 
\end{align*}
which proves Lemma \ref{lem 4.9}.
\end{proof}
\begin{lem}
\label{lem 4.10}
Let $ (t,\eta ) \in\mathbb{P}\times \mathbb{P} $ and set \\
$$  h(t,\eta ) =  \sum_{deg Y\leq k-1}\sum_{deg Z = s-1}E(tYZ)\sum_{deg U = s+m-1}E(\eta YU). $$
  Then
 \begin{equation}
 \label{eq 4.7}
h(t,\eta )  =   \begin{cases}
 2^{2s +m +k-2 - r( D^{\left[\stackrel{s-1}{s-1+m}\right] \times k }(t,\eta  ) ) }  & \text{if }
 r( D^{\left[\stackrel{s-1}{s-1+m}\right] \times k }(t,\eta  ) )    =  r( D^{\left[\stackrel{s-1}{\stackrel{s+m-1}{\alpha_{s -} + \beta_{s+m -} }}\right] \times k }(t,\eta )), \\
  0  & \text{otherwise }.
    \end{cases}
\end{equation} 
\end{lem}
\begin{proof}We have by the proof of lemma \ref{lem 4.6} and equation \eqref{eq 4.6}
\begin{align*}
  h(t,\eta ) &  =  \sum_{deg Y\leq k-1}\sum_{deg Z = s-1}E(tYZ)\sum_{deg U = s+m-1}E(\eta YU)  \\
   & =   \sum_{deg Y\leq k-1}E(Y(tT^{s-1}+\eta T^{s+m-1})) \sum_{deg Z \leq  s-2}E(tYZ)\sum_{deg U \leq  s+m-2}E(\eta YU)  \\
   & = 2^{2s +m-2}\sum_{{deg Y\leq k-1\atop Y\in \ker D^{\left[\stackrel{s-1}{s-1+m}\right] \times k }}} E(Y(tT^{s-1}+\eta T^{s+m-1}))  \\
   & = 2^{2s +m-2}\left[  \sum_{\substack{deg Y\leq k-1  \\
 Y\in \ker D^{\big[\stackrel{s-1}{s-1+m}\big] \times k }(t,\eta ) \\
 \sum_{j=1}^{k-1}(\alpha _{s+j}+\beta _{s+m+j})\delta _{j} = 0  }}1  -   \sum_{\substack{deg Y\leq k-1 \\
 Y\in \ker D^{\big[\stackrel{s-1}{s-1+m}\big] \times k }(t,\eta ) \\
 \sum_{j=1}^{k-1}(\alpha _{s+j}+\beta _{s+m+j})\delta _{j} = 1  }}1 \right]  \\
  & =  2^{2s +m-2}\left[  \sum_{deg Y\leq k-1 \atop
 Y\in \ker D^{\big[\stackrel{s-1}{\stackrel{s+m-1}
{\alpha_{s -} + \beta_{s+m -} }}\big] \times k }(t,\eta  )} 1  - \left( \sum_{ { deg Y\leq k-1 \atop
 Y\in \ker D^{\big[\stackrel{s-1}{s+m-1}\big] \times k }(t,\eta  )}} 1 -  \sum_{deg Y\leq k-1 \atop
 Y\in \ker D^{\big[\stackrel{s-1}{\stackrel{s+m-1}
{\alpha_{s -} + \beta_{s+m -} }}\big] \times k }(t,\eta  )} 1\right) \right]  \\
& =  2^{2s +m-2}\left[2\cdot  \sum_{deg Y\leq k-1 \atop
 Y\in \ker D^{\big[\stackrel{s-1}{\stackrel{s+m-1}
{\alpha_{s -} + \beta_{s+m -} }}\big] \times k }(t,\eta  )} 1  - \sum_{ { deg Y\leq k-1 \atop
 Y\in \ker D^{\big[\stackrel{s-1}{s+m-1}\big] \times k }(t,\eta  )}} 1 \right]  \\
 & = 2^{2s +m-2}\left[2\cdot 2^{k -  r( D^{\left[\stackrel{s-1}{\stackrel{s+m-1}{\alpha_{s -} + \beta_{s+m -} }}\right] \times k }(t,\eta  )  )} 
 - 2^{ k -   r( D^{\left[\stackrel{s-1}{s-1+m}\right] \times k }(t,\eta  ) ) }\right]
 \end{align*}
 which proves Lemma \ref{lem 4.10}.
\end{proof}

\begin{lem}
\label{lem 4.11}
Let $ (t,\eta ) \in\mathbb{P}\times \mathbb{P} $ and set \\

\begin{align*}
  h(t,\eta ) &  =  \sum_{deg Y\leq k-1}\sum_{deg Z = s-1}E(tYZ)\sum_{deg U = s+m-1}E(\eta YU),  \\
   v(t,\eta ) &  =  \sum_{deg Y\leq k-1}\sum_{deg Z \leq  s-2}E(tYZ)\sum_{deg U \leq  s+m-2}E(\eta YU),  \\
 f_{1}(t,\eta ) &  =  \sum_{deg Y\leq k-1}\sum_{deg Z = s-1}E(tYZ)\sum_{deg U \leq  s+m-2}E(\eta YU),  \\
   f_{2}(t,\eta ) &  =  \sum_{deg Y\leq k-1}\sum_{deg Z \leq  s-2}E(tYZ)\sum_{deg U =  s+m-1}E(\eta YU). 
\end{align*}
Then we have
\begin{equation}
\label{eq 4.8}
 h^{q}(t,\eta ) = h(t,\eta )\cdot v^{q-1}(t,\eta )\; for\; every\; rational\; integer\; q\geq 2, 
\end{equation}
\begin{equation}
\label{eq 4.9}
 f_{1}^{q}(t,\eta ) = f_{1}(t,\eta )\cdot v^{q-1}(t,\eta )\; for\; every\; rational\; integer\; q\geq 2, 
\end{equation}
\begin{equation}
\label{eq 4.10}
 f_{2}^{q}(t,\eta ) = f_{2}(t,\eta )\cdot v^{q-1}(t,\eta )\; for\; every\; rational\; integer\; q\geq 2. 
\end{equation}
\end{lem}
\begin{proof} We get obviously\vspace{0.1 cm}
\Small
$$ h^{2}(t,\eta ) 
  =  \big[ 2^{2s +m-2}\sum_{{deg Y_{1}\leq k-1\atop Y_{1}\in \ker D^{\left[\stackrel{s-1}{s-1+m}\right] \times k }}}
  E(Y_{1}(tT^{s-1}+\eta T^{s+m-1}))\big]\cdot \big[ 2^{2s +m-2}\sum_{{deg Y_{2}\leq k-1\atop Y_{2}\in \ker D^{\left[\stackrel{s-1}{s-1+m}\right] \times k }}}
  E(Y_{2}(tT^{s-1}+\eta T^{s+m-1}))\big]. $$ \\
  \normalsize
 We set  \[\left\{\begin{array}{cc}
Y_{1} + Y_{2} = Y_{3}  &  deg Y_{3}\leq k-1, \\
              Y_{2} = Y_{4}   &   deg Y_{4}\leq k-1.
\end{array}\right.\]\\
Then we obtain\\
  \begin{align*} h^{2}(t,\eta ) & = 
  2^{(2s+m-2)2}\sum_{{deg Y_{3}\leq k-1\atop Y_{3}\in \ker D^{\left[\stackrel{s-1}{s-1+m}\right] \times k }}}
  \sum_{{deg Y_{4}\leq k-1\atop Y_{4}\in \ker D^{\left[\stackrel{s-1}{s-1+m}\right] \times k }}}
   E((Y_{4}+ Y_{3}+Y_{4})(tT^{s-1}+\eta T^{s+m-1})) \\
  & = \big[ 2^{2s +m-2}\sum_{{deg Y_{3}\leq k-1\atop Y_{3}\in \ker D^{\left[\stackrel{s-1}{s-1+m}\right] \times k }}}
  E(Y_{3}(tT^{s-1}+\eta T^{s+m-1}))\big]\cdot \big[ 2^{2s +m-2}\sum_{{deg Y_{4}\leq k-1\atop Y_{4}\in \ker D^{\left[\stackrel{s-1}{s-1+m}\right] \times k }}}
  1 \big] \\
  & =  h(t,\eta )\cdot v(t,\eta ).\\
  & 
 \end{align*}
 By recurrence on q we get \quad $ h^{q}(t,\eta ) = h(t,\eta )\cdot v^{q-1}(t,\eta ). $\vspace{0.1 cm}\\
 
  The proof  of \eqref{eq 4.9} and \eqref{eq 4.10} is similar to the proof of \eqref{eq 4.10}, that is 
 $$ f_{1}^{2}(t,\eta ) 
  =  \big[ 2^{2s +m-2}\sum_{{deg Y_{1}\leq k-1\atop Y_{1}\in \ker D^{\left[\stackrel{s-1}{s-1+m}\right] \times k }}}
  E(Y_{1}tT^{s-1})\big]\cdot \big[ 2^{2s +m-2}
  \sum_{{deg Y_{2}\leq k-1\atop Y_{2}\in \ker D^{\left[\stackrel{s-1}{s-1+m}\right] \times k }}}
  E(Y_{2}tT^{s-1})\big]. $$ \\
   We set  \[\left\{\begin{array}{cc}
Y_{1} + Y_{2} = Y_{3}  &  deg Y_{3}\leq k-1, \\
              Y_{2} = Y_{4}   &   deg Y_{4}\leq k-1.
\end{array}\right.\]\\
Then we obtain\\
  \begin{align*} f_{1}^{2}(t,\eta ) & = 
  2^{(2s+m-2)2}\sum_{{deg Y_{3}\leq k-1\atop Y_{3}\in \ker D^{\left[\stackrel{s-1}{s-1+m}\right] \times k }}}
  \sum_{{deg Y_{4}\leq k-1\atop Y_{4}\in \ker D^{\left[\stackrel{s-1}{s-1+m}\right] \times k }}}
   E( Y_{3}tT^{s-1}) \\
  & = \big[ 2^{2s +m-2}\sum_{{deg Y_{3}\leq k-1\atop Y_{3}\in \ker D^{\left[\stackrel{s-1}{s-1+m}\right] \times k }}}
  E(Y_{3}tT^{s-1})\big]\cdot \big[ 2^{2s +m-2}
  \sum_{{deg Y_{4}\leq k-1\atop Y_{4}\in \ker D^{\left[\stackrel{s-1}{s-1+m}\right] \times k }}}
  1 \big] \\
  & =  f_{1}(t,\eta )\cdot v(t,\eta ).
 \end{align*}
  By recurrence on q we get respectively
 \begin{align*}
  f_{1}^{q}(t,\eta ) &  =  f_{1}(t,\eta )\cdot v^{q-1}(t,\eta ), \\
  f_{2}^{q}(t,\eta ) &  =  f_{2}(t,\eta )\cdot v^{q-1}(t,\eta ). 
 \end{align*}
  \end{proof}
\begin{lem}
\label{lem 4.12}
Let $ (t,\eta ) \in\mathbb{P}\times \mathbb{P}\; and\; q\geq 2\;, then\; we\; have\; for \;1\leq i\leq q-1 $
 \begin{equation}
 \label{eq 4.11}
h^{i}f_{1}^{q-i}  = h^{i}f_{2}^{q-i} = \begin{cases}
 v^{q}  & \text{if }    f_{1}f_{2} \neq  0, \\
  0  & \text{otherwise},
    \end{cases}
\end{equation}that is \\
\begin{equation}
 \label{eq 4.12}
 h^{i}(t,\eta )f_{1}^{q-i}(t,\eta )  = h^{i}(t,\eta )f_{2}^{q-i}(t,\eta ) 
\end{equation}
 is equal to \vspace{0.1 cm}\\
  \begin{equation}
 \label{eq 4.13}
 \begin{cases}
 2^{q\cdot(2s +m +k-2 - r( D^{\left[\stackrel{s-1}{s+m-1}\right] \times k)}(t,\eta  ) )) }  & \text{if }
 r( D^{\left[\stackrel{s-1}{s -1+m}\right] \times k }(t,\eta  ) )    =  r( D^{\left[\stackrel{s}{s+(m-1)}\right] \times k }(t,\eta  ) ) = 
  r( D^{\left[\stackrel{s-1}{s+ m}\right] \times k }(t,\eta  ) ), \\
  0  & \text{otherwise}.
    \end{cases}
\end{equation} 
\end{lem}
\begin{proof}We have \\
\Small
\begin{align*}
f_{1}(t,\eta ) h(t,\eta ) 
&  =  \big[ 2^{2s +m-2}\sum_{{deg Y_{1}\leq k-1\atop Y_{1}\in \ker D^{\left[\stackrel{s-1}{s-1+m}\right] \times k }}}
  E(Y_{1}(tT^{s-1}))\big]\cdot \big[ 2^{2s +m-2}\sum_{{deg Y_{2}\leq k-1\atop Y_{2}\in \ker D^{\left[\stackrel{s-1}{s-1+m}\right] \times k }}}
  E(Y_{2}(tT^{s-1}+\eta T^{s+m-1}))\big] \\
  & =  2^{(2s+m-2)2}\sum_{{deg Y_{1}\leq k-1\atop Y_{1}\in \ker D^{\left[\stackrel{s-1}{s-1+m}\right] \times k }}}
  \sum_{{deg Y_{2}\leq k-1\atop Y_{2}\in \ker D^{\left[\stackrel{s-1}{s-1+m}\right] \times k }}}
   E(t(Y_{1}+ Y_{2}) T^{s-1}) E( \eta Y_{2}T^{s+m-1}). 
\end{align*}
  \normalsize
 We set  \[\left\{\begin{array}{cc}
Y_{1} + Y_{2} = Y_{3}  &  deg Y_{3}\leq k-1, \\
              Y_{2} = Y_{4}   &   deg Y_{4}\leq k-1.
\end{array}\right.\]\\
Then we obtain\\
  \begin{align*} f_{1}(t,\eta ) h(t,\eta )  & = 
  2^{(2s+m-2)2}\sum_{{deg Y_{3}\leq k-1\atop Y_{3}\in \ker D^{\left[\stackrel{s-1}{s-1+m}\right] \times k }}}
  \sum_{{deg Y_{4}\leq k-1\atop Y_{4}\in \ker D^{\left[\stackrel{s-1}{s-1+m}\right] \times k }}}
   E(tY_{3}T^{s-1})E(\eta Y_{4}T^{s+m-1}) \\
  & = \big[ 2^{2s +m-2}\sum_{{deg Y_{3}\leq k-1\atop Y_{3}\in \ker D^{\left[\stackrel{s-1}{s-1+m}\right] \times k }}}
  E(tY_{3}T^{s-1}) \big]\cdot \big[ 2^{2s +m-2}\sum_{{deg Y_{4}\leq k-1\atop Y_{4}\in \ker D^{\left[\stackrel{s-1}{s-1+m}\right] \times k }}}
  E(\eta Y_{4}T^{s+m-1}) \big] \\
  & =  f_{1}(t,\eta )\cdot f_{2}(t,\eta ). 
 \end{align*}\vspace{0.1 cm}
 In the same way we get  $$ f_{2}(t,\eta )\cdot h(t,\eta ) =  f_{1}(t,\eta )\cdot f_{2}(t,\eta ). $$ Then 
 
  \begin{equation}
  \label{eq 4.14}
 f_{1}(t,\eta )\cdot h(t,\eta ) = f_{2}(t,\eta )\cdot h(t,\eta ) =  f_{1}(t,\eta )\cdot f_{2}(t,\eta ).
\end{equation}
By respectively \eqref{eq 4.8}, \eqref{eq 4.9}  and  \eqref{eq 4.10} we deduce
 \begin{equation}
 \label{eq 4.15}
h = \begin{cases}
 v  & \text{if }     h \neq  0, \\
  0  & \text{otherwise},
    \end{cases}
\end{equation}
 \begin{equation}
 \label{eq 4.16}
f_{1} = \begin{cases}
 v  & \text{if }     f_{1}\neq  0, \\
  0  & \text{otherwise},
    \end{cases}
\end{equation}
 \begin{equation}
 \label{eq 4.17}
f_{2} = \begin{cases}
 v  & \text{if }     f_{2}\neq  0, \\
  0  & \text{otherwise}.
    \end{cases}
\end{equation}
Now by  \eqref{eq 4.14}, \eqref{eq 4.15}, \eqref{eq 4.16}   and  \eqref{eq 4.17}  we have for $ 1\leq i\leq q-1 $
 \begin{equation}
 \label{eq 4.18}
  h^{i}f_{1}^{q-i}  = h^{i}f_{2}^{q-i} = \begin{cases}
 v^{i}\cdot v^{q-i}= v^{q}  & \text{if   }  h\cdot f_{1} = h\cdot f_{2} =  f_{1}\cdot f_{2}\not\neq 0,    \\
  0  & \text{otherwise}.
    \end{cases}
\end{equation}
Finally  from  \eqref{eq 4.4}, \eqref{eq 4.5} and  \eqref{eq 4.18} we deduce  \eqref{eq 4.13}.
\end{proof}

 \begin{lem}
\label{lem 4.13}
 We denote by  $ R_{q}(k,s,m) $ the number of solutions \\
 $(Y_1,Z_1,U_{1}, \ldots,Y_q,Z_q,U_{q}) $  of the polynomial equations
   \[\left\{\begin{array}{c}
 Y_{1}Z_{1} +Y_{2}Z_{2}+ \ldots + Y_{q}Z_{q} = 0,  \\
   Y_{1}U_{1} + Y_{2}U_{2} + \ldots  + Y_{q}U_{q} = 0,
 \end{array}\right.\]
  satisfying the degree conditions \\
                   $$  degY_i \leq k-1 , \quad degZ_i \leq s-1 ,\quad degU_{i}\leq s+m-1 \quad for \quad 1\leq i \leq q. $$ \\                           
Then \\
\begin{align}
  R_{q}(k,s,m) & =  \int_{\mathbb{P}\times \mathbb{P}} g_{k,s,m}^{q}(t,\eta )dtd\eta  \label{eq 4.19}\\
&   = 2^{(2s+m+k)(q-1)}\cdot 2^{-k+2}\cdot \sum_{i = 0}^{\inf(2s+m,k)} \Gamma _{i}^{\Big[\substack{s  \\ s+m }\Big] \times k}\cdot2^{- qi}. \nonumber 
\end{align}
\end{lem}
\begin{proof}
We have by \eqref{eq 4.1} observing that $ g(t,\eta ) $ is constant on cosets of $ \mathbb{P}_{k+s-1}\times \mathbb{P}_{k+s+m-1} $ \\
 \begin{align*}
\int_{\mathbb{P}^{2}} g^{q}(t,\eta )dtd\eta  & =
 \sum_{(t,\eta )\in \mathbb{P}/\mathbb{P}_{k+s-1}\times \mathbb{P}/\mathbb{P}_{k+s+m-1}}
 2^{q(2s+m+k- r(D^{\big[\stackrel{s}{s+m}\big] \times k }(t,\eta ) )}\int_{\mathbb{P}_{k+s-1}}dt \int_{\mathbb{P}_{k+s+m-1}}dt d\eta \\
 & = \sum_{i = 0}^{\inf(2s+m,k)} \sum_{(t,\eta )\in \mathbb{P}/\mathbb{P}_{k+s-1}\times \mathbb{P}/\mathbb{P}_{k+s+m-1}}
  2^{q(2s+m+k- i) }\int_{\mathbb{P}_{k+s-1}}dt \int_{\mathbb{P}_{k+s+m-1}}dt d\eta \\
  & = \sum_{i = 0}^{\inf(2s+m,k)}\Gamma _{i}^{\Big[\substack{s  \\ s+m }\Big] \times k} 
  \cdot2^{-qi}\cdot2^{q(2s+m+k)}\cdot2^{-(k+s-1)}\cdot2^{-(k+s+m-1)}.
 \end{align*}
\end{proof}

 \section{\textbf{A RECURRENT FORMULA FOR THE NUMBER  $\Gamma_{i}^{\left[s\atop s+m\right]\times k}$ OF RANK i MATRICES OF THE FORM
  $ \left[{A\over B}\right], $ WHERE A IS A   $ s \times k $ PERSYMMETRIC MATRIX  AND  B a  $ (s+m) \times k $ PERSYMMETRIC  MATRIX WITH  ENTRIES  IN  $ \mathbb{F}_{2} $}}
   \label{sec 5}
 In this section we establish a recurrent formula for the number of rank i matrices of the form $\left[A\over B\right], $ where A and B are persymmetric.\\
   
  \begin{lem}
\label{lem 5.1} 
  We have
  \begin{align}
    \int_{\mathbb{P}^{2}}g_{1}^{q}(t,\eta )dtd\eta & =
  \int_{  \left\{(t,\eta ) \in \mathbb{P}^{2}\mid  h(t,\eta )\neq 0 \right\}}v^{q}(t,\eta ) dt d\eta  +
  \int_{  \left\{(t,\eta ) \in \mathbb{P}^{2}\mid  f_{1}(t,\eta )\neq 0 \right\}}v^{q}(t,\eta ) dt d\eta  \label{eq 5.1} \\
   & +  (2^{q} - 2)\cdot  \int_{  \left\{(t,\eta ) \in \mathbb{P}^{2}\mid  f_{1}(t,\eta )\cdot f_{2}(t,\eta )\neq 0 \right\}}v^{q}(t,\eta ) dt d\eta.  \nonumber
   \end{align}
 \end{lem}
\begin{proof} By the binomial theorem we obtain \\
\begin{equation}
\label{eq 5.2}
g_{1}^{q} = (h + f_{1})^{q} = h^{q} + f_{1}^{q} + \sum_{i = 1}^{q -1} \binom{q}{i}h^{i}\cdot f_{1}^{q-i}. 
\end{equation}
By integrating  \eqref{eq 5.2} on the the unit interval of  $ \mathbb{K}^{2} $ and using  \eqref{eq 4.15},\eqref{eq 4.16} and \eqref{eq 4.9} we get \\
\begin{align}
\int_{\mathbb{P}^{2}} g_{1}^{q}(t,\eta )dtd\eta  &  = \int_{\mathbb{P}^{2}} (h + f_{1})^{q}(t,\eta )dtd\eta  \label{eq 5.3} \\
& = \int_{\mathbb{P}^{2}} h^{q}(t,\eta )dtd\eta  +\int_{\mathbb{P}^{2}} f_{1}^{q}(t,\eta )dtd\eta  + \int_{\mathbb{P}^{2}}\sum_{i = 1}^{q -1} \binom{q}{i}h^{i}(t,\eta )\cdot f_{1}^{q-i}(t,\eta )dtd\eta \nonumber \\
& =   \int_{  \left\{(t,\eta ) \in \mathbb{P}^{2}\mid  h(t,\eta )\neq 0 \right\}}v^{q}(t,\eta ) dt d\eta  +
  \int_{  \left\{(t,\eta ) \in \mathbb{P}^{2}\mid  f_{1}(t,\eta )\neq 0 \right\}}v^{q}(t,\eta ) dt d\eta  \nonumber \\
   & +  (2^{q} - 2)\cdot  \int_{  \left\{(t,\eta ) \in \mathbb{P}^{2}\mid  f_{1}(t,\eta )\cdot f_{2}(t,\eta )\neq 0 \right\}}v^{q}(t,\eta ) dt d\eta . \nonumber
 \end{align}
\end{proof}
 \begin{lem}
\label{lem 5.2} 
We have 
\begin{equation}
\label{eq 5.4}
\int_{\mathbb{P}^{2}} g_{1}^{q}(t,\eta )dtd\eta = 2^{q(2s +m+k -1)}\cdot2^{-(2k+2s+m-2)}\cdot
\sum_{i = 0}^{\inf(2s+m-1,k)} \sigma _{i,i}^{\left[\stackrel{s-1}{\stackrel{s+m }{\overline{\alpha_{s -}}}}\right] \times k }\cdot2^{-qi}.
\end{equation}
\end{lem}
\begin{proof}
We have by \eqref{eq 4.2} observing that $ g_{1}(t,\eta ) $ is constant on cosets of $ \mathbb{P}_{k+s-1}\times \mathbb{P}_{k+s+m-1} $ \\
 \begin{align*}
\int_{\mathbb{P}^{2}} g_{1}^{q}(t,\eta )dtd\eta  & =
 \sum_{(t,\eta )\in \mathbb{P}/\mathbb{P}_{k+s-1}\times \mathbb{P}/\mathbb{P}_{k+s+m-1}\atop
 { r(D^{\big[\stackrel{s-1}{s+m}\big] \times k }(t,\eta ))  = r(D^{\big[\stackrel{s}{s+m}\big] \times k }(t,\eta )) }}
 2^{q(2s+m+k-1- r(D^{\big[\stackrel{s}{s+m}\big] \times k }(t,\eta )) )}\int_{\mathbb{P}_{k+s-1}}dt \int_{\mathbb{P}_{k+s+m-1}}d\eta \\
 & = \sum_{i = 0}^{\inf(2s+m-1,k)} \sum_{(t,\eta )\in \mathbb{P}/\mathbb{P}_{k+s-1}\times \mathbb{P}/\mathbb{P}_{k+s+m-1}\atop
 { r(D^{\big[\stackrel{s-1}{s+m}\big] \times k }(t,\eta ))  = r(D^{\big[\stackrel{s}{s+m}\big] \times k }(t,\eta )) = i}}
  2^{q(2s+m+k-1- i) }\int_{\mathbb{P}_{k+s-1}}dt \int_{\mathbb{P}_{k+s+m-1}}d\eta \\
  & = \sum_{i = 0}^{\inf(2s+m-1,k)} \sigma _{i,i}^{\left[\stackrel{s-1}{\stackrel{s+m }{\overline{\alpha_{s -}}}}\right] \times k }
  \cdot2^{-qi}\cdot2^{q(2s+m+k-1)}\cdot2^{-(k+s-1)}\cdot2^{-(k+s+m-1)}.\\
  & 
 \end{align*}
\end{proof}
 \begin{lem}
\label{lem 5.3} 
We have 
\begin{equation}
\label{eq 5.5}
\int_{\mathbb{P}^{2}} f_{1}^{q}(t,\eta )dtd\eta = 2^{q(2s +m+k -2)}\cdot2^{ - (2k+2s+m-3)}\cdot
\sum_{i = 0}^{\inf(2s+m-2,k)} \sigma _{i,i}^{\left[\stackrel{s-1}{\stackrel{s+m-1 }{\overline{\alpha_{s -}}}}\right] \times k }\cdot2^{-qi}.
\end{equation}
\end{lem}
\begin{proof} Similar to the proof of Lemma \ref{lem 5.2} using  \eqref{eq 4.4} and  observing that
 $ f_{1}(t,\eta ) $ is constant on cosets of $ \mathbb{P}_{k+s-1}\times \mathbb{P}_{k+s+m-2}. $ \\
  \end{proof}
 
  \begin{lem}
\label{lem 5.4} We have 
\begin{equation}
\label{eq 5.6}
 \int_{  \left\{(t,\eta ) \in \mathbb{P}^{2}\mid  f_{1}(t,\eta )\cdot f_{2}(t,\eta )\neq 0 \right\}}v^{q}(t,\eta ) dt d\eta  
 =   2^{q(2s +m+k -2)}\cdot2^{ - (2k+2s+m-2)}\cdot
\sum_{i = 0}^{\inf(2s+m-2,k)} \sigma _{i,i,i}^{\left[\stackrel{s-1}{\stackrel{s+m-1 }
{\overline {\stackrel{\alpha_{s -}}{\beta_{s+m-} }}}}\right] \times k } \cdot2^{-qi}.
\end{equation}
 \end{lem}
 \begin{proof}We have by \eqref{eq 4.12} and  \eqref{eq 4.13}\vspace{0.1 cm} \\
   $\displaystyle  \int_{  \left\{(t,\eta ) \in \mathbb{P}^{2}\mid  f_{1}(t,\eta )\cdot f_{2}(t,\eta )\neq 0 \right\}}v^{q}(t,\eta ) dt d\eta $
  \Small
    \begin{align*}
   & = \sum_{(t,\eta )\in \mathbb{P}/\mathbb{P}_{k+s-1}\times \mathbb{P}/\mathbb{P}_{k+s+m-1}\atop
 { r(D^{\big[\stackrel{s-1}{s-1+m}\big] \times k }(t,\eta ))  = r(D^{\big[\stackrel{s}{s+ m-1}\big] \times k }(t,\eta ))=  r(D^{\big[\stackrel{s-1}{s - 1+ m+1}\big] \times k }(t,\eta )) }}
 2^{q(2s+m+k-2- r(D^{\big[\stackrel{s-1}{s-1+m}\big] \times k }(t,\eta )) )}\int_{\mathbb{P}_{k+s-1}}dt \int_{\mathbb{P}_{k+s+m-1}}d\eta \\
 & = \sum_{i = 0}^{\inf(2s+m-2,k)} \sum_{(t,\eta )\in \mathbb{P}/\mathbb{P}_{k+s-1}\times \mathbb{P}/\mathbb{P}_{k+s+m-1}\atop
 { r(D^{\big[\stackrel{s-1}{s-1+m}\big] \times k }(t,\eta ))  = r(D^{\big[\stackrel{s}{s+m-1}\big] \times k }(t,\eta ))=  r(D^{\big[\stackrel{s-1}{s-1+m+1}\big] \times k }(t,\eta )) = i}}
  2^{q(2s+m+k-2- i) }\int_{\mathbb{P}_{k+s-1}}dt \int_{\mathbb{P}_{k+s+m-1}}d\eta \\
  & = \sum_{i = 0}^{\inf(2s+m-2,k)} \sigma _{i,i,i}^{\left[\stackrel{s-1}{\stackrel{s+m-1 }
{\overline {\stackrel{\alpha_{s -}}{\beta_{s+m-} }}}}\right] \times k }
  \cdot2^{-qi}\cdot2^{q(2s+m+k-2)}\cdot2^{-(k+s-1)}\cdot2^{-(k+s+m-1)}.
 \end{align*}\vspace{0.1 cm}
\end{proof}
   \begin{lem}
\label{lem 5.5} We have 
$$  \sigma _{i,i}^{\left[\stackrel{s-1}{\stackrel{s+m -1}{\overline{\alpha_{s -}+ \beta _{s+m-}}}}\right] \times k } =
  2\cdot \sigma _{i,i}^{\left[\stackrel{s-1}{\stackrel{s+m -1}{\overline{\alpha_{s -}}}}\right] \times k }. $$
  \end{lem}
\begin{proof}We consider the following $ (2s+m-1)\times k $ matrix denoted by 
$ D^{\left[\stackrel{s-1}{\stackrel{s+m-1}{\alpha_{s -} + \beta_{s+m -} }}\right] \times k }(t,\eta  ) $ (see Section \ref{sec 1}) \vspace{0.1 cm}\\

\footnotesize
  $  \bordermatrix{%
                    &                   &                      &                      &                            &    \cr
   r_{1} &    \alpha _{1} & \alpha _{2} & \alpha _{3} &  \ldots & \alpha _{k-1}  &  \alpha _{k} \cr
 r_{2} &   \alpha _{2 } & \alpha _{3} & \alpha _{4}&  \ldots  &  \alpha _{k} &  \alpha _{k+1} \cr
  \vdots &   \vdots & \vdots & \vdots    &  \vdots & \vdots  &  \vdots  \cr
 r_{s-1}   &  \alpha _{s-1} & \alpha _{s} & \alpha _{s +1} & \ldots  &  \alpha _{s+k-3} &  \alpha _{s+k-2}  \cr
r_{s}   & \beta  _{1} & \beta  _{2} & \beta  _{3} & \ldots  &  \beta_{k-1} &  \beta _{k}  \cr
 r_{s+1} & \beta  _{2} & \beta  _{3} & \beta  _{4} & \ldots  &  \beta_{k} &  \beta _{k+1}  \cr
\vdots  &  \vdots & \vdots & \vdots    &  \vdots & \vdots  &  \vdots \cr
r_{s+m} &      \beta  _{m+1} & \beta  _{m+2} & \beta  _{m+3} & \ldots  &  \beta_{k+m-1} &  \beta _{k+m}  \cr
\vdots  & \vdots & \vdots & \vdots    &  \vdots & \vdots  &  \vdots \cr
   r_{2s+m-2} & \beta  _{s+m-1} & \beta  _{s+m} & \beta  _{s+m+1} & \ldots  &  \beta_{s+m+k-3} &  \beta _{s+m+k-2}  \cr
  \hline
r_{2s+m-1} & \alpha _{s} + \beta  _{s+m} & \alpha _{s+1} + \beta  _{s+m+1} & \alpha _{s+2} + \beta  _{s+m+2} & \ldots 
& \alpha _{s+k-2} + \beta  _{s+k+m-2} &  \alpha _{s+k-1} + \beta  _{s+k+m-1} \cr
} $.  \vspace{0.5 cm}\\
$$\big\uparrow $$
\normalsize
 We recall  that the rank of a matrix does not change under elementary row operations.\\
   On the above  matrix, we add  the s+m+j-th row to the j-th row for $ 0\leq j\leq s-2 $ obtaining \\

\scriptsize
   $ \bordermatrix{%
                    &                   &                      &                      &                            &    \cr
 r_{1} + r_{s+m} &  \alpha _{1} + \beta  _{m+1} & \alpha _{2} + \beta  _{m+2} & \alpha _{3}  + \beta  _{m+3} & \ldots & 
\alpha _{k-1}  + \beta  _{m+k-1} &  \alpha _{k}  + \beta  _{m+k} \cr
 r_{2} + r_{s+m +1} & \alpha _{2 } + \beta  _{m+2} & \alpha _{3} + \beta  _{m+3} & \alpha _{4} + \beta  _{m+4}&  \ldots  & 
 \alpha _{k} + \beta  _{m+k} &  \alpha _{k+1} + \beta  _{m+k+1} \cr 
\vdots &  \vdots & \vdots & \vdots    &  \vdots & \vdots  &  \vdots \cr
 r_{s-1} + r_{2s+m -2} & \alpha _{s-1} + \beta  _{s-1+m} & \alpha _{s} + \beta  _{s+m} & \alpha _{s +1} + \beta  _{s+m+1} & \ldots  & 
 \alpha _{s+k-3} + \beta  _{s+k+m-3} &  \alpha _{s+k-2}  + \beta  _{s+k+m-2} \cr
r_{s} & \beta  _{1} & \beta  _{2} & \beta  _{3} & \ldots  &  \beta_{k-1} &  \beta _{k}  \cr
r_{s+1} &  \beta  _{2} & \beta  _{3} & \beta  _{4} & \ldots  &  \beta_{k} &  \beta _{k+1}  \cr
\vdots & \vdots & \vdots    &  \vdots & \vdots  &  \vdots \cr
r_{s+m} & \beta  _{m+1} & \beta  _{m+2} & \beta  _{m+3} & \ldots  &  \beta_{k+m-1} &  \beta _{k+m}  \cr
\vdots & \vdots & \vdots & \vdots    &  \vdots & \vdots  &  \vdots \cr
r_{2s+m-2} & \beta  _{s+m-1} & \beta  _{s+m} & \beta  _{s+m+1} & \ldots  &  \beta_{s+m+k-3} &  \beta _{s+m+k-2}  \cr
\hline
r_{2s+m-1} & \alpha _{s} + \beta  _{s+m} & \alpha _{s+1} + 
\beta  _{s+m+1} & \alpha _{s+2} + \beta  _{s+m+2} & \ldots 
& \alpha _{s+k-2} + \beta  _{s+k+m-2} &  \alpha _{s+k-1} + \beta  _{s+k+m-1}\cr 
}. $  \vspace{0.5 cm}\\
\normalsize
$$\big\uparrow $$
In the above  matrix we set 
$\begin{cases}
 \alpha _{i}+\beta _{i+m} & = \mu _{i}\quad \text{for} \quad 1\leq i\leq k+s-2, \\
  \beta _{j} & = \tau _{j}\quad \text{for}\quad  1\leq j \leq k+s+m-2.
\end{cases} $ \vspace{0.1 cm} \\
( Remark that the map $ \kappa:  \mathbb{F}_{2}^{2k+2s+m-4}\longmapsto \mathbb{F}_{2}^{2k+2s+m-4}$
defined by \vspace{0.1 cm} \\
 $ (\alpha _{1},\alpha _{2}, \ldots,\alpha _{k+s-2},\beta _{1},
 \beta _{2},\ldots,\beta _{k+s+m-2})\longmapsto 
 (\mu _{1},\mu _{2},\ldots,\mu _{k+s-2}, \tau _{1}, \tau _{2},\ldots, \tau _{k+s+m-2}) $ \vspace{0.1 cm} \\
 is an isomorphisme). We then obtain  the below  matrix \vspace{0.1 cm} \\
 $$\big\downarrow $$
  $  \bordermatrix{%
                    &                   &                      &                      &                            &    \cr
   r_{1} &    \mu  _{1} & \mu  _{2} & \mu  _{3} &  \ldots & \mu  _{k-1}  &  \mu  _{k} \cr
 r_{2} &   \mu  _{2 } & \mu  _{3} & \mu  _{4}&  \ldots  &  \mu  _{k} &  \mu  _{k+1} \cr
  \vdots &   \vdots & \vdots & \vdots    &  \vdots & \vdots  &  \vdots  \cr
 r_{s-1}   &  \mu  _{s-1} & \mu  _{s} & \mu  _{s +1} & \ldots  &  \mu _{s+k-3} &  \mu _{s+k-2}  \cr
r_{s}   & \tau  _{1} & \tau  _{2} & \tau  _{3} & \ldots  &  \tau _{k-1} &  \tau  _{k}  \cr
 r_{s+1} & \tau   _{2} & \tau  _{3} & \tau   _{4} & \ldots  &  \tau _{k} &  \tau  _{k+1}  \cr
\vdots  &  \vdots & \vdots & \vdots    &  \vdots & \vdots  &  \vdots \cr
r_{s+m} &      \tau   _{m+1} & \tau  _{m+2} & \tau  _{m+3} & \ldots  &  \tau _{k+m-1} &  \tau _{k+m}  \cr
\vdots  & \vdots & \vdots & \vdots    &  \vdots & \vdots  &  \vdots \cr
   r_{2s+m-2} & \tau  _{s+m-1} & \tau  _{s+m} & \tau  _{s+m+1} & \ldots  &  \tau _{s+m+k-3} &  \tau  _{s+m+k-2}  \cr
  \hline
r_{2s+m-1} & \mu _{s} & \mu _{s+1} & \mu _{s+2} & \ldots 
& \mu _{s+k-2} & \alpha _{s+k-1} + \beta  _{s+k+m-1}  \cr
}. $  \vspace{0.5 cm}\\
 Comparing the above  matrix  with the below  matrix $ \bigg\updownarrow   $        \\

 $  \bordermatrix{%
                    &                   &                      &                      &                            &    \cr
   r_{1} &    \mu  _{1} & \mu  _{2} & \mu  _{3} &  \ldots & \mu  _{k-1}  &  \mu  _{k} \cr
 r_{2} &   \mu  _{2 } & \mu  _{3} & \mu  _{4}&  \ldots  &  \mu  _{k} &  \mu  _{k+1} \cr
  \vdots &   \vdots & \vdots & \vdots    &  \vdots & \vdots  &  \vdots  \cr
 r_{s-1}   &  \mu  _{s-1} & \mu  _{s} & \mu  _{s +1} & \ldots  &  \mu _{s+k-3} &  \mu _{s+k-2}  \cr
r_{s}   & \tau  _{1} & \tau  _{2} & \tau  _{3} & \ldots  &  \tau _{k-1} &  \tau  _{k}  \cr
 r_{s+1} & \tau   _{2} & \tau  _{3} & \tau   _{4} & \ldots  &  \tau _{k} &  \tau  _{k+1}  \cr
\vdots  &  \vdots & \vdots & \vdots    &  \vdots & \vdots  &  \vdots \cr
r_{s+m} &      \tau   _{m+1} & \tau  _{m+2} & \tau  _{m+3} & \ldots  &  \tau _{k+m-1} &  \tau _{k+m}  \cr
\vdots  & \vdots & \vdots & \vdots    &  \vdots & \vdots  &  \vdots \cr
   r_{2s+m-2} & \tau  _{s+m-1} & \tau  _{s+m} & \tau  _{s+m+1} & \ldots  &  \tau _{s+m+k-3} &  \tau  _{s+m+k-2}  \cr
  \hline
r_{2s+m-1} & \mu _{s} & \mu _{s+1} & \mu _{s+2} & \ldots 
& \mu _{s+k-2} & \mu  _{s+k-1}   \cr
} $  \vspace{0.5 cm}\\
and  observing that for all $ \mu _{s+k-1}\in \mathbb{F}_{2} $\\
 $$ Card\left\{(\alpha _{s+k-1},\beta _{s+k+m-1})\in \mathbb{F}_{2}^{2} \vert \quad
\alpha _{s+k-1} +  \beta _{s+k+m-1}= \mu _{s+k-1}\right\} = 2 $$ \\
we obtain the following equality 
\begin{align*}
&  Card \left\{(t,\eta )\in \mathbb{P}/\mathbb{P}_{k+s-1}\times \mathbb{P}/\mathbb{P}_{k+s+m-1}
\mid  r(D^{\big[\stackrel{s-1}{s-1+m}\big] \times k }(t,\eta ))  =  r( D^{\left[\stackrel{s-1}{\stackrel{s+m-1}{\alpha_{s -} + \beta_{s+m -} }}\right] \times k }(t,\eta  )) = i
  \right\} \\
 &  = 2\cdot Card \left\{(t,\eta )\in \mathbb{P}/\mathbb{P}_{k+s-1}\times \mathbb{P}/\mathbb{P}_{k+s+m-2}
\mid  r(D^{\big[\stackrel{s-1}{s-1+m}\big] \times k }(t,\eta ))  = r(D^{\big[\stackrel{s}{s-1+m}\big] \times k }(t,\eta )) = i
 \right\}. 
\end{align*}
 Alternatively consider the following equivalences 
  \begin{align*}
&  r \left( \begin{smallmatrix}
 \alpha _{1} & \alpha _{2} & \alpha _{3} &  \ldots & \alpha _{k-1}  &  \alpha _{k} \\
\alpha _{2 } & \alpha _{3} & \alpha _{4}&  \ldots  &  \alpha _{k} &  \alpha _{k+1} \\
\vdots & \vdots & \vdots    &  \vdots & \vdots  &  \vdots \\
\alpha _{s-1} & \alpha _{s} & \alpha _{s +1} & \ldots  &  \alpha _{s+k-3} &  \alpha _{s+k-2}  \\
\beta  _{1} & \beta  _{2} & \beta  _{3} & \ldots  &  \beta_{k-1} &  \beta _{k}  \\
\beta  _{2} & \beta  _{3} & \beta  _{4} & \ldots  &  \beta_{k} &  \beta _{k+1}  \\
\vdots & \vdots & \vdots    &  \vdots & \vdots  &  \vdots \\
\beta  _{m+1} & \beta  _{m+2} & \beta  _{m+3} & \ldots  &  \beta_{k+m-1} &  \beta _{k+m}  \\
\vdots & \vdots & \vdots    &  \vdots & \vdots  &  \vdots \\
\beta  _{s+m-1} & \beta  _{s+m} & \beta  _{s+m+1} & \ldots  &  \beta_{s+m+k-3} &  \beta _{s+m+k-2}  \\
\hline
 \alpha _{s}+ \beta  _{s+m}  & \alpha _{s+1}+ \beta  _{s+m+1}  & \alpha _{s +2}+\beta  _{s+m+2} & \ldots  &  \alpha _{s+k-2}+ \beta_{s+m+k-2} &  \alpha _{s+k-1}+   \beta _{s+m+k-1}
  \end{smallmatrix}\right) \\
 & \\
  &  =  r \left( \begin{smallmatrix}
 \alpha _{1} & \alpha _{2} & \alpha _{3} &  \ldots & \alpha _{k-1}  &  \alpha _{k} \\
\alpha _{2 } & \alpha _{3} & \alpha _{4}&  \ldots  &  \alpha _{k} &  \alpha _{k+1} \\
\vdots & \vdots & \vdots    &  \vdots & \vdots  &  \vdots \\
\alpha _{s-1} & \alpha _{s} & \alpha _{s +1} & \ldots  &  \alpha _{s+k-3} &  \alpha _{s+k-2}  \\
\beta  _{1} & \beta  _{2} & \beta  _{3} & \ldots  &  \beta_{k-1} &  \beta _{k}  \\
\beta  _{2} & \beta  _{3} & \beta  _{4} & \ldots  &  \beta_{k} &  \beta _{k+1}  \\
\vdots & \vdots & \vdots    &  \vdots & \vdots  &  \vdots \\
\beta  _{m+1} & \beta  _{m+2} & \beta  _{m+3} & \ldots  &  \beta_{k+m-1} &  \beta _{k+m}  \\
\vdots & \vdots & \vdots    &  \vdots & \vdots  &  \vdots \\
\beta  _{s+m-1} & \beta  _{s+m} & \beta  _{s+m+1} & \ldots  &  \beta_{s+m+k-3} &  \beta _{s+m+k-2}  
   \end{smallmatrix}\right) \vspace{2 cm}\\
&  \Longleftrightarrow  \\
&   r \left( \begin{smallmatrix}
  \alpha _{1} + \beta  _{m+1} & \alpha _{2} + \beta  _{m+2} & \alpha _{3}  + \beta  _{m+3} & \ldots & 
\alpha _{k-1}  + \beta  _{m+k-1} &  \alpha _{k}  + \beta  _{m+k} \\
 \alpha _{2 } + \beta  _{m+2} & \alpha _{3} + \beta  _{m+3} & \alpha _{4} + \beta  _{m+4}&  \ldots  & 
 \alpha _{k} + \beta  _{m+k} &  \alpha _{k+1} + \beta  _{m+k+1}\\
  \vdots &  \vdots & \vdots & \vdots    &  \vdots & \vdots  \\
   \alpha _{s-1} + \beta  _{s-1+m} & \alpha _{s} + \beta  _{s+m} & \alpha _{s +1} + \beta  _{s+m+1} & \ldots  & 
 \alpha _{s+k-3} + \beta  _{s+k+m-3} &  \alpha _{s+k-2}  + \beta  _{s+k+m-2} \\
 \beta  _{1} & \beta  _{2} & \beta  _{3} & \ldots  &  \beta_{k-1} &  \beta _{k}  \\
\beta  _{2} & \beta  _{3} & \beta  _{4} & \ldots  &  \beta_{k} &  \beta _{k+1}  \\
\vdots & \vdots & \vdots    &  \vdots & \vdots  &  \vdots \\
\beta  _{m+1} & \beta  _{m+2} & \beta  _{m+3} & \ldots  &  \beta_{k+m-1} &  \beta _{k+m}  \\
\vdots & \vdots & \vdots    &  \vdots & \vdots  &  \vdots \\
\beta  _{s+m-1} & \beta  _{s+m} & \beta  _{s+m+1} & \ldots  &  \beta_{s+m+k-3} &  \beta _{s+m+k-2}  \\
\hline
 \alpha _{s}+ \beta  _{s+m}  & \alpha _{s+1}+ \beta  _{s+m+1}  & \alpha _{s +2}+\beta  _{s+m+2} & \ldots  &  \alpha _{s+k-2}+ \beta_{s+m+k-2} &  \alpha _{s+k-1}+   \beta _{s+m+k-1}
  \end{smallmatrix}\right) \\
  & \\
    & =   r \left( \begin{smallmatrix}
 \alpha _{1} & \alpha _{2} & \alpha _{3} &  \ldots & \alpha _{k-1}  &  \alpha _{k} \\
\alpha _{2 } & \alpha _{3} & \alpha _{4}&  \ldots  &  \alpha _{k} &  \alpha _{k+1} \\
\vdots & \vdots & \vdots    &  \vdots & \vdots  &  \vdots \\
\alpha _{s-1} & \alpha _{s} & \alpha _{s +1} & \ldots  &  \alpha _{s+k-3} &  \alpha _{s+k-2}  \\
\beta  _{1} & \beta  _{2} & \beta  _{3} & \ldots  &  \beta_{k-1} &  \beta _{k}  \\
\beta  _{2} & \beta  _{3} & \beta  _{4} & \ldots  &  \beta_{k} &  \beta _{k+1}  \\
\vdots & \vdots & \vdots    &  \vdots & \vdots  &  \vdots \\
\beta  _{m+1} & \beta  _{m+2} & \beta  _{m+3} & \ldots  &  \beta_{k+m-1} &  \beta _{k+m}  \\
\vdots & \vdots & \vdots    &  \vdots & \vdots  &  \vdots \\
\beta  _{s+m-1} & \beta  _{s+m} & \beta  _{s+m+1} & \ldots  &  \beta_{s+m+k-3} &  \beta _{s+m+k-2}  
   \end{smallmatrix}\right) \vspace{2 cm}\\
& \Longleftrightarrow   \vspace{2 cm}\\ 
 &   r \left( \begin{smallmatrix}                   
     \mu  _{1} & \mu  _{2} & \mu  _{3} &  \ldots & \mu  _{k-1}  &  \mu  _{k} \\
       \mu  _{2 } & \mu  _{3} & \mu  _{4}&  \ldots  &  \mu  _{k} &  \mu  _{k+1}\\
      \vdots &   \vdots & \vdots & \vdots    &  \vdots & \vdots     \\
    \mu  _{s-1} & \mu  _{s} & \mu  _{s +1} & \ldots  &  \mu _{s+k-3} &  \mu _{s+k-2} \\
     \tau  _{1} & \tau  _{2} & \tau  _{3} & \ldots  &  \tau _{k-1} &  \tau  _{k} \\
   \tau   _{2} & \tau  _{3} & \tau   _{4} & \ldots  &  \tau _{k} &  \tau  _{k+1}\\
    \vdots  &  \vdots & \vdots & \vdots    &  \vdots & \vdots     \\
  \tau   _{m+1} & \tau  _{m+2} & \tau  _{m+3} & \ldots  &  \tau _{k+m-1} &  \tau _{k+m}\\
    \vdots  & \vdots & \vdots & \vdots    &  \vdots & \vdots    \\
 \tau  _{s+m-1} & \tau  _{s+m} & \tau  _{s+m+1} & \ldots  &  \tau _{s+m+k-3} &  \tau  _{s+m+k-2}  \\
   \hline 
    \mu _{s} & \mu _{s+1} & \mu _{s+2} & \ldots & \mu _{s+k-2} & \alpha _{s+k-1} + \beta  _{s+k+m-1}  
 \end{smallmatrix}\right) 
 =  r \left( \begin{smallmatrix}          
      \mu  _{1} & \mu  _{2} & \mu  _{3} &  \ldots & \mu  _{k-1}  &  \mu  _{k}\\
        \mu  _{2 } & \mu  _{3} & \mu  _{4}&  \ldots  &  \mu  _{k} &  \mu  _{k+1}\\
     \vdots &   \vdots & \vdots & \vdots    &  \vdots & \vdots    \\
   \mu  _{s-1} & \mu  _{s} & \mu  _{s +1} & \ldots  &  \mu _{s+k-3} &  \mu _{s+k-2}\\
     \tau  _{1} & \tau  _{2} & \tau  _{3} & \ldots  &  \tau _{k-1} &  \tau  _{k} \\
   \tau   _{2} & \tau  _{3} & \tau   _{4} & \ldots  &  \tau _{k} &  \tau  _{k+1} \\
 \vdots  &  \vdots & \vdots & \vdots    &  \vdots & \vdots  &   \\
  \tau   _{m+1} & \tau  _{m+2} & \tau  _{m+3} & \ldots  &  \tau _{k+m-1} &  \tau _{k+m} \\
   \vdots  & \vdots & \vdots & \vdots    &  \vdots & \vdots   \\
  \tau  _{s+m-1} & \tau  _{s+m} & \tau  _{s+m+1} & \ldots  &  \tau _{s+m+k-3} &  \tau  _{s+m+k-2}
    \end{smallmatrix}\right). 
     \end{align*}

  We then get \\
 
 $$  \sigma _{i,i}^{\left[\stackrel{s-1}{\stackrel{s+m -1}
 {\overline{\alpha_{s -} + \beta_{s+m-} }}}\right] \times k } =
  2\cdot \sigma _{i,i}^{\left[\stackrel{s-1}{\stackrel{s+m -1}
 {\overline{\alpha_{s -}}}}\right] \times k }. $$
 \end{proof}

   \begin{lem}
\label{lem 5.6} We have 
\begin{equation}
\label{eq 5.7}
\int_{\mathbb{P}^{2}} h^{q}(t,\eta )dtd\eta = 2^{q(2s +m+k -2)}\cdot2^{-(2k+2s+m-2)}\cdot
\sum_{i = 0}^{\inf(2s+m-2,k)} \sigma _{i,i}^{\left[\stackrel{s-1}{\stackrel{s+m -1}{\overline{\alpha_{s -}+ \beta _{s+m-}}}}\right] \times k } \cdot2^{-qi}.
\end{equation}
 \end{lem}
 \begin{proof}
 We have by \eqref{eq 4.7} observing that $ h(t,\eta ) $ is constant on cosets of $ \mathbb{P}_{k+s-1}\times \mathbb{P}_{k+s+m-1} $ \\
 \begin{align*}
\int_{\mathbb{P}^{2}} h^{q}(t,\eta )dtd\eta  & =
 \sum_{(t,\eta )\in \mathbb{P}/\mathbb{P}_{k+s-1}\times \mathbb{P}/\mathbb{P}_{k+s+m-1}\atop
 { r(D^{\big[\stackrel{s-1}{s-1+m}\big] \times k }(t,\eta ))  =    r( D^{\left[\stackrel{s-1}{\stackrel{s+m-1}{\alpha_{s -} + \beta_{s+m -} }}\right] \times k }(t,\eta  ))  }}
 2^{q(2s+m+k-2- r(D^{\big[\stackrel{s-1}{s-1+m}\big] \times k }(t,\eta )) )}\int_{\mathbb{P}_{k+s-1}}dt \int_{\mathbb{P}_{k+s+m-1}}d\eta \\
 & = \sum_{i = 0}^{\inf(2s+m-2,k)} \sum_{(t,\eta )\in \mathbb{P}/\mathbb{P}_{k+s-1}\times \mathbb{P}/\mathbb{P}_{k+s+m-1}\atop
 { r(D^{\big[\stackrel{s-1}{s-1+m}\big] \times k }(t,\eta ))  =    r( D^{\left[\stackrel{s-1}{\stackrel{s+m-1}{\alpha_{s -} + \beta_{s+m -} }}\right] \times k }(t,\eta  )) = i}}
  2^{q(2s+m+k-2- i) }\int_{\mathbb{P}_{k+s-1}}dt \int_{\mathbb{P}_{k+s+m-1}}d\eta \\
  & = \sum_{i = 0}^{\inf(2s+m-2,k)}  \sigma _{i,i}^{\left[\stackrel{s-1}{\stackrel{s+m -1}{\overline{\alpha_{s -}+ \beta _{s+m-}}}}\right] \times k }
  \cdot2^{-qi}\cdot2^{q(2s+m+k-2)}\cdot2^{-(k+s-1)}\cdot2^{-(k+s+m-1)}.
 \end{align*}
  \end{proof}
     \begin{lem}
\label{lem 5.7} We have 
\begin{equation}
\label{eq 5.8}
  \int_{\mathbb{P}^{2}} h^{q}(t,\eta )dtd\eta = \int_{\mathbb{P}^{2}} f_{1}^{q}(t,\eta )dtd\eta. 
 \end{equation} 
\end{lem}  
  \begin{proof}
 Immediately obtained by comparing  \eqref{eq 5.5} and \eqref{eq 5.7} using  Lemma \ref{lem 5.5}.
  \end{proof}
      \begin{lem}
\label{lem 5.8} We have for  $ 1\leq i\leq \inf(2s+m-1, k) $
\begin{equation}
\label{eq 5.9}
  \sigma _{i-1,i,i}^{\left[\stackrel{s-1}{\stackrel{s+m-1 }
{\overline {\stackrel{\beta _{s+m -}}{\alpha _{s -} }}}}\right] \times k } = 
2\cdot \sigma _{i-1,i-1,i}^{\left[\stackrel{s-1}{\stackrel{s+m-1 }
{\overline {\stackrel{\alpha_{s -}}{\beta_{s+m-} }}}}\right] \times k }. 
  \end{equation} 
\end{lem}
 \begin{proof}
From \eqref{eq 5.1}, \eqref{eq 5.4},\eqref{eq 5.5},\eqref{eq 5.7} and  \eqref{eq 5.8} we have
\begin{align}
\int_{\mathbb{P}^{2}}g_{1}^{q} & 
= \int_{\mathbb{P}^{2}}h^{q}+ \int_{\mathbb{P}^{2}}f_{1}^{q}
+ (2^{q} - 2)\cdot\int_{f_{1}f_{2}\not\neq  0}v^{q} \nonumber \\
& = 2\cdot\int_{\mathbb{P}^{2}}f_{1}^{q} + (2^{q} - 2)\cdot\int_{f_{1}f_{2}\not\neq  0}v^{q}   \nonumber  \\
& \Leftrightarrow 2^{q(2s +m+k -1)}\cdot2^{-(2k+2s+m-2)}\cdot
\sum_{i = 0}^{\inf(2s+m-1,k)} \sigma _{i,i}^{\left[\stackrel{s-1}{\stackrel{s+m }{\overline{\alpha_{s -}}}}\right] \times k }\cdot2^{-qi} \nonumber \\
& = 2\cdot 2^{q(2s +m+k -2)}\cdot2^{ - (2k+2s+m-3)}\cdot
\sum_{i = 0}^{\inf(2s+m-2,k)} \sigma _{i,i}^{\left[\stackrel{s-1}{\stackrel{s+m-1 }{\overline{\alpha_{s -}}}}\right] \times k }\cdot2^{-qi} \nonumber \\
& + (2^{q} - 2 )\cdot2^{q(2s +m+k -2)}\cdot2^{ - (2k+2s+m-2)}\cdot
\sum_{i = 0}^{\inf(2s+m-2,k)} \sigma _{i,i,i}^{\left[\stackrel{s-1}{\stackrel{s+m-1 }
{\overline {\stackrel{\alpha_{s -}}{\beta_{s+m-} }}}}\right] \times k } \cdot2^{-qi} \nonumber \\
& \Leftrightarrow \sum_{i = 0}^{\inf(2s+m-1,k)} \sigma _{i,i}^{\left[\stackrel{s-1}{\stackrel{s+m }{\overline{\alpha_{s -}}}}\right] \times k }\cdot2^{-qi} \nonumber \\
& = 2^{2-q}\cdot\sum_{i = 0}^{\inf(2s+m-2,k)} \sigma _{i,i}^{\left[\stackrel{s-1}{\stackrel{s+m-1 }{\overline{\alpha_{s -}}}}\right] \times k }\cdot2^{-qi}
+  (2^{q} - 2 )\cdot2^{-q}\sum_{i = 0}^{\inf(2s+m-2,k)} \sigma _{i,i,i}^{\left[\stackrel{s-1}{\stackrel{s+m-1 }
{\overline {\stackrel{\alpha_{s -}}{\beta_{s+m-} }}}}\right] \times k } \cdot2^{-qi} \nonumber \\
& \Leftrightarrow \sum_{i = 0}^{\inf(2s+m-1,k)} \sigma _{i,i}^{\left[\stackrel{s-1}{\stackrel{s+m }{\overline{\alpha_{s -}}}}\right] \times k }\cdot2^{-qi} \nonumber \\
& = \sum_{i = 0}^{\inf(2s+m-2,k)}4\cdot \sigma _{i,i}^{\left[\stackrel{s-1}{\stackrel{s+m-1 }{\overline{\alpha_{s -}}}}\right] \times k }\cdot2^{-q(i+1)} \nonumber \\
& + \sum_{i = 0}^{\inf(2s+m-2,k)}  \sigma _{i,i,i}^{\left[\stackrel{s-1}{\stackrel{s+m-1 }
{\overline {\stackrel{\alpha_{s -}}{\beta_{s+m-} }}}}\right] \times k } \cdot2^{-qi}  - \sum_{i = 0}^{\inf(2s+m-2,k)}  2\cdot \sigma _{i,i,i}^{\left[\stackrel{s-1}{\stackrel{s+m-1 }
{\overline {\stackrel{\alpha_{s -}}{\beta_{s+m-} }}}}\right] \times k } \cdot2^{-q(i+1)}  \nonumber \\
& \Leftrightarrow \sum_{i = 0}^{\inf(2s+m-1,k)} \sigma _{i,i}^{\left[\stackrel{s-1}{\stackrel{s+m }{\overline{\alpha_{s -}}}}\right] \times k }\cdot2^{-qi} \label{eq 5.10}\\
& = \sum_{i = 1}^{\inf(2s+m-2,k) + 1}4\cdot \sigma _{i-1,i-1}^{\left[\stackrel{s-1}{\stackrel{s+m-1 }{\overline{\alpha_{s -}}}}\right] \times k }\cdot2^{-qi}  \nonumber \\
& + \sum_{i = 0}^{\inf(2s+m-2,k)}  \sigma _{i,i,i}^{\left[\stackrel{s-1}{\stackrel{s+m-1 }
{\overline {\stackrel{\alpha_{s -}}{\beta_{s+m-} }}}}\right] \times k } \cdot2^{-qi}  - \sum_{i = 1}^{\inf(2s+m-2,k)+ 1}  2\cdot \sigma _{i-1, i-1, i-1}^{\left[\stackrel{s-1}{\stackrel{s+m-1 }
{\overline {\stackrel{\alpha_{s -}}{\beta_{s+m-} }}}}\right] \times k } \cdot2^{-qi}. \nonumber
\end{align}
\underline{The case  $ k\leq  2s+m-2 $}\\
Observing that $ 4\cdot \sigma _{k,k}^{\left[\stackrel{s-1}{\stackrel{s+m-1 }{\overline{\alpha_{s -}}}}\right] \times k }= 
2\cdot \sigma _{k,k,k}^{\left[\stackrel{s-1}{\stackrel{s+m-1 }{\overline {\stackrel{\alpha_{s -}}{\beta_{s+m-} }}}}\right] \times k }$
we deduce from \eqref{eq 5.10}
\begin{equation}
\label{eq 5.11}
\sum_{i=1}^{k}\left[  \sigma _{i,i}^{\left[\stackrel{s-1}{\stackrel{s+m }{\overline{\alpha_{s -}}}}\right] \times k }
-  4\cdot \sigma _{i-1,i-1}^{\left[\stackrel{s-1}{\stackrel{s+m-1 }{\overline{\alpha_{s -}}}}\right] \times k } 
 - \sigma _{i,i,i}^{\left[\stackrel{s-1}{\stackrel{s+m-1 }{\overline {\stackrel{\alpha_{s -}}{\beta_{s+m-} }}}}\right] \times k } +  
   2\cdot \sigma _{i-1, i-1, i-1}^{\left[\stackrel{s-1}{\stackrel{s+m-1 }
{\overline {\stackrel{\alpha_{s -}}{\beta_{s+m-} }}}}\right] \times k }  \right]\cdot2^{-qi} = 0\; for\; all\; q\geq 2.
\end{equation}
Now we have obviously
\begin{align}
 \sigma _{i,i}^{\left[\stackrel{s-1}{\stackrel{s+m }{\overline{\alpha_{s -}}}}\right] \times k } & = 
 \sigma _{i-1,i,i}^{\left[\stackrel{s-1}{\stackrel{s+m-1 }{\overline {\stackrel{\beta _{s +m -}}{\alpha _{s -} }}}}\right] \times k } + 
 \sigma _{i,i,i}^{\left[\stackrel{s-1}{\stackrel{s+m-1 }{\overline {\stackrel{\beta _{s +m -}}{\alpha _{s -} }}}}\right] \times k }\label{eq 5.12}, \\
2\cdot \sigma _{i-1,i-1}^{\left[\stackrel{s-1}{\stackrel{s-1+ m }{\overline{\alpha_{s -}}}}\right] \times k } & =
   \sigma _{i-1,i-1,i-1}^{\left[\stackrel{s-1}{\stackrel{s+m-1 }{\overline {\stackrel{\alpha_{s -}}{\beta_{s+m-} }}}}\right] \times k } +
  \sigma _{i-1,i-1,i}^{\left[\stackrel{s-1}{\stackrel{s+m-1 }{\overline {\stackrel{\alpha_{s -}}{\beta_{s+m-} }}}}\right] \times k }. \label{eq 5.13}
  \end{align} \vspace{0.1 cm}\\
  By  \eqref{eq 5.11}, (5.12) and (5.13) we obtain for all $ q\geq 2 $\vspace{0.1 cm}\\
  \small
  \begin{equation}
  \label{eq 5.14}
\sum_{i=1}^{k}\left[  \sigma _{i-1,i,i}^{\left[\stackrel{s-1}{\stackrel{s+m-1 }{\overline {\stackrel{\beta _{s +m -}}{\alpha _{s -} }}}}\right] \times k } + 
 \sigma _{i,i,i}^{\left[\stackrel{s-1}{\stackrel{s+m-1 }{\overline {\stackrel{\beta _{s +m -}}{\alpha _{s -} }}}}\right] \times k }
 -   2\cdot \sigma _{i-1,i-1,i-1}^{\left[\stackrel{s-1}{\stackrel{s+m-1 }{\overline {\stackrel{\alpha_{s -}}{\beta_{s+m-} }}}}\right] \times k } 
-2\cdot \sigma _{i-1,i-1,i}^{\left[\stackrel{s-1}{\stackrel{s+m-1 }{\overline {\stackrel{\alpha_{s -}}{\beta_{s+m-} }}}}\right] \times k } 
 - \sigma _{i,i,i}^{\left[\stackrel{s-1}{\stackrel{s+m-1 }{\overline {\stackrel{\alpha_{s -}}{\beta_{s+m-} }}}}\right] \times k } +  
   2\cdot \sigma _{i-1, i-1, i-1}^{\left[\stackrel{s-1}{\stackrel{s+m-1 }
{\overline {\stackrel{\alpha_{s -}}{\beta_{s+m-} }}}}\right] \times k }  \right]\cdot2^{-qi} =  0.  
\end{equation}

From \eqref{eq 5.14} we deduce 
 \begin{equation}
  \label{eq 5.15}
\sum_{i=1}^{k}\left[  \sigma _{i-1,i,i}^{\left[\stackrel{s-1}{\stackrel{s+m-1 }{\overline {\stackrel{\beta _{s +m -}}{\alpha _{s -} }}}}\right] \times k } 
-2\cdot \sigma _{i-1,i-1,i}^{\left[\stackrel{s-1}{\stackrel{s+m-1 }{\overline {\stackrel{\alpha_{s -}}{\beta_{s+m-} }}}}\right] \times k }  \right]\cdot2^{-qi} =  0  \; for \;
all \; q\geq 2.
\end{equation}

\underline{The case  $ k\geq  2s+m-1 $}\vspace{0.1 cm} \\
\begin{align}
\sum_{i = 0}^{2s+m-1} \sigma _{i,i}^{\left[\stackrel{s-1}{\stackrel{s+m }{\overline{\alpha_{s -}}}}\right] \times k }\cdot2^{-qi} & = 
\sum_{i = 1}^{2s+m-1} 4\cdot \sigma _{i-1,i-1}^{\left[\stackrel{s-1}{\stackrel{s+m-1 }{\overline{\alpha_{s -}}}}\right] \times k }\cdot2^{-qi} \label{eq 5.16}\\
& + \sum_{i = 0}^{2s+m-2}  \sigma _{i,i,i}^{\left[\stackrel{s-1}{\stackrel{s+m-1 }
{\overline {\stackrel{\alpha_{s -}}{\beta_{s+m-} }}}}\right] \times k } \cdot2^{-qi}  - \sum_{i = 1}^{2s+m-1}  2\cdot \sigma _{i-1, i-1, i-1}^{\left[\stackrel{s-1}{\stackrel{s+m-1 }
{\overline {\stackrel{\alpha_{s -}}{\beta_{s+m-} }}}}\right] \times k } \cdot2^{-qi}. \nonumber
\end{align}\vspace{0.1 cm}\\
Observing that $ \sigma _{2s+m-1,2s+m-1,2s+m-1}^{\left[\stackrel{s-1}{\stackrel{s+m-1 }
{\overline {\stackrel{\alpha_{s -}}{\beta_{s+m-} }}}}\right] \times k } = 0 $ we deduce from \eqref{eq 5.16}, (5.12) and (5.13) that for all  $ q\geq 2 $ \vspace{0.1 cm}\\
  \begin{equation}
  \label{eq 5.17}
\sum_{i=1}^{2s+m-1}\left[  \sigma _{i-1,i,i}^{\left[\stackrel{s-1}{\stackrel{s+m-1 }{\overline {\stackrel{\beta _{s +m -}}{\alpha _{s -} }}}}\right] \times k } + 
 \sigma _{i,i,i}^{\left[\stackrel{s-1}{\stackrel{s+m-1 }{\overline {\stackrel{\beta _{s +m -}}{\alpha _{s -} }}}}\right] \times k }
 -   2\cdot \sigma _{i-1,i-1,i-1}^{\left[\stackrel{s-1}{\stackrel{s+m-1 }{\overline {\stackrel{\alpha_{s -}}{\beta_{s+m-} }}}}\right] \times k } 
-2\cdot \sigma _{i-1,i-1,i}^{\left[\stackrel{s-1}{\stackrel{s+m-1 }{\overline {\stackrel{\alpha_{s -}}{\beta_{s+m-} }}}}\right] \times k } 
 - \sigma _{i,i,i}^{\left[\stackrel{s-1}{\stackrel{s+m-1 }{\overline {\stackrel{\alpha_{s -}}{\beta_{s+m-} }}}}\right] \times k } +  
   2\cdot \sigma _{i-1, i-1, i-1}^{\left[\stackrel{s-1}{\stackrel{s+m-1 }
{\overline {\stackrel{\alpha_{s -}}{\beta_{s+m-} }}}}\right] \times k }  \right]\cdot2^{-qi} =  0.  
\end{equation}\vspace{0.1 cm}\\
From \eqref{eq 5.17} we obtain
 \begin{equation}
  \label{eq 5.18}
\sum_{i=1}^{2s+m-1}\left[  \sigma _{i-1,i,i}^{\left[\stackrel{s-1}{\stackrel{s+m-1 }{\overline {\stackrel{\beta _{s +m -}}{\alpha _{s -} }}}}\right] \times k } 
-2\cdot \sigma _{i-1,i-1,i}^{\left[\stackrel{s-1}{\stackrel{s+m-1 }{\overline {\stackrel{\alpha_{s -}}{\beta_{s+m-} }}}}\right] \times k }  \right]\cdot2^{-qi} =  0  \; for \;
all \; q\geq 2
\end{equation}
and Lemma \ref{lem 5.8}  is proved.
 \end{proof}
 
\begin{lem}
\label{lem 5.9} 
 \begin{align}
    \int_{\mathbb{P}^{2}}g_{2}^{q}(t,\eta )dtd\eta & =
  \int_{  \left\{(t,\eta ) \in \mathbb{P}^{2}\mid  h(t,\eta )\neq 0 \right\}}v^{q}(t,\eta ) dt d\eta  +
  \int_{  \left\{(t,\eta ) \in \mathbb{P}^{2}\mid  f_{2}(t,\eta )\neq 0 \right\}}v^{q}(t,\eta ) dt d\eta \label{eq 5.19}  \\
   & +  (2^{q} - 2)\cdot  \int_{  \left\{(t,\eta ) \in \mathbb{P}^{2}\mid  f_{1}(t,\eta )\cdot f_{2}(t,\eta )\neq 0 \right\}}v^{q}(t,\eta ) dt d\eta.  \nonumber
   \end{align}
  \end{lem}
 \begin{proof} Similar to the proof of Lemma \ref{lem 5.1}.
 \end{proof}
  
\begin{lem}
\label{lem 5.10}
We have for $1\leq i\leq \inf(2s+m-1,k) $ \vspace{0.1 cm}\\
\begin{equation}
\label{eq 5.20}
 \sigma _{i-1,i,i}^{\left[\stackrel{s-1}{\stackrel{s+m-1 }{\overline {\stackrel{\beta _{s+m -}}{\alpha _{s -} }}}}\right] \times k } +  
   \sigma _{i-1,i-1,i}^{\left[\stackrel{s-1}{\stackrel{s+m-1 }{\overline {\stackrel{\beta _{s+m -}}{\alpha _{s -} }}}}\right] \times k }=
    \sigma _{i-1,i,i}^{\left[\stackrel{s-1}{\stackrel{s+m-1 }{\overline {\stackrel{\alpha_{s -}}{\beta_{s+m-} }}}}\right] \times k } +
     \sigma _{i-1,i-1,i}^{\left[\stackrel{s-1}{\stackrel{s+m-1 }{\overline {\stackrel{\alpha_{s -}}{\beta_{s+m-} }}}}\right] \times k }. 
 \end{equation}
  \end{lem} 
  
 \begin{proof}From \eqref{eq 5.1} and \eqref{eq 5.19} we get for all  $ q\geq 2 $
 \begin{equation}
 \label{eq 5.21}
  \int_{\mathbb{P}^{2}}g_{1}^{q}(t,\eta )dtd\eta - \int_{\mathbb{P}^{2}}g_{2}^{q}(t,\eta )dtd\eta =
   \int_{  \left\{(t,\eta ) \in \mathbb{P}^{2}\mid  f_{1}(t,\eta )\neq 0 \right\}}v^{q}(t,\eta ) dt d\eta -  \int_{  \left\{(t,\eta ) \in \mathbb{P}^{2}\mid  f_{2}(t,\eta )\neq 0 \right\}}v^{q}(t,\eta ) dt d\eta.
 \end{equation}
 By \eqref{eq 5.21} we get for all  $ q\geq 2 $ the following equivalences \\
 
 \begin{align}
&  2^{q(2s +m+k -1)}\cdot2^{-(2k+2s+m-2)}\cdot
\sum_{i = 0}^{\inf(2s+m-1,k)}\left[ \sigma _{i,i}^{\left[\stackrel{s-1}{\stackrel{s+m }{\overline{\alpha_{s -}}}}\right] \times k } - 
 \sigma _{i,i}^{\left[\stackrel{s}{\stackrel{s+m -1}{\overline{\beta _{s+m -}}}}\right] \times k }\right]\cdot2^{-qi} \nonumber \\
 & =  2^{q(2s +m+k -2)}\cdot2^{ - (2k+2s+m-3)}\cdot
\sum_{i = 0}^{\inf(2s+m-2,k)}\left[ \sigma _{i,i}^{\left[\stackrel{s-1}{\stackrel{s+m-1 }{\overline{\alpha_{s -}}}}\right] \times k } -
  \sigma _{i,i}^{\left[\stackrel{s-1}{\stackrel{s+m-1 }{\overline{\beta_{s+m -}}}}\right] \times k }\right] \cdot 2^{-qi}\nonumber \\
  & \Leftrightarrow  \sum_{i = 0}^{\inf(2s+m-1,k)}\left[ \sigma _{i,i}^{\left[\stackrel{s-1}{\stackrel{s+m }{\overline{\alpha_{s -}}}}\right] \times k } - 
 \sigma _{i,i}^{\left[\stackrel{s}{\stackrel{s+m -1}{\overline{\beta _{s+m -}}}}\right] \times k }\right]\cdot2^{-qi} 
   =  \sum_{i = 0}^{\inf(2s+m-2,k)}2\cdot\left[ \sigma _{i,i}^{\left[\stackrel{s-1}{\stackrel{s+m-1 }{\overline{\alpha_{s -}}}}\right] \times k } -
  \sigma _{i,i}^{\left[\stackrel{s-1}{\stackrel{s+m-1 }{\overline{\beta_{s+m -}}}}\right] \times k }\right] \cdot 2^{-q(i+1)}\nonumber \\
  &\Leftrightarrow   \sum_{i = 0}^{\inf(2s+m-1,k)}
\left[   \sigma _{i-1,i,i}^{\left[\stackrel{s-1}{\stackrel{s+m-1 }{\overline {\stackrel{\beta _{s+m -}}{\alpha _{s -} }}}}\right] \times k }
- \sigma _{i-1,i,i}^{\left[\stackrel{s-1}{\stackrel{s+m-1 }{\overline {\stackrel{\alpha_{s -}}{\beta_{s+m-} }}}}\right] \times k }   \right]\cdot 2^{-qi}
 =  \sum_{i = 1}^{\inf(2s+m-2,k) +1}\left[ \sigma _{i-1,i-1,i}^{\left[\stackrel{s-1}{\stackrel{s+m-1 }{\overline {\stackrel{\alpha_{s -}}{\beta_{s+m-} }}}}\right] \times k }-  
  \sigma _{i-1,i-1,i}^{\left[\stackrel{s-1}{\stackrel{s+m-1 }{\overline {\stackrel{\beta _{s+m -}}{\alpha _{s -} }}}}\right] \times k } \right]\cdot 2^{-qi}\nonumber \\
  &\Leftrightarrow  \sum_{i = 1}^{\inf(2s+m-1,k)}\left[  \sigma _{i-1,i,i}^{\left[\stackrel{s-1}{\stackrel{s+m-1 }{\overline {\stackrel{\beta _{s+m -}}{\alpha _{s -} }}}}\right] \times k }
- \sigma _{i-1,i,i}^{\left[\stackrel{s-1}{\stackrel{s+m-1 }{\overline {\stackrel{\alpha_{s -}}{\beta_{s+m-} }}}}\right] \times k } - \sigma _{i-1,i-1,i}^{\left[\stackrel{s-1}{\stackrel{s+m-1 }{\overline {\stackrel{\alpha_{s -}}{\beta_{s+m-} }}}}\right] \times k }+ 
 \sigma _{i-1,i-1,i}^{\left[\stackrel{s-1}{\stackrel{s+m-1 }{\overline {\stackrel{\beta _{s+m -}}{\alpha _{s -} }}}}\right] \times k }       \right]\cdot 2^{-qi} = 0.  \label{eq 5.22}
  \end{align}
Now \eqref{eq 5.22} holds for all $ q\geq 2 $ which deduce Lemma \ref{lem 5.10}.\\
\end{proof}
\begin{lem}
\label{lem 5.11}We have 
\begin{align}
 \Gamma _{i}^{\Big[\substack{s \\ s+m }\Big] \times k}
& = \sigma _{i,i,i}^{\left[\stackrel{s-1}{\stackrel{s+m-1 }{\overline {\stackrel{\alpha_{s -}}{\beta_{s+m-} }}}}\right] \times k } 
 + \sigma _{i-1,i,i}^{\left[\stackrel{s-1}{\stackrel{s+m-1 }{\overline {\stackrel{\alpha_{s -}}{\beta_{s+m-} }}}}\right] \times k }
  +\sigma _{i-1,i-1,i}^{\left[\stackrel{s-1}{\stackrel{s+m-1 }{\overline {\stackrel{\alpha_{s -}}{\beta_{s+m-} }}}}\right] \times k } 
  + \sigma _{i-2,i-1,i}^{\left[\stackrel{s-1}{\stackrel{s+m-1 }{\overline {\stackrel{\alpha_{s -}}{\beta_{s+m-} }}}}\right] \times k }\label{eq 5.23} \\
  & = \sigma _{i,i,i}^{\left[\stackrel{s-1}{\stackrel{s+m-1 }{\overline {\stackrel{\beta_{s+m -}}{\alpha_{s -} }}}}\right] \times k } 
 + \sigma _{i-1,i,i}^{\left[\stackrel{s-1}{\stackrel{s+m-1 }{\overline {\stackrel{\beta_{s+m -}}{\alpha_{s -} }}}}\right] \times k }
  +\sigma _{i-1,i-1,i}^{\left[\stackrel{s-1}{\stackrel{s+m-1 }{\overline {\stackrel{\beta_{s+m -}}{\alpha_{s -} }}}}\right] \times k } 
  + \sigma _{i-2,i-1,i}^{\left[\stackrel{s-1}{\stackrel{s+m-1 }{\overline {\stackrel{\beta_{s+m -}}{\alpha_{s -} }}}}\right] \times k }.\nonumber
   \end{align}
  \end{lem} 
\begin{proof} The proof is obvious. \\
\end{proof}
 
\begin{lem}
\label{lem 5.12} We have 
\begin{equation}
\label{eq 5.24}
 \sigma _{i-2,i-1,i}^{\left[\stackrel{s-1}{\stackrel{s+m-1 }{\overline {\stackrel{\alpha_{s -}}{\beta_{s+m-} }}}}\right] \times k }=
\sigma _{i-2,i-1,i}^{\left[\stackrel{s-1}{\stackrel{s+m-1 }{\overline {\stackrel{\beta_{s+m -}}{\alpha_{s -} }}}}\right] \times k }.
\end{equation}
 \end{lem}
 \begin{proof} 
 \eqref{eq 5.24} follows from  \eqref{eq 5.23} and  \eqref{eq 5.20}.  \vspace{0.1 cm}\\
 \end{proof}
 \begin{lem}
 \label{lem 5.13}
 The following holds :
 \begin{itemize}
\item $$ \sigma _{i-1,i,i}^{\left[\stackrel{s-1}{\stackrel{s+m-1 }
{\overline {\stackrel{\beta _{s+m -}}{\alpha _{s -} }}}}\right] \times k } = 
2\cdot \sigma _{i-1,i-1,i}^{\left[\stackrel{s-1}{\stackrel{s+m-1 }
{\overline {\stackrel{\alpha_{s -}}{\beta_{s+m-} }}}}\right] \times k } $$
\item
$$ \sigma _{i-1,i,i}^{\left[\stackrel{s-1}{\stackrel{s+m-1 }{\overline {\stackrel{\beta _{s+m -}}{\alpha _{s -} }}}}\right] \times k } +  
   \sigma _{i-1,i-1,i}^{\left[\stackrel{s-1}{\stackrel{s+m-1 }{\overline {\stackrel{\beta _{s+m -}}{\alpha _{s -} }}}}\right] \times k }=
    \sigma _{i-1,i,i}^{\left[\stackrel{s-1}{\stackrel{s+m-1 }{\overline {\stackrel{\alpha_{s -}}{\beta_{s+m-} }}}}\right] \times k } +
     \sigma _{i-1,i-1,i}^{\left[\stackrel{s-1}{\stackrel{s+m-1 }{\overline {\stackrel{\alpha_{s -}}{\beta_{s+m-} }}}}\right] \times k } $$
 \item  $$  \sigma _{i-2,i-1,i}^{\left[\stackrel{s-1}{\stackrel{s+m-1 }{\overline {\stackrel{\alpha_{s -}}{\beta_{s+m-} }}}}\right] \times k }=
\sigma _{i-2,i-1,i}^{\left[\stackrel{s-1}{\stackrel{s+m-1 }{\overline {\stackrel{\beta_{s+m -}}{\alpha_{s -} }}}}\right] \times k } $$
\item
 $$  \sigma _{i,i,i}^{\left[\stackrel{s-1}{\stackrel{s+m-1 }{\overline {\stackrel{\alpha_{s -}}{\beta_{s+m-} }}}}\right] \times k }=
\sigma _{i,i,i}^{\left[\stackrel{s-1}{\stackrel{s+m-1 }{\overline {\stackrel{\beta_{s+m -}}{\alpha_{s -} }}}}\right] \times k } $$
\item 
$$ \Gamma _{i}^{\Big[\substack{s \\ s+m }\Big] \times k}
 = \sigma _{i,i,i}^{\left[\stackrel{s-1}{\stackrel{s+m-1 }{\overline {\stackrel{\alpha_{s -}}{\beta_{s+m-} }}}}\right] \times k } 
 + \sigma _{i-1,i,i}^{\left[\stackrel{s-1}{\stackrel{s+m-1 }{\overline {\stackrel{\alpha_{s -}}{\beta_{s+m-} }}}}\right] \times k }
  +\sigma _{i-1,i-1,i}^{\left[\stackrel{s-1}{\stackrel{s+m-1 }{\overline {\stackrel{\alpha_{s -}}{\beta_{s+m-} }}}}\right] \times k } 
  + \sigma _{i-2,i-1,i}^{\left[\stackrel{s-1}{\stackrel{s+m-1 }{\overline {\stackrel{\alpha_{s -}}{\beta_{s+m-} }}}}\right] \times k } $$
 \item 
$$  \Gamma _{i}^{\Big[\substack{s \\ s+m }\Big] \times k}
      = \sigma _{i,i,i}^{\left[\stackrel{s-1}{\stackrel{s+m-1 }{\overline {\stackrel{\beta_{s+m -}}{\alpha_{s -} }}}}\right] \times k } 
 + \sigma _{i-1,i,i}^{\left[\stackrel{s-1}{\stackrel{s+m-1 }{\overline {\stackrel{\beta_{s+m -}}{\alpha_{s -} }}}}\right] \times k }
  +\sigma _{i-1,i-1,i}^{\left[\stackrel{s-1}{\stackrel{s+m-1 }{\overline {\stackrel{\beta_{s+m -}}{\alpha_{s -} }}}}\right] \times k } 
  + \sigma _{i-2,i-1,i}^{\left[\stackrel{s-1}{\stackrel{s+m-1 }{\overline {\stackrel{\beta_{s+m -}}{\alpha_{s -} }}}}\right] \times k }. $$
    \end{itemize}
 \end{lem}
 \begin{proof}
  The Theorem follows from  Lemmas \ref{lem 5.8}, \ref{lem 5.10}, \ref{lem 5.11} and \ref{lem 5.12}.\\
   \end{proof}
 
  \begin{lem}
 \label{lem 5.14} Let $ s\geq 2, \; m\geq 0, \; k\geq 1 $ and  $ 0\leq i\leq \inf{(2s+m,k)}. $ Then we have the following recurrent formula for the number of rank i matrices 
 of the form $ \left[{A\over B}\right], $ such that A is a $ s\times k $ persymmetric matrix and B a  $ (s+m)\times k $ persymmetric matrix  with entries in  $ \mathbb{F}_{2} $\\
   \begin{align}
  \Gamma _{i}^{\Big[\substack{s \\ s+m }\Big] \times k} & = 2\cdot \Gamma _{i-1}^{\Big[\substack{s -1\\ s-1+(m+1) }\Big] \times k}
+ 4\cdot \Gamma _{i-1}^{\Big[\substack{s \\ s+(m-1) }\Big] \times k} - 8\cdot \Gamma _{i-2}^{\Big[\substack{s-1 \\ s-1+m }\Big] \times k}
 + \Delta _{i}^{\Big[\substack{s \\ s+m }\Big] \times k} \label{eq 5.25} \\
 & \text{where the remainder $\Delta _{i}^{\Big[\substack{s \\ s+m }\Big] \times k}$ is equal to} \nonumber \\
  & \sigma _{i,i,i}^{\left[\stackrel{s-1}{\stackrel{s+m-1 }{\overline {\stackrel{\alpha_{s -}}{\beta_{s+m-} }}}}\right] \times k }
  - 3\cdot \sigma _{i-1,i-1,i-1}^{\left[\stackrel{s-1}{\stackrel{s+m-1 }{\overline {\stackrel{\alpha_{s -}}{\beta_{s+m-} }}}}\right] \times k } 
  + 2\cdot \sigma _{i-2,i-2,i-2}^{\left[\stackrel{s-1}{\stackrel{s+m-1 }{\overline {\stackrel{\alpha_{s -}}{\beta_{s+m-} }}}}\right] \times k }.\label{eq 5.26} 
  \end{align}
   \end{lem} 
 \begin{proof}
Consider the  following two matrices \vspace{0.2 cm}\\
$ \left( \begin{smallmatrix}
 \alpha _{1} & \alpha _{2} & \alpha _{3} &  \ldots & \alpha _{k-1}  &  \alpha _{k} \\
\alpha _{2 } & \alpha _{3} & \alpha _{4}&  \ldots  &  \alpha _{k} &  \alpha _{k+1} \\
\vdots & \vdots & \vdots    &  \vdots & \vdots  &  \vdots \\
\alpha _{s-1} & \alpha _{s} & \alpha _{s +1} & \ldots  &  \alpha _{s+k-3} &  \alpha _{s+k-2}  \\
\beta  _{1} & \beta  _{2} & \beta  _{3} & \ldots  &  \beta_{k-1} &  \beta _{k}  \\
\beta  _{2} & \beta  _{3} & \beta  _{4} & \ldots  &  \beta_{k} &  \beta _{k+1}  \\
\vdots & \vdots & \vdots    &  \vdots & \vdots  &  \vdots \\
\beta  _{m+1} & \beta  _{m+2} & \beta  _{m+3} & \ldots  &  \beta_{k+m-1} &  \beta _{k+m}  \\
\vdots & \vdots & \vdots    &  \vdots & \vdots  &  \vdots \\
\beta  _{s+m-1} & \beta  _{s+m} & \beta  _{s+m+1} & \ldots  &  \beta_{s+m+k-3} &  \beta _{s+m+k-2}  \\
\hline
 \alpha _{s} & \alpha _{s+1} & \alpha _{s +2} & \ldots  &  \alpha _{s+k-2} &  \alpha _{s+k-1}  \\
 \beta  _{s+m} & \beta  _{s+m+1} & \beta  _{s+m+2} & \ldots  &  \beta_{s+m+k-2} &  \beta _{s+m+k-1}
  \end{smallmatrix}\right) $  and
 $ \left( \begin{smallmatrix}
 \alpha _{1} & \alpha _{2} & \alpha _{3} &  \ldots & \alpha _{k-1}  &  \alpha _{k} \\
\alpha _{2 } & \alpha _{3} & \alpha _{4}&  \ldots  &  \alpha _{k} &  \alpha _{k+1} \\
\vdots & \vdots & \vdots    &  \vdots & \vdots  &  \vdots \\
\alpha _{s-1} & \alpha _{s} & \alpha _{s +1} & \ldots  &  \alpha _{s+k-3} &  \alpha _{s+k-2}  \\
\beta  _{1} & \beta  _{2} & \beta  _{3} & \ldots  &  \beta_{k-1} &  \beta _{k}  \\
\beta  _{2} & \beta  _{3} & \beta  _{4} & \ldots  &  \beta_{k} &  \beta _{k+1}  \\
\vdots & \vdots & \vdots    &  \vdots & \vdots  &  \vdots \\
\beta  _{m+1} & \beta  _{m+2} & \beta  _{m+3} & \ldots  &  \beta_{k+m-1} &  \beta _{k+m}  \\
\vdots & \vdots & \vdots    &  \vdots & \vdots  &  \vdots \\
\beta  _{s+m-1} & \beta  _{s+m} & \beta  _{s+m+1} & \ldots  &  \beta_{s+m+k-3} &  \beta _{s+m+k-2}  \\
\hline
 \beta  _{s+m} & \beta  _{s+m+1} & \beta  _{s+m+2} & \ldots  &  \beta_{s+m+k-2} &  \beta _{s+m+k-1}\\
  \alpha _{s} & \alpha _{s+1} & \alpha _{s +2} & \ldots  &  \alpha _{s+k-2} &  \alpha _{s+k-1}  
   \end{smallmatrix}\right). $ \vspace{0.2 cm}\\
   
   We have by Lemma \ref{lem 5.13} \vspace{0.2 cm} \\
   \begin{align*}
&   2\cdot \Gamma _{i-1}^{\Big[\substack{s -1\\ s-1+(m+1) }\Big] \times k}
+ 4\cdot \Gamma _{i-1}^{\Big[\substack{s \\ s+(m-1) }\Big] \times k} - 8\cdot \Gamma _{i-2}^{\Big[\substack{s-1 \\ s-1+m }\Big] \times k}
 + \Delta _{i}^{\Big[\substack{s \\ s+m }\Big] \times k} \\
 & =    \sigma _{i-1,i-1,i-1}^{\left[\stackrel{s-1}{\stackrel{s+m-1 }{\overline {\stackrel{\beta_{s+m -}}{\alpha_{s -} }}}}\right] \times k } 
 + \sigma _{i-1,i-1,i}^{\left[\stackrel{s-1}{\stackrel{s+m-1 }{\overline {\stackrel{\beta_{s+m -}}{\alpha_{s -} }}}}\right] \times k }
  +\sigma _{i-2,i-1,i-1}^{\left[\stackrel{s-1}{\stackrel{s+m-1 }{\overline {\stackrel{\beta_{s+m -}}{\alpha_{s -} }}}}\right] \times k } 
  + \sigma _{i-2,i-1,i}^{\left[\stackrel{s-1}{\stackrel{s+m-1 }{\overline {\stackrel{\beta_{s+m -}}{\alpha_{s -} }}}}\right] \times k } \\
  & + 2\cdot \sigma _{i-1,i-1,i-1}^{\left[\stackrel{s-1}{\stackrel{s+m-1 }{\overline {\stackrel{\alpha_{s -}}{\beta_{s+m-} }}}}\right] \times k } 
 + 2\cdot\sigma _{i-2,i-1,i-1}^{\left[\stackrel{s-1}{\stackrel{s+m-1 }{\overline {\stackrel{\alpha_{s -}}{\beta_{s+m-} }}}}\right] \times k }
  + 2\cdot\sigma _{i-1,i-1,i}^{\left[\stackrel{s-1}{\stackrel{s+m-1 }{\overline {\stackrel{\alpha_{s -}}{\beta_{s+m-} }}}}\right] \times k } 
  + 2\cdot \sigma _{i-2,i-1,i}^{\left[\stackrel{s-1}{\stackrel{s+m-1 }{\overline {\stackrel{\alpha_{s -}}{\beta_{s+m-} }}}}\right] \times k }\\
  & -  2\cdot \sigma _{i-2,i-2,i-2}^{\left[\stackrel{s-1}{\stackrel{s+m-1 }{\overline {\stackrel{\alpha_{s -}}{\beta_{s+m-} }}}}\right] \times k } 
 - 2\cdot\sigma _{i-2,i-2,i-1}^{\left[\stackrel{s-1}{\stackrel{s+m-1 }{\overline {\stackrel{\alpha_{s -}}{\beta_{s+m-} }}}}\right] \times k }
  - 2\cdot\sigma _{i-2,i-1,i-1}^{\left[\stackrel{s-1}{\stackrel{s+m-1 }{\overline {\stackrel{\alpha_{s -}}{\beta_{s+m-} }}}}\right] \times k } 
  - 2\cdot \sigma _{i-2,i-1,i}^{\left[\stackrel{s-1}{\stackrel{s+m-1 }{\overline {\stackrel{\alpha_{s -}}{\beta_{s+m-} }}}}\right] \times k }\\
  &  +  \sigma _{i,i,i}^{\left[\stackrel{s-1}{\stackrel{s+m-1 }{\overline {\stackrel{\alpha_{s -}}{\beta_{s+m-} }}}}\right] \times k }
  - 3\cdot \sigma _{i-1,i-1,i-1}^{\left[\stackrel{s-1}{\stackrel{s+m-1 }{\overline {\stackrel{\alpha_{s -}}{\beta_{s+m-} }}}}\right] \times k } 
  + 2\cdot \sigma _{i-2,i-2,i-2}^{\left[\stackrel{s-1}{\stackrel{s+m-1 }{\overline {\stackrel{\alpha_{s -}}{\beta_{s+m-} }}}}\right] \times k } \\
  & =   \sigma _{i-1,i-1,i}^{\left[\stackrel{s-1}{\stackrel{s+m-1 }{\overline {\stackrel{\beta_{s+m -}}{\alpha_{s -} }}}}\right] \times k }
  +\sigma _{i-2,i-1,i-1}^{\left[\stackrel{s-1}{\stackrel{s+m-1 }{\overline {\stackrel{\beta_{s+m -}}{\alpha_{s -} }}}}\right] \times k } 
  + \sigma _{i-2,i-1,i}^{\left[\stackrel{s-1}{\stackrel{s+m-1 }{\overline {\stackrel{\beta_{s+m -}}{\alpha_{s -} }}}}\right] \times k } 
   + 2\cdot\sigma _{i-1,i-1,i}^{\left[\stackrel{s-1}{\stackrel{s+m-1 }{\overline {\stackrel{\alpha_{s -}}{\beta_{s+m-} }}}}\right] \times k } \\
   & - 2\cdot\sigma _{i-2,i-2,i-1}^{\left[\stackrel{s-1}{\stackrel{s+m-1 }{\overline {\stackrel{\alpha_{s -}}{\beta_{s+m-} }}}}\right] \times k }
+  \sigma _{i,i,i}^{\left[\stackrel{s-1}{\stackrel{s+m-1 }{\overline {\stackrel{\alpha_{s -}}{\beta_{s+m-} }}}}\right] \times k }  \\
& = \sigma _{i,i,i}^{\left[\stackrel{s-1}{\stackrel{s+m-1 }{\overline {\stackrel{\beta_{s+m -}}{\alpha_{s -} }}}}\right] \times k } 
 + \sigma _{i-1,i,i}^{\left[\stackrel{s-1}{\stackrel{s+m-1 }{\overline {\stackrel{\beta_{s+m -}}{\alpha_{s -} }}}}\right] \times k }
  +\sigma _{i-1,i-1,i}^{\left[\stackrel{s-1}{\stackrel{s+m-1 }{\overline {\stackrel{\beta_{s+m -}}{\alpha_{s -} }}}}\right] \times k } 
  + \sigma _{i-2,i-1,i}^{\left[\stackrel{s-1}{\stackrel{s+m-1 }{\overline {\stackrel{\beta_{s+m -}}{\alpha_{s -} }}}}\right] \times k } \\
  & =  \Gamma _{i}^{\Big[\substack{s \\ s+m }\Big] \times k}.
 \end{align*}
  \end{proof}

   \section{\textbf{RANK PROPERTIES OF A PARTITION OF DOUBLE PERSYMMETRIC MATRICES} }
   \label{sec 6}
    Consider the following partition of the matrix $$ D^{\left[\stackrel{s-1}{\stackrel{s+m-1 }
{\overline {\stackrel{\alpha  _{s -}}{\beta  _{s +m -} }}}}\right] \times k } $$
     $$   \left ( \begin{array} {ccccc|c}
\alpha _{1} & \alpha _{2} & \alpha _{3} &  \ldots & \alpha _{k-1}  &  \alpha _{k} \\
\alpha _{2 } & \alpha _{3} & \alpha _{4}&  \ldots  &  \alpha _{k} &  \alpha _{k+1} \\
\vdots & \vdots & \vdots    &  \vdots & \vdots  &  \vdots \\
\alpha _{s-1} & \alpha _{s} & \alpha _{s +1} & \ldots  &  \alpha _{s+k-3} &  \alpha _{s+k-2}  \\
 \beta  _{1} & \beta  _{2} & \beta  _{3} & \ldots  &  \beta_{k-1} &  \beta _{k}  \\
\beta  _{2} & \beta  _{3} & \beta  _{4} & \ldots  &  \beta_{k} &  \beta _{k+1}  \\
\vdots & \vdots & \vdots    &  \vdots & \vdots  &  \vdots \\
\beta  _{m+1} & \beta  _{m+2} & \beta  _{m+3} & \ldots  &  \beta_{k+m-1} &  \beta _{k+m}  \\
\vdots & \vdots & \vdots    &  \vdots & \vdots  &  \vdots \\
\beta  _{s+m-1} & \beta  _{s+m} & \beta  _{s+m+1} & \ldots  &  \beta_{s+m+k-3} &  \beta _{s+m+k-2}  \\
\hline
\alpha _{s} & \alpha _{s+1} & \alpha _{s +2} & \ldots  &  \alpha _{s+k-2} &  \alpha _{s+k-1} \\
\hline
\beta  _{s+m} & \beta  _{s+m+1} & \beta  _{s+m+2} & \ldots  &  \beta_{s+m+k-2} &  \beta _{s+m+k-1}
\end{array}  \right). $$ \vspace{0.1 cm}

    By integrating some appropriate exponential sums on the unit interval of $\mathbb{K}^2$ with integral equal to zero,\vspace{0.1 cm}\\
    we deduce the following rank formulas  for all $ j\in [0,\inf(2s+m-2,k-1)]: $ \vspace{0.2 cm} \\
\SMALL
    \begin{align*}
  &   {}^{\#}\left(\begin{array}{c | c}
           j & j \\
           \hline
           j & j +1 \\
           \hline
            j & j +1 
           \end{array} \right)_{\mathbb{P}/\mathbb{P}_{k+s -1}\times
           \mathbb{P}/\mathbb{P}_{k+s+m-1} }^{{\alpha  \over \beta }} = 
2\cdot{}^{\#}\left(\begin{array}{c | c}
           j & j \\
           \hline
           j & j\\
           \hline
            j & j 
           \end{array} \right)_{\mathbb{P}/\mathbb{P}_{k+s -1}\times
           \mathbb{P}/\mathbb{P}_{k+s+m-1} }^{{\alpha  \over \beta }} = 
        2\cdot{}^{\#}\left(\begin{array}{c | c}
           j & j \\
           \hline
           j & j\\
           \hline
            j & j +1
           \end{array} \right)_{\mathbb{P}/\mathbb{P}_{k+s -1}\times
           \mathbb{P}/\mathbb{P}_{k+s+m-1} }^{{\alpha  \over \beta }}.\\
            &  
          \end{align*} 
   \normalsize
     \begin{lem}
\label{lem 6.1}
Let $ (t,\eta ) \in\mathbb{P}\times \mathbb{P} $ and set \\
$$  \psi  (t,\eta ) =  \sum_{deg Y = k-1}\sum_{deg Z \leq s-1}E(tYZ)\sum_{deg U \leq  s -1 +m}E(\eta YU). $$
  Then
 \begin{equation}
 \label{eq 6.1}
 \psi (t,\eta )  =   \begin{cases}
 2^{2s +m +k-1 - r( D^{\left[\stackrel{s}{s+m}\right] \times k }(t,\eta  ) ) }  & \text{if }
 r( D^{\left[\stackrel{s}{s+m}\right] \times k }(t,\eta  ) )    =  r( D^{\left[\stackrel{s}{s+m}\right] \times (k-1) }(t,\eta)), \\
  0  & \text{otherwise}.
    \end{cases}
\end{equation} 
\end{lem}
  \begin{proof}
  We have \\
    \begin{align*}
 &  \psi  (t,\eta ) =  \sum_{deg Y = k-1}\sum_{deg Z \leq s-1}E(tYZ)\sum_{deg U\leq  s+m-1}E(\eta YU) \\
 & =  \sum_{deg Y\leq k-1}\sum_{deg Z \leq  s-1}E(tYZ)\sum_{deg U \leq s+m-1}E(\eta YU) - 
\sum_{deg Y\leq k-2}\sum_{deg Z \leq  s-1}E(tYZ)\sum_{deg U \leq s+m-1}E(\eta YU) \\
&   = 2^{2s + m }\cdot\sum_{{deg Y\leq k-1\atop Y\in \ker D^{\left[\stackrel{s}{s+m}\right] \times k }(t,\eta  ) }}1
 -2^{2s + m }\cdot\sum_{{deg Y\leq k-1\atop Y\in \ker D^{\left[\stackrel{s}{s+m}\right] \times (k-1) }(t,\eta  ) }}1 \\
 &  = 2^{2s + m }\cdot2^{k - r( D^{\left[\stackrel{s}{s+m}\right] \times k }(t,\eta  ) ) } -  2^{2s + m }\cdot2^{k-1 - r( D^{\left[\stackrel{s}{s+m}\right] \times (k-1) }(t,\eta  ) ) }. \\
\end{align*} 
 \end{proof} 
  \begin{lem}
\label{lem 6.2} We have
\begin{align}
 \phi^{2} (t,\eta ) & =  \phi_{1} (t,\eta ) \cdot \phi_{2} (t,\eta ), where  \label{eq 6.2} \\
 \phi (t,\eta ) & =   \sum_{deg Y = k-1}\sum_{deg Z \leq s-1}E(tYZ)\sum_{deg U =  s+m-1}E(\eta YU), \nonumber \\
   \phi_{1} (t,\eta ) & =  \sum_{deg Y = k-1}\sum_{deg Z \leq s-1}E(tYZ)\sum_{deg U \leq  s+m-2}E(\eta YU) \quad  and \nonumber\\
    \phi_{2} (t,\eta ) &  =  \sum_{deg Y \leq k-2}\sum_{deg Z \leq s-1}E(tYZ)\sum_{deg U = s+m-1}E(\eta YU). \nonumber
  \end{align}
  \end{lem}
  
 \begin{proof}
 We obtain \\
  \begin{align*}
& \phi^{2} (t,\eta ) \\
  & = \left[ 2^{2s + m -1}\cdot\sum_{{deg Y_{1} = k-1\atop Y_{1}\in \ker D^{\left[\stackrel{s}{s+m-1}\right] \times k }(t,\eta  ) }}E(\eta Y_{1}T^{s+m-1})\right]
 \cdot \left[ 2^{2s + m -1 }\cdot\sum_{{deg Y_{2} = k-1\atop Y_{2}\in \ker D^{\left[\stackrel{s}{s+m-1}\right] \times k }(t,\eta  ) }}E(\eta Y_{2}T^{s+m-1})\right].
 \end{align*}
   We set  \[\left\{\begin{array}{cc}
Y_{1} + Y_{2} = Y_{3}  &  deg Y_{3}\leq k-2, \\
              Y_{1} = Y_{4}   &   deg Y_{4} = k-1.
\end{array}\right.\]\\
Then we get obviously 
  \[\left\{\begin{array}{cc}
 Y_{1} \in \ker D^{\left[\stackrel{s}{s+m-1}\right] \times k }(t,\eta  )   &  deg Y_{1} = k-1, \\
   Y_{2} \in \ker D^{\left[\stackrel{s}{s+m-1}\right] \times k }(t,\eta  )   &  deg Y_{2} = k-1,            
\end{array}\right.\]

  \[\Leftrightarrow \left\{\begin{array}{cc}
 Y_{3} \in \ker D^{\left[\stackrel{s}{s+m-1}\right] \times (k-1) }(t,\eta  )   &  deg Y_{3} \leq k-2, \\
   Y_{4} \in \ker D^{\left[\stackrel{s}{s+m-1}\right] \times k }(t,\eta  )   &  deg Y_{4} = k-1.            
\end{array}\right.\]
 And we obtain  \\
  \begin{align*} \phi^{2} (t,\eta ) & = 
  2^{(2s+m-1)2}\sum_{{deg Y_{4} = k-1\atop Y_{4}\in \ker D^{\left[\stackrel{s}{s+m-1}\right] \times k }(t,\eta  ) }}
  \sum_{{deg Y_{3}\leq k-2\atop Y_{3}\in \ker D^{\left[\stackrel{s}{s +m-1}\right] \times (k-1) }(t,\eta  ) }}
   E(\eta  Y_{3}T^{s+m-1}) \\
  & = \big[ 2^{2s +m-1}\sum_{{deg Y_{4} = k-1\atop Y_{4}\in \ker D^{\left[\stackrel{s}{s +m-1}\right] \times k }(t,\eta  ) }}
   1 \big]\cdot \big[ 2^{2s +m-1}
  \sum_{{deg Y_{3}\leq k-2 \atop Y_{3}\in \ker D^{\left[\stackrel{s}{s+m-1}\right] \times (k-1) }(t,\eta  ) }}
  E(\eta Y_{3}T^{s+m-1})  \big] \\
  & =  \phi_{1} (t,\eta ) \cdot \phi_{2} (t,\eta ).
 \end{align*}
  \end{proof}  
 
   \begin{lem}
\label{lem 6.3} We have the following equivalences 
\begin{align*}
 \phi (t,\eta ) \not\neq 0 &  \Leftrightarrow  \phi^{2} (t,\eta ) \not\neq 0 \\
 & \Leftrightarrow  \phi_{1} (t,\eta ) \cdot \phi_{2} (t,\eta )\not\neq 0 \\
 & \Leftrightarrow  r( D^{\left[\stackrel{s}{s+m-1}\right] \times (k-1) }(t,\eta  ) )    =  r( D^{\left[\stackrel{s}{s+m-1}\right] \times k }(t,\eta  )  )
  =r( D^{\left[\stackrel{s}{s+m}\right] \times (k-1) }(t,\eta  ) ). 
\end{align*}
\end{lem}
\begin{proof} By Lemma \ref{lem 6.1} with $ m \rightarrow m-1 $ we have
 \begin{equation}
 \label{eq 6.3}
 \phi_{1} (t,\eta )    =   \begin{cases}
 2^{2s +m +k-2 - r( D^{\left[\stackrel{s}{s+(m-1)}\right] \times k }(t,\eta  ) ) }  & \text{if }
 r( D^{\left[\stackrel{s}{s+(m-1)}\right] \times k }(t,\eta  ) )    =  r( D^{\left[\stackrel{s}{s+(m-1)}\right] \times (k-1) }(t,\eta)), \\
  0  & \text{otherwise }
    \end{cases}
  \end{equation}  
  and by   \eqref{eq 4.5} with $ k \rightarrow k-1,\; s  \rightarrow s +1 \;\text{and} \; m \rightarrow m-1\; \text{we have}$          
   \begin{equation}
 \label{eq 6.4}   
\phi_{2} (t,\eta )  =   \begin{cases}
 2^{2s +m +k-2 - r( D^{\left[\stackrel{s}{s+m}\right] \times (k-1) }(t,\eta  ) ) }  & \text{if }
 r( D^{\left[\stackrel{s}{s+m}\right] \times (k-1) }(t,\eta  ) )    =  r( D^{\left[\stackrel{s}{s+(m-1)}\right] \times (k-1) }(t,\eta)), \\
  0  & \text{otherwise}.
    \end{cases}
\end{equation}
Lemma \ref{lem 6.3} follows then from \eqref{eq 6.3} and  \eqref{eq 6.4}.\\
 \end{proof}
 
 \begin{lem}
\label{lem 6.4} Let $ (t,\eta )\in \mathbb{P}\times \mathbb{P}$, then 
$$ \phi (t,\eta )  =   \sum_{deg Y = k-1}\sum_{deg Z \leq s-1}E(tYZ)\sum_{deg U =  s+m-1}E(\eta YU) $$ is given by
\begin{equation*}
 \begin{cases}
   2^{2s+m+k-j-2} &  \text{if       }\quad (t,\eta ) \in  
\begin{vmatrix}
k-1 & \vline & k  & \vline \\
\hline
\cdot &\vline & \cdot & \vline & D^{\left[s-1\atop s+m-1 \right]\times \cdot}  \\
\hline
j &\vline & j &  \vline  & \alpha _{s}-\\
\hline 
j &\vline & j & \vline  & \beta _{s+m}- 
\end{vmatrix} \\
&  \\
  -  2^{2s+m+k-j-2} &  \text{if     } \quad (t,\eta )\in    
\begin{vmatrix}
k-1 & \vline & k  & \vline \\
\hline
\cdot &\vline & \cdot & \vline & D^{\left[s-1\atop s+m-1 \right]\times \cdot}  \\
\hline
j &\vline & j &  \vline  & \alpha _{s}-\\
\hline 
j &\vline & j+1 & \vline  & \beta _{s+m}- 
\end{vmatrix} \\
   0    &   \text{otherwise}.
\end{cases}
\end{equation*}
 \end{lem}

 \begin{proof}
 We consider the following two cases in which by Lemma \ref{lem 6.3} $  \phi (t,\eta ) $ is different from zero.\vspace{0.1 cm} \\
 
  \small
 First case : \vspace{0.1 cm}\\
  $ r( D^{\left[\stackrel{s}{s+m-1}\right] \times (k-1) }(t,\eta  ) )    =  r( D^{\left[\stackrel{s}{s+m-1}\right] \times k }(t,\eta  )  )
  = r( D^{\left[\stackrel{s}{s+m}\right] \times (k-1) }(t,\eta  ) ) =   r( D^{\left[\stackrel{s}{s+m}\right] \times (k-1) }(t,\eta  ) ) = j. $ \vspace{0.1 cm}\\
  
 Second case : \vspace{0.1 cm}\\
   $ r( D^{\left[\stackrel{s}{s+m-1}\right] \times (k-1) }(t,\eta  ) )    =  r( D^{\left[\stackrel{s}{s+m-1}\right] \times k }(t,\eta  )  )
  = r( D^{\left[\stackrel{s}{s+m}\right] \times (k-1) }(t,\eta  ) ) = j \quad  \text{and} \quad    r( D^{\left[\stackrel{s}{s+m}\right] \times  k }(t,\eta  ) ) = j+1. $
\normalsize \vspace{0.1 cm}\\
We obtain  using \eqref{eq 4.3} :\vspace{0.1 cm}\\

In the first case \vspace{0.1 cm}\\
\begin{align*}
\phi (t,\eta ) &  =   \sum_{deg Y = k-1}\sum_{deg Z \leq s-1}E(tYZ)\sum_{deg U =  s+m-1}E(\eta YU) \\
& =  \sum_{deg Y \leq  k-1}\sum_{deg Z \leq s-1}E(tYZ)\sum_{deg U =  s+m-1}E(\eta YU) -
\sum_{deg Y \leq  k-2}\sum_{deg Z \leq s-1}E(tYZ)\sum_{deg U =  s+m-1}E(\eta YU)\\
& = 2^{2s+m+k-1-j} - 2^{2s +m +k -2 -j} =  2^{2s +m +k -2 -j}.
\end{align*}\vspace{0.1 cm}\\
In the second case \vspace{0.1 cm}\\
\begin{align*}
\phi (t,\eta ) &  =   \sum_{deg Y = k-1}\sum_{deg Z \leq s-1}E(tYZ)\sum_{deg U =  s+m-1}E(\eta YU) \\
& =  \sum_{deg Y \leq  k-1}\sum_{deg Z \leq s-1}E(tYZ)\sum_{deg U =  s+m-1}E(\eta YU) -
\sum_{deg Y \leq  k-2}\sum_{deg Z \leq s-1}E(tYZ)\sum_{deg U =  s+m-1}E(\eta YU)\\
& =  0  - 2^{2s +m +k -2 -j} = - 2^{2s +m +k -2 -j}.
\end{align*}\vspace{0.1 cm}\\
And otherwise  $ \phi (t,\eta ) $  is equal to zero.\\
\end{proof} 
\begin{lem}
\label{lem 6.5}
We have 
$$\int_{\mathbb{P}\times \mathbb{P}} \phi^{2q+1} (t,\eta )dtd\eta =  0. $$
\end{lem}
\begin{proof}The integral  $ \int_{\mathbb{P}\times \mathbb{P}} \phi^{2q+1} (t,\eta )dtd\eta  $ is equal to the 
number of solutions \\
 $(Y_{1},Z_{1},U_{1}, Y_{2},Z_{2},U_{2},\ldots, Y_{2q+1},Z_{2q+1},U_{2q+1}) $
of the polynomial equations
    \[\left\{\begin{array}{cc}
Y_{1}Z_{1} + Y_{2}Z_{2} + \ldots + Y_{2q+1}Z_{2q+1}= 0, \\
Y_{1}U_{1} + Y_{2}U_{2} + \ldots + Y_{2q+1}U_{2q+1}= 0,
    \end{array}\right.\]\\ 
 satisfying the degree conditions \\
$ \deg Y_{i} = k-1,\quad  \deg Z_{i} \leq s-1  \quad  \deg U_{i}= s+m-1 \quad for \quad 1\leq i\leq 2q+1. $\vspace{0.1 cm}\\
Now 2q+1 is odd so  $ \deg\sum_{i = 1}^{2q+1}Y_{i}U_{i} $ is equal to k+s+m-2. The Lemma follows.\\
\end{proof}
   \begin{lem}
\label{lem 6.6}
 Let $ j \in\mathbb{N} \; such\; that \; 0\leq j\leq\inf(2s+m-1,k-1), $ then \\
 \begin{align*}
{}^{\#}\left(\begin{array}{c | c}
           j & j \\
           \hline
           j & j 
           \end{array} \right)_{\mathbb{P}/\mathbb{P}_{k+s -1}\times
           \mathbb{P}/\mathbb{P}_{k+s+m-1} }^{{\alpha  \over \beta }}=
   {}^{\#}\left(\begin{array}{c | c}
           j & j \\
           \hline
           j & j+1 
           \end{array} \right)_{\mathbb{P}/\mathbb{P}_{k+s -1}\times
           \mathbb{P}/\mathbb{P}_{k+s+m-1} }^{{\alpha \over \beta }}.        
           \end{align*}
 \end{lem}
 \begin{proof} 
   We define  
  $$\begin{array}{l} \mathbb{A} =  \{(t,\eta ) \in \mathbb{P}^{2}
\mid  r( D^{\left[\stackrel{s}{s+m-1}\right] \times (k-1) }(t,\eta  ) )    =  r( D^{\left[\stackrel{s}{s+m-1}\right] \times k }(t,\eta  )  ) \\
 = r( D^{\left[\stackrel{s}{s+m}\right] \times (k-1) }(t,\eta  ) ) =   r( D^{\left[\stackrel{s}{s+m}\right] \times k }(t,\eta  ) ) \} \end{array} $$\\
 and 
  $$\begin{array}{l}  \mathbb{B} =  \{(t,\eta ) \in \mathbb{P}^{2}
\mid  r( D^{\left[\stackrel{s}{s+m-1}\right] \times (k-1) }(t,\eta  ) )    =  r( D^{\left[\stackrel{s}{s+m-1}\right] \times k }(t,\eta  )  )= r( D^{\left[\stackrel{s}{s+m}\right] \times (k-1) }(t,\eta  ) ), \\
 r( D^{\left[\stackrel{s}{s+m}\right] \times k }(t,\eta  ) ) =   r( D^{\left[\stackrel{s}{s+m-1}\right] \times (k-1) }(t,\eta  ) ) +1 \}. \end{array} $$\\
 
  By lemma \ref{lem 6.4} we have by observing that $ \phi (t,\eta ) $ is constant on cosets of $ \mathbb{P}_{k+s-1}\times \mathbb{P}_{k+s+m-1} $
$$ \int_{\mathbb{P}\times \mathbb{P}} \phi^{2q+1} (t,\eta ) dt\eta 
 = \int_{  \mathbb{A}} \phi^{2q+1} (t,\eta ) dt\eta 
    +  \int_{ \mathbb{B}}\phi^{2q+1} (t,\eta ) dt\eta   $$
$$ =\sum_{(t,\eta ) \in \mathbb{A}\bigcap\big(\mathbb{P}/\mathbb{P}_{k+s-1}\times\mathbb{P}/\mathbb{P}_{k+s+m-1}\big) }
 2^{(2s+m+k-2-  r( D^{\left[\stackrel{s}{s+m-1}\right] \times (k-1) }(t,\eta  ) )(2q+1) }\int_{\mathbb{P}_{s+k-1}}dt\int_{\mathbb{P}_{s+k+m-1}}d\eta  $$
$$  + \sum_{(t,\eta ) \in \mathbb{B}\bigcap\big(\mathbb{P}/\mathbb{P}_{k+s-1}\times\mathbb{P}/\mathbb{P}_{k+s+m-1}\big) }
- 2^{(2s+m+k-2-  r( D^{\left[\stackrel{s}{s+m-1}\right] \times (k-1) }(t,\eta  ))(2q+1)  }\int_{\mathbb{P}_{s+k-1}}dt\int_{\mathbb{P}_{s+k+m-1}}d\eta  $$ 
=  \begin{align*}
 \sum_{j=0}^{\inf(2s+m-1,k-1)}  2^{(2s+m+k-2-j)(2q+1)}\cdot {}^{\#}\left(\begin{array}{c | c}
           j & j \\
           \hline
           j & j 
           \end{array} \right)_{\mathbb{P}/\mathbb{P}_{k+s -1}\times
           \mathbb{P}/\mathbb{P}_{k+s+m-1} }^{{\alpha \over \beta  }}\int_{\mathbb{P}_{s+k-1}}dt\int_{\mathbb{P}_{s+k+m-1}}d\eta  \\
           - \sum_{j=0}^{\inf(2s+m-1,k-1)} 2^{(2s+m+k-2-j)(2q+1)}\cdot {}^{\#}\left(\begin{array}{c | c}
           j & j \\
           \hline
           j & j+1 
           \end{array} \right)_{\mathbb{P}/\mathbb{P}_{k+s -1}\times
           \mathbb{P}/\mathbb{P}_{k+s+m-1} }^{{\alpha \over \beta  }} \int_{\mathbb{P}_{s+k-1}}dt\int_{\mathbb{P}_{s+k+m-1}}d\eta.  \\
            \end{align*}
 Now by Lemma \ref{lem 6.5} we get  for all  $ q\in \mathbb{N} $ \\
  \begin{align*}
  \sum_{j=0}^{\inf(2s+m-1,k-1)}2^{(-j)(2q+1)}\cdot\Big( {}^{\#}\left(\begin{array}{c | c}
           j & j \\
           \hline
           j & j 
           \end{array} \right)_{\mathbb{P}/\mathbb{P}_{k+s -1}\times \mathbb{P}/\mathbb{P}_{k+s+m-1} }^{{\alpha \over \beta }} 
               -   {}^{\#}\left( \begin{array}{c | c}
           j & j \\
           \hline
           j & j +1
           \end{array} \right)_{\mathbb{P}/\mathbb{P}_{k+s -1}\times
           \mathbb{P}/\mathbb{P}_{k+s+m-1} }^{{\alpha \over \beta }} \Big)  = 0,
            \end{align*}
            which proves Lemma \ref{lem 6.6}.\\
             \end{proof}

       \begin{lem}
\label{lem 6.7} Let $ (t,\eta)\in\mathbb{P}^{2}, $  q be a rational integer $ \geq 2 $ and 
\begin{align*}
 \theta _{1} (t,\eta ) & =   \sum_{deg Y = k-1}\sum_{deg Z = s-1}E(tYZ)\sum_{deg U \leq  s+m-2}E(\eta YU), \\
 \theta _{2} (t,\eta ) & =  \sum_{deg Y \leq  k-2}\sum_{deg Z \leq s-1}E(tYZ)\sum_{deg U = s+m-1}E(\eta YU) 
   \end{align*}
Then  $ \theta _{1} (t,\eta ) $ is given by 
\begin{equation}
\label{eq 6.5}
 \begin{cases}
   2^{2s+m+k-j-3} &  \text{if       }\quad (t,\eta ) \in    \begin{vmatrix}
k-1 & \vline & k  & \vline \\
\hline
j &\vline & j & \vline & D^{\left[s-1\atop s+m-1 \right]\times \cdot}  \\
\hline
j &\vline & j &  \vline  & \alpha _{s}-\\
\hline 
\cdot &\vline & \cdot & \vline  & \beta _{s+m}- 
\end{vmatrix} \\
 & \\
  -  2^{2s+m+k-j-3} &  \text{if     } \quad (t,\eta ) \in   \begin{vmatrix}
k-1 & \vline & k  & \vline \\
\hline
j &\vline & j & \vline & D^{\left[s-1\atop s+m-1 \right]\times \cdot}  \\
\hline
j &\vline & j+1 &  \vline  & \alpha _{s}-\\
\hline 
\cdot &\vline & \cdot & \vline  & \beta _{s+m}- 
\end{vmatrix} \\
     0    &   \text{otherwise},
\end{cases}
\end{equation} 
and  $ \theta _{2} (t,\eta ) $ is given by 
 \begin{equation}
\label{eq 6.6}
  \begin{cases}
 2^{2s +m +k-2 - j }  & \text{if }\quad  (t,\eta ) \in    \begin{vmatrix}
k-1 & \vline & k  & \vline \\
\hline
\cdot &\vline & \cdot & \vline & D^{\left[s-1\atop s+m-1 \right]\times \cdot}  \\
\hline
j &\vline & \cdot &  \vline  & \alpha _{s}-\\
\hline 
j &\vline & \cdot & \vline  & \beta _{s+m}- 
\end{vmatrix} \\
& \\
  0  & \text{otherwise},
    \end{cases}
\end{equation} 
\end{lem}

\begin{proof}
  The proofs of \eqref{eq 6.5}, \eqref{eq 6.6} are respectively  similar to the proofs of  Lemmas \ref{lem 6.4},
  \ref{lem 4.7}.
   \end{proof}
   
     \begin{lem}
\label{lem 6.8}
 Let $ (t,\eta)\in\mathbb{P}^{2} $ and  q be a rational integer $ \geq 2, $  then
$  \theta _{1} (t,\eta )\cdot  \theta _{2}^{q-1} (t,\eta )  $ is equal to 
\begin{equation}
\label{eq 6.7}
  \begin{cases}
 2^{q-1}\cdot2^{(2s+m+k-3-j)q}     &  \text{if       }\quad  (t,\eta ) \in    \begin{vmatrix}
k-1 & \vline & k  & \vline \\
\hline
j &\vline & j & \vline & D^{\left[s-1\atop s+m-1 \right]\times \cdot}  \\
\hline
j &\vline & j &  \vline  & \alpha _{s}-\\
\hline 
j &\vline & \cdot & \vline  & \beta _{s+m}- 
\end{vmatrix} \\
   & \\
   -    2^{q-1}\cdot2^{(2s+m+k-3-j)q}  &  \text{if     } \quad  (t,\eta ) \in    \begin{vmatrix}
k-1 & \vline & k  & \vline \\
\hline
j &\vline & j & \vline & D^{\left[s-1\atop s+m-1 \right]\times \cdot}  \\
\hline
j &\vline & j+1 &  \vline  & \alpha _{s}-\\
\hline 
j &\vline & j+1 & \vline  & \beta _{s+m}- 
\end{vmatrix} \\
   0    &   \text{otherwise}.
\end{cases}
\end{equation} 
\end{lem}

\begin{proof}
     We consider the following partition of the matrix $$ D^{\left[\stackrel{s-1}{\stackrel{s+m-1 }
{\overline {\stackrel{\alpha  _{s -}}{\beta  _{s +m -} }}}}\right] \times k } $$
 
    $$   \left ( \begin{array} {ccccc|c}
\alpha _{1} & \alpha _{2} & \alpha _{3} &  \ldots & \alpha _{k-1}  &  \alpha _{k} \\
\alpha _{2 } & \alpha _{3} & \alpha _{4}&  \ldots  &  \alpha _{k} &  \alpha _{k+1} \\
\vdots & \vdots & \vdots    &  \vdots & \vdots  &  \vdots \\
\alpha _{s-1} & \alpha _{s} & \alpha _{s +1} & \ldots  &  \alpha _{s+k-3} &  \alpha _{s+k-2}  \\
\beta  _{1} & \beta  _{2} & \beta  _{3} & \ldots  &  \beta_{k-1} &  \beta _{k}  \\
\beta  _{2} & \beta  _{3} & \beta  _{4} & \ldots  &  \beta_{k} &  \beta _{k+1}  \\
\vdots & \vdots & \vdots    &  \vdots & \vdots  &  \vdots \\
\beta  _{m+1} & \beta  _{m+2} & \beta  _{m+3} & \ldots  &  \beta_{k+m-1} &  \beta _{k+m}  \\
\vdots & \vdots & \vdots    &  \vdots & \vdots  &  \vdots \\
\beta  _{s+m-1} & \beta  _{s+m} & \beta  _{s+m+1} & \ldots  &  \beta_{s+m+k-3} &  \beta _{s+m+k-2}  \\
\hline
\alpha _{s} & \alpha _{s+1} & \alpha _{s +2} & \ldots  &  \alpha _{s+k-2} &  \alpha _{s+k-1}\\
\hline
\beta  _{s+m} & \beta  _{s+m+1} & \beta  _{s+m+2} & \ldots  &  \beta_{s+m+k-2} &  \beta _{s+m+k-1}
  \end{array}  \right). $$ \vspace{0.5 cm}\\
  
Obviously by \eqref{eq 6.6} we have \\
  \begin{equation}
  \label{eq 6.8}
 \theta _{2}^{q-1} (t,\eta ) =  \begin{cases}
    2^{(2s+m+k-2-j)(q-1)}  & \text{if }
 r( D^{\left[\stackrel{s}{s-1+m}\right] \times (k-1) }(t,\eta  ) )    =  r( D^{\left[\stackrel{s}{s+m}\right] \times (k-1) }(t,\eta  ) ) = j,  \\
  0  & \text{otherwise }.
    \end{cases}
\end{equation}

Then by considering the above matrix,  using   \eqref{eq 6.5},  \eqref{eq 6.8} and elementary rank properties we get \\

$ \theta _{1} (t,\eta )\cdot  \theta _{2}^{q-1} (t,\eta )  $\; is  equal to\; $  2^{q-1}\cdot2^{(2s+m+k-3-j)q}$ \quad if and only if \\

\begin{align*}
 &  r( D^{\left[\stackrel{s-1}{s-1+m}\right] \times (k-1) }(t,\eta  ) )  
     =  r( D^{\left[\stackrel{s-1}{s+m-1}\right] \times k }(t,\eta  )  )  = r( D^{\left[\stackrel{s}{s+m-1}\right] \times (k-1) }(t,\eta  ) ) \\
      &   =  r( D^{\left[\stackrel{s}{s+m-1}\right] \times k }(t,\eta  ) ) =   r( D^{\left[\stackrel{s}{s+m}\right] \times (k-1) }(t,\eta  )  )   = j \\
       & \text{and} \\
       &   r( D^{\left[\stackrel{s}{s+m}\right] \times k }(t,\eta  ) )    =  j \; \text{or}\; j+1\\
       & \\
 \end{align*}       
       
   $ \theta _{1} (t,\eta )\cdot  \theta _{2}^{q-1} (t,\eta )  $\; is  equal to\; $ - 2^{q-1}\cdot2^{(2s+m+k-3-j)q}$ \quad if and only if \\

\begin{align*}
 &  r( D^{\left[\stackrel{s-1}{s-1+m}\right] \times (k-1) }(t,\eta  ) )  
     =  r( D^{\left[\stackrel{s}{s+m-1}\right] \times (k-1) }(t,\eta  )  )  = r( D^{\left[\stackrel{s}{s+m}\right] \times (k-1) }(t,\eta  ) ) 
         =  r( D^{\left[\stackrel{s-1}{s+m-1}\right] \times k }(t,\eta  ) )   = j \\
       & \text{and} \\
       &   r( D^{\left[\stackrel{s-1}{s+m-1}\right] \times k }(t,\eta  ) ) =  r( D^{\left[\stackrel{s}{s+m}\right] \times k }(t,\eta  ) )   =  j +1\\
       & \\
 \end{align*}       
And in all the other cases $ \theta _{1} (t,\eta )\cdot  \theta _{2}^{q-1} (t,\eta )  $ is equal to zero,
which proves Lemma \ref{lem 6.8}.         
  \end{proof}

   \begin{lem}
\label{lem 6.9} 
We have 
$$\int_{\mathbb{P}\times \mathbb{P}}\theta _{1} (t,\eta )\cdot  \theta _{2}^{q-1} (t,\eta )  dtd\eta =  0. $$
\end{lem}
\begin{proof}The integral  $\int_{\mathbb{P}\times \mathbb{P}}\theta _{1} (t,\eta )\cdot  \theta _{2}^{q-1} (t,\eta )  dt d\eta  $ is equal to the 
number of solutions \\
 $ (Y, Z,U, Y_{1},Z_{1},U_{1},\ldots, Y_{q-1},Z_{q-1},U_{q-1}) $
of the polynomial equations
    \[\left\{\begin{array}{cc}
YZ + Y_{1}Z_{1} + \ldots + Y_{q-1}Z_{q-1}= 0, \\
YU + Y_{1}U_{1} + \ldots + Y_{q-1}U_{q-1}= 0,
    \end{array}\right.\]\\ 
 satisfying the degree conditions \\
$ \deg Y = k-1,\quad  \deg Z = s-1  \quad  \deg U \leq s+m-2  $ \vspace{0.1 cm}\\
$ \deg Y_{i} \leq  k-2, \quad  \deg Z_{i} \leq  s-1  \quad  \deg U_{i} =  s+m-1  \quad for 1\leq i\leq q-1.  $ \vspace{0.1 cm}\\
By degree considerations  we have that $ \deg \left[YZ + Y_{1}Z_{1} + \ldots + Y_{q-1}Z_{q-1} \right] $ is equal to s + k - 2.\\
 Lemma \ref{lem 6.9} follows.
 \end{proof}

   \begin{lem}
\label{lem 6.10}
We have for all $ j\in \mathbb{Z} $ such that $ 0\leq j\leq \inf(2s+m-2, k-1), $
\Small
\begin{align*}
{}^{\#}\left(\begin{array}{c | c}
           j & j \\
           \hline
           j & j\\
           \hline
            j & j 
           \end{array} \right)_{\mathbb{P}/\mathbb{P}_{k+s -1}\times
           \mathbb{P}/\mathbb{P}_{k+s+m-1} }^{{\alpha  \over \beta }}
            +    
{}^{\#}\left(\begin{array}{c | c}
           j & j \\
           \hline
           j & j\\
           \hline
            j & j +1
           \end{array} \right)_{\mathbb{P}/\mathbb{P}_{k+s -1}\times
           \mathbb{P}/\mathbb{P}_{k+s+m-1} }^{{\alpha  \over \beta }}
          -  
{}^{\#}\left(\begin{array}{c | c}
           j & j \\
           \hline
           j & j+1 \\
           \hline
            j & j +1
           \end{array} \right)_{\mathbb{P}/\mathbb{P}_{k+s -1}\times
           \mathbb{P}/\mathbb{P}_{k+s+m-1} }^{{\alpha  \over \beta }} = 0.
            \end{align*}\vspace{0.1 cm}\\
            
\normalsize

  Recall that   
\begin{align*}
{}^{\#}\left(\begin{array}{c | c}
           j_{1} & j_{2} \\
           \hline
           j_{3} & j_{4}\\
           \hline
            j_{5} & j_{6} 
           \end{array} \right)_{\mathbb{P}/\mathbb{P}_{k+s -1}\times
           \mathbb{P}/\mathbb{P}_{k+s+m-1} }^{{\alpha  \over \beta }}
            \end{align*}
    denotes  the cardinality of the following set
    $$\begin{array}{l}\Big\{ (t,\eta ) \in \mathbb{P}/\mathbb{P}_{k+s -1}\times
           \mathbb{P}/\mathbb{P}_{k+s+m-1} 
\mid r(  D^{\left[\stackrel{s-1}{s+m-1}\right] \times (k-1) }(t,\eta  ) ) = j_{1}, \quad 
r( D^{\left[\stackrel{s-1}{s+m-1}\right] \times k }(t,\eta  ) ) = j_{2},  \\
 r(  D^{\left[\stackrel{s}{s+m-1}\right] \times (k-1) }(t,\eta  )  ) = j_{3},\quad  
  r( D^{\left[\stackrel{s}{s+m-1}\right] \times k }(t,\eta  ) ) = j_{4}, \\
   r(  D^{\left[\stackrel{s}{s+m}\right] \times (k-1) }(t,\eta  )  ) = j_{5}, \quad
  r( D^{\left[\stackrel{s}{s+m}\right] \times k }(t,\eta  ) ) = j_{6} \Big\}
    \end{array}$$\\
 for $  (j_{1},j_{2},j_{3},j_{4},j_{5},j_{6}) \in \mathbb{N}^{6}. $ \\
 \end{lem} 
 \begin{proof} 
    We define  
  $$\begin{array}{l} \mathbb{A} =  \{(t,\eta ) \in \mathbb{P}^{2}
\mid  r( D^{\left[\stackrel{s-1}{s+m-1}\right] \times (k-1) }(t,\eta  ) )  
  =  r( D^{\left[\stackrel{s-1}{s+m-1}\right] \times k }(t,\eta  )  ) \\
 = r( D^{\left[\stackrel{s}{s+m-1}\right] \times (k-1) }(t,\eta  ) ) =   r( D^{\left[\stackrel{s}{s+m-1}\right] \times k }(t,\eta  ) ) 
 =  r( D^{\left[\stackrel{s}{s+m}\right] \times (k-1) }(t,\eta  ) )  \} \end{array} $$\\
 and 
  $$\begin{array}{l}  \mathbb{B} =  \{(t,\eta ) \in \mathbb{P}^{2}
\mid  r( D^{\left[\stackrel{s-1}{s+m-1}\right] \times (k-1) }(t,\eta  ) ) 
  =  r( D^{\left[\stackrel{s-1}{s+m-1}\right] \times k }(t,\eta  )  )= 
  r( D^{\left[\stackrel{s}{s+m-1}\right] \times (k-1) }(t,\eta  ) ) \\
= r( D^{\left[\stackrel{s}{s+m}\right] \times (k-1) }(t,\eta  ) ), \quad r( D^{\left[\stackrel{s}{s+m-1}\right] \times k }(t,\eta  ) )
 =  r( D^{\left[\stackrel{s-1}{s+m-1}\right] \times (k-1) }(t,\eta  ) ) +1   \}. \end{array} $$
 
    By \eqref{eq 6.7} we have,  observing that   $ \theta _{1} (t,\eta )\cdot  \theta _{2}^{q-1} (t,\eta )  $ is constant on cosets of $ \mathbb{P}_{k+s-1}\times \mathbb{P}_{k+s+m-1} $
$$ \int_{\mathbb{P}\times \mathbb{P}}\theta _{1} (t,\eta )\cdot  \theta _{2}^{q-1} (t,\eta )  dt\eta 
 = \int_{  \mathbb{A}}\theta _{1} (t,\eta )\cdot  \theta _{2}^{q-1} (t,\eta )  dt\eta 
    +  \int_{ \mathbb{B}} \theta _{1} (t,\eta )\cdot  \theta _{2}^{q-1} (t,\eta )  dt\eta  $$
$$ =\sum_{(t,\eta ) \in \mathbb{A}\bigcap\big(\mathbb{P}/\mathbb{P}_{k+s-1}\times\mathbb{P}/\mathbb{P}_{k+s+m-1}\big) }
2^{q-1}\cdot2^{(2s+m+k-3-  r( D^{\left[\stackrel{s-1}{s+m-1}\right] \times (k-1) }(t,\eta  ) ) q }\int_{\mathbb{P}_{s+k-1}}dt\int_{\mathbb{P}_{s+k+m-1}}d\eta  $$
$$  + \sum_{(t,\eta ) \in \mathbb{B}\bigcap\big(\mathbb{P}/\mathbb{P}_{k+s-1}\times\mathbb{P}/\mathbb{P}_{k+s+m-1}\big) }
-2^{q-1}\cdot2^{(2s+m+k-3-  r( D^{\left[\stackrel{s-1}{s+m-1}\right] \times (k-1) }(t,\eta  )) q  }\int_{\mathbb{P}_{s+k-1}}dt\int_{\mathbb{P}_{s+k+m-1}}d\eta  $$ 
=  \begin{align*}
 \sum_{j=0}^{\inf(2s+m-2,k-1)} 2^{q-1}\cdot 2^{(2s+m+k-3-j) q}\cdot 
{}^{\#}\left(\begin{array}{c | c}
           j & j \\
           \hline
           j & j\\
           \hline
            j & j 
           \end{array} \right)_{\mathbb{P}/\mathbb{P}_{k+s -1}\times
           \mathbb{P}/\mathbb{P}_{k+s+m-1} }^{{\alpha  \over \beta }}\cdot2^{-(2k+2s+m-2)} \\
               +  \sum_{j=0}^{\inf(2s+m-2,k-1)}2^{q-1}\cdot 2^{(2s+m+k-3-j)q}\cdot {}^{\#}\left(\begin{array}{c | c}
           j & j \\
           \hline
           j & j\\
           \hline
            j & j +1
           \end{array} \right)_{\mathbb{P}/\mathbb{P}_{k+s -1}\times
           \mathbb{P}/\mathbb{P}_{k+s+m-1} }^{{\alpha  \over \beta }}\cdot2^{-(2k+2s+m-2)} \\
           -  \sum_{j=0}^{\inf(2s+m-2,k-1)}2^{q-1}\cdot 2^{(2s+m+k-3-j) q}\cdot {}^{\#}\left(\begin{array}{c | c}
           j & j \\
           \hline
           j & j+1 \\
           \hline
            j & j +1
           \end{array} \right)_{\mathbb{P}/\mathbb{P}_{k+s -1}\times
           \mathbb{P}/\mathbb{P}_{k+s+m-1} }^{{\alpha  \over \beta }}\cdot2^{-(2k+2s+m-2)}. \\
           &   
           \end{align*}

 Now by Lemma \ref{lem 6.9} we get  for all  $ q \geq 2 $ \\
 
  \begin{align*}
  \sum_{j=0}^{\inf(2s+m-2,k-1)}2^{-j q}\cdot\Bigg[{}^{\#}\left(\begin{array}{c | c}
           j & j \\
           \hline
           j & j\\
           \hline
            j & j 
           \end{array} \right)_{\mathbb{P}/\mathbb{P}_{k+s -1}\times
           \mathbb{P}/\mathbb{P}_{k+s+m-1} }^{{\alpha  \over \beta }}
            +    
{}^{\#}\left(\begin{array}{c | c}
           j & j \\
           \hline
           j & j\\
           \hline
            j & j +1
           \end{array} \right)_{\mathbb{P}/\mathbb{P}_{k+s -1}\times
           \mathbb{P}/\mathbb{P}_{k+s+m-1} }^{{\alpha  \over \beta }} \\
           - {}^{\#}\left(\begin{array}{c | c}
           j & j \\
           \hline
           j & j+1 \\
           \hline
            j & j +1
           \end{array} \right)_{\mathbb{P}/\mathbb{P}_{k+s -1}\times
           \mathbb{P}/\mathbb{P}_{k+s+m-1} }^{{\alpha  \over \beta }}\Bigg] = 0.
            \end{align*}\\
      \end{proof}

 \begin{lem}
\label{lem 6.11}
Let $ (t,\eta )\in \mathbb{P}\times \mathbb{P}, $and q be a rational integer $ \geq 2. $ \\
 Set 
\begin{align*}
  \phi_{3} (t,\eta ) &  =  \sum_{deg Y \leq k-2}\sum_{deg Z  = s-1}E(tYZ)\sum_{deg U \leq s+m-2}E(\eta YU).
\end{align*}
Then we have 
  \begin{equation}
 \label{eq 6.9}
   \phi_{3}^{q-1} (t,\eta )      =   \begin{cases}
 2^{(2s +m +k - 3 - j)\cdot(q-1) }  & \text{if }\quad  (t,\eta ) \in    \begin{vmatrix}
k-1 & \vline & k  & \vline \\
\hline
j &\vline & \cdot & \vline & D^{\left[s-1\atop s+m-1 \right]\times \cdot}  \\
\hline
j &\vline & \cdot &  \vline  & \alpha _{s}-\\
\hline 
\cdot &\vline & \cdot & \vline  & \beta _{s+m}- 
\end{vmatrix} \\
& \\
  0  & \text{otherwise }.
    \end{cases}           
  \end{equation}
  \end{lem}  
\begin{proof}
Lemma \ref{lem 6.11} follows from  \eqref{eq 4.4} with  $ k \rightarrow  k-1. $
\end{proof}
 \begin{lem}
\label{lem 6.12}
 Let $ (t,\eta)\in\mathbb{P}^{2}, $ and  q be a rational integer $ \geq 2, $  then
$  \phi (t,\eta )\cdot  \theta _{3}^{q-1} (t,\eta )  $ is equal to 
\begin{equation}
\label{eq 6.10}
  \begin{cases}
 2\cdot2^{(2s+m+k-3-j)q}    &  \text{if       }\quad (t,\eta ) \in    \begin{vmatrix}
k-1 & \vline & k  & \vline \\
\hline
j &\vline & j & \vline & D^{\left[s-1\atop s+m-1 \right]\times \cdot}  \\
\hline
j &\vline & j &  \vline  & \alpha _{s}-\\
\hline 
j &\vline & j & \vline  & \beta _{s+m}- 
\end{vmatrix} \\
 & \\
 -  2\cdot2^{(2s+m+k-3-j)q}     &  \text{if     } \quad (t,\eta ) \in    \begin{vmatrix}
k-1 & \vline & k  & \vline \\
\hline
j &\vline & j & \vline & D^{\left[s-1\atop s+m-1 \right]\times \cdot}  \\
\hline
j &\vline & j &  \vline  & \alpha _{s}-\\
\hline 
j &\vline & j +1 & \vline  & \beta _{s+m}- 
\end{vmatrix} \\
      0    &   \text{otherwise}.
\end{cases}
\end{equation}
  \end{lem}
  
\begin{proof}Recall that by Lemma \ref{lem 6.4} $ \phi (t,\eta ) $ is equal to 
\begin{equation*}
 \begin{cases}
   2^{2s+m+k-j-2} &  \text{if       }\quad (t,\eta ) \in    \begin{vmatrix}
k-1 & \vline & k  & \vline \\
\hline
\cdot &\vline & \cdot & \vline & D^{\left[s-1\atop s+m-1 \right]\times \cdot}  \\
\hline
j &\vline & j &  \vline  & \alpha _{s}-\\
\hline 
j &\vline & j  & \vline  & \beta _{s+m}- 
\end{vmatrix} \\
 &  \\
      -  2^{2s+m+k-j-2} &  \text{if     } \quad (t,\eta ) \in    \begin{vmatrix}
k-1 & \vline & k  & \vline \\
\hline
\cdot &\vline & \cdot & \vline & D^{\left[s-1\atop s+m-1 \right]\times \cdot}  \\
\hline
j &\vline & j &  \vline  & \alpha _{s}-\\
\hline 
j &\vline & j +1 & \vline  & \beta _{s+m}- 
\end{vmatrix} \\
     0    &   \text{otherwise}.
     \end{cases}
\end{equation*}
Then by \eqref{eq 6.9} we obtain by elementary rank properties  \eqref{eq 6.10}  and Lemma  \ref{lem 6.12}  is proved.\\
\end{proof}

 \begin{lem}
\label{lem 6.13} 
We have 
$$\int_{\mathbb{P}\times \mathbb{P}} \phi (t,\eta )\cdot  \theta _{3}^{q-1} (t,\eta )  dtd\eta =  0. $$
\end{lem}
\begin{proof}The integral  $\int_{\mathbb{P}\times \mathbb{P}} (t,\eta )  \phi (t,\eta )\cdot  \theta _{3}^{q-1} (t,\eta )  dtd\eta $ is equal to the 
number of solutions \\
 $ (Y, Z,U, Y_{1},Z_{1},U_{1},\ldots, Y_{q-1},Z_{q-1},U_{q-1}) $
of the polynomial equations
    \[\left\{\begin{array}{cc}
YZ + Y_{1}Z_{1} + \ldots + Y_{q-1}Z_{q-1}= 0, \\
YU + Y_{1}U_{1} + \ldots + Y_{q-1}U_{q-1}= 0,
    \end{array}\right.\]\\ 
 satisfying the degree conditions \\
$ \deg Y = k-1,\quad  \deg Z \leq  s-1,  \quad  \deg U = s+m-1,  $ \vspace{0.1 cm}\\
$ \deg Y_{i} \leq  k-2, \quad  \deg Z_{i} = s-1,  \quad  \deg U_{i} \leq  s+m-2  \quad for 1\leq i\leq q-1.  $ \vspace{0.1 cm}\\
By degree considerations  we have that  $ \deg \left[YU + Y_{1}U_{1} + \ldots + Y_{q-1}U_{q-1} \right] $ is equal to k+s+m-2. \\
 Lemma \ref{lem 6.13} follows.\\
 \end{proof}
   \begin{lem}
\label{lem 6.14}
We have for all $ j\in \mathbb{Z} $ such that $ 0\leq j\leq \inf(2s+m-2, k-1). $
\begin{align*}
{}^{\#}\left(\begin{array}{c | c}
           j & j \\
           \hline
           j & j\\
           \hline
            j & j 
           \end{array} \right)_{\mathbb{P}/\mathbb{P}_{k+s -1}\times
           \mathbb{P}/\mathbb{P}_{k+s+m-1} }^{{\alpha  \over \beta }}
            -   
{}^{\#}\left(\begin{array}{c | c}
           j & j \\
           \hline
           j & j\\
           \hline
            j & j +1
           \end{array} \right)_{\mathbb{P}/\mathbb{P}_{k+s -1}\times
           \mathbb{P}/\mathbb{P}_{k+s+m-1} }^{{\alpha  \over \beta }}  = 0
 \end{align*}\vspace{0.1 cm}\\
    
  Recall that   
\begin{align*}
{}^{\#}\left(\begin{array}{c | c}
           j_{1} & j_{2} \\
           \hline
           j_{3} & j_{4}\\
           \hline
            j_{5} & j_{6} 
           \end{array} \right)_{\mathbb{P}/\mathbb{P}_{k+s -1}\times
           \mathbb{P}/\mathbb{P}_{k+s+m-1} }^{{\alpha  \over \beta }}
            \end{align*}
    denotes  the cardinality of the following set
    $$\begin{array}{l}\Big\{ (t,\eta ) \in \mathbb{P}/\mathbb{P}_{k+s -1}\times
           \mathbb{P}/\mathbb{P}_{k+s+m-1} 
\mid r(  D^{\left[\stackrel{s-1}{s+m-1}\right] \times (k-1) }(t,\eta  ) ) = j_{1}, \quad 
r( D^{\left[\stackrel{s-1}{s+m-1}\right] \times k }(t,\eta  ) ) = j_{2},  \\
 r(  D^{\left[\stackrel{s}{s+m-1}\right] \times (k-1) }(t,\eta  )  ) = j_{3},\quad  
  r( D^{\left[\stackrel{s}{s+m-1}\right] \times k }(t,\eta  ) ) = j_{4}, \\
   r(  D^{\left[\stackrel{s}{s+m}\right] \times (k-1) }(t,\eta  )  ) = j_{5}, \quad
  r( D^{\left[\stackrel{s}{s+m}\right] \times k }(t,\eta  ) ) = j_{6} \Big\}
    \end{array}$$\\
 for $  (j_{1},j_{2},j_{3},j_{4},j_{5},j_{6}) \in \mathbb{N}^{6}. $ 
  \end{lem}
   
 \begin{proof} 
    We define  
  $$\begin{array}{l} \mathbb{A} =  \{(t,\eta ) \in \mathbb{P}^{2}
\mid  r( D^{\left[\stackrel{s-1}{s+m-1}\right] \times (k-1) }(t,\eta  ) )  
  =  r( D^{\left[\stackrel{s-1}{s+m-1}\right] \times k }(t,\eta  )  ) = r( D^{\left[\stackrel{s}{s+m-1}\right] \times (k-1) }(t,\eta  ) )  \\
 =   r( D^{\left[\stackrel{s}{s+m-1}\right] \times k }(t,\eta  ) )  =  r( D^{\left[\stackrel{s}{s+m}\right] \times (k-1) }(t,\eta  ) )
  = r( D^{\left[\stackrel{s}{s+m}\right] \times k }(t,\eta  ) ) \}. \end{array} $$\\
 and 
  $$\begin{array}{l}  \mathbb{B} =  \{(t,\eta ) \in \mathbb{P}^{2}
\mid  r( D^{\left[\stackrel{s-1}{s+m-1}\right] \times (k-1) }(t,\eta  ) ) 
  =  r( D^{\left[\stackrel{s-1}{s+m-1}\right] \times k }(t,\eta  )  )= 
  r( D^{\left[\stackrel{s}{s+m-1}\right] \times (k-1) }(t,\eta  ) ) \\
  = r( D^{\left[\stackrel{s}{s+m-1}\right] \times k }(t,\eta  ) ) =  r( D^{\left[\stackrel{s}{s+m}\right] \times (k-1) }(t,\eta  ) ),   \\
 r( D^{\left[\stackrel{s}{s+m}\right] \times k }(t,\eta  ) )  =  r( D^{\left[\stackrel{s-1}{s+m-1}\right] \times (k-1) }(t,\eta  ) ) +1  \}. \end{array} $$\\
 
  By \eqref{eq 6.10} we have,  observing that   $    \phi (t,\eta )\cdot  \theta _{3}^{q-1} (t,\eta ) $ is constant on cosets of $ \mathbb{P}_{k+s-1}\times \mathbb{P}_{k+s+m-1} $
$$ \int_{\mathbb{P}\times \mathbb{P}} \phi (t,\eta )\cdot  \theta _{3}^{q-1} (t,\eta )dt\eta 
 = \int_{  \mathbb{A}}\phi (t,\eta )\cdot  \theta _{3}^{q-1} (t,\eta )  dt\eta 
    +  \int_{ \mathbb{B}}\phi (t,\eta )\cdot  \theta _{3}^{q-1} (t,\eta )  dt\eta  $$
$$ =\sum_{(t,\eta ) \in \mathbb{A}\bigcap\big(\mathbb{P}/\mathbb{P}_{k+s-1}\times\mathbb{P}/\mathbb{P}_{k+s+m-1}\big) }
2\cdot 2^{(2s+m+k-3-  r( D^{\left[\stackrel{s-1}{s+m-1}\right] \times (k-1) }(t,\eta  ) ) q }\int_{\mathbb{P}_{s+k-1}}dt\int_{\mathbb{P}_{s+k+m-1}}d\eta  $$
$$  + \sum_{(t,\eta ) \in \mathbb{B}\bigcap\big(\mathbb{P}/\mathbb{P}_{k+s-1}\times\mathbb{P}/\mathbb{P}_{k+s+m-1}\big) }
- 2\cdot2^{(2s+m+k-3-  r( D^{\left[\stackrel{s-1}{s+m-1}\right] \times (k-1) }(t,\eta  )) q  }\int_{\mathbb{P}_{s+k-1}}dt\int_{\mathbb{P}_{s+k+m-1}}d\eta  $$ 
=  \begin{align*}
 \sum_{j=0}^{\inf(2s+m-2,k-1)} 2\cdot 2^{(2s+m+k-3-j) q}\cdot 
{}^{\#}\left(\begin{array}{c | c}
           j & j \\
           \hline
           j & j\\
           \hline
            j & j 
           \end{array} \right)_{\mathbb{P}/\mathbb{P}_{k+s -1}\times
           \mathbb{P}/\mathbb{P}_{k+s+m-1} }^{{\alpha  \over \beta }}\cdot2^{-(2k+2s+m-2)} \\
               -  \sum_{j=0}^{\inf(2s+m-2,k-1)}2\cdot 2^{(2s+m+k-3-j)q}\cdot {}^{\#}\left(\begin{array}{c | c}
           j & j \\
           \hline
           j & j\\
           \hline
            j & j +1
           \end{array} \right)_{\mathbb{P}/\mathbb{P}_{k+s -1}\times
           \mathbb{P}/\mathbb{P}_{k+s+m-1} }^{{\alpha  \over \beta }}\cdot2^{-(2k+2s+m-2)}. 
          \end{align*}

 Now by Lemma \ref{lem 6.13} we get  for all  $ q \geq 2 $ \\
 
  \begin{align*}
  \sum_{j=0}^{\inf(2s+m-2,k-1)}2^{-j q}\cdot\Bigg[{}^{\#}\left(\begin{array}{c | c}
           j & j \\
           \hline
           j & j\\
           \hline
            j & j 
           \end{array} \right)_{\mathbb{P}/\mathbb{P}_{k+s -1}\times
           \mathbb{P}/\mathbb{P}_{k+s+m-1} }^{{\alpha  \over \beta }}
            -   
{}^{\#}\left(\begin{array}{c | c}
           j & j \\
           \hline
           j & j\\
           \hline
            j & j +1
           \end{array} \right)_{\mathbb{P}/\mathbb{P}_{k+s -1}\times
           \mathbb{P}/\mathbb{P}_{k+s+m-1} }^{{\alpha  \over \beta }}\Bigg] = 0. \\
           \end{align*}\vspace{0.1 cm}\\
      \end{proof}

  \section{\textbf{RANK PROPERTIES OF SUBMATRICES OF DOUBLE PERSYMMETRIC MATRICES}}
   \label{sec 7}
   Consider the following partition of the matrix $  D^{\left[\stackrel{s-1}{\stackrel{s+m-1 }
{\overline {\stackrel{\alpha  _{s -}}{\beta  _{s +m -} }}}}\right] \times k }(t,\eta )  $ \\

         $$   \left ( \begin{array} {ccccc|c}
\alpha _{1} & \alpha _{2} & \alpha _{3} &  \ldots & \alpha _{k-1}  &  \alpha _{k} \\
\alpha _{2 } & \alpha _{3} & \alpha _{4}&  \ldots  &  \alpha _{k} &  \alpha _{k+1} \\
\vdots & \vdots & \vdots    &  \vdots & \vdots  &  \vdots \\
\alpha _{s-1} & \alpha _{s} & \alpha _{s +1} & \ldots  &  \alpha _{s+k-3} &  \alpha _{s+k-2}  \\
 \beta  _{1} & \beta  _{2} & \beta  _{3} & \ldots  &  \beta_{k-1} &  \beta _{k}  \\
\beta  _{2} & \beta  _{3} & \beta  _{4} & \ldots  &  \beta_{k} &  \beta _{k+1}  \\
\vdots & \vdots & \vdots    &  \vdots & \vdots  &  \vdots \\
\beta  _{m+1} & \beta  _{m+2} & \beta  _{m+3} & \ldots  &  \beta_{k+m-1} &  \beta _{k+m}  \\
\vdots & \vdots & \vdots    &  \vdots & \vdots  &  \vdots \\
\beta  _{s+m-1} & \beta  _{s+m} & \beta  _{s+m+1} & \ldots  &  \beta_{s+m+k-3} &  \beta _{s+m+k-2}  \\
\hline
\alpha _{s} & \alpha _{s+1} & \alpha _{s +2} & \ldots  &  \alpha _{s+k-2} &  \alpha _{s+k-1} \\
\hline
\beta  _{s+m} & \beta  _{s+m+1} & \beta  _{s+m+2} & \ldots  &  \beta_{s+m+k-2} &  \beta _{s+m+k-1}
\end{array}  \right). $$ \vspace{0.2 cm} \\

 By studying rank properties of submatrices of the above  double persymmetric matrix, we deduce by contradiction  the following rank formula  for all $ i \in [0,\inf(2s+m-3,k-2)] $ \vspace{0.2 cm}
   \begin{align*}
 {}^{\#}\left(\begin{array}{c | c}
           i & i +1 \\
           \hline
           i & i +1 \\
           \hline
            i & i +1
           \end{array} \right)_{\mathbb{P}/\mathbb{P}_{k+s -1}\times
           \mathbb{P}/\mathbb{P}_{k+s+m-1} }^{{\alpha  \over \beta }} = 0. 
\end{align*}\vspace{0.2 cm}\\
  Recall that   
\begin{align*}
{}^{\#}\left(\begin{array}{c | c}
           i & i+1 \\
           \hline
           i & i+1\\
           \hline
            i & i+1 
           \end{array} \right)_{\mathbb{P}/\mathbb{P}_{k+s -1}\times
           \mathbb{P}/\mathbb{P}_{k+s+m-1} }^{{\alpha  \over \beta }}
            \end{align*}
    denotes  the cardinality of the following set
    $$\begin{array}{l}\Big\{ (t,\eta ) \in \mathbb{P}/\mathbb{P}_{k+s -1}\times
           \mathbb{P}/\mathbb{P}_{k+s+m-1} 
\mid r(  D^{\left[\stackrel{s-1}{s+m-1}\right] \times (k-1) }(t,\eta  ) ) = i, \quad 
r( D^{\left[\stackrel{s-1}{s+m-1}\right] \times k }(t,\eta  ) ) = i+1,  \\
 r(  D^{\left[\stackrel{s}{s+m-1}\right] \times (k-1) }(t,\eta  )  ) = i,\quad  
  r( D^{\left[\stackrel{s}{s+m-1}\right] \times k }(t,\eta  ) ) = i+1, \\
   r(  D^{\left[\stackrel{s}{s+m}\right] \times (k-1) }(t,\eta  )  ) = i, \quad
  r( D^{\left[\stackrel{s}{s+m}\right] \times k }(t,\eta  ) ) = i+1 \Big\}. 
    \end{array} $$
     \begin{lem}
   \label{lem 7.1}
  For all i such that $ 0\leq i\leq \inf(2s+m-3, k-2) $ we have  \begin{align*}
 {}^{\#}\left(\begin{array}{c | c}
           i & i +1 \\
           \hline
           i & i +1 \\
           \hline
            i & i +1
           \end{array} \right)_{\mathbb{P}/\mathbb{P}_{k+s -1}\times
           \mathbb{P}/\mathbb{P}_{k+s+m-1} }^{{\alpha  \over \beta }} = 0.
\end{align*}
 \end{lem}
  \begin{proof}   
    $Set \;(t,\eta )= (\sum_{i\geq 1}\alpha _{i}T^{-i},\sum_{i\geq 1}\beta  _{i}T^{-i})
\in \mathbb{P}\times\mathbb{P} \; and\; let \; (s,k,m,j)\in \mathbb{N}^{*}\times\mathbb{N}^{*}
\times\mathbb{N}\times\mathbb{N}^{*},\\
where  \quad 1\leq j\leq k-1. $\\
 We denote by  $  D_{j}^{\left[\stackrel{s}{s+m}\right] \times (k-j+1) }(t,\eta  )  $ the following 
$ (2s +m)\times (k-j+1) $ matrix 
 $$   \left ( \begin{array} {ccccc|c}
\alpha _{j} & \alpha _{j+1} & \alpha _{j+2} &  \ldots & \alpha _{k-1}  &  \alpha _{k} \\
 \beta  _{j} & \beta  _{j+1} & \beta  _{j+2} & \ldots  &  \beta_{k-1} &  \beta _{k}  \\
\alpha _{j+1 } & \alpha _{j+2} & \alpha _{j+3}&  \ldots  &  \alpha _{k} &  \alpha _{k+1} \\
\beta  _{j+1} & \beta  _{j+2} & \beta  _{j+3} & \ldots  &  \beta_{k} &  \beta _{k+1}  \\
\vdots & \vdots & \vdots    &  \vdots & \vdots  &  \vdots \\
\alpha _{j+s-2} & \alpha _{j+s-1} & \alpha _{j+s} & \ldots  &  \alpha _{k+s-3} &  \alpha _{k+s-2}  \\
\beta  _{j+s-2} & \beta  _{j+s-1} & \beta  _{j+s} & \ldots  &  \beta_{k+s-3} &  \beta _{k+s-2}  \\
\beta  _{j+s-1} & \beta  _{j+s} & \beta  _{j+s+1} & \ldots  &  \beta_{k+s-2} &  \beta _{k+s-1}  \\
\beta  _{j+s} & \beta  _{j+s+1} & \beta  _{j+s+2} & \ldots  &  \beta_{k+s-1} &  \beta _{k+s}  \\
\vdots & \vdots & \vdots    &  \vdots & \vdots  &  \vdots \\
\beta  _{j+s+m-3} & \beta  _{j+s+m-2} & \beta  _{j+s+m-1} & \ldots  &  \beta_{k+s+m -4} &  \beta _{k+s+m-3}  \\
\beta  _{j+s+m-2} & \beta  _{j+s+m-1} & \beta  _{j+s+m} & \ldots  &  \beta_{k+s+m-3} &  \beta _{s+m+k-2}  \\
\hline \\
\alpha _{j+s-1} & \alpha _{j+s} & \alpha _{j+s+1} & \ldots  &  \alpha _{k+s-2} &  \alpha _{k+s-1} \\
\beta  _{j+s+m-1} & \beta  _{j+s+m} & \beta  _{j+s+m+1} & \ldots  &  \beta_{k+s+m-2} &  \beta _{k+s+m-1}
\end{array}  \right). $$

Remark that after a rearrangement of the rows in the above matrix we obtain the following double persymmetric matrix
where the first s rows form a $ s\times (k-j+1) $ persymmetric matrix and the last s+m rows form a $ (s+m)\times (k-j+1) $
persymmetric matrix with entries in  $ \mathbb{F}_{2} $

    $$   \left ( \begin{array} {ccccc|c}
\alpha _{j} & \alpha _{j+1} & \alpha _{j+2} &  \ldots & \alpha _{k-1}  &  \alpha _{k} \\
\alpha _{j+1 } & \alpha _{j+2} & \alpha _{j+3}&  \ldots  &  \alpha _{k} &  \alpha _{k+1} \\
\vdots & \vdots & \vdots    &  \vdots & \vdots  &  \vdots \\
\alpha _{j+s-2} & \alpha _{j+s-1} & \alpha _{j+s} & \ldots  &  \alpha _{k+s-3} &  \alpha _{k+s-2}  \\
\alpha _{j+s-1} & \alpha _{j+s} & \alpha _{j+s+1} & \ldots  &  \alpha _{k+s-2} &  \alpha _{k+s-1}\\
\beta  _{j} & \beta  _{j+1} & \beta  _{j+2} & \ldots  &  \beta_{k-1} &  \beta _{k}  \\
\beta  _{j+1} & \beta  _{j+2} & \beta  _{j+3} & \ldots  &  \beta_{k} &  \beta _{k+1}  \\
\vdots & \vdots & \vdots    &  \vdots & \vdots  &  \vdots \\
\beta  _{j+s-2} & \beta  _{j+s-1} & \beta  _{j+s} & \ldots  &  \beta_{k+s-3} &  \beta _{k+s-2}  \\
\beta  _{j+s-1} & \beta  _{j+s} & \beta  _{j+s+1} & \ldots  &  \beta_{k+s-2} &  \beta _{k+s-1}  \\
\beta  _{j+s} & \beta  _{j+s+1} & \beta  _{j+s+2} & \ldots  &  \beta_{k+s-1} &  \beta _{k+s}  \\
\vdots & \vdots & \vdots    &  \vdots & \vdots  &  \vdots \\
\beta  _{j+s+m-3} & \beta  _{j+s+m-2} & \beta  _{j+s+m-1} & \ldots  &  \beta_{k+s+m -4} &  \beta _{k+s+m-3}  \\
\beta  _{j+s+m-2} & \beta  _{j+s+m-1} & \beta  _{j+s+m} & \ldots  &  \beta_{k+s+m-3} &  \beta _{k+s+m-2}  \\
\beta  _{j+s+m-1} & \beta  _{j+s+m} & \beta  _{j+s+m+1} & \ldots  &  \beta_{k+s+m-2} &  \beta _{k+s+m-1}
\end{array}  \right). $$ \\

Proof by contradiction.\\
Assume on the contrary  that there exists $ i_{0} \in [0,\inf(2s+m-3, k-2)] $  such that
 \begin{equation}
\label{eq 7.1}
{}^{\#}\left(\begin{array}{c | c}
           i_{0} & i_{0}+1 \\
           \hline
           i_{0} & i_{0}+1\\
           \hline
            i_{0} & i_{0}+1 
           \end{array} \right)_{\mathbb{P}/\mathbb{P}_{k+s -1}\times
           \mathbb{P}/\mathbb{P}_{k+s+m-1} }^{{\alpha  \over \beta }} > 0.
 \end{equation}
 We are going to show that\\

 \begin{align*}
& {}^{\#}\left(\begin{array}{c | c}
           i_{0} & i_{0}+1 \\
           \hline
           i_{0} & i_{0}+1\\
           \hline
            i_{0} & i_{0}+1 
           \end{array} \right)_{\mathbb{P}/\mathbb{P}_{k+s -1}\times
           \mathbb{P}/\mathbb{P}_{k+s+m-1} }^{{\alpha  \over \beta }} > 0
 \Longrightarrow 
 {}^{\#}\left(\begin{array}{c | c}
           i_{0}-1 & i_{0} \\
           \hline
           i_{0}-1 & i_{0}\\
           \hline
            i_{0}-1 & i_{0}
           \end{array} \right)_{\mathbb{P}_{1}/\mathbb{P}_{k+s -1}\times
           \mathbb{P}_{1}/\mathbb{P}_{k+s+m-1} }^{{\alpha  \over \beta }} > 0   \\
&  \Longrightarrow  \ldots  \Longrightarrow
  {}^{\#}\left(\begin{array}{c | c}
           0 & 1 \\
           \hline
           0 & 1 \\
           \hline
            0 & 1 
           \end{array} \right)_{\mathbb{P}_{i_{0}}/\mathbb{P}_{k+s -1}\times
           \mathbb{P}_{i_{0}}/\mathbb{P}_{k+s+m-1} }^{{\alpha  \over \beta }} > 0,
\end{align*}
  which obviously contradicts
 \begin{equation*}
{}^{\#}\left(\begin{array}{c | c}
           0 & 1 \\
           \hline
           0 & 1 \\
           \hline
            0 & 1 
           \end{array} \right)_{\mathbb{P}_{i_{0}}/\mathbb{P}_{k+s -1}\times
           \mathbb{P}_{i_{0}}/\mathbb{P}_{k+s+m-1} }^{{\alpha  \over \beta }} = 0.
 \end{equation*} 
 
  By   \eqref{eq 7.1} there exists $( t_{0},\eta_{0}) \in \mathbb{P}/\mathbb{P}_{k +s -1}\times
           \mathbb{P}/\mathbb{P}_{k+s+m-1}  $ such that \\
    \begin{align*}
&  r(  D^{\left[\stackrel{s-1}{s+m-1}\right] \times (k-1) }( t_{0},\eta_{0})) = 
      r(  D^{\left[\stackrel{s}{s+m-1}\right] \times (k-1) }( t_{0},\eta_{0})  )       
  =     r(  D^{\left[\stackrel{s}{s+m}\right] \times (k-1) } ( t_{0},\eta_{0}) ) = i_{0},  \\
  &    r(  D^{\left[\stackrel{s-1}{s+m-1}\right] \times k }( t_{0},\eta_{0})) = 
      r(  D^{\left[\stackrel{s}{s+m-1}\right] \times k }( t_{0},\eta_{0})  )       
  =     r(  D^{\left[\stackrel{s}{s+m}\right] \times k } ( t_{0},\eta_{0}) ) = i_{0} +1.    
       \end{align*}
 Set   $ (t_{0},\eta_{0} )= (\sum_{i\geq 1}\alpha _{i}T^{-i},\sum_{i\geq 1}\beta  _{i}T^{-i}), $   and 
 consider  the following partition of the matrix  $  D^{\left[\stackrel{s}{s+m}\right] \times k }(t_{0},\eta_{0}  )  $ 

 $$   \left ( \begin{array} {c|cccc|c}
\alpha _{1} & \alpha _{2} & \alpha _{3} &  \ldots & \alpha _{k-1}  &  \alpha _{k} \\
 \beta  _{1} & \beta  _{2} & \beta  _{3} & \ldots  &  \beta_{k-1} &  \beta _{k}  \\
 \hline
\alpha _{2 } & \alpha _{3} & \alpha _{4}&  \ldots  &  \alpha _{k} &  \alpha _{k+1} \\
\beta  _{2} & \beta  _{3} & \beta  _{4} & \ldots  &  \beta_{k} &  \beta _{k+1}  \\
\vdots & \vdots & \vdots    &  \vdots & \vdots  &  \vdots \\
\alpha _{s-1} & \alpha _{s} & \alpha _{s+1} & \ldots  &  \alpha _{k+s-3} &  \alpha _{k+s-2}  \\
\beta  _{s-1} & \beta  _{s} & \beta  _{s+1} & \ldots  &  \beta_{k+s-3} &  \beta _{k+s-2}  \\
\beta  _{s} & \beta  _{s+1} & \beta  _{s+2} & \ldots  &  \beta_{k+s-2} &  \beta _{k+s-1}  \\
\beta  _{s+1} & \beta  _{s+2} & \beta  _{s+3} & \ldots  &  \beta_{k+s-1} &  \beta _{k+s}  \\
\vdots & \vdots & \vdots    &  \vdots & \vdots  &  \vdots \\
\beta  _{s+m-2} & \beta  _{s+m-1} & \beta  _{s+m} & \ldots  &  \beta_{k+s+m -4} &  \beta _{k+s+m-3}  \\
\beta  _{s+m-1} & \beta  _{s+m} & \beta  _{s+m+1} & \ldots  &  \beta_{k+s+m-3} &  \beta _{k+s+m-2}  \\
\hline \\
\alpha _{s} & \alpha _{s+1} & \alpha _{s+2} & \ldots  &  \alpha _{k+s-2} &  \alpha _{k+s-1} \\
\beta  _{s+m} & \beta  _{s+m+1} & \beta  _{s+m+2} & \ldots  &  \beta_{k+s+m-2} &  \beta _{k+s+m-1}
\end{array}  \right). $$ \\
         
   By \eqref{eq 7.1} we have
     $ r(  D_{2}^{\left[\stackrel{s-1}{s+m-1}\right] \times (k-1) }( t_{0},\eta_{0})) \leq i_{0}\quad and \quad
   r(  D^{\left[\stackrel{s-1}{s+m-1}\right] \times k }( t_{0},\eta_{0})) =  i_{0}+1. $\\
    It follows that $  r(  D_{2}^{\left[\stackrel{s-1}{s+m-1}\right] \times (k-1) }( t_{0},\eta_{0})) = i_{0}. $\\
     Let $ a_{1},a_{2},\ldots,a_{k}$ denote the  columns of \quad $ D^{\left[\stackrel{s-1}{s+m-1}\right] \times k }( t_{0},\eta_{0})  $  that is,\\
     
\begin{displaymath}
a_{j} = \left(\begin{array}{c}
\alpha _{j} \\
 \beta  _{j}  \\
\alpha _{j+1 } \\
\beta  _{j+1}\\
\vdots \\
\alpha _{j+s-2}  \\
\beta  _{j+s-2} \\
\beta  _{j+s-1}  \\
\beta  _{j+s}  \\
\vdots \\
\beta  _{j+s+m-3}  \\
\beta  _{j+s+m-2} 
\end{array}\right).
\end{displaymath} \\

Since $ r(  D_{2}^{\left[\stackrel{s-1}{s+m-1}\right] \times (k-1) }( t_{0},\eta_{0})) = i_{0}\;
and \;  r(  D^{\left[\stackrel{s-1}{s+m-1}\right] \times k }( t_{0},\eta_{0})) =  i_{0}+1, $\\
 we  have 
$ a_{1}\notin span\left\{a_{2},a_{3},\ldots,a_{k}\right\},
    \text{therefore  }\;  r(  D_{2}^{\left[\stackrel{s-1}{s+m-1}\right] \times (k-2) }( t_{0},\eta_{0})) = i_{0}-1. $\\
  We have then\\
  \begin{align}
  r(  D_{2}^{\left[\stackrel{s-1}{s+m-1}\right] \times (k-2) }( t_{0},\eta_{0})) & = i_{0}-1,  \label{eq 7.2} \\
  & \nonumber \\
   r(  D^{\left[\stackrel{s-1}{s+m-1}\right] \times (k-1) }( t_{0},\eta_{0})) & =  i_{0}, \label{eq 7.3} \\
    & \nonumber \\
    r(  D^{\left[\stackrel{s}{s+m-1}\right] \times (k-1) }( t_{0},\eta_{0})) &  =  i_{0}, \label{eq 7.4}\\
    & \nonumber \\ 
     r(  D^{\left[\stackrel{s}{s+m}\right] \times (k-1) }( t_{0},\eta_{0})) &  =  i_{0} \label{eq 7.5} & \text{and, see above,}\\
      & \nonumber \\
      r(  D_{2}^{\left[\stackrel{s-1}{s+m-1}\right] \times (k-1) }( t_{0},\eta_{0})) &  = i_{0}. \label{eq 7.6}
 \end{align}
Now we consider the matrix obtained from the matrix  $  D^{\left[\stackrel{s}{s+m}\right] \times k }(t_{0},\eta_{0}  )  $ 
by delating the first column and replacing the last column by the first one :  
  $$   \left ( \begin{array} {cccc|c}
 \alpha _{2} & \alpha _{3} &  \ldots & \alpha _{k-1}  &  \alpha _{1} \\
  \beta  _{2} & \beta  _{3} & \ldots  &  \beta_{k-1} &  \beta _{1}  \\
  \alpha _{3} & \alpha _{4}&  \ldots  &  \alpha _{k} &  \alpha _{2} \\
\beta  _{3} & \beta  _{4} & \ldots  &  \beta_{k} &  \beta _{2}  \\
 \vdots & \vdots    &  \vdots & \vdots  &  \vdots \\
\alpha _{s} & \alpha _{s+1} & \ldots  &  \alpha _{k+s-3} &  \alpha _{s-1}  \\
 \beta  _{s} & \beta  _{s+1} & \ldots  &  \beta_{k+s-3} &  \beta _{s-1}  \\
 \beta  _{s+1} & \beta  _{s+2} & \ldots  &  \beta_{k+s-2} &  \beta _{s}  \\
 \beta  _{s+2} & \beta  _{s+3} & \ldots  &  \beta_{k+s-1} &  \beta _{s+1}  \\
\vdots & \vdots    &  \vdots & \vdots  &  \vdots \\
 \beta  _{s+m-1} & \beta  _{s+m} & \ldots  &  \beta_{k+s+m -4} &  \beta _{s+m-2}  \\
 \beta  _{s+m} & \beta  _{s+m+1} & \ldots  &  \beta_{k+s+m-3} &  \beta _{s+m-1}  \\
\hline \\
 \alpha _{s+1} & \alpha _{s+2} & \ldots  &  \alpha _{k+s-2} &  \alpha _{s} \\
 \beta  _{s+m+1} & \beta  _{s+m+2} & \ldots  &  \beta_{k+s+m-2} &  \beta _{s+m}
\end{array}  \right). $$   
From \eqref{eq 7.2}, \eqref{eq 7.3}, \eqref{eq 7.4} and \eqref{eq 7.5} we obtain by elementary rank considerations 
 \begin{align}
  r(  D_{2}^{\left[\stackrel{s}{s+m-1}\right] \times (k-2) }( t_{0},\eta_{0})) & = i_{0}-1,  \label{eq 7.7} \\
   & \nonumber \\
   r(  D_{2}^{\left[\stackrel{s}{s+m}\right] \times (k-2) }( t_{0},\eta_{0})) & =  i_{0}-1. \label{eq 7.8} \\
    & \nonumber 
   \end{align}
Consider now the matrix $ D_{2}^{\left[\stackrel{s}{s+m}\right] \times (k-1) }(t_{0},\eta_{0}  ) $ obtained by the matrix 
$ D^{\left[\stackrel{s}{s+m}\right] \times k }(t_{0},\eta_{0}  ) $ by delating the first column, that is \\
 
   $$   \left ( \begin{array} {cccc|c}
 \alpha _{2} & \alpha _{3} &  \ldots & \alpha _{k-1}  &  \alpha _{k} \\
  \beta  _{2} & \beta  _{3} & \ldots  &  \beta_{k-1} &  \beta _{k}  \\
  \alpha _{3} & \alpha _{4}&  \ldots  &  \alpha _{k} &  \alpha _{k+1} \\
\beta  _{3} & \beta  _{4} & \ldots  &  \beta_{k} &  \beta _{k+1}  \\
 \vdots & \vdots    &  \vdots & \vdots  &  \vdots \\
\alpha _{s} & \alpha _{s+1} & \ldots  &  \alpha _{k+s-3} &  \alpha _{k+s-2}  \\
 \beta  _{s} & \beta  _{s+1} & \ldots  &  \beta_{k+s-3} &  \beta _{k+s-2}  \\
 \beta  _{s+1} & \beta  _{s+2} & \ldots  &  \beta_{k+s-2} &  \beta _{k+s-1}  \\
 \beta  _{s+2} & \beta  _{s+3} & \ldots  &  \beta_{k+s-1} &  \beta _{k+s}  \\
\vdots & \vdots    &  \vdots & \vdots  &  \vdots \\
 \beta  _{s+m-1} & \beta  _{s+m} & \ldots  &  \beta_{k+s+m -4} &  \beta _{k+s+m-3}  \\
 \beta  _{s+m} & \beta  _{s+m+1} & \ldots  &  \beta_{k+s+m-3} &  \beta _{k+s+m-2}  \\
\hline \\
 \alpha _{s+1} & \alpha _{s+2} & \ldots  &  \alpha _{k+s-2} &  \alpha _{k+s-1} \\
 \beta  _{s+m+1} & \beta  _{s+m+2} & \ldots  &  \beta_{k+s+m-2} &  \beta _{k+s+m-1}
\end{array}  \right). $$ \\

 From \eqref{eq 7.2}, \eqref{eq 7.6}, \eqref{eq 7.7} and   \eqref{eq 7.8}, we obtain \\
$$  r(  D_{2}^{\left[\stackrel{s}{s+m-1}\right] \times (k-1) }( t_{0},\eta_{0}))  =  
    r(  D_{2}^{\left[\stackrel{s}{s+m}\right] \times (k-1) }( t_{0},\eta_{0})) =  i_{0}. $$
We get 
\begin{equation*}
 {}^{\#}\left(\begin{array}{c | c}
           i_{0}-1 & i_{0} \\
           \hline
           i_{0}-1 & i_{0}\\
           \hline
            i_{0}-1 & i_{0}
           \end{array} \right)_{\mathbb{P}_{1}/\mathbb{P}_{k+s -1}\times
           \mathbb{P}_{1}/\mathbb{P}_{k+s+m-1} }^{{\alpha  \over \beta }} > 0.  
\end{equation*}
  Recall that   
\begin{align*}
{}^{\#}\left(\begin{array}{c | c}
            i_{0}-1 & i_{0} \\
           \hline
           i_{0}-1 & i_{0}\\
           \hline
            i_{0}-1 & i_{0}
           \end{array} \right)_{\mathbb{P}_{1}/\mathbb{P}_{k+s -1}\times
           \mathbb{P}_{1}/\mathbb{P}_{k+s+m-1} }^{{\alpha  \over \beta }}
            \end{align*}
    denotes  the cardinality of the following set
    $$\begin{array}{l}\Big\{ (t,\eta ) \in \mathbb{P}_{1}/\mathbb{P}_{k+s -1}\times
           \mathbb{P}_{1}/\mathbb{P}_{k+s+m-1} 
\mid r(  D_{2}^{\left[\stackrel{s-1}{s+m-1}\right] \times (k-2) }(t,\eta  ) ) = i_{0}-1, \quad 
r( D_{2}^{\left[\stackrel{s-1}{s+m-1}\right] \times (k-1) }(t,\eta  ) ) = i_{0},  \\
 r(  D_{2}^{\left[\stackrel{s}{s+m-1}\right] \times (k-2) }(t,\eta  )  ) = i_{0}-1,\quad  
  r( D_{2}^{\left[\stackrel{s}{s+m-1}\right] \times (k-1) }(t,\eta  ) ) = i_{0}, \\
   r(  D_{2}^{\left[\stackrel{s}{s+m}\right] \times (k-2) }(t,\eta  )  ) = i_{0}-1, \quad
  r( D_{2}^{\left[\stackrel{s}{s+m}\right] \times (k-1) }(t,\eta  ) ) = i_{0} \Big\}.
    \end{array}$$\\
We have now proved that \\
 \begin{align*}
& {}^{\#}\left(\begin{array}{c | c}
           i_{0} & i_{0}+1 \\
           \hline
           i_{0} & i_{0}+1\\
           \hline
            i_{0} & i_{0}+1 
           \end{array} \right)_{\mathbb{P}/\mathbb{P}_{k+s -1}\times
           \mathbb{P}/\mathbb{P}_{k+s+m-1} }^{{\alpha  \over \beta }} > 0
 \Longrightarrow 
 {}^{\#}\left(\begin{array}{c | c}
           i_{0}-1 & i_{0} \\
           \hline
           i_{0}-1 & i_{0}\\
           \hline
            i_{0}-1 & i_{0}
           \end{array} \right)_{\mathbb{P}_{1}/\mathbb{P}_{k+s -1}\times
           \mathbb{P}_{1}/\mathbb{P}_{k+s+m-1} }^{{\alpha  \over \beta }} > 0.   
     \end{align*}
  
 We repeat this procedure and obtain after finitely many steps \\
  \begin{equation}
  \label{eq 7.9}
 {}^{\#}\left(\begin{array}{c | c}
           0 & 1 \\
           \hline
           0 & 1 \\
           \hline
            0 & 1 
           \end{array} \right)_{\mathbb{P}_{i_{0}}/\mathbb{P}_{k+s -1}\times
           \mathbb{P}_{i_{0}}/\mathbb{P}_{k+s+m-1} }^{{\alpha  \over \beta }} > 0.
 \end{equation} 
  From  \begin{displaymath}
   D_{i_{0}+1}^{\left[\stackrel{s}{s+m}\right] \times (k-i_{0}) }(t_{0},\eta_{0}  ) = \left ( \begin{array} {ccccc|c}
\alpha _{i_{0}+1} & \alpha _{i_{0}+2} & \alpha _{i_{0}+3} &  \ldots & \alpha _{k-1}  &  \alpha _{k} \\
 \beta  _{i_{0}+1} & \beta  _{i_{0}+2} & \beta  _{i_{0}+3} & \ldots  &  \beta_{k-1} &  \beta _{k}  \\
\alpha _{ i_{0}+2} & \alpha _{i_{0}+3} & \alpha _{i_{0}+4}&  \ldots  &  \alpha _{k} &  \alpha _{k+1} \\
\beta  _{i_{0}+2} & \beta  _{i_{0}+3} & \beta  _{i_{0}+4} & \ldots  &  \beta_{k} &  \beta _{k+1}  \\
\vdots & \vdots & \vdots    &  \vdots & \vdots  &  \vdots \\
\alpha _{i_{0}+s-1} & \alpha _{i_{0}+s} & \alpha _{i_{0}+s+1} & \ldots  &  \alpha _{k+s-3} &  \alpha _{k+s-2}  \\
\beta  _{i_{0}+s} & \beta  _{i_{0}+s +1} & \beta  _{i_{0}+s +2} & \ldots  &  \beta_{k+s-3} &  \beta _{k+s-2}  \\
\beta  _{i_{0}+s +1} & \beta  _{i_{0}+s +2} & \beta  _{i_{0}+s +3} & \ldots  &  \beta_{k+s-2} &  \beta _{k+s-1}  \\
\beta  _{i_{0}+s +2} & \beta  _{i_{0}+s +3} & \beta  _{i_{0}+s +4} & \ldots  &  \beta_{k+s-1} &  \beta _{k+s}  \\
\vdots & \vdots & \vdots    &  \vdots & \vdots  &  \vdots \\
\beta  _{i_{0}+s+m-2} & \beta  _{i_{0}+s+m-1} & \beta  _{i_{0}+s+m} & \ldots  &  \beta_{k+s+m -4} &  \beta _{k+s+m-3}  \\
\beta  _{i_{0}+s+m-1} & \beta  _{i_{0}+s+m} & \beta  _{i_{0}+s+m+1} & \ldots  &  \beta_{k+s+m-3} &  \beta _{s+m+k-2}  \\
\hline \\
\alpha _{i_{0}+s} & \alpha _{i_{0}+s +1} & \alpha _{i_{0}+s+2} & \ldots  &  \alpha _{k+s-2} &  \alpha _{k+s-1} \\
\beta  _{i_{0}+s+m} & \beta  _{i_{0}+s+m+1} & \beta  _{i_{0}+s+m+2} & \ldots  &  \beta_{k+s+m-2} &  \beta _{k+s+m-1},
\end{array}  \right). 
 \end{displaymath} \\
 
 we  obviously get 
 \begin{equation*}
   {}^{\#}\left(\begin{array}{c | c}
           0 & 1 \\
           \hline
           0 & 1 \\
           \hline
            0 & 1 
           \end{array} \right)_{\mathbb{P}_{i_{0}}/\mathbb{P}_{k+s -1}\times
           \mathbb{P}_{i_{0}}/\mathbb{P}_{k+s+m-1} }^{{\alpha  \over \beta }} = 0 
 \end{equation*}
  which clearly contradicts \eqref{eq 7.9}.\\
  
   Recall that   
\begin{align*}
{}^{\#}\left(\begin{array}{c | c}
           0 & 1 \\
           \hline
           0 & 1\\
           \hline
            0 & 1 
           \end{array} \right)_{\mathbb{P}_{i_{0}}/\mathbb{P}_{k+s -1}\times
           \mathbb{P}_{i_{0}}/\mathbb{P}_{k+s+m-1} }^{{\alpha  \over \beta }}
            \end{align*}
    denotes  the cardinality of the following set
    $$\begin{array}{l}\Big\{ (t,\eta ) \in \mathbb{P}_{i_{0}}/\mathbb{P}_{k+s -1}\times
           \mathbb{P}_{i_{0}}/\mathbb{P}_{k+s+m-1} 
\mid r(  D^{\left[\stackrel{s-1}{s+m-1}\right] \times (k- i_{0}-1) }(t,\eta  ) ) = 0, \quad 
r( D^{\left[\stackrel{s-1}{s+m-1}\right] \times (k-i_{0}) }(t,\eta  ) ) =  1,  \\
 r(  D^{\left[\stackrel{s}{s+m-1}\right] \times (k-i_{0}-1) }(t,\eta  )  ) = 0,\quad  
  r( D^{\left[\stackrel{s}{s+m-1}\right] \times (k-i_{0}) }(t,\eta  ) ) =  1, \\
   r(  D^{\left[\stackrel{s}{s+m}\right] \times (k-i_{0}-1) }(t,\eta  )  ) = 0, \quad
  r( D^{\left[\stackrel{s}{s+m}\right] \times (k-i_{0}) }(t,\eta  ) ) =  1 \Big\}.
    \end{array}$$
  \end{proof} 
 \begin{lem}
 \label{lem 7.2}We have \\
 \begin{equation}\\
 \label{eq 7.10}
 {}^{\#}\left(\begin{array}{c | c}
           i & i +1 \\
           \hline
           i & i +1 \\
           \hline
            i & i +1
           \end{array} \right)_{\mathbb{P}/\mathbb{P}_{k+s -1}\times
           \mathbb{P}/\mathbb{P}_{k+s+m-1} }^{{\alpha  \over \beta }} = 0 \quad   if \;  0\leq i\leq \inf(2s+m-3, k-2),  
           \end{equation}\\
           
\begin{equation}
\label{eq 7.11}
{}^{\#}\left(\begin{array}{c | c}
           i & i \\
           \hline
           i & i\\
           \hline
            i & i 
           \end{array} \right)_{\mathbb{P}/\mathbb{P}_{k+s -1}\times
           \mathbb{P}/\mathbb{P}_{k+s+m-1} }^{{\alpha  \over \beta }}
            = 
{}^{\#}\left(\begin{array}{c | c}
           i & i \\
           \hline
           i & i\\
           \hline
            i & i +1
           \end{array} \right)_{\mathbb{P}/\mathbb{P}_{k+s -1}\times
           \mathbb{P}/\mathbb{P}_{k+s+m-1} }^{{\alpha  \over \beta }}  \quad   if \;  0\leq i\leq \inf(2s+m-2, k-1),  
 \end{equation}\\
 
\begin{equation}
\label{eq 7.12}
{}^{\#}\left(\begin{array}{c | c}
           i & i \\
           \hline
           i & i+1\\
           \hline
            i & i +1
           \end{array} \right)_{\mathbb{P}/\mathbb{P}_{k+s -1}\times
           \mathbb{P}/\mathbb{P}_{k+s+m-1} }^{{\alpha  \over \beta }} =
               2\cdot{}^{\#}\left(\begin{array}{c | c}
           i & i \\
           \hline
           i & i\\
           \hline
            i & i
           \end{array} \right)_{\mathbb{P}/\mathbb{P}_{k+s -1}\times
           \mathbb{P}/\mathbb{P}_{k+s+m-1} }^{{\alpha  \over \beta }} \quad if \; 0\leq i\leq \inf(2s+m-2, k-1).
\end{equation}\\
 \end{lem}
 
 \begin{proof}          
 The Lemma  follows respectively from Lemmas \ref{lem 7.1}, \ref{lem 6.14} and \ref{lem 6.10}.
  \end{proof}

  \section{\textbf{ STUDY OF THE REMAINDER   $ \Delta _{i}^{\Big[{s \atop s+m}\Big] \times k}$ IN THE RECURRENT FORMULA }}  
\label{sec 8}
From the rank formulas established in sections 6 and 7, we deduce by elementary rank considerations the following formula  for $ 1\leq i \leq 2s+m,\; k\geq i+1 $\vspace{0.1 cm}\\
$$ \Delta _{i}^{\Big[{s \atop s+m}\Big] \times k} = \sum_{j = i-2}^{i+1}a_{j}\cdot \Gamma _{j}^{\Big[\substack{s -1\\ s-1 +m }\Big] \times j}$$
 where  the $ a_{j}\in \mathbb{Z}$    are  explicitely determined.\vspace{0.1 cm}\\
We get $$\Delta _{i}^{\Big[{s \atop s+m}\Big] \times k} = \Delta _{i}^{\Big[{s \atop s+m}\Big] \times (i+1)}\quad   whenever  \;  k \geq i+1. $$\vspace{0.1 cm}\\

  Let  $ ( j_{1}, j_{2},   j_{3})  \in \mathbb{N}^{3}$ recall that we define $$  \sigma _{j_{1},j_{2},j_{3}}^{\left[\stackrel{s-1}{\stackrel{s+m-1 }
{\overline {\stackrel{\alpha_{s -}}{\beta_{s+m-} }}}}\right] \times k } $$\\ to be the cardinality of the following set \\
 \small
 $$\begin{array}{l}\Big\{ (t,\eta ) \in \mathbb{P}/\mathbb{P}_{k+s -1}\times
           \mathbb{P}/\mathbb{P}_{k+s+m-1} 
\mid r(  D^{\left[\stackrel{s-1}{s+m-1}\right] \times k }(t,\eta  ) ) = j_{1}, \quad 
r( D^{\left[\stackrel{s}{s+m-1}\right] \times k }(t,\eta  ) ) = j_{2},  \\
 r(  D^{\left[\stackrel{s}{s+m}\right] \times k }(t,\eta  )  ) = j_{3}  \Big\}.
    \end{array}$$\\
  \begin{lem}
   \label{lem 8.1}
   \label{eq 8.1}
     For $ 1\leq i\leq 2s+m-3,\quad k \geq  i+1  $ we have
    \begin{equation}
   \sigma _{i,i,i}^{\left[\stackrel{s-1}{\stackrel{s+m-1 }
{\overline {\stackrel{\alpha_{s -}}{\beta_{s+m-} }}}}\right] \times k } = 
      \sigma _{i,i,i}^{\left[\stackrel{s-1}{\stackrel{s+m-1 }
{\overline {\stackrel{\alpha_{s -}}{\beta_{s+m-} }}}}\right] \times (i+1) }. 
 \end{equation}
\end{lem} 
\begin{proof}
The formula \eqref{eq 8.1} is obvious for k = i+1.\\

   We consider the following partition of the matrix $$ D^{\left[\stackrel{s-1}{\stackrel{s+m-1 }
{\overline {\stackrel{\alpha  _{s -}}{\beta  _{s +m -} }}}}\right] \times k }(t,\eta ) $$
 
    $$   \left ( \begin{array} {ccccc|c}
\alpha _{1} & \alpha _{2} & \alpha _{3} &  \ldots & \alpha _{k-1}  &  \alpha _{k} \\
\alpha _{2 } & \alpha _{3} & \alpha _{4}&  \ldots  &  \alpha _{k} &  \alpha _{k+1} \\
\vdots & \vdots & \vdots    &  \vdots & \vdots  &  \vdots \\
\alpha _{s-1} & \alpha _{s} & \alpha _{s +1} & \ldots  &  \alpha _{s+k-3} &  \alpha _{s+k-2}  \\
\beta  _{1} & \beta  _{2} & \beta  _{3} & \ldots  &  \beta_{k-1} &  \beta _{k}  \\
\beta  _{2} & \beta  _{3} & \beta  _{4} & \ldots  &  \beta_{k} &  \beta _{k+1}  \\
\vdots & \vdots & \vdots    &  \vdots & \vdots  &  \vdots \\
\beta  _{m+1} & \beta  _{m+2} & \beta  _{m+3} & \ldots  &  \beta_{k+m-1} &  \beta _{k+m}  \\
\vdots & \vdots & \vdots    &  \vdots & \vdots  &  \vdots \\
\beta  _{s+m-1} & \beta  _{s+m} & \beta  _{s+m+1} & \ldots  &  \beta_{s+m+k-3} &  \beta _{s+m+k-2}  \\
\hline
\alpha _{s} & \alpha _{s+1} & \alpha _{s +2} & \ldots  &  \alpha _{s+k-2} &  \alpha _{s+k-1}\\
\hline
\beta  _{s+m} & \beta  _{s+m+1} & \beta  _{s+m+2} & \ldots  &  \beta_{s+m+k-2} &  \beta _{s+m+k-1}
  \end{array}  \right). $$ \vspace{0.5 cm}

Let $ k \geq  i+2. $  Using Lemma \ref{lem 7.1}, we obtain by elementary rank considerations  \\
\begin{align}
\label{eq 8.2}
  \sigma _{i,i,i}^{\left[\stackrel{s-1}{\stackrel{s+m-1 }
{\overline {\stackrel{\alpha_{s -}}{\beta_{s+m-} }}}}\right] \times k } 
&  = 
 {}^{\#}\left(\begin{array}{c | c}
           i & i \\
           \hline
           i & i \\
           \hline
            i & i 
           \end{array} \right)_{\mathbb{P}/\mathbb{P}_{k+s -1}\times
           \mathbb{P}/\mathbb{P}_{k+s+m-1} }^{{\alpha  \over \beta }} +
        {}^{\#}\left(\begin{array}{c | c}
             i -1 & i \\
           \hline
           i & i \\
           \hline
            i & i 
           \end{array} \right)_{\mathbb{P}/\mathbb{P}_{k+s -1}\times
           \mathbb{P}/\mathbb{P}_{k+s+m-1} }^{{\alpha  \over \beta }}\\
           &  +
           {}^{\#}\left(\begin{array}{c | c}
           i-1 & i \\
           \hline
           i-1 & i \\
           \hline
            i & i 
           \end{array} \right)_{\mathbb{P}/\mathbb{P}_{k+s -1}\times
           \mathbb{P}/\mathbb{P}_{k+s+m-1} }^{{\alpha  \over \beta  }}  + 
             {}^{\#}\left(\begin{array}{c | c}
            i -1 & i \\
           \hline
           i-1 & i \\
           \hline
            i-1 & i 
           \end{array} \right)_{\mathbb{P}/\mathbb{P}_{k+s -1}\times
           \mathbb{P}/\mathbb{P}_{k+s+m-1} }^{{\alpha  \over \beta  }}  \nonumber  \\
           & =  {}^{\#}\left(\begin{array}{c | c}
           i & i \\
           \hline
           i & i \\
           \hline
            i & i 
           \end{array} \right)_{\mathbb{P}/\mathbb{P}_{k+s -1}\times
           \mathbb{P}/\mathbb{P}_{k+s+m-1} }^{{\alpha  \over \beta  }}. \nonumber
            \end{align}
    Further,  by  Lemma \ref{lem 7.2} we have  \\
         
         \begin{align}
\label{eq 8.3}
4\cdot \sigma _{i,i,i}^{\left[\stackrel{s-1}{\stackrel{s+m-1 }
{\overline {\stackrel{\alpha_{s -}}{\beta_{s+m-} }}}}\right] \times (k-1) } 
&  = 
 {}^{\#}\left(\begin{array}{c | c}
           i & i \\
           \hline
           i & i \\
           \hline
            i & i 
           \end{array} \right)_{\mathbb{P}/\mathbb{P}_{k+s -1}\times
           \mathbb{P}/\mathbb{P}_{k+s+m-1} }^{{\alpha \over \beta  }} +
        {}^{\#}\left(\begin{array}{c | c}
             i  & i \\
           \hline
           i & i+1 \\
           \hline
            i & i +1
           \end{array} \right)_{\mathbb{P}/\mathbb{P}_{k+s -1}\times
           \mathbb{P}/\mathbb{P}_{k+s+m-1} }^{{\alpha \over \beta }}\\
           &  +
           {}^{\#}\left(\begin{array}{c | c}
           i & i \\
           \hline
           i & i \\
           \hline
            i & i +1
           \end{array} \right)_{\mathbb{P}/\mathbb{P}_{k+s -1}\times
           \mathbb{P}/\mathbb{P}_{k+s+m-1} }^{{\alpha  \over \beta }}  + 
             {}^{\#}\left(\begin{array}{c | c}
            i  & i +1\\
           \hline
           i & i +1 \\
           \hline
            i & i +1
           \end{array} \right)_{\mathbb{P}/\mathbb{P}_{k+s -1}\times
           \mathbb{P}/\mathbb{P}_{k+s+m-1} }^{{\alpha \over \beta }}  \nonumber  \\
           & = 4\cdot {}^{\#}\left(\begin{array}{c | c}
           i & i \\
           \hline
           i & i \\
           \hline
            i & i 
           \end{array} \right)_{\mathbb{P}/\mathbb{P}_{k+s -1}\times
           \mathbb{P}/\mathbb{P}_{k+s+m-1} }^{{\alpha  \over \beta }}. \nonumber
            \end{align}
  By \eqref{eq 8.2},  \eqref{eq 8.3} we deduce \\
    \begin{equation}
  \label{eq 8.4}
  \sigma _{i,i,i}^{\left[\stackrel{s-1}{\stackrel{s+m-1 }{\overline {\stackrel{\alpha_{s -}}{\beta_{s+m-} }}}}\right] \times (k-1) } = 
    \sigma _{i,i,i}^{\left[\stackrel{s-1}{\stackrel{s+m-1 }{\overline {\stackrel{\alpha_{s -}}{\beta_{s+m-} }}}}\right] \times k }\quad  for \; all \; k\geq i+2. 
    \end{equation}
Hence by \eqref{eq 8.4} we obtain successively\\

$$  \sigma _{i,i,i}^{\left[\stackrel{s-1}{\stackrel{s+m-1 }{\overline {\stackrel{\alpha_{s -}}{\beta_{s+m-} }}}}\right] \times k }
= \sigma _{i,i,i}^{\left[\stackrel{s-1}{\stackrel{s+m-1 }{\overline {\stackrel{\alpha_{s -}}{\beta_{s+m-} }}}}\right] \times (k-1) }
= \sigma _{i,i,i}^{\left[\stackrel{s-1}{\stackrel{s+m-1 }{\overline {\stackrel{\alpha_{s -}}{\beta_{s+m-} }}}}\right] \times (k-2)}
= \ldots
= \sigma _{i,i,i}^{\left[\stackrel{s-1}{\stackrel{s+m-1 }{\overline {\stackrel{\alpha_{s -}}{\beta_{s+m-} }}}}\right] \times (i+1) }. $$
\end{proof}

   \begin{lem}
   \label{lem 8.2} For all i such that $ 1 \leq i \leq 2s+m-3  $ we have \\
   \begin{equation}
\label{eq 8.5}
    \sigma _{i,i,i}^{\left[\stackrel{s-1}{\stackrel{s+m-1 }{\overline {\stackrel{\alpha_{s -}}{\beta_{s+m-} }}}}\right] \times (i+1)} = 
   4 \cdot \Gamma _{i}^{\Big[\substack{s -1\\ s-1 +m }\Big] \times i} -  \Gamma _{i+1}^{\Big[\substack{s-1 \\ s-1+m }\Big] \times (i+1)}.
   \end{equation}
    \end{lem}
 \begin{proof}
  Consider the following partition of the matrix $$ D^{\left[\stackrel{s-1}{\stackrel{s+m-1 }
{\overline {\stackrel{\alpha  _{s -}}{\beta  _{s +m -} }}}}\right] \times (i+1) }(t,\eta ) $$
 
    $$   \left ( \begin{array} {ccccc|c}
\alpha _{1} & \alpha _{2} & \alpha _{3} &  \ldots & \alpha _{i}  &  \alpha _{i+1} \\
\alpha _{2 } & \alpha _{3} & \alpha _{4}&  \ldots  &  \alpha _{1+1} &  \alpha _{i+2} \\
\vdots & \vdots & \vdots    &  \vdots & \vdots  &  \vdots \\
\alpha _{s-1} & \alpha _{s} & \alpha _{s +1} & \ldots  &  \alpha _{i+s-2} &  \alpha _{i+s-1}  \\
\beta  _{1} & \beta  _{2} & \beta  _{3} & \ldots  &  \beta_{i} &  \beta _{i+1}  \\
\beta  _{2} & \beta  _{3} & \beta  _{4} & \ldots  &  \beta_{i+1} &  \beta _{i+2}  \\
\vdots & \vdots & \vdots    &  \vdots & \vdots  &  \vdots \\
\beta  _{m+1} & \beta  _{m+2} & \beta  _{m+3} & \ldots  &  \beta_{i+m} &  \beta _{i+m+1}  \\
\vdots & \vdots & \vdots    &  \vdots & \vdots  &  \vdots \\
\beta  _{s+m-1} & \beta  _{s+m} & \beta  _{s+m+1} & \ldots  &  \beta_{s+m+i-2} &  \beta _{s+m+i-1}  \\
\hline
\alpha _{s} & \alpha _{s+1} & \alpha _{s +2} & \ldots  &  \alpha _{s+i-1} &  \alpha _{s+i}\\
\hline
\beta  _{s+m} & \beta  _{s+m+1} & \beta  _{s+m+2} & \ldots  &  \beta_{s+m+i-1} &  \beta _{s+m+i}
  \end{array}  \right). $$ \vspace{0.5 cm}
  
   By elementary rank considerations and using Lemma \ref{lem 7.1} we obtain \\
\begin{align}
\label{eq 8.6}
  \sigma _{i,i,i}^{\left[\stackrel{s-1}{\stackrel{s+m-1 }
{\overline {\stackrel{\alpha_{s -}}{\beta_{s+m-} }}}}\right] \times (i+1) } 
&  = 
 {}^{\#}\left(\begin{array}{c | c}
           i & i \\
           \hline
           i & i \\
           \hline
            i & i 
           \end{array} \right)_{\mathbb{P}/\mathbb{P}_{i+s }\times
           \mathbb{P}/\mathbb{P}_{i+s+m} }^{{\alpha  \over \beta  }} +
        {}^{\#}\left(\begin{array}{c | c}
             i -1 & i \\
           \hline
           i & i \\
           \hline
            i & i 
           \end{array} \right)_{\mathbb{P}/\mathbb{P}_{i+s }\times
           \mathbb{P}/\mathbb{P}_{i+s+m} }^{{\alpha  \over \beta  }}\\
           &  +
           {}^{\#}\left(\begin{array}{c | c}
           i-1 & i \\
           \hline
           i-1 & i \\
           \hline
            i & i 
           \end{array} \right)_{\mathbb{P}/\mathbb{P}_{i+s }\times
           \mathbb{P}/\mathbb{P}_{i+s+m} }^{{\alpha \over \beta }}  + 
             {}^{\#}\left(\begin{array}{c | c}
            i -1 & i \\
           \hline
           i-1 & i \\
           \hline
            i-1 & i 
           \end{array} \right)_{\mathbb{P}/\mathbb{P}_{i+s }\times
           \mathbb{P}/\mathbb{P}_{i+s+m} }^{{\alpha \over \beta  }}  \nonumber  \\
           & =  {}^{\#}\left(\begin{array}{c | c}
           i & i \\
           \hline
           i & i \\
           \hline
            i & i 
           \end{array} \right)_{\mathbb{P}/\mathbb{P}_{i+s }\times
           \mathbb{P}/\mathbb{P}_{i+s+m} }^{{\alpha  \over \beta }}. \nonumber
            \end{align}
    Further,  by  Lemma \ref{lem 7.2} we have  \\
         
         \begin{align}
\label{eq 8.7}
4\cdot \sigma _{i,i,i}^{\left[\stackrel{s-1}{\stackrel{s+m-1 }
{\overline {\stackrel{\alpha_{s -}}{\beta_{s+m-} }}}}\right] \times i } 
&  =  4 \cdot 4\cdot \Gamma _{i}^{\Big[\substack{s -1\\ s-1 +m }\Big] \times i} \\
& = 
 {}^{\#}\left(\begin{array}{c | c}
           i & i \\
           \hline
           i & i \\
           \hline
            i & i 
           \end{array} \right)_{\mathbb{P}/\mathbb{P}_{i+s }\times
           \mathbb{P}/\mathbb{P}_{i+s+m} }^{{\alpha \over \beta }} +
        {}^{\#}\left(\begin{array}{c | c}
             i  & i \\
           \hline
           i & i+1 \\
           \hline
            i & i +1
           \end{array} \right)_{\mathbb{P}/\mathbb{P}_{i+s }\times
           \mathbb{P}/\mathbb{P}_{i+s+m} }^{{\alpha \over \beta  }}\nonumber \\
           &  +
           {}^{\#}\left(\begin{array}{c | c}
           i & i \\
           \hline
           i & i \\
           \hline
            i & i +1
           \end{array} \right)_{\mathbb{P}/\mathbb{P}_{i+s }\times
           \mathbb{P}/\mathbb{P}_{i+s+m} }^{{\alpha  \over \beta }}  + 
             {}^{\#}\left(\begin{array}{c | c}
            i  & i +1\\
           \hline
           i & i +1 \\
           \hline
            i & i +1
           \end{array} \right)_{\mathbb{P}/\mathbb{P}_{i+s }\times
           \mathbb{P}/\mathbb{P}_{i+s+m} }^{{\alpha  \over \beta }}  \nonumber  \\
           & = 4\cdot {}^{\#}\left(\begin{array}{c | c}
           i & i \\
           \hline
           i & i \\
           \hline
            i & i 
           \end{array} \right)_{\mathbb{P}/\mathbb{P}_{i+s }\times
           \mathbb{P}/\mathbb{P}_{i+s+m} }^{{\alpha  \over \beta }} +    
            {}^{\#}\left(\begin{array}{c | c}
            i  & i +1\\
           \hline
           i & i +1 \\
           \hline
            i & i +1
           \end{array} \right)_{\mathbb{P}/\mathbb{P}_{i+s }\times
           \mathbb{P}/\mathbb{P}_{i+s+m} }^{{\alpha  \over \beta }}  \nonumber  \\
           & =   4\cdot {}^{\#}\left(\begin{array}{c | c}
           i & i \\
           \hline
           i & i \\
           \hline
            i & i 
           \end{array} \right)_{\mathbb{P}/\mathbb{P}_{i+s }\times
           \mathbb{P}/\mathbb{P}_{i+s+m} }^{{\alpha  \over \beta }} +  
              4 \cdot \Gamma _{i+1}^{\Big[\substack{s -1\\ s-1 +m }\Big] \times (i+1)}. \nonumber        
            \end{align}
   Combining \eqref{eq 8.6},  \eqref{eq 8.7} we deduce \\
  
  \begin{equation}
  \label{eq 8.8}
   16 \cdot \Gamma _{i}^{\Big[\substack{s -1\\ s-1 +m }\Big] \times i} =  
  4\cdot \sigma _{i,i,i}^{\left[\stackrel{s-1}{\stackrel{s+m-1 }
{\overline {\stackrel{\alpha_{s -}}{\beta_{s+m-} }}}}\right] \times (i+1)} +
    4\cdot \Gamma _{i+1}^{\Big[\substack{s -1\\ s-1 +m }\Big] \times (i+1)}
  \end{equation}
and \eqref{eq 8.5} is established.
 \end{proof}
 
\begin{lem}
\label{lem 8.3}Let $ s\geq 2 $ and $ m\geq 0, $ we have in the following two cases : \vspace{0.1 cm}\\

\underline {The case $1 \leq  k \leq 2s+m-2 $}\\
\begin{equation}
\label{eq 8.9}
 \sigma _{i,i,i}^{\left[\stackrel{s-1}{\stackrel{s+m -1}{\overline {\stackrel{\alpha_{s -}}{\beta_{s+m-} }}}}\right] \times k }= 
  \begin{cases}
1 & \text{if  } i = 0,\quad k\geq 1, \\
   4 \cdot \Gamma _{i}^{\Big[\substack{s -1\\ s-1 +m }\Big] \times i} -  \Gamma _{i+1}^{\Big[\substack{s-1 \\ s-1+m }\Big] \times (i+1)} 
      & \text{if   } 1 \leq i \leq k-1, \\
   4 \cdot \Gamma _{k}^{\Big[\substack{s -1\\ s-1 +m }\Big] \times k}    & \text{if   } i = k.
\end{cases}
\end{equation}

\underline {The case $ k\geq  2s+m-2 $ }\\
\begin{equation}
\label{eq 8.10}
 \sigma _{i,i,i}^{\left[\stackrel{s-1}{\stackrel{s+m-1 }{\overline {\stackrel{\alpha_{s -}}{\beta_{s+m-} }}}}\right] \times k }= 
  \begin{cases}
1 & \text{if  } i = 0,          \\
   4 \cdot \Gamma _{i}^{\Big[\substack{s -1\\ s-1 +m }\Big] \times i} -  \Gamma _{i+1}^{\Big[\substack{s-1 \\ s-1+m }\Big] \times (i+1)} 
      & \text{if   } 1 \leq i \leq 2s+m-3,  \\
   4 \cdot \Gamma _{2s+m-2}^{\Big[\substack{s -1\\ s-1 +m }\Big] \times (2s+m-2)}    & \text{if   } i = 2s+m-2. 
\end{cases}
\end{equation}
\end{lem}
\begin{proof}
 Combining \eqref{eq 8.1} and \eqref{eq 8.5} and recalling the definition of 
$ \sigma _{i,i,i}^{\left[\stackrel{s}{\stackrel{s+m }{\overline {\stackrel{\alpha_{s -}}{\beta_{s+m-} }}}}\right] \times k } $ we deduce easily Lemma \ref{lem  8.3}.
\end{proof}
\begin{lem}
\label{lem 8.4} The remainder   $ \displaystyle  \Delta _{i}^{\Big[\substack{s \\ s+m }\Big] \times k}   $ in  the recurrent formula 
is equal to 

   \begin{equation}
\label{eq 8.11}
  \begin{cases}
1 & \text{if  } i = 0,\; k \geq 1,         \\
   4 \cdot \Gamma _{1}^{\Big[\substack{s -1\\ s-1 +m }\Big] \times 1} -  \Gamma _{2}^{\Big[\substack{s-1 \\ s-1+m }\Big] \times 2} -3
      & \text{if   } i = 1,\; k\geq 2, \\
    4 \cdot \Gamma _{1}^{\Big[\substack{s -1\\ s-1 +m }\Big] \times 1} - 3   & \text{if   } i = 1,\; k =1, \\
     7 \cdot \Gamma _{2}^{\Big[\substack{s -1\\ s-1 +m }\Big] \times 2} -  12\cdot \Gamma _{1}^{\Big[\substack{s-1 \\ s-1+m }\Big] \times 1}
     -  \Gamma _{3}^{\Big[\substack{s -1\\ s-1 +m }\Big] \times 3} + 2  & \text{if   } i = 2,\; k\geq 3, \\
        7 \cdot \Gamma _{2}^{\Big[\substack{s -1\\ s-1 +m }\Big] \times 2} -  12\cdot \Gamma _{1}^{\Big[\substack{s-1 \\ s-1+m }\Big] \times 1} +2
         & \text{if   } i = 2,\; k = 2, \\
         7 \cdot \Gamma _{i}^{\Big[\substack{s -1\\ s-1 +m }\Big] \times i} -  14\cdot \Gamma _{i-1}^{\Big[\substack{s-1 \\ s-1+m }\Big] \times (i-1)}
     +8\cdot \Gamma _{i-2}^{\Big[\substack{s -1\\ s-1 +m }\Big] \times (i-2)}  -  \Gamma _{i+1}^{\Big[\substack{s -1\\ s-1 +m }\Big] \times (i+1)}   & \text{if } 3\leq i\leq 2s+m-3,  k\geq i+1, \\
        7 \cdot \Gamma _{i}^{\Big[\substack{s -1\\ s-1 +m }\Big] \times i} -  14\cdot \Gamma _{i-1}^{\Big[\substack{s-1 \\ s-1+m }\Big] \times (i-1)}
     +8\cdot \Gamma _{i-2}^{\Big[\substack{s -1\\ s-1 +m }\Big] \times (i-2)}    & \text{if   }  3\leq i\leq 2s+m-3, \; k = i, \\
       7 \cdot \Gamma _{2s+m-2}^{\Big[\substack{s -1\\ s-1 +m }\Big] \times (2s+m-2)} -  14\cdot \Gamma _{2s+m-3}^{\Big[\substack{s-1 \\ s-1+m }\Big] \times (2s+m-3)}
     +8\cdot \Gamma _{2s+m-4}^{\Big[\substack{s -1\\ s-1 +m }\Big] \times (2s+m-4)}    & \text{if   }  i = 2s+m-2, \; k \geq  i, \\
     - 14\cdot \Gamma _{2s+m-2}^{\Big[\substack{s-1 \\ s-1+m }\Big] \times (2s+m-2)}
     +8\cdot \Gamma _{2s+m-3}^{\Big[\substack{s -1\\ s-1 +m }\Big] \times (2s+m-3)}    & \text{if   }  i = 2s+m-1, \; k \geq  i, \\
      8\cdot \Gamma _{2s+m-2}^{\Big[\substack{s -1\\ s-1 +m }\Big] \times (2s+m-2)}    & \text{if   }  i = 2s+m , \; k \geq  i. 
\end{cases}
\end{equation}
  \end{lem} 
  
  \begin{proof}
From  \eqref{eq 8.9} and \eqref{eq 8.10} and recalling that $ \Delta _{i}^{\Big[\substack{s \\ s+m }\Big] \times k} $ is equal to
$$ \sigma _{i,i,i}^{\left[\stackrel{s-1}{\stackrel{s+m-1 }{\overline {\stackrel{\alpha_{s -}}{\beta_{s+m-} }}}}\right] \times k }
  - 3\cdot \sigma _{i-1,i-1,i-1}^{\left[\stackrel{s-1}{\stackrel{s+m-1 }{\overline {\stackrel{\alpha_{s -}}{\beta_{s+m-} }}}}\right] \times k } 
  + 2\cdot \sigma _{i-2,i-2,i-2}^{\left[\stackrel{s-1}{\stackrel{s+m-1 }{\overline {\stackrel{\alpha_{s -}}{\beta_{s+m-} }}}}\right] \times k }  $$  we deduce easily Lemma \ref{lem 8.4}.
\end{proof}
 \begin{lem}
\label{lem  8.5}We have for all i such that $ 0\leq i\leq 2s+m-3 $ and for all $ k\geq i+1 $
 \begin{equation}
 \label{eq 8.12}
 \Delta _{i}^{\Big[\substack{s \\ s+m }\Big] \times k}  =  \Delta _{i}^{\Big[\substack{s \\ s+m }\Big] \times (i+1)},
\end{equation}
 \begin{equation}
 \label{eq 8.13}
 \Delta _{i}^{\Big[\substack{s \\ s+m }\Big] \times k}  =  \Delta _{i}^{\Big[\substack{s \\ s+m }\Big] \times i} \;  for \; i\in \left\{2s+m-2, 2s+m-1,2s+m \right\}, \; k\geq i.
\end{equation}
 \end{lem}
  \begin{proof} 
Follows immediately from Lemma \ref{lem 8.4}. \vspace{0.05 cm}\\
 \end{proof}

  \section{\textbf{ A RECURRENT FORMULA FOR THE DIFFERENCE $\Gamma_{i}^{\left[s\atop s+m\right]\times (k+1)} - \Gamma_{i}^{\left[s\atop s+ m\right]\times k},  1\leq i\leq 2s+m,\; k\geq i+1 $}}
     \label{sec 9}
In this section we deduce from the recurrent formula in Lemma \ref{lem 5.14} and from the fact that the remainder $  \Delta_{i}^{\Big[\substack{s\\ s +m }\Big] \times k} $ is independent of k if  $ k\geq i+1, $
the following recurrent formula  \vspace{0.1 cm }\\
\begin{align*}
&  \Gamma _{i}^{\Big[\substack{s \\ s +m }\Big] \times (k+1)} -  \Gamma _{i}^{\Big[\substack{s \\ s +m }\Big] \times k}
   =  4\cdot \left[\Gamma _{i-1}^{\Big[\substack{s \\ s +(m-1) }\Big] \times (k+1)} - \Gamma _{i-1}^{\Big[\substack{s \\ s +(m-1) }\Big] \times  k}\right] \\
 &  + R(i,s,m,k) \quad \text{if  }  0 \leq i\leq \inf(2s+m,k-1 ),\quad \text{ where R(i,s,m,k)  is  equal to } 
       \end{align*}
 
  \begin{equation*}
  \begin{cases}
    0  & \text{if  }\quad 1\leq i\leq s-1,\; k\geq i+1,  \\
    2^{s-1}\left[ \Gamma _{j+1}^{\Big[\substack{1 \\ 1+ (s+m-1)}\Big] \times (k+1)} -  \Gamma _{j+1}^{\Big[\substack{1 \\ 1+(s+m-1) }\Big] \times k}\right] \\
      -   2^{s+1}\left[ \Gamma _{j}^{\Big[\substack{1 \\ 1+ (s+m-2)}\Big] \times (k+1)} -  \Gamma _{j}^{\Big[\substack{1 \\ 1+(s+m-2) }\Big] \times k}\right] & \text{if }\quad i = s+j,\quad 0\leq j\leq \inf(s+m,k-s-1).
   \end{cases}
  \end{equation*}\vspace{0.1 cm }\\
If i = s+j we set  R(i,s,m,k ) =  R(j,s,m,k)  $ \; \text{where}\;  0\leq j\leq \inf(s+m,k-s-1).$ \\

\begin{lem}
\label{lem  9.1} $ Let \; s \geq 2,  m\geq 0 \; and \; k \geq 1,$ then  $ \Gamma _{i}^{\Big[\substack{s \\ s +m }\Big] \times k} - 4\cdot \Gamma _{i-1}^{\Big[\substack{s \\ s +(m-1) }\Big] \times k} $ is given by \\

 \begin{equation}
\label{eq 9.1} 
 \begin{cases} \displaystyle 
 5\cdot2^{i-1} +   \sum_{j = 0}^{i-2}2^{j}\Delta_{i-j}^{\Big[\substack{s -j\\ s-j +(m+j) }\Big] \times k}  & \text{if  }  2\leq i\leq \inf(s-1,k),   \\
 \displaystyle   2^{s-1}\left[ \Gamma _{i-(s-1)}^{\Big[\substack{1 \\ 1+m+(s-1)}\Big] \times k} - 4\cdot \Gamma _{i-s}^{\Big[\substack{1 \\ 1+m+(s-2) }\Big] \times k}\right] 
  +     \sum_{j = 0}^{s-2}2^{j}\Delta_{i-j}^{\Big[\substack{s -j\\ s-j +(m+j) }\Big] \times k} & \text{if  } s\leq i\leq \inf(k,2s+m).
 \end{cases}
\end{equation}
\end{lem}

\begin{proof}
\underline{The first case $ 2\leq i\leq \inf(s-1,k) $}\vspace{0.1 cm} \\

Set for $ 1\leq i\leq \inf(2s+m,k) $ $$\Omega _{i}(s,m,k)= \Gamma _{i}^{\Big[\substack{s \\ s +m }\Big] \times k} - 4\cdot \Gamma _{i-1}^{\Big[\substack{s \\ s +(m-1) }\Big] \times k} $$
 so  \\
 $$ \Omega _{i-1}(s-1,m+1,k)= \Gamma _{i-1}^{\Big[\substack{s-1 \\ s-1 +(m+1)}\Big] \times k} - 4\cdot \Gamma _{i-2}^{\Big[\substack{s -1\\ s-1 + m }\Big] \times k}. $$ \\
 We obtain then by Lemma \ref{lem 5.14}\\
 
  \begin{align}
  \Gamma _{i}^{\Big[\substack{s \\ s+m }\Big] \times k} & = 2\cdot \Gamma _{i-1}^{\Big[\substack{s -1\\ s-1+(m+1) }\Big] \times k}
+ 4\cdot \Gamma _{i-1}^{\Big[\substack{s \\ s+(m-1) }\Big] \times k} - 8\cdot \Gamma _{i-2}^{\Big[\substack{s-1 \\ s-1+m }\Big] \times k}
 + \Delta _{i}^{\Big[\substack{s \\ s+m }\Big] \times k}, \label{eq 9.2} \\
 & \nonumber \\
  & \Longleftrightarrow  \Omega _{i}(s,m,k)= 2\cdot \Omega _{i-1}(s-1,m+1,k) +  \Delta _{i}^{\Big[\substack{s \\ s+m }\Big] \times k}.\label{eq 9.3}\\
  & \nonumber 
 \end{align}
 By \eqref{eq 9.3} we obtain succesively
 \Small
  \begin{align*}
 &  \Omega _{i}(s,m,k)  = 2\cdot \Omega _{i-1}(s-1,m+1,k) +  \Delta _{i}^{\Big[\substack{s \\ s+m }\Big] \times k}\\
 &   2\cdot \Omega _{i-1}(s-1,m+1,k)  =2\cdot\left( 2\cdot \Omega _{i-2}(s-2,m+2,k) +  \Delta _{i-1}^{\Big[\substack{s-1 \\ s-1+(m+1) }\Big] \times k}\right)\\
   &                                                      \vdots       \\
    &  2^{j}\cdot \Omega _{i-j}(s-j,m+j,k)  =2^{j}\cdot\left( 2\cdot \Omega _{i-(j+1)}(s-(j+1),m+(j+1),k) +  \Delta _{i-j}^{\Big[\substack{s-j \\ s-j+(m+j) }\Big] \times k}\right) \\
     &                                                      \vdots       \\
   &    2^{i-2}\cdot \Omega _{i-(i-2)}(s-(i-2),m+(i-2),k)  =2^{i-2}\cdot\left( 2\cdot \Omega _{i-(i-1)}(s-(i-1),m+(i-1),k) +  \Delta _{i-(i-2)}^{\Big[\substack{s-(i-2) \\ s-(i-2)+(m+(i-2)) }\Big] \times k}\right) \\
       &                                                              
  \end{align*}
  \normalsize
  By summing the above equations we get \\
   \begin{align}
\displaystyle   \sum_{j = 0}^{i-2}2^{j}\cdot \Omega _{i-j}(s-j,m+j,k) &  =   \sum_{j = 0}^{i-2}2^{j+1}\cdot \Omega _{i-(j+1)}(s-(j+1),m+(j+1),k) 
  +  \sum_{j = 0}^{i-2}2^{j}\cdot \Delta _{i-j}^{\Big[\substack{s-j \\ s-j+(m+j) }\Big] \times k}  \label{eq 9.4} \\
&   =  \sum_{j = 1}^{i-1}2^{j}\cdot \Omega _{i-j}(s-j,m+ j,k) 
  +  \sum_{j = 0}^{i-2}2^{j}\cdot \Delta _{i-j}^{\Big[\substack{s-j \\ s-j+(m+j) }\Big] \times k}. \nonumber \\
  & \nonumber
\end{align}  
By \eqref{eq 9.4}  we get after some simplifications \\
\begin{equation}
\label{eq 9.5}
  \Omega _{i}(s,m,k) = 2^{i-1}\cdot \Omega _{1}(s-(i-1),m+(i-1),k) +  \sum_{j = 0}^{i-2}2^{j}\cdot \Delta _{i-j}^{\Big[\substack{s-j \\ s-j+(m+j) }\Big] \times k}. \\
   \end{equation}
   
 By the definition of  $ \Omega _{i}(s,m,k) \;\text{we have}\; \Omega _{1}(s-(i-1),m+(i-1),k) = 
 \Gamma _{1}^{\Big[\substack{s -i +1\\ s +m }\Big] \times k} - 4\cdot \Gamma _{0}^{\Big[\substack{s - i+1\\ s +(m-1) }\Big] \times k}. $ \\
 
 Recall that  $  D^{\left[\stackrel{s-i+1}{s+m}\right] \times k }(t,\eta  ) = 
\left[{D_{(s-i+1)\times k}(t)\over D_{(s+m)\times k}(\eta )}\right] $ denotes 
the following $(2s+m-i+1)\times k $ matrix where the first (s -i +1) rows form a  $(s-i+1)\times k $
persymmetric matrix and the last ( s+m) rows form a $( s+m)\times k $ persymmetric matrix with entries in  $\mathbb{F}_{2} $\\

   $$   \left ( \begin{array} {cccccc}
\alpha _{1} & \alpha _{2} & \alpha _{3} &  \ldots & \alpha _{k-1}  &  \alpha _{k} \\
\alpha _{2 } & \alpha _{3} & \alpha _{4}&  \ldots  &  \alpha _{k} &  \alpha _{k+1} \\
\vdots & \vdots & \vdots    &  \vdots & \vdots  &  \vdots \\
\alpha _{s-i+1} & \alpha _{s-i+2} & \alpha _{s-i +3} & \ldots  &  \alpha _{s-i+k-1} &  \alpha _{s-i+k}  \\
\hline
\beta  _{1} & \beta  _{2} & \beta  _{3} & \ldots  &  \beta_{k-1} &  \beta _{k}  \\
\beta  _{2} & \beta  _{3} & \beta  _{4} & \ldots  &  \beta_{k} &  \beta _{k+1}  \\
\vdots & \vdots & \vdots    &  \vdots & \vdots  &  \vdots \\
\beta  _{m+1} & \beta  _{m+2} & \beta  _{m+3} & \ldots  &  \beta_{k+m-1} &  \beta _{k+m}  \\
\vdots & \vdots & \vdots    &  \vdots & \vdots  &  \vdots \\
\beta  _{s+m-1} & \beta  _{s+m} & \beta  _{s+m+1} & \ldots  &  \beta_{s+m+k-3} &  \beta _{s+m+k-2}  \\
\beta  _{s+m} & \beta  _{s+m+1} & \beta  _{s+m+2} & \ldots  &  \beta_{s+m+k-2} &  \beta _{s+m+k-1}\\
\end{array}  \right). $$ \vspace{0.5 cm}\\
Then obviously we have  \\ 
\small
\begin{align*}
&  \Gamma _{1}^{\Big[\substack{s -i+1 \\ s+m }\Big] \times k} \\
  &  =   Card
  \left\{(t,\eta )\in \mathbb{P}/\mathbb{P}_{k+s-i-1}\times \mathbb{P}/\mathbb{P}_{k+s+m-1}
\mid   r(D^{\big[\stackrel{s-i+1}{s+m}\big] \times k }(t,\eta )) = 1 \right\} \\
& =   Card
  \left\{(t,\eta )\in \mathbb{P}/\mathbb{P}_{k+s-i-1}\times \mathbb{P}/\mathbb{P}_{k+s+m-1}
\mid r(D_{(s-i+1)\times k}(t)) = 0,\quad    r(D_{(s+m)\times k}(\eta )) = 1 \right\}\\
& +  Card
  \left\{(t,\eta )\in \mathbb{P}/\mathbb{P}_{k+s-i-1}\times \mathbb{P}/\mathbb{P}_{k+s+m-1}
\mid r(D_{(s-i+1)\times k}(t)) = 1,\quad    r(D_{(s+m)\times k}(\eta )) = 0 \right\}\\
& + Card
  \left\{(t,\eta )\in \mathbb{P}/\mathbb{P}_{k+s-i-1}\times \mathbb{P}/\mathbb{P}_{k+s+m-1}
\mid r(D_{(s-i+1)\times k}(t)) = 1,\; r(D_{(s+m)\times k}(\eta )) = 1,\; 
  r(D^{\big[\stackrel{s-i+1}{s+m}\big] \times k }(t,\eta )) = 1  \right\}\\
  & = 3 + 3 + 3 = 9.
\end{align*}
\normalsize
And \eqref{eq 9.1} follows in the first case.\\

\underline{The second  case $ s\leq i\leq \inf(2s+m,k) $}\vspace{0.1 cm} \\

We proceed as in the first case.\\

Again  by \eqref{eq 9.3} we obtain succesively
\small
  \begin{align*}
 &  \Omega _{i}(s,m,k)  = 2\cdot \Omega _{i-1}(s-1,m+1,k) +  \Delta _{i}^{\Big[\substack{s \\ s+m }\Big] \times k}\\
 &   2\cdot \Omega _{i-1}(s-1,m+1,k)  =2\cdot\left( 2\cdot \Omega _{i-2}(s-2,m+2,k) +  \Delta _{i-1}^{\Big[\substack{s-1 \\ s-1+(m+1) }\Big] \times k}\right)\\
                                                       &  \vdots       \\
  &    2^{j}\cdot \Omega _{i-j}(s-j,m+j,k)  = 2^{j}\cdot\left( 2\cdot \Omega _{i-(j+1)}(s-(j+1),m+(j+1),k) +  \Delta _{i-j}^{\Big[\substack{s-j \\ s-j+(m+j) }\Big] \times k}\right) \\
 &  \vdots       \\
   &    2^{s-2}\cdot \Omega _{i-(s-2)}(s-(s-2),m+(s-2),k)  =2^{s-2}\cdot\left( 2\cdot \Omega _{i-(s-1)}(s-(s-1),m+(s-1),k) +  \Delta _{i-(s-2)}^{\Big[\substack{s-(s-2) \\ s+m }\Big] \times k}\right).                                                               
  \end{align*}
  By summing the above equations we obtain \\
  \begin{align}
\displaystyle   \sum_{j = 0}^{s -2}2^{j}\cdot \Omega _{i-j}(s-j,m+j,k) &  =   \sum_{j = 0}^{s-2}2^{j+1}\cdot \Omega _{i-(j+1)}(s-(j+1),m+(j+1),k) 
  +  \sum_{j = 0}^{s-2}2^{j}\cdot \Delta _{i-j}^{\Big[\substack{s-j \\ s-j+(m+j) }\Big] \times k}  \label{eq 9.6} \\
&   =  \sum_{j = 1}^{s-1}2^{j}\cdot \Omega _{i-j}(s-j,m+ j,k) 
  +  \sum_{j = 0}^{s-2}2^{j}\cdot \Delta _{i-j}^{\Big[\substack{s-j \\ s-j+(m+j) }\Big] \times k}. \nonumber
\end{align}  
By \eqref{eq 9.6}  we get after some simplifications \\
\begin{equation}
\label{eq 9.7}
  \Omega _{i}(s,m,k) = 2^{s-1}\cdot \Omega _{i-(s-1)}(1,m+s-1,k) +  \sum_{j = 0}^{s-2}2^{j}\cdot \Delta _{i-j}^{\Big[\substack{s-j \\ s-j+(m+j) }\Big] \times k}. \\
   \end{equation}
 By the definition of  $ \Omega _{i}(s,m,k) \; \text{we have} \; \Omega _{i-(s-1)}(1,m+s-1,k) = 
 \Gamma _{i-(s-1)}^{\Big[\substack{1 \\ 1+m+(s-1)}\Big] \times k} - 4\cdot \Gamma _{i-s}^{\Big[\substack{ 1 \\ 1+m+(s-2) }\Big] \times k}. $ \\
 And \eqref{eq 9.1} follows in the second case.
\end{proof}
\begin{lem}
\label{lem 9.2}  We have  \vspace{0.1 cm}\\
 \begin{equation}
\label{eq 9.8}
 \begin{cases} 
 \Gamma _{i}^{\Big[\substack{s \\ s +m }\Big] \times (k+1)} -  \Gamma _{i}^{\Big[\substack{s \\ s +m }\Big] \times k} =  4\cdot \left[\Gamma _{i-1}^{\Big[\substack{s \\ s +(m-1) }\Big] \times (k+1)} - \Gamma _{i-1}^{\Big[\substack{s \\ s +(m-1) }\Big] \times  k}\right]  
 & \text{if  }  2 \leq i\leq s-1,\quad k\geq i+1,   \\
  \Gamma _{i}^{\Big[\substack{s \\ s +m }\Big] \times (k+1)} -  \Gamma _{i}^{\Big[\substack{s \\ s +m }\Big] \times k} =  4\cdot \left[\Gamma _{i-1}^{\Big[\substack{s \\ s +(m-1) }\Big] \times (k+1)} - \Gamma _{i-1}^{\Big[\substack{s \\ s +(m-1) }\Big] \times  k}\right]  \\
  +   2^{s-1}\left[ \Gamma _{i-(s-1)}^{\Big[\substack{1 \\ 1+ (s+m-1)}\Big] \times (k+1)} -  \Gamma _{i-(s-1)}^{\Big[\substack{1 \\ 1+(s+m-1) }\Big] \times k}\right] \\
   -   2^{s+1}\left[ \Gamma _{i- s}^{\Big[\substack{1 \\ 1+ (s+m-2)}\Big] \times (k+1)} -  \Gamma _{i- s}^{\Big[\substack{1 \\ 1+(s+m-2) }\Big] \times k}\right] 
    & \text{if  }  s \leq i\leq 2s+m ,\quad k\geq i+1. 
 \end{cases}
\end{equation}
\end{lem}

\begin{proof}
By  \eqref{eq 8.12} we obtain 
$$ \sum_{j = 0}^{i-2}2^{j}\Delta_{i-j}^{\Big[\substack{s -j\\ s-j +(m+j) }\Big] \times (k+1)} =  \sum_{j = 0}^{i-2}2^{j}\Delta_{i-j}^{\Big[\substack{s -j\\ s-j +(m+j) }\Big] \times k}\;for \; 2 \leq i\leq s-1,\quad k\geq i+1  $$ 
and equally 
$$ \sum_{j = 0}^{s-2}2^{j}\Delta_{i-j}^{\Big[\substack{s -j\\ s-j +(m+j) }\Big] \times (k+1)} =  \sum_{j = 0}^{s-2}2^{j}\Delta_{i-j}^{\Big[\substack{s -j\\ s-j +(m+j) }\Big] \times k}\;for \;2 \leq  s \leq i\leq 2s+m ,\quad k\geq i+1 . $$ 
 Lemma  \ref{lem 9.2} now  follows  from Theorem \ref{lem  9.1} using the above equations.\\
 \end{proof}

  \section{\textbf{COMPUTATION OF THE REMAINDER  R(j,s,m,k)  IN THE RECURRENT FORMULA FOR THE DIFFERENCE
  $\Gamma_{s+j}^{\left[s\atop s+m\right]\times (k+1)} - \Gamma_{s+j}^{\left[s\atop s+m\right]\times k}\quad for \;  0\leq j\leq s+m,\; k > s+j  $  }}
\label{sec 10}
From a formula of $ \Gamma _{i}^{\Big[\substack{1 \\ 1+ m}\Big] \times k} $ obtained in [4] we compute in this section the 
remainder R(j,s,m,k) in the recurrent formula for the difference $ \Gamma _{i}^{\Big[\substack{s \\ s+ m}\Big] \times (k+1)} -  \Gamma _{i}^{\Big[\substack{s \\ s+ m}\Big] \times k} $\vspace{0.1 cm}.\\
 By  definition  we have for $ s\leq i\leq \inf(2s+m,k) $ \\
$$  \Omega _{i-(s-1)}(1,m+s-1,k) = 
 \Gamma _{i-(s-1)}^{\Big[\substack{1 \\ 1+m+(s-1)}\Big] \times k} - 4\cdot \Gamma _{i-s}^{\Big[\substack{ 1 \\ 1+m+(s-2) }\Big] \times k}. $$ \vspace{0.1 cm} \\
Set i - s = j, then $$\Omega _{j+1}(1,m+s-1,k) =  \Gamma _{j+1}^{\Big[\substack{1 \\ 1+m+(s-1)}\Big] \times k} - 4\cdot \Gamma _{j}^{\Big[\substack{ 1 \\ 1+m+(s-2) }\Big] \times k}. $$ \vspace{0.1 cm} \\

\begin{lem}
\label{lem 10.1}We have by Theorem 3.8 [see  [4] ]  with $ m \rightarrow s+m-1,  m \rightarrow s+m-2  $ \\

\underline {The case  s+m -1  = 0, $ k\geq 2 $}\\
\begin{equation*}
 \Gamma _{j}^{\Big[\substack{1 \\ 1+(s+m-1)}\Big] \times k} = \Gamma _{j}^{\Big[\substack{1 \\ 1}\Big] \times k}= \begin{cases}
 1  & \text{if  } j = 0, \\
 3\cdot(2^k-1)   &  \text{if  }    j =1, \\
2^{2k}- 3\cdot2^k + 2      & \text{if   } j  = 2.
\end{cases}
\end{equation*}

\underline {The case  s+m -1 = 1, $ k\geq 3 $}\\
\begin{equation*}
 \Gamma _{j}^{\Big[\substack{1 \\ 1+ (s+m-1)}\Big] \times k} = \Gamma _{j}^{\Big[\substack{1 \\ 1+ 1}\Big] \times k}= \begin{cases}
 1  & \text{if  } j = 0, \\
2^k + 5  &  \text{if  }    j =1, \\
11\cdot(2^k - 1)   & \text{if   } j = 2, \\
2^{2k+1}- 3\cdot2^{k+2} + 2^4    & \text{if   } j = 3.
\end{cases}
\end{equation*}
\underline {The case $3\leq  k \leq s+m $}\\
\begin{equation*}
 \Gamma _{j+1}^{\Big[\substack{1 \\ 1+(s+m-1) }\Big] \times k}= \begin{cases}
2^k +5  &  \text{if  }    j = 0, \\
3\cdot2^{k+ 2j -2 } + 21\cdot2^{3j -2} & \text{if   } 1\leq j\leq k-2, \\
2^{2k + s+m-1} - 5\cdot2^{3k -5} & \text{if   } j = k -1.
\end{cases}
\end{equation*}
 \underline {The case $ 4 \leq s+m\leq k $}
\begin{equation*}
 \Gamma _{j+1}^{\Big[\substack{1 \\ 1+(s+m-1) }\Big] \times k}= \begin{cases}
2^k +5  &  \text{if  }    j = 0, \\
3\cdot2^{k+ 2j -2 } + 21\cdot2^{3j -2} & \text{if   } 1 \leq j\leq s+m-2, \\
11[2^{k+2s+2m-4} - 2^{3s+3m-5}]    & \text{if   } j = s+m-1, \\
2^{2k+s+m-1}  - 3\cdot2^{k+2s+2m-2} +2^{3s+3m-2}  & \text{if   } j = s+m.
\end{cases}
\end{equation*}
\underline {The case $3\leq  k \leq s+m $}\\
\begin{equation*}
 \Gamma _{j}^{\Big[\substack{1 \\ 1+(s+m-2) }\Big] \times k}= \begin{cases}
1 & \text{if  } j = 0, \\
2^k +5  &  \text{if  }    j=1, \\
3\cdot2^{k+ 2j -4 } + 21\cdot2^{3j -5} & \text{if   } 2\leq j\leq k-1, \\
2^{2k + s+m-2} - 5\cdot2^{3k -5} & \text{if   } j = k.
\end{cases}
\end{equation*}
 \underline {The case $ 4 \leq s+m \leq k $}
\begin{equation*}
 \Gamma _{j}^{\Big[\substack{1 \\ 1+(s+m-2) }\Big] \times k}= \begin{cases}
1 & \text{if  } j = 0, \\
2^k +5  &  \text{if  }    j =1, \\
3\cdot2^{k+ 2j -4 } + 21\cdot2^{3j -5} & \text{if   } 2\leq j\leq s+m-2, \\
11[2^{k+2s+2m-6} - 2^{3s-3m-8}]    & \text{if   } j = s+m -1,\\
2^{2k+s+m-2}  - 3\cdot2^{k+2s+2m-4} +2^{3s+3m-5}  & \text{if   } j = s+m.
\end{cases}
\end{equation*}
\end{lem}
\begin{align*}
&
\end{align*}
\begin{lem}
\label{lem 10.2}
 For $  0\leq j \leq s+m, k\geq s+j+1 $\\
 
  $$   2^{s-1}\left[ \Gamma _{j+1}^{\Big[\substack{1 \\ 1+ (s+m-1)}\Big] \times (k+1)} -  \Gamma _{j+1}^{\Big[\substack{1 \\ 1+(s+m-1) }\Big] \times k}\right] 
   -   2^{s+1}\left[ \Gamma _{j}^{\Big[\substack{1 \\ 1+ (s+m-2)}\Big] \times (k+1)} -  \Gamma _{j}^{\Big[\substack{1 \\ 1+(s+m-2) }\Big] \times k}\right] $$
        is equal to
        \begin{equation*}
  \begin{cases}
2^{k+s-1} & \text{if  } j = 0 ,\quad k\geq s+1, \\
- 2^{k+s-1} &  \text{if  }    j =1, \quad k\geq s+2, \\
 0  & \text{if   } 2\leq j\leq s+m-1, \quad k \geq s+j+1, \\
  -3\cdot 2^{2k+2s+m-2}         & \text{if   } j = s+m, \quad k \geq 2s+m +1. 
\end{cases}
\end{equation*}
  \end{lem}
\begin{proof}

From Lemma \ref{lem 10.1} we deduce \\

 \underline {In the case $j = 0,\; k\geq s+1 $}   \vspace{0.1 cm} \\
  \begin{align*}
&  2^{s-1}\left[ \Gamma _{1}^{\Big[\substack{1 \\ 1+ (s+m-1)}\Big] \times (k+1)} -  \Gamma _{1}^{\Big[\substack{1 \\ 1+(s+m-1) }\Big] \times k}\right] 
   -   2^{s+1}\left[ \Gamma _{0}^{\Big[\substack{1 \\ 1+ (s+m-2)}\Big] \times (k+1)} -  \Gamma _{0}^{\Big[\substack{1 \\ 1+(s+m-2) }\Big] \times k}\right] \\
   & =  2^{s-1}\left[(2^{k+1} + 5) -  (2^{k} + 5)  \right] - 2^{s+1}\left[ 1 - 1\right] = 2^{k+s-1}.\\
   & 
 \end{align*}
 \underline {In the case $j = 1,\; k\geq s+2 $}   \vspace{0.1 cm} \\
   \begin{align*}
&  2^{s-1}\left[ \Gamma _{2}^{\Big[\substack{1 \\ 1+ (s+m-1)}\Big] \times (k+1)} -  \Gamma _{2}^{\Big[\substack{1 \\ 1+(s+m-1) }\Big] \times k}\right] 
   -   2^{s+1}\left[ \Gamma _{1}^{\Big[\substack{1 \\ 1+ (s+m-2)}\Big] \times (k+1)} -  \Gamma _{1}^{\Big[\substack{1 \\ 1+(s+m-2) }\Big] \times k}\right] \\
   & =  2^{s-1}\left[(3\cdot 2^{k+1} + 42) -  (3\cdot 2^{k} + 42)  \right] - 2^{s+1}\left[  (2^{k+1} + 5) -  (2^{k} + 5)  \right] \\
   & = 3\cdot 2^{k+s-1} - 2^{k+s+1} = -2^{k+s-1}.\\
   & 
 \end{align*}
 \underline {In the case $ 2\leq j\leq m-1,\;s+j+1\leq k\leq s+m $}   \vspace{0.1 cm} \\
   \begin{align*}
&  2^{s-1}\left[ \Gamma _{j+1}^{\Big[\substack{1 \\ 1+ (s+m-1)}\Big] \times (k+1)} -  \Gamma _{j+1}^{\Big[\substack{1 \\ 1+(s+m-1) }\Big] \times k}\right] 
   -   2^{s+1}\left[ \Gamma _{j}^{\Big[\substack{1 \\ 1+ (s+m-2)}\Big] \times (k+1)} -  \Gamma _{j}^{\Big[\substack{1 \\ 1+(s+m-2) }\Big] \times k}\right] \\
   & =  2^{s-1}\left[(3\cdot 2^{k+2j-1} + 21\cdot2^{3j-2}) -  (3\cdot 2^{k+2j-2} + 21\cdot2^{3j-2})   \right]\\
   &  - 2^{s+1}\left[ (3\cdot 2^{k+2j-3} + 21\cdot2^{3j-5}) -  (3\cdot 2^{k+2j-4} + 21\cdot2^{3j-5})   \right] \\
   & = 3\cdot 2^{k+2j+s-3} -3\cdot 2^{k+2j+s-3} = 0 \\
   &
 \end{align*}
 \underline {In the case $ 2\leq j\leq m-1\leq s+m-2,\;k\geq s+m $}   \vspace{0.1 cm} \\
   \begin{align*}
&  2^{s-1}\left[ \Gamma _{j+1}^{\Big[\substack{1 \\ 1+ (s+m-1)}\Big] \times (k+1)} -  \Gamma _{j+1}^{\Big[\substack{1 \\ 1+(s+m-1) }\Big] \times k}\right] 
   -   2^{s+1}\left[ \Gamma _{j}^{\Big[\substack{1 \\ 1+ (s+m-2)}\Big] \times (k+1)} -  \Gamma _{j}^{\Big[\substack{1 \\ 1+(s+m-2) }\Big] \times k}\right] \\
   & =  2^{s-1}\left[(3\cdot 2^{k+2j-1} + 21\cdot2^{3j-2}) -  (3\cdot 2^{k+2j-2} + 21\cdot2^{3j-2})   \right]\\
   &  - 2^{s+1}\left[ (3\cdot 2^{k+2j-3} + 21\cdot2^{3j-5}) -  (3\cdot 2^{k+2j-4} + 21\cdot2^{3j-5})   \right] \\
   & = 3\cdot 2^{k+2j+s-3} -3\cdot 2^{k+2j+s-3} = 0.
 \end{align*}
\underline {In the case $2 \leq m \leq j \leq s+m-2 ,\;k\geq s+j+1\geq s+m+1 $}   \vspace{0.1 cm} \\
   \begin{align*}
&  2^{s-1}\left[ \Gamma _{j+1}^{\Big[\substack{1 \\ 1+ (s+m-1)}\Big] \times (k+1)} -  \Gamma _{j+1}^{\Big[\substack{1 \\ 1+(s+m-1) }\Big] \times k}\right] 
   -   2^{s+1}\left[ \Gamma _{j}^{\Big[\substack{1 \\ 1+ (s+m-2)}\Big] \times (k+1)} -  \Gamma _{j}^{\Big[\substack{1 \\ 1+(s+m-2) }\Big] \times k}\right] \\
   & =  2^{s-1}\left[(3\cdot 2^{k+2j-1} + 21\cdot2^{3j-2}) -  (3\cdot 2^{k+2j-2} + 21\cdot2^{3j-2})   \right]\\
   &  - 2^{s+1}\left[ (3\cdot 2^{k+2j-3} + 21\cdot2^{3j-5}) -  (3\cdot 2^{k+2j-4} + 21\cdot2^{3j-5})   \right] \\
   & = 3\cdot 2^{k+2j+s-3} -3\cdot 2^{k+2j+s-3} = 0.
 \end{align*}
\underline {In the case $ j = s+m-1 ,\;k\geq s+j+1=2s+m $}   \vspace{0.1 cm} \\
   \begin{align*}
&  2^{s-1}\left[ \Gamma _{s+m}^{\Big[\substack{1 \\ 1+ (s+m-1)}\Big] \times (k+1)} -  \Gamma _{s+m}^{\Big[\substack{1 \\ 1+(s+m-1) }\Big] \times k}\right] 
   -   2^{s+1}\left[ \Gamma _{s+m-1}^{\Big[\substack{1 \\ 1+ (s+m-2)}\Big] \times (k+1)} -  \Gamma _{s+m-1}^{\Big[\substack{1 \\ 1+(s+m-2) }\Big] \times k}\right] \\
   & =  2^{s-1}\left[11(2^{k+2s+2m-3}- 2^{3s+3m-5}) -  11(2^{k+2s+2m-4}- 2^{3s+3m-5})  \right] \\
   &  - 2^{s+1}\left[11(2^{k+2s+2m-5}- 2^{3s+3m-8}) -  11(2^{k+2s+2m-6}- 2^{3s+3m-8})  \right] \\
   & = 11\cdot2^{k+3s+2m-5} - 11\cdot2^{k+3s+2m-5} = 0.
 \end{align*}
\underline {In the case $ j = s+m ,\;k\geq s+j+1=2s+m+1 $}   \vspace{0.1 cm} \\
   \begin{align*}
&  2^{s-1}\left[ \Gamma _{s+m+1}^{\Big[\substack{1 \\ 1+ (s+m-1)}\Big] \times (k+1)} -  \Gamma _{s+m+1}^{\Big[\substack{1 \\ 1+(s+m-1) }\Big] \times k}\right] 
   -   2^{s+1}\left[ \Gamma _{s+m}^{\Big[\substack{1 \\ 1+ (s+m-2)}\Big] \times (k+1)} -  \Gamma _{s+m}^{\Big[\substack{1 \\ 1+(s+m-2) }\Big] \times k}\right] \\
   & =  2^{s-1}\left[(2^{2k+s+m +1}  - 3\cdot2^{k+2s+2m-1} +2^{3s+3m-2}) -(2^{2k+s+m-1}  - 3\cdot2^{k+2s+2m-2} +2^{3s+3m-2})  \right] \\
   &  - 2^{s+1}\left[(2^{2k+s+m}  - 3\cdot2^{k+2s+2m-3} +2^{3s+3m-5}) -(2^{2k+s+m-2}  - 3\cdot2^{k+2s+2m-4} +2^{3s+3m-5}) \right] \\
   & =  2^{s-1}\left[3\cdot2^{2k+s+m-1} -3\cdot2^{k+2s+2m-2}\right] -  2^{s+1}\left[3\cdot2^{2k+s+m-2} -3\cdot2^{k+2s+2m-4}\right] \\
   & = -3\cdot2^{2k+2s+m-2}.
 \end{align*}
\underline {In all the others cases the proofs are similar.}   
\end{proof}
\begin{lem}
\label{lem 10.3}We have for $ s\geq 2, m\geq 0 $ the following recurrent formula for the difference 
$  \Gamma _{i}^{\Big[\substack{s \\ s +m }\Big] \times (k+1)} -  \Gamma _{i}^{\Big[\substack{s \\ s +m }\Big] \times k} $ where  $ 0\leq i\leq 2s+m, k >  i.$\vspace{0.1 cm} \\

 \begin{equation}
\label{eq 10.1}
 \begin{cases} 
 \Gamma _{j}^{\Big[\substack{s \\ s +m }\Big] \times (k+1)} -  \Gamma _{j}^{\Big[\substack{s \\ s +m }\Big] \times k} =  4\cdot \left[\Gamma _{j-1}^{\Big[\substack{s \\ s +(m-1) }\Big] \times (k+1)} - \Gamma _{j-1}^{\Big[\substack{s \\ s +(m-1) }\Big] \times  k}\right]  
 & \text{if  }  1 \leq j\leq s-1,\quad k\geq j +1,   \\
   \Gamma _{s+j}^{\Big[\substack{s \\ s +m }\Big] \times (k+1)} -  \Gamma _{s+j}^{\Big[\substack{s \\ s +m }\Big] \times k} =  4\cdot \left[\Gamma _{s+j-1}^{\Big[\substack{s \\ s +(m-1) }\Big] \times (k+1)} - \Gamma _{s+j-1}^{\Big[\substack{s \\ s +(m-1) }\Big] \times  k}\right]  \\
+   \begin{cases}
2^{k+s-1} & \text{if  } j = 0 ,\quad k\geq s+1, \\
- 2^{k+s-1} &  \text{if  }    j =1, \quad k\geq s+2, \\
 0  & \text{if   } 2\leq j\leq s+m-1, \quad k \geq s+j+1, \\
  -3\cdot 2^{2k+2s+m-2}         & \text{if   } j = s+m, \quad k \geq 2s+m +1. 
\end{cases}
 \end{cases}
    \end{equation}
    \end{lem}
 \begin{align*}
         &
\end{align*}
\begin{proof}
Lemma  \ref{lem 10.3} follows from Lemma  \ref{lem 9.2} and lemma   \ref{lem 10.2}.\\
\end{proof}

  \section{\textbf{ COMPUTATION OF  $\Gamma_{i}^{\left[s\atop s+m\right]\times k}\;  for\; i\leq s-1,\;i\leq k-1,\;m\geq 0 $}}   
\label{sec 11}
In this section we apply successively the recurrent formula \eqref{eq 10.1} for the difference $\Gamma_{i}^{\left[s\atop s+m\right]\times (k+1)}- \Gamma_{i}^{\left[s\atop s+m\right]\times k}$ 
to compute  $\Gamma_{i}^{\left[s\atop s+m\right]\times k}\;\text{for} \; i\leq s-1,\;i\leq k-1,\;m\geq 0. $\\

\begin{lem}
\label{lem 11.1}Consider the matrix\\
 
  $$   \left ( \begin{array} {cccccc}
\alpha _{1} & \alpha _{2} & \alpha _{3} &  \ldots & \alpha _{k-1}  &  \alpha _{k} \\
\alpha _{2 } & \alpha _{3} & \alpha _{4}&  \ldots  &  \alpha _{k} &  \alpha _{k+1} \\
\vdots & \vdots & \vdots    &  \vdots & \vdots  &  \vdots \\
\alpha _{s-1} & \alpha _{s} & \alpha _{s +1} & \ldots  &  \alpha _{s+k-3} &  \alpha _{s+k-2}  \\
\alpha _{s} & \alpha _{s+1} & \alpha _{s +2} & \ldots  &  \alpha _{s+k-2} &  \alpha _{s+k-1} \\
 \beta  _{1} & \beta  _{2} & \beta  _{3} & \ldots  &  \beta_{k-1} &  \beta _{k}  \\
\beta  _{2} & \beta  _{3} & \beta  _{4} & \ldots  &  \beta_{k} &  \beta _{k+1}  \\
\vdots & \vdots & \vdots    &  \vdots & \vdots  &  \vdots \\
\beta  _{m+1} & \beta  _{m+2} & \beta  _{m+3} & \ldots  &  \beta_{k+m-1} &  \beta _{k+m}  \\
\vdots & \vdots & \vdots    &  \vdots & \vdots  &  \vdots \\
\beta  _{s+m-1} & \beta  _{s+m} & \beta  _{s+m+1} & \ldots  &  \beta_{s+m+k-3} &  \beta _{s+m+k-2}  \\
\beta  _{s+m} & \beta  _{s+m+1} & \beta  _{s+m+2} & \ldots  &  \beta_{s+m+k-2} &  \beta _{s+m+k-1}
\end{array}  \right). $$ \vspace{0.5 cm}\\

Let $ g_{k,s,m}(t,\eta )  $ be the quadratic exponential sum in  $ \mathbb{P}\times \mathbb{P} $ defined by 
  $$ (t,\eta ) \in  \mathbb{P}\times \mathbb{P}\longmapsto  
  \sum_{deg Y\leq k-1}\sum_{deg Z \leq s-1}E(tYZ)\sum_{deg U \leq  s+m-1}E(\eta YU) \in \mathbb{Z}.  $$\vspace{0.5 cm}
   We have
   \begin{align}
  \sum_{ i = 0}^{\inf(2s+m,k)} \Gamma _{i}^{\Big[\substack{s  \\ s+m }\Big] \times k} & = 2^{2k+2s+m-2} \label{eq 11.1}\\
 &  \text{and}\nonumber \\
  \sum_{i = 0}^{\inf(2s+m,k)} \Gamma _{i}^{\Big[\substack{s  \\ s+m }\Big] \times k}\cdot2^{-i} & =  2^{k+2s+m-2} + 2^{2k-2} - 2^{k-2}. \label{eq 11.2}
\end{align}
\end{lem}
\begin{proof}
The proof of \eqref{eq 11.1} is obvious.\\

By Lemma  \ref{lem 4.6} we obtain, observing that  $ g_{k,s,m}(t,\eta )  $ is constant on cosets of $ \mathbb{P}_{k+s-1}\times \mathbb{P}_{k+s+m-1} $ \\
 \begin{align*}
\int_{\mathbb{P}\times \mathbb{P}} g_{k,s,m}(t,\eta )dtd\eta &  =\sum_{(t,\eta )\in \mathbb{P}/\mathbb{P}_{k+s-1}\times \mathbb{P}/\mathbb{P}_{k+s +m-1}}
 2^{2s+m+k- r( D^{\left[\stackrel{s}{s+m}\right] \times k }(t,\eta  ) )}\int_{\mathbb{P}_{k+s-1}}dt \int_{\mathbb{P}_{k+s+m-1}}d\eta  \\
 & = \sum_{i=0}^{\inf(2s+m,k)} \Gamma _{i}^{\Big[\substack{s  \\ s+m }\Big] \times k}\cdot 2^{2s+m+k-i}\int_{\mathbb{P}_{k+s-1}}dt \int_{\mathbb{P}_{k+s+m-1}}d\eta \\
 & = 2^{-k+2}\cdot \sum_{i = 0}^{\inf(2s+m,k)} \Gamma _{i}^{\Big[\substack{s  \\ s+m }\Big] \times k}\cdot2^{-i}. 
\end{align*}
On the other hand \\
 \begin{align*}
\int_{\mathbb{P}\times \mathbb{P}} g_{k,s,m}(t,\eta ) dt d \eta & =
 Card \left\{(Y,Z,U), deg Y \leq k-1, deg Z \leq s-1,  deg U \leq s+m-1 \mid Y \cdot Z = Y \cdot U = 0 \right\}\\
 & = 2^{2s+m} + 2^{k} -1. 
 \end{align*}
 The above equations imply  \eqref{eq 11.2}.
\end{proof}
\begin{lem}
\label{lem 11.2}  We have  for all $ j\geq 1 $\vspace{0.1 cm}\\
 \begin{equation}
 \label{eq 11.3}
 (H_{j}) \hspace{3 cm}  \Gamma _{j}^{\Big[\substack{s  \\ s+m }\Big] \times k} =  \Gamma _{j}^{\Big[\substack{j  \\ j}\Big] \times (j+1)}
  \overset{ \text{def}}{=} \gamma _{j}
    \quad   \text{ for   \quad   $ s\geq j+1,\;k\geq j+1,\; m\geq 0.  $}\\
 \end{equation}
\end{lem}
\begin{proof}The proof is by strong induction, that is:\vspace{0.1 cm}\\
If
\begin{itemize}
\item \text{$ (H_{1}) $ is true, and  }
\item \text{for all $ j\geq 1, \quad (H_{1})\wedge (H_{2})\wedge\ldots\wedge (H_{j})\quad implies \; (H_{j+1}) $}
\end{itemize}
then $ (H_{j}) $ is true for all $ j \geq 1. $\vspace{0.1 cm}\\

\underline {$ (H_{1})\, is \;true\; $}\vspace{0.1 cm}\\

Indeed from  \eqref{eq 10.1}  with   j = 1   \vspace{0.1 cm} \\
$ \Gamma _{1}^{\Big[\substack{s \\ s +m }\Big] \times (k+1)} -  \Gamma _{1}^{\Big[\substack{s \\ s +m }\Big] \times k}
 =  4\cdot \left[\Gamma _{0}^{\Big[\substack{s \\ s +(m-1) }\Big] \times (k+1)} - \Gamma _{0}^{\Big[\substack{s \\ s +(m-1) }\Big] \times  k}\right] $
  for all $ s \geq 2,\; k\geq 2,\; m\geq  0 $ \vspace{0.1 cm}\\
  which implies that $ (H_{j})$ holds for j = 1, that is
  \begin{equation}
\label{eq 11.4}
   \Gamma _{1}^{\Big[\substack{s \\ s +m }\Big] \times k}=   \Gamma _{1}^{\Big[\substack{s \\ s +m }\Big] \times 2} = \Gamma _{1}^{\Big[\substack{1 \\ 1}\Big] \times 2}= \gamma _{1}= 9 
 \quad  \forall\; s\geq 2, \;  \forall\; k \geq 2, \; m\geq  0. 
\end{equation}
\underline { $ (H_{1}) \Rightarrow  (H_{2}) $}\vspace{0.1 cm}\\

  From  \eqref{eq 10.1}  with   j = 2   \vspace{0.1 cm} \\
$ \Gamma _{2}^{\Big[\substack{s \\ s +m }\Big] \times (k+1)} -  \Gamma _{2}^{\Big[\substack{s \\ s +m }\Big] \times k}
 =  4\cdot \left[\Gamma _{1}^{\Big[\substack{s \\ s +(m-1) }\Big] \times (k+1)} - \Gamma _{1}^{\Big[\substack{s \\ s +(m-1) }\Big] \times  k}\right] $
  for all $ s\geq 3,\;k\geq 3,\ m\geq 0. $ \vspace{0.1 cm}\\
  Hence by \eqref{eq 11.4} we have 
    \begin{equation}
    \label{eq 11.5}
   \Gamma _{2}^{\Big[\substack{s \\ s +m }\Big] \times k } = \Gamma _{2}^{\Big[\substack{s \\ s +m }\Big] \times 3 } \quad \text{for all}\; s \geq 3,\;k\geq 3,\; m\geq 0.
\end{equation}
Consider the matrix\\
   $$   \left ( \begin{array} {ccc}
\alpha _{1} & \alpha _{2} & \alpha _{3}  \\
\alpha _{2 } & \alpha _{3} & \alpha _{4} \\
\vdots & \vdots & \vdots     \\
\alpha _{s} & \alpha _{s+1} & \alpha _{s +2}   \\
\beta  _{1} & \beta  _{2} & \beta  _{3}   \\
\beta  _{2} & \beta  _{3} & \beta  _{4}   \\
\vdots & \vdots & \vdots   \\
\beta  _{m+1} & \beta  _{m+2}  & \beta  _{m+3}  \\
\vdots & \vdots & \vdots  \\
\beta  _{s+m} & \beta  _{s+m+1} & \beta  _{s+m+2}   
\end{array}  \right). $$ \vspace{0.1 cm}\\
 We have  respectively  by \eqref{eq 11.1},  \eqref{eq 11.2} with k = 3 \\
 \begin{align}
  \sum_{i=0}^{3} \Gamma _{i}^{\Big[\substack{s \\ s +m }\Big] \times 3 } &  = 2^{2s+m+4}, \label{eq 11.6} \\
  & \nonumber \\
    \sum_{i=0}^{3} \Gamma _{i}^{\Big[\substack{s \\ s +m }\Big] \times 3 }\cdot 2^{-i} &  =  2^{2s+m+1} + 2^{4} - 2.  \label{eq 11.7}
\end{align}
From \eqref{eq 11.6}  and  \eqref{eq 11.7} using   \eqref{eq 11.4}  we deduce \vspace{0.1 cm}\\
$$ \Gamma _{2}^{\Big[\substack{s \\ s +m }\Big] \times 3 } = \Gamma _{2}^{\Big[\substack{2 \\ 2}\Big] \times 3 }=  78 \quad \text{and} \quad \Gamma _{3}^{\Big[\substack{s \\ s +m }\Big] \times 3 }=  2^{2s+m+4} - 88 $$ \vspace{0.1 cm}\\
which implies, using \eqref{eq 11.5}, that $ (H_{j})$ holds for j = 2, that is
  \begin{equation}
\label{eq 11.8}
   \Gamma _{2}^{\Big[\substack{s \\ s +m }\Big] \times k}=  \Gamma _{2}^{\Big[\substack{2 \\ 2}\Big] \times 3}= 78 
 \quad  \forall\; s\geq 3,  \; \forall\; k \geq 3, \; \forall\; m \geq 0. 
\end{equation}\vspace{0.1 cm}\\
 \underline {$ \quad (H_{1})\wedge (H_{2})\; implies \; (H_{3}) $}\vspace{0.1 cm}\\
 
    From  \eqref{eq 10.1}  with   j = 3   \vspace{0.1 cm} \\
$ \Gamma _{3}^{\Big[\substack{s \\ s +m }\Big] \times (k+1)} -  \Gamma _{3}^{\Big[\substack{s \\ s +m }\Big] \times k}
 =  4\cdot \left[\Gamma _{2}^{\Big[\substack{s \\ s +(m-1) }\Big] \times (k+1)} - \Gamma _{2}^{\Big[\substack{s \\ s +(m-1) }\Big] \times  k}\right] $
  for all $ s\geq 4, \;k\geq 4,\; m\geq 0. $ \vspace{0.1 cm}\\
  Hence by \eqref{eq 11.8} we have 
    \begin{equation}
    \label{eq 11.9}
   \Gamma _{3}^{\Big[\substack{s \\ s +m }\Big] \times k } = \Gamma _{3}^{\Big[\substack{s \\ s +m }\Big] \times 4 } \quad \text{for all}\; s \geq 4, \; k \geq 4,\; m\geq 0.
\end{equation}
Consider the matrix\\
   $$   \left ( \begin{array} {cccc}
\alpha _{1} & \alpha _{2} & \alpha _{3} & \alpha _{4}  \\
\alpha _{2 } & \alpha _{3} & \alpha _{4}& \alpha _{5}  \\
\vdots & \vdots & \vdots & \vdots    \\
\alpha _{s} & \alpha _{s+1} & \alpha _{s +2} & \alpha _{s +3}   \\
\beta  _{1} & \beta  _{2} & \beta  _{3} & \beta  _{4}   \\
\beta  _{2} & \beta  _{3} & \beta  _{4}   & \beta  _{5}  \\
\vdots & \vdots & \vdots & \vdots     \\
\beta  _{m+1} & \beta  _{m+2}  & \beta  _{m+3}  & \beta  _{m+4} \\
\vdots & \vdots & \vdots  & \vdots \\
\beta  _{s+m} & \beta  _{s+m+1} & \beta  _{s+m+2}  & \beta  _{s+m+3}     
\end{array}  \right). $$ \vspace{0.1 cm}\\
 We have  respectively  by \eqref{eq 11.1},  \eqref{eq 11.2} with k = 4 \\
 \begin{align}
  \sum_{i=0}^{4} \Gamma _{i}^{\Big[\substack{s \\ s +m }\Big] \times 4 } &  = 2^{2s+m+6}, \label{eq 11.10} \\
  & \nonumber \\
   \sum_{i=0}^{4} \Gamma _{i}^{\Big[\substack{s \\ s +m }\Big] \times 4 }\cdot 2^{-i} &  =  2^{2s+m+2} + 2^{6} - 2^{2}.  \label{eq 11.11}
\end{align}
From \eqref{eq 11.10}  and  \eqref{eq 11.11} using   \eqref{eq 11.5}  and  \eqref{eq 11.9} we deduce \vspace{0.1 cm}\\
$$ \Gamma _{3}^{\Big[\substack{s \\ s +m }\Big] \times 4 } = \Gamma _{3}^{\Big[\substack{3 \\ 3}\Big] \times 4 }= 648  \quad \text{and} \quad \Gamma _{4}^{\Big[\substack{s \\ s +m }\Big] \times 4 }=  2^{2s+m+6}- 736 $$
which implies, using \eqref{eq 11.9}, that $ (H_{j})$ holds for j = 3, that is
  \begin{equation}
\label{eq 11.12}
   \Gamma _{3}^{\Big[\substack{s \\ s +m }\Big] \times k}=  \Gamma _{3}^{\Big[\substack{3 \\ 3}\Big] \times 4}= 648 = \gamma _{3}
 \quad \text{for all}\; s\geq 4, \; k \geq 4, \; m\geq 0.
\end{equation}\vspace{0.05 cm}\\

 \underline { $  (H_{1})\wedge (H_{2})\wedge\ldots\wedge (H_{j})\; implies \; (H_{j+1})  $}\vspace{0.1 cm}\\
 
 From  \eqref{eq 10.1} with  j $\rightarrow $ j+1   \vspace{0.1 cm} \\
 
$ \Gamma _{j+1}^{\Big[\substack{s \\ s +m }\Big] \times (k+1)} -  \Gamma _{j+1}^{\Big[\substack{s \\ s +m }\Big] \times k}
 =  4\cdot \left[\Gamma _{j}^{\Big[\substack{s \\ s +(m-1) }\Big] \times (k+1)} - \Gamma _{j}^{\Big[\substack{s \\ s +(m-1) }\Big] \times  k}\right] $
  for all $ s\geq j+2,\;k\geq j+2,\; m\geq 0. $ \vspace{0.05 cm}\\
  By $ (H_{j}) $ we obtain  $$ \Gamma _{j}^{\Big[\substack{s \\ s +(m-1) }\Big] \times (k+1)} -  \Gamma _{j}^{\Big[\substack{s \\ s +(m-1) }\Big] \times k} = \gamma _{j} - \gamma _{j} = 0.  $$
  Indeed if m = 0 \vspace{0.1 cm}\\
$$  \Gamma _{j}^{\Big[\substack{s \\ s +(m-1) }\Big] \times k}=  \Gamma _{j}^{\Big[\substack{s -1 \\ s -1 +1 }\Big] \times k} = 
   \Gamma _{j}^{\Big[\substack{j \\ j }\Big] \times (j+1)} = 
   \gamma _{j} $$ \vspace{0.1 cm}\\
    since $$  s-1\leq j+1,\quad k \geq j+2 \geq j+1, \quad m\rightarrow 1\geq 0. $$\vspace{0.05 cm}\\
 If $ m\geq 1 $  \vspace{0.1 cm}\\
 $$ \Gamma _{j}^{\Big[\substack{s \\ s +(m-1) }\Big] \times k} = \Gamma _{j}^{\Big[\substack{j \\ j }\Big] \times (j+1)} = \gamma _{j} $$ \vspace{0.1 cm}\\
   since $$ s\geq j+2\geq j+1,\quad k \geq j+2 \geq j+1 , \quad m \rightarrow (m-1)\geq 0. $$\vspace{0.05 cm}\\
 Thus, we obtain \vspace{0.1 cm}\\
 \begin{equation}
\label{eq 11.13}
\Gamma _{j+1}^{\Big[\substack{s \\ s +m }\Big] \times k} = \Gamma _{j+1}^{\Big[\substack{s \\ s +m }\Big] \times (j+2)}\quad \text{for}\; s \geq j+2,\;k\geq j+2,\;m \geq 0.
\end{equation}
  Consider the matrix\\
   
   $$   \left ( \begin{array} {cccccc}
\alpha _{1} & \alpha _{2} & \alpha _{3} &  \ldots & \alpha _{j+1}  &  \alpha _{j+2} \\
\alpha _{2 } & \alpha _{3} & \alpha _{4}&  \ldots  &  \alpha _{j+2} &  \alpha _{j+3} \\
\vdots & \vdots & \vdots    &  \vdots & \vdots  &  \vdots \\
\alpha _{s-1} & \alpha _{s} & \alpha _{s +1} & \ldots  &  \alpha _{s+j-1} &  \alpha _{s+j}  \\
\alpha _{s} & \alpha _{s+1} & \alpha _{s +2} & \ldots  &  \alpha _{s+ j} &  \alpha _{s+j+1} \\
 \beta  _{1} & \beta  _{2} & \beta  _{3} & \ldots  &  \beta_{j+1} &  \beta _{j+2}  \\
\beta  _{2} & \beta  _{3} & \beta  _{4} & \ldots  &  \beta_{j+2} &  \beta _{j+3}  \\
\vdots & \vdots & \vdots    &  \vdots & \vdots  &  \vdots \\
\beta  _{m+1} & \beta  _{m+2} & \beta  _{m+3} & \ldots  &  \beta_{j+m+1} &  \beta _{j+m+2}  \\
\vdots & \vdots & \vdots    &  \vdots & \vdots  &  \vdots \\
\beta  _{s+m-1} & \beta  _{s+m} & \beta  _{s+m+1} & \ldots  &  \beta_{s+m+j-1} &  \beta _{s+m+j}  \\
\beta  _{s+m} & \beta  _{s+m+1} & \beta  _{s+m+2} & \ldots  &  \beta_{s+m+ j} &  \beta _{s+m+j+1}
\end{array}  \right). $$ \\

 We have  respectively  by \eqref{eq 11.1},  \eqref{eq 11.2} with k = j+2 \\
   \begin{align}
  \sum_{ i = 0}^{j+2} \Gamma _{i}^{\Big[\substack{s  \\ s+m }\Big] \times (j+2)} & = 2^{2j+2s+m +2} \label{eq 11.14}\\
 &  \text{and}\nonumber \\
  \sum_{i = 0}^{j+2} \Gamma _{i}^{\Big[\substack{s  \\ s+m }\Big] \times (j+2)}\cdot2^{-i} & =  2^{j+2s+m} + 2^{2j +2} - 2^{j}. \label{eq 11.15}
\end{align}
  Setting s = j+1, m = 0   in  \eqref{eq 11.14},  \eqref{eq 11.15} we have \vspace{0.5 cm}  \\
     \begin{align}
  \sum_{ i = 0}^{j+2} \Gamma _{i}^{\Big[\substack{j+1 \\ j+1}\Big] \times (j+2)} & = 2^{4j +4} \label{eq 11.16}\\
 &  \text{and}\nonumber \\
  \sum_{i = 0}^{j+2} \Gamma _{i}^{\Big[\substack{j+1 \\ j+1}\Big] \times (j+2)}\cdot2^{- i} & =  2^{3j +2} + 2^{2j +2} - 2^{j}. \label{eq 11.17}
\end{align}
From $(H_{1})\wedge (H_{2})\wedge\ldots\wedge (H_{j}), $ it follows that\vspace{0.5 cm}  \\
\begin{equation}
\label{eq 11.18}
 \Gamma _{i}^{\Big[\substack{s  \\ s+m }\Big] \times (j+2)} =  \Gamma _{i}^{\Big[\substack{j+1  \\ j+1 }\Big] \times (j+2)}= 
  \Gamma _{i}^{\Big[\substack{i \\ i}\Big] \times (i+1)}= \gamma _{i}\quad \text{for}\quad  0\leq i\leq j.
\end{equation}\vspace{0.5 cm}  \\
We deduce from  \eqref{eq 11.14},  \eqref{eq 11.15} \vspace{0.5 cm}  \\
\begin{equation}
\label{eq 11.19}
 \Gamma _{j+1}^{\Big[\substack{s  \\ s+m }\Big] \times (j+2)} =  2^{3j+4} - 2^{2j+2} + \sum_{i = 0}^{j}\gamma _{i}(1-2^{j+2-i})
\end{equation}
and from  \eqref{eq 11.16},  \eqref{eq 11.17} \vspace{0.5 cm}  \\
\begin{equation}
\label{eq 11.20}
 \Gamma _{j+1}^{\Big[\substack{j+1  \\ j+1}\Big] \times (j+2)} =  2^{3j+4} - 2^{2j +2} + \sum_{i = 0}^{j}\gamma _{i}(1-2^{j+2-i}).
\end{equation}\vspace{0.5 cm}  \\
By  \eqref{eq 11.13},  \eqref{eq 11.19} and   \eqref{eq 11.20}  it follows that  \vspace{0.5 cm}  \\ 
\begin{equation*}
 \Gamma _{j+1}^{\Big[\substack{s  \\ s+m }\Big] \times k} = 
  \Gamma _{j+1}^{\Big[\substack{j+1 \\ j+1 }\Big] \times (j+2)}= \gamma _{j+1}\quad \text{for all}\;  s\geq j+2, \;k\geq j+2, \;m\geq 0.
\end{equation*}
\end{proof}
\begin{lem}
\label{lem 11.3}We have \vspace{0.5 cm}  \\
\begin{equation}
\label{eq 11.21}
 \Gamma _{i}^{\Big[\substack{s \\ s+m }\Big] \times k}= \Gamma _{i}^{\Big[\substack{i \\ i }\Big] \times (i+1)} = \begin{cases}
 1  & \text{if  } i = 0, \\
 21\cdot2^{3i-4} - 3\cdot2^{2i -3}  &  \text{if  }  1\leq i\leq s-1,\;k\geq i+1,\;m\geq 0. 
\end{cases}
\end{equation}\vspace{0.1 cm}  \\
\end{lem}
\begin{proof}
We have respectively by \eqref{eq 11.1}, \eqref{eq 11.2} with k = i+1, s = i  and m = 0,  using \eqref{eq 11.3}  \vspace{0.1 cm}  \\
  \begin{align}
  \sum_{ j = 0}^{i+1} \Gamma _{j}^{\Big[\substack{i  \\ i }\Big] \times (i+1)}& = \sum_{j = 0}^{i}\gamma _{j} + \Gamma _{i+1}^{\Big[\substack{i  \\ i }\Big] \times (i+1)}   = 2^{4i} \label{eq 11.22}\\
 &  \text{and}\nonumber \\
  \sum_{j = 0}^{i+1} \Gamma _{j}^{\Big[\substack{i  \\ i }\Big] \times (i+1)}\cdot2^{-j} & = 
    \sum_{j = 0}^{i}\gamma _{j}\cdot2^{-j} + \Gamma _{i+1}^{\Big[\substack{i  \\ i }\Big] \times (i+1)}\cdot2^{-(i+1)} = 2^{3i-1} + 2^{2i} - 2^{i-1}.  \label{eq 11.23}
\end{align}\vspace{0.1 cm}  \\
From \eqref{eq 11.22} and \eqref{eq 11.23} we get \vspace{0.1 cm}  \\
\begin{align}
 & 2^{i+1}\cdot \left( \sum_{j = 0}^{i}\gamma _{j}\cdot2^{-j} + \Gamma _{i+1}^{\Big[\substack{i  \\ i }\Big] \times (i+1)}\cdot2^{-(i+1)}\right) - 
\left(\sum_{j = 0}^{i}\gamma _{j} + \Gamma _{i+1}^{\Big[\substack{i  \\ i }\Big] \times (i+1)} \right) = 2^{3i+1} - 2^{2i} \nonumber    \\ 
&\Leftrightarrow \sum_{j = 0}^{i}\gamma _{j}\cdot(2^{i+1-j} - 1) =  2^{3i+1} - 2^{2i}. \label{eq 11.24}
\end{align}\vspace{0.1 cm}  \\
Hence by \eqref{eq 11.24} we deduce \vspace{0.1 cm}  \\
\begin{align}
& \sum_{j = 0}^{i}\gamma _{j}\cdot(2^{i+1-j} - 1) - 2\cdot\left(\sum_{j = 0}^{i-1}\gamma _{j}\cdot(2^{i -j} - 1)\right) 
 = ( 2^{3i+1} - 2^{2i})  - 2\cdot( 2^{3i -2} - 2^{2i-2}) \nonumber \\
 & \Leftrightarrow \sum_{j = 0}^{i}\gamma _{j}= 3\cdot2^{3i-1} - 2^{2i-1}. \label{eq 11.25}
\end{align}\vspace{0.1 cm}  \\
From \eqref{eq 11.25} we get for $ i\geq 2 $ \vspace{0.1 cm}  \\
\begin{align*}
& \sum_{j = 0}^{i}\gamma _{j} - \sum_{j = 0}^{i-1}\gamma _{j}  = (3\cdot2^{3i-1} - 2^{2i-1}) - (3\cdot2^{3i-4} - 2^{2i-3}) \\
& \Leftrightarrow  \Gamma _{i}^{\Big[\substack{i \\ i }\Big] \times (i+1)}= \gamma _{i} = 21\cdot2^{3i-4} - 3\cdot 2^{2i -3}.
\end{align*}
\end{proof}
\begin{lem}
\label{lem  11.4}We have
\begin{align}
 \Gamma _{k}^{\Big[\substack{s  \\ s+m }\Big] \times k} &  = 2^{2k+2s+m-2} -3\cdot2^{3k-4} + 2^{2k-3}\quad for \; 1\leq k\leq s, \label{eq 11.26}\\
 & \nonumber \\
  \sum_{j = 0}^{s -1} \Gamma _{j}^{\Big[\substack{s  \\ s+m }\Big] \times k} & = 3\cdot2^{3s-4} - 2^{2s-3}\quad for \; k\geq s, \label{eq 11.27}\\
  & \nonumber  \\
 \sum_{j = 0}^{s -1} \Gamma _{j}^{\Big[\substack{s  \\ s+m }\Big] \times k}\cdot 2^{-j}  & = 7\cdot2^{2s - 4} - 3\cdot2^{s-3} \quad for \; k\geq s. \label{eq 11.28}
\end{align}
\end{lem}
\begin{proof}
Consider the matrix \\
   $$   \left ( \begin{array} {cccccc}
\alpha _{1} & \alpha _{2} & \alpha _{3} &  \ldots & \alpha _{k-1}  &  \alpha _{k} \\
\alpha _{2 } & \alpha _{3} & \alpha _{4}&  \ldots  &  \alpha _{k} &  \alpha _{k+1} \\
\vdots & \vdots & \vdots    &  \vdots & \vdots  &  \vdots \\
\alpha _{s-1} & \alpha _{s} & \alpha _{s +1} & \ldots  &  \alpha _{s+k-3} &  \alpha _{s+k-2}  \\
\alpha _{s} & \alpha _{s+1} & \alpha _{s +2} & \ldots  &  \alpha _{s+k-2} &  \alpha _{s+k-1} \\
 \beta  _{1} & \beta  _{2} & \beta  _{3} & \ldots  &  \beta_{k-1} &  \beta _{k}  \\
\beta  _{2} & \beta  _{3} & \beta  _{4} & \ldots  &  \beta_{k} &  \beta _{k+1}  \\
\vdots & \vdots & \vdots    &  \vdots & \vdots  &  \vdots \\
\beta  _{m+1} & \beta  _{m+2} & \beta  _{m+3} & \ldots  &  \beta_{k+m -1} &  \beta _{k+m}  \\
\vdots & \vdots & \vdots    &  \vdots & \vdots  &  \vdots \\
\beta  _{s+m-1} & \beta  _{s+m} & \beta  _{s+m+1} & \ldots  &  \beta_{s+m+k-3} &  \beta _{s+m+k-2}  \\
\beta  _{s+m} & \beta  _{s+m+1} & \beta  _{s+m+2} & \ldots  &  \beta_{s+m+k-2} &  \beta _{s+m+k-1}
\end{array}  \right). $$ \\
 
Obviously by \eqref{eq 11.25} we have  \vspace{0.1 cm}\\
\begin{align*}
 \Gamma _{k}^{\Big[\substack{s  \\ s+m }\Big] \times k}  
& =  2^{2k+2s+m-2} - \sum_{j = 0}^{k-1} \Gamma _{j}^{\Big[\substack{s  \\ s+m }\Big] \times k} \\
& =   2^{2k+2s+m-2} - \sum_{j = 0}^{k-1} \gamma _{j} \\
&  =  2^{2k+2s+m-2}- (3\cdot 2^{3(k-1)-1} - 2^{2(k-1) -1}) = 2^{2k+2s+m-2} -3\cdot2^{3k-4} + 2^{2k-3}.  
\end{align*} \vspace{0.1 cm}\\
We deduce \eqref{eq 11.27} from \eqref{eq 11.25} with i = s -1. \vspace{0.1 cm}\\
From \eqref{eq 11.21} we obtain after some calculations \vspace{0.1 cm}\\
$$ \sum_{j = 0}^{s -1} \Gamma _{j}^{\Big[\substack{s  \\ s+m }\Big] \times k} \cdot 2^{-j} =
 1 +  \sum_{j = 1}^{s -1} ( 21\cdot2^{2j-4} - 3\cdot 2^{j-3}) = 7\cdot2^{2s-4} - 3\cdot 2^{s-3}. $$\\
 
  \end{proof}

 \section{\textbf{ COMPUTATION OF  $\quad \Gamma_{s+j}^{\left[s \atop s+m\right]\times (k+1)} - 
 \Gamma_{s+j}^{\left[s \atop s+m\right]\times k}\;  for \; 0 \leq j\leq s+m,\;k\geq s+j +1,\; m\geq 0 $}}
  \label{sec 12}
  In this section we apply successively the recurrent formula \eqref{eq 10.1} to compute explicitly the difference $\Gamma_{s+j}^{\left[s\atop s+m\right]\times (k+1)}- \Gamma_{s+j}^{\left[s\atop s+m\right]\times k}.$ \\

 \begin{lem}
\label{lem 12.1} We have
\begin{equation}
\label{eq 12.1}
 \Gamma _{s+j}^{\Big[\substack{s \\ s}\Big] \times (k+1)}
  - \Gamma _{s+j}^{\Big[\substack{s \\ s}\Big] \times k}=
  \begin{cases}
 3\cdot 2^{k+s-1}  & \text{if  } j = 0,\; k >s, \\
 21\cdot 2^{k+s+3j-4}   &  \text{if  }    1\leq j\leq s-1,\; k>s+j, \\
3\cdot2^{2k+2s-2}- 3\cdot2^{k+4s - 4}       & \text{if   } j  = s,  \; k>2s.
\end{cases}
\end{equation}
\end{lem}
\begin{proof}
From \eqref{eq 10.1} we have the following formula \vspace{0.1 cm}\\

 \begin{align}
 \label{eq 12.2}
 \Gamma _{s+j}^{\Big[\substack{s \\ s +m }\Big] \times (k+1)} -  \Gamma _{s+j}^{\Big[\substack{s \\ s +m }\Big] \times k} & =  4\cdot \left[\Gamma _{s+j-1}^{\Big[\substack{s  \\ s +(m-1)) }\Big] \times (k+1)} - \Gamma _{s+j-1}^{\Big[\substack{s   \\ s +(m-1))  }\Big] \times  k}\right] \\
 &  +   \begin{cases}
2^{k+s-1} & \text{if  } j = 0 ,\quad k\geq s+1, \\
- 2^{k+s-1} &  \text{if  }    j =1, \quad k\geq s+2, \\
 0  & \text{if   } 2\leq j\leq s +m  -1, \quad k \geq s+j+1, \\
  -3\cdot 2^{2k+2s+ m  -2}         & \text{if   } j = s +m ,\;  k \geq 2s + m+1. \nonumber
\end{cases}
  \end{align}\vspace{0.1 cm}\\
  
\underline {The case  j  = 0,  $ k > s $}\vspace{0.1 cm}\\

Applying formula  \eqref{eq 12.2} we obtain using \eqref{eq 11.21} \vspace{0.1 cm}\\
\begin{align}
 \Gamma _{s}^{\Big[\substack{s \\ s  }\Big] \times (k+1)} -  \Gamma _{s}^{\Big[\substack{s \\ s  }\Big] \times k} &  = 
 4\cdot \left[\Gamma _{s-1}^{\Big[\substack{s -1 \\ s -1 +1 }\Big] \times (k+1)} - \Gamma _{s-1}^{\Big[\substack{s -1  \\ s -1+1  }\Big] \times  k}\right]  + 2^{k+s-1}, \label{eq 12.3}\\
 & \nonumber \\
  \Gamma _{s -1}^{\Big[\substack{s -1 \\ s-1+1  }\Big] \times (k+1)} -  \Gamma _{s-1}^{\Big[\substack{s-1 \\ s -1+1 }\Big] \times k} & = 
 4\cdot \left[\Gamma _{s-2}^{\Big[\substack{s -1 \\ s -1  }\Big] \times (k+1)} - \Gamma _{s-2 }^{\Big[\substack{s -1  \\ s -1  }\Big] \times  k}\right]  + 2^{k+s-2} =  2^{k+s-2}. \label{eq 12.4}
 \end{align}\vspace{0.1 cm}\\
From \eqref{eq 12.3},  \eqref{eq 12.4} we get \vspace{0.1 cm}\\
\begin{equation}
\label{eq 12.5}
 \Gamma _{s}^{\Big[\substack{s \\ s  }\Big] \times (k+1)} -  \Gamma _{s}^{\Big[\substack{s \\ s  }\Big] \times k}  = 
 4\cdot2^{k+s-2} + 2^{k+s-1} = 3\cdot2^{k+s-1}.
\end{equation}\\

 \underline {The case  j  = 1, $ k > s+1 $}\vspace{0.1 cm}\\
 
We proceed as in the case j = 0 using \eqref{eq 12.5} with $ s\rightarrow s-1  $\vspace{0.1 cm}\\
\begin{align}
 \Gamma _{s +1}^{\Big[\substack{s \\ s  }\Big] \times (k+1)} -  \Gamma _{s +1}^{\Big[\substack{s \\ s  }\Big] \times k} &  = 
 4\cdot \left[\Gamma _{s}^{\Big[\substack{s -1 \\ s -1 +1 }\Big] \times (k+1)} - \Gamma _{s}^{\Big[\substack{s -1  \\ s -1+1  }\Big] \times  k}\right]  - 2^{k+s-1}, \label{eq 12.6}\\
 & \nonumber \\
  \Gamma _{s }^{\Big[\substack{s -1 \\ s-1+1  }\Big] \times (k+1)} -  \Gamma _{s}^{\Big[\substack{s-1 \\ s -1+1 }\Big] \times k} &  = 
 4\cdot \left[\Gamma _{s-1}^{\Big[\substack{s -1 \\ s -1  }\Big] \times (k+1)} - \Gamma _{s-1 }^{\Big[\substack{s -1  \\ s -1  }\Big] \times  k}\right]  - 2^{k+s-2}    \label{eq 12.7}\\
& = 4\cdot3\cdot2^{k+s-2} - 2^{k+s-2} = 11\cdot2^{k+s-2}. \nonumber
 \end{align}\vspace{0.1 cm}\\
From \eqref{eq 12.6},  \eqref{eq 12.7} we get  \vspace{0.1 cm}\\
\begin{equation}
\label{eq 12.8}
 \Gamma _{s+1}^{\Big[\substack{s \\ s  }\Big] \times (k+1)} -  \Gamma _{s+1}^{\Big[\substack{s \\ s  }\Big] \times k}  = 
 4\cdot11 \cdot2^{k+s-2} - 2^{k+s-1} = 21\cdot2^{k+s-1}. 
\end{equation}\\

 \underline {The case  j  = 2, $ k > s+2 $}\vspace{0.1 cm}\\
 
 Proceeding as before, using \eqref{eq 12.8} with $ s \rightarrow s-1 $\vspace{0.1 cm}\\
\begin{align}
 \Gamma _{s +2}^{\Big[\substack{s \\ s  }\Big] \times (k+1)} -  \Gamma _{s +2}^{\Big[\substack{s \\ s  }\Big] \times k} &  = 
 4\cdot \left[\Gamma _{s+1}^{\Big[\substack{s -1 \\ s -1 +1 }\Big] \times (k+1)} - \Gamma _{s+1}^{\Big[\substack{s -1  \\ s -1+1  }\Big] \times  k}\right],   \label{eq 12.9}\\
 & \nonumber \\
  \Gamma _{s -1+2 }^{\Big[\substack{s -1 \\ s-1+1  }\Big] \times (k+1)} -  \Gamma _{s-1+2}^{\Big[\substack{s-1 \\ s -1+1 }\Big] \times k} &  = 
 4\cdot \left[\Gamma _{s -1+1}^{\Big[\substack{s -1 \\ s -1  }\Big] \times (k+1)} - \Gamma _{s-1+1 }^{\Big[\substack{s -1  \\ s -1  }\Big] \times  k}\right]     \label{eq 12.10}\\
& = 4\cdot21\cdot2^{k+s-2}  = 21 \cdot2^{k+s}. \nonumber
 \end{align}\vspace{0.1 cm}\\
From \eqref{eq 12.9},  \eqref{eq 12.10} we get  \vspace{0.1 cm}\\
\begin{equation}
\label{eq 12.11}
 \Gamma _{s+2}^{\Big[\substack{s \\ s  }\Big] \times (k+1)} -  \Gamma _{s+2}^{\Big[\substack{s \\ s  }\Big] \times k}  = 
  4 \cdot\cdot 21 \cdot 2^{k+s} = 21\cdot2^{k+s +2}. 
\end{equation} \vspace{0.1 cm}\\
 \underline {The case $ 2 \leq j\leq s-1, \; k > s+j $} \vspace{0.1 cm}\\
 
 By \eqref{eq 12.2} we obtain  \vspace{0.1 cm}\\
 
\begin{align*}
 \Gamma _{s +j}^{\Big[\substack{s \\ s  }\Big] \times (k+1)} -  \Gamma _{s +j}^{\Big[\substack{s \\ s  }\Big] \times k} &  = 
 4\cdot \left[\Gamma _{s +j-1}^{\Big[\substack{s -1 \\ (s-1) +1 }\Big] \times (k+1)} - 
 \Gamma _{s -1+j}^{\Big[\substack{s -1  \\ (s-1)+1  }\Big] \times  k}\right]  && \text{if $ 2\leq j\leq s-1 $}, \\
 & \\
  \Gamma _{s -1+j}^{\Big[\substack{s -1\\ (s-1) +1}\Big] \times (k+1)} -  \Gamma _{s -1+j}^{\Big[\substack{s-1 \\ (s-1)+1 }\Big] \times k} &  = 
 4\cdot \left[\Gamma _{s -1+(j-1)}^{\Big[\substack{s -1 \\  s-1 }\Big] \times (k+1)} - 
 \Gamma _{s -1+(j-1)}^{\Big[\substack{s -1  \\ s-1  }\Big] \times  k}\right] && \text{if $ 2\leq j\leq s-1+1-1= s-1 $}.  \\
\end{align*}\\
From the above equations we get \vspace{0.1 cm}\\
\begin{equation}
\label{eq 12.12}
 \Gamma _{s +j}^{\Big[\substack{s \\ s  }\Big] \times (k+1)} -  \Gamma _{s +j}^{\Big[\substack{s \\ s  }\Big] \times k}   = 
 4^{2}\cdot \left[\Gamma _{s -1 +(j-1)}^{\Big[\substack{s -1 \\ s-1 }\Big] \times (k+1)} - 
 \Gamma _{s -1+(j-1)}^{\Big[\substack{s -1  \\ s-1  }\Big] \times  k}\right]  \quad \text{if} \;  2\leq j\leq s-1. 
\end{equation}\vspace{0.1 cm} \\
Using successively \eqref{eq 12.12} we get \vspace{0.1 cm}\\
\begin{align*}
& \Gamma _{s +j}^{\Big[\substack{s \\ s  }\Big] \times (k+1)} -  \Gamma _{s +j}^{\Big[\substack{s \\ s  }\Big] \times k}\\
& =   4^{2}\cdot \left[\Gamma _{s -1 +(j-1)}^{\Big[\substack{s -1 \\ s-1 }\Big] \times (k+1)} - 
 \Gamma _{s -1+(j-1)}^{\Big[\substack{s -1  \\ s-1  }\Big] \times  k}\right]  && \text{if $ 2\leq j\leq s-1 $} \\
 & =    4^{4}\cdot \left[\Gamma _{s -2 +(j-2)}^{\Big[\substack{s -2 \\ s-2 }\Big] \times (k+1)} - 
 \Gamma _{s -2+(j-2)}^{\Big[\substack{s -2  \\ s-2  }\Big] \times  k}\right]  && \text{if $ 2\leq j-1\leq s-1-1 $}\\
 & \vdots \\
  & = 4^{2(l-1)}\cdot \left[\Gamma _{s -(l-1) +(j- (l-1))}^{\Big[\substack{s -(l-1) \\ s-(l-1) }\Big] \times (k+1)} - 
 \Gamma _{s -(l-1) +(j-(l-1))}^{\Big[\substack{s -(l-1)  \\ s-(l-1)  }\Big] \times  k}\right]  && \text{if $ 2\leq j-(l-2)\leq s-(l-2) -1 $}\\
  & = 4^{2l}\cdot \left[\Gamma _{s -l +(j- l)}^{\Big[\substack{s -l \\ s-l }\Big] \times (k+1)} - 
 \Gamma _{s -l +(j-l)}^{\Big[\substack{s -l  \\ s-l  }\Big] \times  k}\right]  && \text{if $ 2\leq j-(l-1)\leq s- (l-1) -1 $}\\
  & \vdots \\
  & =  4^{2(j-1)}\cdot \left[\Gamma _{s -(j-1) +(j- (j-1))}^{\Big[\substack{s -(j-1) \\ s-(j-1) }\Big] \times (k+1)} - 
 \Gamma _{s -(j-1) +(j-(j-1))}^{\Big[\substack{s -(j-1)  \\ s-(j-1)  }\Big] \times  k}\right]  && \text{if $ 2\leq j-((j-1)-1)\leq s- ((j-1)-1) -1 $}.
\end{align*}\vspace{0.1 cm}\\
From the above equations we get, using \eqref{eq 12.8} with  $ s\rightarrow s-j+1  $   \vspace{0.1 cm}\\
\begin{align}
\label{eq 12.13}
\Gamma _{s +j}^{\Big[\substack{s \\ s  }\Big] \times (k+1)} -  \Gamma _{s +j}^{\Big[\substack{s \\ s  }\Big] \times k} & =
 4^{2(j-1)}\cdot \left[\Gamma _{s -j +2}^{\Big[\substack{s -j+1) \\ s -j+1) }\Big] \times (k+1)} - 
 \Gamma _{s -j+2}^{\Big[\substack{s -j+1  \\ s-j+1  }\Big] \times  k}\right] \quad \text{if} \; 2\leq j\leq s-1, \\
 & =  4^{2(j-1)}\cdot21\cdot2^{k+s-j+1-1}= 21\cdot2^{k+s+3j-4}.\nonumber
\end{align}\\

 \underline {The case  j = s, \;$ k > 2s $} \vspace{0.1 cm} \\
 
 By \eqref{eq 12.2} we obtain  \vspace{0.1 cm}\\
\begin{align*}
 \Gamma _{2s}^{\Big[\substack{s \\ s  }\Big] \times (k+1)} -  \Gamma _{2s}^{\Big[\substack{s \\ s  }\Big] \times k} &  = 
 4\cdot \left[\Gamma _{2s-1}^{\Big[\substack{s -1 \\ (s-1) +1 }\Big] \times (k+1)} - 
 \Gamma _{2s-1}^{\Big[\substack{s -1  \\ (s-1)+1  }\Big] \times  k}\right] - 3\cdot2^{2k +2s -2}, \\
    \Gamma _{(s-1) + s}^{\Big[\substack{s -1\\ (s-1) +1}\Big] \times (k+1)} -  \Gamma _{(s-1) + s}^{\Big[\substack{s-1 \\ (s-1)+1 }\Big] \times k} &  = 
 4\cdot \left[\Gamma _{2(s-1)}^{\Big[\substack{s -1 \\  s-1 }\Big] \times (k+1)} - 
 \Gamma _{2(s-1)}^{\Big[\substack{s -1  \\ s-1  }\Big] \times  k}\right] - 3\cdot 2^{2k + 2(s-1) +1 -2}. \\
\end{align*} \\
From the above equations we get \vspace{0.1 cm}\\
\begin{align*}
 \Gamma _{2s}^{\Big[\substack{s \\ s  }\Big] \times (k+1)} -  \Gamma _{2s}^{\Big[\substack{s \\ s  }\Big] \times k} &  = 
 4^{2}\cdot \left[\Gamma _{2(s-1)}^{\Big[\substack{s -1 \\ s-1 }\Big] \times (k+1)} - 
 \Gamma _{2(s-1)}^{\Big[\substack{s -1  \\ s-1  }\Big] \times  k}\right] - 9\cdot 2^{2k+2s-2} \\
  \Gamma _{2(s-1)}^{\Big[\substack{s-1 \\ s -1 }\Big] \times (k+1)} -  \Gamma _{2(s-1)}^{\Big[\substack{s-1 \\ s -1 }\Big] \times k} &  = 
 4^{2}\cdot \left[\Gamma _{2(s-2)}^{\Big[\substack{s -2 \\ s-2 }\Big] \times (k+1)} - 
 \Gamma _{2(s-2)}^{\Big[\substack{s -2  \\ s-2  }\Big] \times  k}\right] - 9\cdot 2^{2k+2(s-1)-2} \\
 & = \\
 & \vdots \\
   \Gamma _{2(s-l)}^{\Big[\substack{s -l\\ s-l  }\Big] \times (k+1)} -  \Gamma _{2(s-l)}^{\Big[\substack{s -l \\ s -l }\Big] \times k} &  = 
 4^{2}\cdot \left[\Gamma _{2(s-l-1)}^{\Big[\substack{s -l-1 \\ s-l -1}\Big] \times (k+1)} - 
 \Gamma _{2(s-1-1)}^{\Big[\substack{s -l -1 \\ s-l-1  }\Big] \times  k}\right] - 9\cdot 2^{2k+2(s-l)-2} \\
  & \vdots \\
   \Gamma _{2(s-(s-2))}^{\Big[\substack{s -(s-2)\\ s-(s-2)  }\Big] \times (k+1)} -  \Gamma _{2(s-(s-2))}^{\Big[\substack{s -(s-2) \\ s - (s-2) }\Big] \times k} &  = 
 4^{2}\cdot \left[\Gamma _{2(s-(s-2)-1)}^{\Big[\substack{s -(s-2)-1 \\ s - (s-2) -1}\Big] \times (k+1)} - 
 \Gamma _{2(s-(s-2)-1)}^{\Big[\substack{s -(s-2) -1 \\ s-(s-2)-1  }\Big] \times  k}\right] - 9\cdot 2^{2k+2(s-(s-2))-2}. 
   \end{align*}\vspace{0.1 cm}\\
We deduce from the above equations \vspace{0.1 cm}\\
\begin{align}
 \displaystyle
&  \sum_{l = 0}^{s-2} 4^{2l}\cdot\left[\Gamma _{2(s-l)}^{\Big[\substack{s - l \\ s-l  }\Big] \times (k+1)} -  \Gamma _{2(s-l)}^{\Big[\substack{s -l \\ s -l }\Big] \times k} \right] \label{eq 12.14} \\
&  = \sum_{l = 0}^{s-2}\left( 4^{2(l+1)}\cdot\left[\Gamma _{2(s-l-1)}^{\Big[\substack{s -l-1 \\ s- l -1}\Big] \times (k+1)} - 
 \Gamma _{2(s-l-1)}^{\Big[\substack{s -l -1 \\ s-l-1  }\Big] \times  k}\right] - 9\cdot 2^{2k+2(s-l)-2}\cdot2^{4l}\right)  \nonumber \\
 & = \sum_{l = 1}^{s-1} 4^{2l}\cdot\left[\Gamma _{2(s-l)}^{\Big[\substack{s - l \\ s- l }\Big] \times (k+1)} - 
 \Gamma _{2(s-l)}^{\Big[\substack{s -l  \\ s-l  }\Big] \times  k}\right] - \sum_{l = 1}^{s-1} 9\cdot 2^{2k+2s + 2l-4}. \nonumber
\end{align}\vspace{0.1 cm}\\
We get from \eqref{eq 12.14} after some simplications \vspace{0.1 cm}\\
\begin{align*}
 \Gamma _{2s}^{\Big[\substack{s \\ s  }\Big] \times (k+1)} -  \Gamma _{2s}^{\Big[\substack{s \\ s  }\Big] \times k} & = 
 4^{2(s-1)}\left[  \Gamma _{2}^{\Big[\substack{1 \\ 1 }\Big] \times (k+1)} -  \Gamma _{2}^{\Big[\substack{1 \\ 1 }\Big] \times k}   \right] -
 9\cdot2^{2k+2s-4}\cdot\sum_{l = 1}^{s-1}2^{2l} \\
 & = 4^{2(s-1)}\left[(2^{2(k+1)} -3\cdot2^{k+1} +2) - (2^{2k} -3\cdot2^{k} +2)  \right] - 3\cdot2^{2k+2s-2}\cdot(2^{2s-2}-1) \\
 & = 3\cdot2^{2k+2s-2}- 3\cdot2^{k+4s -4}.  
\end{align*}

\end{proof}

\begin{lem}
\label{lem 12.2} We have
\begin{equation}
\label{eq 12.15}
 \Gamma _{s+j}^{\Big[\substack{s \\ s+1}\Big] \times (k+1)}
  - \Gamma _{s+j}^{\Big[\substack{s \\ s+1}\Big] \times k}=
  \begin{cases}
  2^{k+s-1}  & \text{if  } j = 0,\; k >s, \\
  11\cdot2^{k+s-1}& \text{if  } j = 1,\; k > s +1, \\
21 \cdot 2^{k+s+3j-5}   &  \text{if  }    2\leq j\leq s ,\; k>s+j, \\
3\cdot2^{2k+2s-1}- 3\cdot2^{k+4s -2}       & \text{if   } j  = s+1 ,  \; k>2s +1.
\end{cases}
\end{equation}
\end{lem}
\begin{proof}We proceed as in the proof of Lemma \ref{lem 12.1}.\vspace{0.1 cm}\\

\underline {The case  j  = 0, $ k > s $}\vspace{0.1 cm}\\

Applying formula  \eqref{eq 12.2} we obtain using \eqref{eq 11.21} \vspace{0.1 cm}\\
\begin{equation}
 \Gamma _{s}^{\Big[\substack{s \\ s +1 }\Big] \times (k+1)} -  \Gamma _{s}^{\Big[\substack{s \\ s +1 }\Big] \times k}  = 
 4\cdot \left[\Gamma _{s-1}^{\Big[\substack{s  \\ s }\Big] \times (k+1)} - \Gamma _{s-1}^{\Big[\substack{s   \\ s   }\Big] \times  k}\right]  + 2^{k+s-1} =  2^{k+s-1}. \label{eq 12.16}\\
\end{equation}
\underline {The case  j  = 1, $ k > s +1  $}\vspace{0.1 cm}\\

From  \eqref{eq 12.2} we obtain using \eqref{eq 12.5} \vspace{0.1 cm}\\
 \begin{align}
 \Gamma _{s+1}^{\Big[\substack{s \\ s +1 }\Big] \times (k+1)} -  \Gamma _{s+1}^{\Big[\substack{s \\ s +1 }\Big] \times k} &  = 
 4\cdot \left[\Gamma _{s}^{\Big[\substack{s  \\ s }\Big] \times (k+1)} - \Gamma _{s}^{\Big[\substack{s   \\ s   }\Big] \times  k}\right]  - 2^{k+s-1}  \label{eq 12.17}\\
 & = 4\cdot(3\cdot 2^{k+s-1}) - 2^{k+s-1} = 11\cdot 2^{k+s-1}. \nonumber
\end{align}
\underline {The case  $ 2\leq j\leq s,\;  k > s + j  $}\vspace{0.1 cm}\\

From  \eqref{eq 12.2} we obtain using \eqref{eq 12.1} with $ j \rightarrow j-1 $ (observing that $1 \leq j-1 \leq s-1 $)\vspace{0.1 cm}\\
 \begin{align}
 \Gamma _{s+j}^{\Big[\substack{s \\ s +1 }\Big] \times (k+1)} -  \Gamma _{s+j}^{\Big[\substack{s \\ s +1 }\Big] \times k} &  = 
 4\cdot \left[\Gamma _{s+j-1}^{\Big[\substack{s  \\ s }\Big] \times (k+1)} - \Gamma _{s+j-1}^{\Big[\substack{s   \\ s   }\Big] \times  k}\right]    \label{eq 12.18}\\
 & = 4\cdot(21\cdot 2^{k+s +3(j-1) -4})  = 21 \cdot 2^{k+s +3j-5}. \nonumber
\end{align}
\underline {The case   j = s+1 }\vspace{0.1 cm}\\
 \begin{align}
 \Gamma _{2s +1}^{\Big[\substack{s \\ s +1 }\Big] \times (k+1)} -  \Gamma _{2s+1}^{\Big[\substack{s \\ s +1 }\Big] \times k} &  = 
 4\cdot \left[\Gamma _{2s}^{\Big[\substack{s  \\ s }\Big] \times (k+1)} - \Gamma _{2s}^{\Big[\substack{s   \\ s   }\Big] \times  k}\right]  - 3\cdot2^{2k+2s-1}  \label{eq 12.19}\\
 & = 4\cdot[3\cdot2^{2k+2s-2} -3\cdot2^{k+4s-4}] -  3\cdot2^{2k+2s-1}  \nonumber\\
 & = 3\cdot2^{2k+2s-1} -3\cdot2^{k+4s-2}.   \nonumber
\end{align}
\end{proof}

\begin{lem}
\label{lem 12.3} We have in the case $ m\geq 2 $\vspace{0.1 cm}\\
 
\begin{equation}
\label{eq 12.20}
 \Gamma _{s+j}^{\Big[\substack{s \\ s+m}\Big] \times (k+1)}
  - \Gamma _{s+j}^{\Big[\substack{s \\ s+m}\Big] \times k}=
  \begin{cases}
  2^{k+s-1}  & \text{if  } j = 0,\; k >s, \\
  3\cdot2^{k+s +2j -3}& \text{if  } 1\leq j\leq m-1,\; k > s+j,  \\
11\cdot2^{k+ s + 2m -3}              &  \text{if  }  j = m, \; k > s+m,
\end{cases}
\end{equation}

\begin{equation}
\label{eq 12.21}
 \Gamma _{s+m+j}^{\Big[\substack{s \\ s+m}\Big] \times (k+1)}
  - \Gamma _{s+m+j}^{\Big[\substack{s \\ s+m}\Big] \times k}=
  \begin{cases}
  21\cdot2^{k+s +2m +3j -4}& \text{if  } 1\leq j\leq s-1,\; k > s+m+ j,  \\
3\cdot2^{2k+2s+m-2} -  3\cdot2^{k+4s+2m-4}            &  \text{if  }  j = s, \; k > 2s+m.
\end{cases}
\end{equation}
\end{lem}

\begin{proof}
 \underline {The case s+j with  $ 0\leq j\leq m,\;  k > s+ j $}\vspace{0.01 cm}\\
 
  \underline {j = 0}\vspace{0.01 cm}\\
  
Applying formula  \eqref{eq 12.2} we obtain using \eqref{eq 11.21}, observing that $ m -1 \geq 1\geq 0  $ \vspace{0.01 cm} \\
\begin{equation}
 \Gamma _{s}^{\Big[\substack{s \\ s +m }\Big] \times (k+1)} -  \Gamma _{s}^{\Big[\substack{s \\ s +m }\Big] \times k}  = 
 4\cdot \left[\Gamma _{s-1}^{\Big[\substack{s  \\ s +(m-1)}\Big] \times (k+1)} - \Gamma _{s-1}^{\Big[\substack{s   \\ s +(m-1)  }\Big] \times  k}\right]  + 2^{k+s-1} =  2^{k+s-1}. \label{eq 12.22}\\
\end{equation}\vspace{0.1 cm}\\
 \underline {j = 1} \vspace{0.01 cm}\\
 
 Applying formula  \eqref{eq 12.2} we obtain using \eqref{eq 11.21}, observing that $ m -2 \geq 0  $ \vspace{0.01 cm} \\
 \begin{align}
& \Gamma _{s  +1}^{\Big[\substack{s \\ s +m }\Big] \times (k+1)} -  \Gamma _{s +1}^{\Big[\substack{s \\ s +m }\Big] \times k} \label{eq 12.23} \\
& =  4\cdot \left[\Gamma _{s}^{\Big[\substack{s  \\ s +(m-1)}\Big] \times (k+1)} - \Gamma _{s}^{\Big[\substack{s   \\ s +(m-1)  }\Big] \times  k}\right]  - 2^{k+s-1} \nonumber \\
& = 4\cdot \left( 4\cdot \left[\Gamma _{s-1}^{\Big[\substack{s  \\ s +(m-2)}\Big] \times (k+1)} - \Gamma _{s -1}^{\Big[\substack{s   \\ s +(m-2)  }\Big] \times  k}\right]  + 2^{k+s-1} \right ) - 2^{k+s-1}\nonumber \\
& = 3\cdot2^{k+s-1}. \nonumber
\end{align}
 \underline {$2\leq j\leq m-1 $} \vspace{0.01 cm}\\
 
Assume $ 0\leq l\leq j-2,$  then by  \eqref{eq 12.2} with $ j\rightarrow j-l,\; m\rightarrow m-l,$ we have
 
 \begin{equation}
 \label{eq 12.24}
 \Gamma _{s+j-l}^{\Big[\substack{s \\ s +(m-l) }\Big] \times (k+1)} -  \Gamma _{s+j-l}^{\Big[\substack{s \\ s +(m-l) }\Big] \times k} 
  =  4\cdot \left[\Gamma _{s+j- (l+1)}^{\Big[\substack{s  \\ s +(m-(l+1)) }\Big] \times (k+1)} - \Gamma _{s+j-(l+1)}^{\Big[\substack{s   \\ s +(m-(l+1))  }\Big] \times  k}\right] 
 \quad \text{for} \; l = 0,1,2,\ldots,j-2.
\end{equation}

From \eqref{eq 12.24} we deduce \vspace{0.01 cm}\\

\begin{align}
 \prod_{l=0}^{j-2}\left[ \Gamma _{s+j-l}^{\Big[\substack{s \\ s +(m-l) }\Big] \times (k+1)} -  \Gamma _{s+j-l}^{\Big[\substack{s \\ s +(m-l) }\Big] \times k}\right ]
&  = \prod_{l=0}^{j-2} 4\cdot \left[\Gamma _{s+j- (l+1)}^{\Big[\substack{s  \\ s +(m-(l+1)) }\Big] \times (k+1)} - \Gamma _{s+j-(l+1)}^{\Big[\substack{s   \\ s +(m-(l+1))  }\Big] \times  k}\right]  \label{eq 12.25} \\
& = \prod_{l=1}^{j-1}4\cdot \left[ \Gamma _{s+j-l}^{\Big[\substack{s \\ s +(m-l) }\Big] \times (k+1)} -  \Gamma _{s+j-l}^{\Big[\substack{s \\ s +(m-l) }\Big] \times k}\right ]. \nonumber
 \end{align} \vspace{0.01 cm}\\
 We then obtain by  \eqref{eq 12.25},  \eqref{eq 12.23}, observing that $ m-(j-1)\geq  2 $ \vspace{0.01 cm}\\
  \begin{align*}
  \Gamma _{s+j}^{\Big[\substack{s \\ s +m  }\Big] \times (k+1)} -  \Gamma _{s+j}^{\Big[\substack{s \\ s + m}\Big] \times k} 
&  =  4^{j-1}\cdot \left[\Gamma _{s+1}^{\Big[\substack{s  \\ s +(m-(j-1)) }\Big] \times (k+1)} - \Gamma _{s+1}^{\Big[\substack{s   \\ s +(m-(j-1))  }\Big] \times  k}\right] \\
& = 2^{2j-2}\cdot3\cdot2^{k+s-1}= 3\cdot2^{k+s+2j-3}.
\end{align*}\\

 \underline {j = m } \vspace{0.01 cm}\\
  
 Assume $2\leq l\leq m \;(\leq s+m-1) $  then by  \eqref{eq 12.2} with $ j \rightarrow l,\; m\rightarrow  l,$ we get \\
 
  \begin{equation}
 \label{eq 12.26}
 \Gamma _{s+l}^{\Big[\substack{s \\ s + l }\Big] \times (k+1)} -  \Gamma _{s+l}^{\Big[\substack{s \\ s + l }\Big] \times k} 
  =  4\cdot \left[\Gamma _{s+(l-1)}^{\Big[\substack{s  \\ s + (l-1) }\Big] \times (k+1)} - \Gamma _{s+ (l-1)}^{\Big[\substack{s   \\ s + (l-1)  }\Big] \times  k}\right] 
 \quad for \; l = 2,3, \ldots, m.
\end{equation}
 From \eqref{eq 12.26} we deduce \vspace{0.01 cm}\\

\begin{align}
 \prod_{l = 2}^{m}\left[ \Gamma _{s+l}^{\Big[\substack{s \\ s + l }\Big] \times (k+1)} -  \Gamma _{s+l}^{\Big[\substack{s \\ s + l }\Big] \times k}\right ]
&  = \prod_{l=2}^{m} 4\cdot \left[\Gamma _{s+ (l-1)}^{\Big[\substack{s  \\ s + (l-1) }\Big] \times (k+1)} - \Gamma _{s+(l-1)}^{\Big[\substack{s   \\ s + (l-1)  }\Big] \times  k}\right]  \label{eq 12.27} \\
& = \prod_{l=1}^{m-1}4\cdot \left[ \Gamma _{s+l}^{\Big[\substack{s \\ s + l }\Big] \times (k+1)} -  \Gamma _{s+ l}^{\Big[\substack{s \\ s + l}\Big] \times k}\right ]. \nonumber
 \end{align} \vspace{0.01 cm}\\
 We then obtain by  \eqref{eq 12.27} and  \eqref{eq 12.15} with j = 1 \vspace{0.01 cm}\\
  \begin{align*}
  \Gamma _{s+m}^{\Big[\substack{s \\ s +m  }\Big] \times (k+1)} -  \Gamma _{s+m}^{\Big[\substack{s \\ s + m}\Big] \times k} 
&  =  4^{m-1}\cdot \left[\Gamma _{s+1}^{\Big[\substack{s  \\ s + 1 }\Big] \times (k+1)} - \Gamma _{s+1}^{\Big[\substack{s   \\ s + 1  }\Big] \times  k}\right] \\
& = 2^{2m-2}\cdot11\cdot2^{k+s-1}= 11\cdot2^{k+s +2m-3}.
\end{align*}\\

  \underline { The case s+m+j with  $1 \leq j\leq s -1 $} \vspace{0.01 cm}\\
  
  Assume $ 0\leq l\leq m -1,  $  then by  \eqref{eq 12.2} with $ j \rightarrow  (m-l) +j ,\; m\rightarrow m- l,$\\
   (observing that $ 2\leq m-l+j\leq s+(m-l)-1 $) we have \vspace{0.01 cm}\\
  \begin{equation}
 \label{eq 12.28}
 \Gamma _{s+(m-l)+j}^{\Big[\substack{s \\ s + (m-l) }\Big] \times (k+1)} -  \Gamma _{s+l}^{\Big[\substack{s \\ s + (m-l) }\Big] \times k} 
  =  4\cdot \left[\Gamma _{s+ (m-(l +1))+j}^{\Big[\substack{s  \\ s + (m-(l+1)) }\Big] \times (k+1)} - \Gamma _{s+ (m-(l+1))+j}^{\Big[\substack{s   \\ s + (m-(l+1))  }\Big] \times  k}\right] 
 \quad \text{for} \; l = 0,1, \ldots, m-1.
\end{equation}
 From \eqref{eq 12.28} we deduce \vspace{0.01 cm}\\

\begin{align}
 \prod_{l = 0}^{m-1}\left[ \Gamma _{s+(m-l)+j}^{\Big[\substack{s \\ s + (m-l) }\Big] \times (k+1)} -  \Gamma _{s+(m-l)+j}^{\Big[\substack{s \\ s + (m-l) }\Big] \times k}\right ]
&  = \prod_{l=0}^{m-1} 4\cdot \left[\Gamma _{s+ (m-(l+1))+j }^{\Big[\substack{s  \\ s + (m-(l+1)) }\Big] \times (k+1)} - \Gamma _{s+(m-(l+1))+j}^{\Big[\substack{s   \\ s + (m-(l+1))  }\Big] \times  k}\right]  \label{eq 12.29} \\
& = \prod_{l=1}^{m}4\cdot \left[ \Gamma _{s+(m-l)+j}^{\Big[\substack{s \\ s + (m-l) }\Big] \times (k+1)} -  \Gamma _{s+ (m-l)}^{\Big[\substack{s \\ s + (m-l)}\Big] \times k}\right ]. \nonumber
 \end{align} \vspace{0.01 cm}\\
 We then obtain by  \eqref{eq 12.29} and  \eqref{eq 12.1}  \vspace{0.01 cm}\\
  \begin{align*}
  \Gamma _{s+m+j}^{\Big[\substack{s \\ s +m  }\Big] \times (k+1)} -  \Gamma _{s+m+j}^{\Big[\substack{s \\ s + m}\Big] \times k} 
&  =  4^{m}\cdot \left[\Gamma _{s+j}^{\Big[\substack{s  \\ s  }\Big] \times (k+1)} - \Gamma _{s+j}^{\Big[\substack{s   \\ s   }\Big] \times  k}\right] \\
& = 2^{2m}\cdot21\cdot2^{k+s +3j-4}= 21\cdot2^{k+s +2m +3j -4}.
\end{align*} \vspace{0.01 cm}\\
 \underline { The case s+m+j with  j = s} \vspace{0.01 cm}\\
 
   Assume $ 0\leq l\leq m -1,  $  then by  \eqref{eq 12.2} with $ j \rightarrow s+  (m-l)  ,\; m\rightarrow m- l,$ we have  \vspace{0.01 cm}\\
  \begin{align}
 \label{eq 12.30}
&  \Gamma _{2s+(m-l)}^{\Big[\substack{s \\ s + (m-l) }\Big] \times (k+1)} -  \Gamma _{2s+(m-l)}^{\Big[\substack{s \\ s + (m-l) }\Big] \times k} \\
&  =  4\cdot \left[\Gamma _{2s+ (m-(l +1))}^{\Big[\substack{s  \\ s + (m-(l+1)) }\Big] \times (k+1)} - \Gamma _{2s+ (m-(l+1))}^{\Big[\substack{s   \\ s + (m-(l+1))  }\Big] \times  k}\right] 
  -3\cdot2^{2k+2s+(m-l)-2} \quad \text{for} \; l = 0,1, \ldots, m-1. \nonumber
\end{align}
 From \eqref{eq 12.30} we deduce \vspace{0.01 cm}\\
\begin{align}
&  \sum_{l = 0}^{m-1}4^{l}\cdot\left[ \Gamma _{2s+(m-l)}^{\Big[\substack{s \\ s + (m-l) }\Big] \times (k+1)} -  \Gamma _{2s+(m-l)}^{\Big[\substack{s \\ s + (m-l) }\Big] \times k}\right ] \label{eq 12.31}\\
&  = \sum_{l=0}^{m-1}4^{l}\cdot\left( 4\cdot \left[\Gamma _{2s+ (m-(l+1)) }^{\Big[\substack{s  \\ s + (m-(l+1)) }\Big] \times (k+1)} - \Gamma _{2s+(m-(l+1))}^{\Big[\substack{s   \\ s + (m-(l+1))  }\Big] \times  k}\right] - 3\cdot2^{2k+2s+(m-l)-2}\right) \nonumber \\
& = \sum_{l=1}^{m}4^{l}\cdot \left[ \Gamma _{2s+(m-l)}^{\Big[\substack{s \\ s + (m-l) }\Big] \times (k+1)} -  \Gamma _{2s+ (m-l)}^{\Big[\substack{s \\ s + (m-l)}\Big] \times k}\right ]- \sum_{l=1}^{m}3\cdot2^{2k+2s+m+l-3}. \nonumber
 \end{align} \vspace{0.01 cm}\\
 We then obtain by  \eqref{eq 12.31} and  \eqref{eq 12.1} after some simplifications  \vspace{0.01 cm}\\
  \begin{align*}
  \Gamma _{2s+m}^{\Big[\substack{s \\ s +m  }\Big] \times (k+1)} -  \Gamma _{2s+m}^{\Big[\substack{s \\ s + m}\Big] \times k} 
&  =  4^{m}\cdot \left[\Gamma _{2s}^{\Big[\substack{s  \\ s  }\Big] \times (k+1)} - \Gamma_{2s}^{\Big[\substack{s   \\ s   }\Big] \times  k}\right] \\
& = 2^{2m}\cdot(3\cdot2^{2k+2s-2} - 3\cdot2^{k+4s-4}) - (3\cdot2^{2k+2s+2m-2} - 3\cdot2^{2k+2s+m-2})\\
& =  3\cdot2^{2k+2s+m-2} - 3\cdot2^{k+4s+2m-4}.
\end{align*} \vspace{0.01 cm}\\
\end{proof}
\begin{lem} 
 \label{lem 12.4}We have in the case $ m \in \left\{0,1\right\} $ \vspace{0.01 cm}\\
 \begin{equation}
\label{eq 12.32}
 \Gamma _{s+j}^{\Big[\substack{s \\ s}\Big] \times (k+1)}
  - \Gamma _{s+j}^{\Big[\substack{s \\ s}\Big] \times k}=
  \begin{cases}
 3\cdot 2^{k+s-1}  & \text{if  } j = 0,\; k >s, \\
 21\cdot 2^{k+s+3j-4}   &  \text{if  }    1\leq j\leq s-1,\; k>s+j, \\
3\cdot2^{2k+2s-2}- 3\cdot2^{k+4s - 4}       & \text{if   } j  = s,  \; k>2s,
\end{cases}
\end{equation}
 \begin{equation}
\label{eq 12.33}
 \Gamma _{s+j}^{\Big[\substack{s \\ s+1}\Big] \times (k+1)}
  - \Gamma _{s+j}^{\Big[\substack{s \\ s+1}\Big] \times k}=
  \begin{cases}
  2^{k+s-1}  & \text{if  } j = 0,\; k >s, \\
  11\cdot2^{k+s-1}& \text{if  } j = 1,\; k > s +1, \\
21 \cdot 2^{k+s+3j-5}   &  \text{if  }    2\leq j\leq s ,\; k>s+j, \\
3\cdot2^{2k+2s-1}- 3\cdot2^{k+4s -2}       & \text{if   } j  = s+1 ,  \; k>2s +1.
\end{cases}
\end{equation}
 In the case $ m\geq 2,  $ we get \vspace{0.1 cm}\\
 \begin{equation}
\label{eq 12.34}
 \Gamma _{s+j}^{\Big[\substack{s \\ s+m}\Big] \times (k+1)}
  - \Gamma _{s+j}^{\Big[\substack{s \\ s+m}\Big] \times k}=
  \begin{cases}
  2^{k+s-1}  & \text{if  } j = 0,\; k >s, \\
  3\cdot2^{k+s +2j -3}& \text{if  } 1\leq j\leq m-1,\; k > s+j,  \\
11\cdot2^{k+ s + 2m -3}              &  \text{if  }  j = m, \; k > s+m,
\end{cases}
\end{equation}
\begin{equation}
\label{eq 12.35}
 \Gamma _{s+m+j}^{\Big[\substack{s \\ s+m}\Big] \times (k+1)}
  - \Gamma _{s+m+j}^{\Big[\substack{s \\ s+m}\Big] \times k}=
  \begin{cases}
  21\cdot2^{k+s +2m +3j -4}& \text{if  } 1\leq j\leq s-1,\; k > s+m+ j,  \\
3\cdot2^{2k+2s+m-2} -  3\cdot2^{k+4s+2m-4}            &  \text{if  }  j = s, \; k > 2s+m.
\end{cases}
\end{equation}
 \end{lem}
 \begin{proof}
The assertions follow  from Lemmas \ref{lem 12.1}, \ref{lem 12.2}  and \ref{lem 12.3}. \\

  \end{proof}

 \section{\textbf{ COMPUTATION OF  $\Gamma_{s+j}^{\left[s\atop s+m\right]\times k}
\;for\; j\in\left\{0,1\right\},\; k \geq s+j,\;m\geq 0 $}}
\label{sec 13}

From the equations \eqref{eq 11.1} and \eqref{eq 11.2} with k = s+1 we deduce $\Gamma_{s}^{\left[s\atop s+m\right]\times (s+1)}$  and $ \Gamma_{s+1}^{\left[s\atop s+m\right]\times (s+1)} $
 since by \eqref{eq 11.21} the terms  $\Gamma_{i}^{\left[s\atop s+m\right]\times (s+1)}$ are known for $ i\leq s-1. $
Then from the recurrent formula \eqref{eq 12.34} with j = 0 we compute  $\Gamma_{s}^{\left[s\atop s+m\right]\times k}$ for $ k\geq s+2. $
 The other results are obtained in a similar way.\\
 
\begin{lem}
\label{lem 13.1}We have in the case m = 0 
\begin{align}
 \Gamma _{s}^{\Big[\substack{s \\ s  }\Big] \times s} & = 2^{4s-2} - 3\cdot2^{3s-4} + 2^{2s-3}, \label{eq 13.1} \\
 & \nonumber \\
  \Gamma _{s}^{\Big[\substack{s \\ s  }\Big] \times k} & = 3\cdot 2^{k+s-1}+ 21\cdot2^{3s-4} -27\cdot 2^{2s-3}\quad \text{if $ k > s $}, \label{eq 13.2} \\
  & \nonumber \\
   \Gamma _{s+1}^{\Big[\substack{s \\ s  }\Big] \times (s+1)} & = 2^{4s} - 3\cdot2^{3s- 1} + 2^{2s- 1}, \label{eq 13.3} \\
   & \nonumber \\
    \Gamma _{s+1}^{\Big[\substack{s \\ s  }\Big] \times k} & = 21\cdot[2^{k+s-1} + 2^{3s-1} - 5\cdot2^{2s-1}]\quad \text{if $ k> s+1 $}.\label{eq 13.4}\\
    & \nonumber 
\end{align}
\end{lem}
\begin{proof}
We get \eqref{eq 13.1} from \eqref{eq 11.26} with k = s, m = 0.\vspace{0.01 cm}\\

To prove \eqref{eq 13.2} and  \eqref{eq 13.3}  we proceed as follows :  \vspace{0.01 cm}\\
Consider the matrix

  $$   \left ( \begin{array} {cccccc}
\alpha _{1} & \alpha _{2} & \alpha _{3} &  \ldots & \alpha _{s}  &  \alpha _{s+1} \\
\alpha _{2 } & \alpha _{3} & \alpha _{4}&  \ldots  &  \alpha _{s+1} &  \alpha _{s+2} \\
\vdots & \vdots & \vdots    &  \vdots & \vdots  &  \vdots \\
\alpha _{s-1} & \alpha _{s} & \alpha _{s +1} & \ldots  &  \alpha _{2s-2} &  \alpha _{2s-1}  \\
\alpha _{s} & \alpha _{s+1} & \alpha _{s +2} & \ldots  &  \alpha _{2s-1} &  \alpha _{2s} \\
 \beta  _{1} & \beta  _{2} & \beta  _{3} & \ldots  &  \beta_{s} &  \beta _{s+1}  \\
\beta  _{2} & \beta  _{3} & \beta  _{4} & \ldots  &  \beta_{s+1} &  \beta _{s+2}  \\
\vdots & \vdots & \vdots    &  \vdots & \vdots  &  \vdots \\
\beta  _{s-1} & \beta  _{s} & \beta  _{s+1} & \ldots  &  \beta_{2s-2} &  \beta _{2s-1}  \\
\beta  _{s} & \beta  _{s+1} & \beta  _{s+2} & \ldots  &  \beta_{2s-1} &  \beta _{2s}
\end{array}  \right). $$ \vspace{0.01 cm}\\

We have respectively  by \eqref{eq 11.1}, \eqref{eq 11.2}, \eqref{eq 11.27} and \eqref{eq 11.28}  with m = 0, k = s+1  \vspace{0.01 cm}\\
  \begin{align}
  \sum_{ i = 0}^{s+1} \Gamma _{i}^{\Big[\substack{s  \\ s }\Big] \times (s+1)} & = 2^{4s},  \label{eq 13.5}\\
  & \nonumber \\
   \sum_{i = 0}^{s+1} \Gamma _{i}^{\Big[\substack{s  \\ s }\Big] \times (s+1)}\cdot2^{-i} & =  2^{3s -1} + 2^{2s} - 2^{s-1},  \label{eq 13.6}\\
   & \nonumber \\
    \sum_{i = 0}^{s -1} \Gamma _{i}^{\Big[\substack{s  \\ s }\Big] \times (s+1)} & = 3\cdot2^{3s-4} - 2^{2s-3} \label{eq 13.7}\\
    &  \text{and}\nonumber \\
   \sum_{i = 0}^{s -1} \Gamma _{i}^{\Big[\substack{s  \\ s }\Big] \times (s+1)}\cdot2^{-i} & =  7\cdot2^{2s-4} - 3\cdot2^{s-3}.  \label{eq 13.8}
  \end{align} \vspace{0.01 cm}\\  
  From \eqref{eq 13.5}, \eqref{eq 13.6},  \eqref{eq 13.7} and \eqref{eq 13.8} we deduce after some calculations  \vspace{0.01 cm}\\
  \begin{align}
\Gamma _{s}^{\Big[\substack{s  \\ s }\Big] \times (s+1)} &  =  21\cdot2^{3s-4} - 3\cdot2^{2s-3},  \label{eq 13.9} \\
& \nonumber \\
\Gamma _{s+1}^{\Big[\substack{s  \\ s }\Big] \times (s+1)} & = 2^{4s}-3\cdot2^{3s-1} +2^{2s-1}. \label{eq 13.10}\\
& \nonumber
\end{align}
By \eqref{eq 12.32} with j = 0 we get \vspace{0.01 cm}\\
\begin{equation}
\label{eq 13.11}
\Gamma _{s}^{\Big[\substack{s  \\ s }\Big] \times (k+1)} - \Gamma _{s}^{\Big[\substack{s  \\ s }\Big] \times k }= 3\cdot2^{k+s-1}\quad \text{if}\quad k > s.
\end{equation}

From \eqref{eq 13.11}, \eqref{eq 13.8} we deduce \vspace{0.01 cm}\\
\begin{align}
& \sum_{j = s+1}^{k}\left( \Gamma _{s}^{\Big[\substack{s  \\ s }\Big] \times (j+1)} - \Gamma _{s}^{\Big[\substack{s  \\ s }\Big] \times j}   \right)  = \sum_{j = s+1}^{k}  3 \cdot 2^{j+s-1} \nonumber \\
& \Leftrightarrow  \sum_{j = s+2}^{k+1}\Gamma _{s}^{\Big[\substack{s  \\ s }\Big] \times j} - \sum_{j = s+1}^{k} \Gamma _{s}^{\Big[\substack{s  \\ s }\Big] \times j} = 3\cdot2^{k+s}- 3\cdot2^{2s} \nonumber \\
& \Leftrightarrow \Gamma _{s}^{\Big[\substack{s  \\ s }\Big] \times (k+1)} - \Gamma _{s}^{\Big[\substack{s  \\ s }\Big] \times (s+1)} = 3\cdot2^{k+s}- 3\cdot2^{2s} \nonumber \\
& \Leftrightarrow \Gamma _{s}^{\Big[\substack{s  \\ s }\Big] \times (k+1)}= 3\cdot2^{k+s}- 3\cdot2^{2s} +  21\cdot2^{3s-4} - 3\cdot2^{2s-3} \nonumber \\
& \Leftrightarrow \Gamma _{s}^{\Big[\substack{s  \\ s }\Big] \times (k+1)}= 3\cdot2^{k+s} +  21\cdot2^{3s-4} - 27\cdot2^{2s-3}  \quad \text{ if $ k > s $}. \label{eq 13.12}\\
& \nonumber
\end{align}
By \eqref{eq 13.8} we see that \eqref{eq 13.12} holds for $ k = s, $ then  \eqref{eq 13.12} holds for $ k \geq s.$\vspace{0.01 cm}\\

To prove \eqref{eq 13.4} we proceed as follows : \vspace{0.01 cm}\\
Consider the matrix\\

  $$   \left ( \begin{array} {cccccc}
\alpha _{1} & \alpha _{2} & \alpha _{3} &  \ldots & \alpha _{s+1}  &  \alpha _{s+2} \\
\alpha _{2 } & \alpha _{3} & \alpha _{4}&  \ldots  &  \alpha _{s+2} &  \alpha _{s+3} \\
\vdots & \vdots & \vdots    &  \vdots & \vdots  &  \vdots \\
\alpha _{s-1} & \alpha _{s} & \alpha _{s +1} & \ldots  &  \alpha _{2s-1} &  \alpha _{2s}  \\
\alpha _{s} & \alpha _{s+1} & \alpha _{s +2} & \ldots  &  \alpha _{2s} &  \alpha _{2s+1} \\
 \beta  _{1} & \beta  _{2} & \beta  _{3} & \ldots  &  \beta_{s+1} &  \beta _{s+2}  \\
\beta  _{2} & \beta  _{3} & \beta  _{4} & \ldots  &  \beta_{s+2} &  \beta _{s+3}  \\
\vdots & \vdots & \vdots    &  \vdots & \vdots  &  \vdots \\
\beta  _{s-1} & \beta  _{s} & \beta  _{s+1} & \ldots  &  \beta_{2s-1} &  \beta _{2s}  \\
\beta  _{s} & \beta  _{s+1} & \beta  _{s+2} & \ldots  &  \beta_{2s} &  \beta _{2s +1}
\end{array}  \right). $$ \vspace{0.01 cm}\\

We have respectively  by \eqref{eq 11.1}, \eqref{eq 11.2}, \eqref{eq 11.27} and \eqref{eq 11.28} with m = 0, k = s+2,   \vspace{0.01 cm}\\
  \begin{align}
  \sum_{ i = 0}^{s+2} \Gamma _{i}^{\Big[\substack{s  \\ s }\Big] \times (s+2)} & = 2^{4s+2},  \label{eq 13.13}\\
  & \nonumber \\
   \sum_{i = 0}^{s+2} \Gamma _{i}^{\Big[\substack{s  \\ s }\Big] \times (s+2)}\cdot2^{-i} & =  2^{3s } + 2^{2s+2} - 2^{s},  \label{eq 13.14} \\
   & \nonumber \\
    \sum_{i = 0}^{s -1} \Gamma _{i}^{\Big[\substack{s  \\ s }\Big] \times (s+1)} & = 3\cdot2^{3s-4} - 2^{2s-3} \label{eq 13.15} \\
      &  \text{and}\nonumber \\
   \sum_{i = 0}^{s -1} \Gamma _{i}^{\Big[\substack{s  \\ s }\Big] \times (s+2)}\cdot2^{-i} & =  7\cdot2^{2s-4} - 3\cdot2^{s-3}.  \label{eq 13.16}
  \end{align} \vspace{0.01 cm}\\
  From \eqref{eq 13.13}, \eqref{eq 13.14}, \eqref{eq 13.15}, \eqref{eq 13.16} and \eqref{eq 13.2} with k = s+2  we obtain  \vspace{0.01 cm}\\
  \begin{align}
& \Gamma _{s+1}^{\Big[\substack{s  \\ s }\Big] \times (s+2)} + \Gamma _{s+2}^{\Big[\substack{s  \\ s }\Big] \times (s+2)} \nonumber \\
& = 2^{4s+2} - (3\cdot2^{3s-4} - 2^{2s-3}) - (3\cdot2^{s+2 +s-1} + 21\cdot2^{3s-4} -27\cdot2^{2s-3}), \label{eq 13.17}\\
& \nonumber \\
& \Gamma _{s+1}^{\Big[\substack{s  \\ s }\Big] \times (s+2)}\cdot2^{-(s+1)} + \Gamma _{s+2}^{\Big[\substack{s  \\ s }\Big] \times (s+2)}\cdot 2^{-(s+2)} \nonumber \\
& =   2^{3s } + 2^{2s+2} - 2^{s} -(7\cdot2^{2s-4} - 3\cdot2^{s-3} ) -  (3\cdot2^{s+2 +s-1} + 21\cdot2^{3s-4} -27\cdot2^{2s-3})\cdot2^{-s}. \label{eq 13.18} 
\end{align}\vspace{0.01 cm}\\
Hence by \eqref{eq 13.17}, \eqref{eq 13.18} we deduce after some calculations \vspace{0.01 cm}\\
\begin{align}
\Gamma _{s+1}^{\Big[\substack{s  \\ s }\Big] \times (s+2)}& = 21\cdot(2^{3s-1}- 2^{2s-1}), \label{eq 13.19}\\
& \nonumber \\
\Gamma _{s+2}^{\Big[\substack{s  \\ s }\Big] \times (s+2)}& = 2^{4s+2} -3\cdot2^{3s+2}+ 2^{2s+3}. \label{eq 13.20}\\
& \nonumber
\end{align}
By \eqref{eq 12.32}, with j = 1, we get \vspace{0.01 cm}\\
\begin{equation}
\label{eq 13.21}
\Gamma _{s+1}^{\Big[\substack{s  \\ s }\Big] \times (k+1)} - \Gamma _{s+1}^{\Big[\substack{s  \\ s }\Big] \times k }= 21\cdot2^{k+s-1}\quad \text{if}\quad k > s +1.
\end{equation}

From \eqref{eq 13.19}, \eqref{eq 13.21} we deduce \vspace{0.01 cm}\\
\begin{align}
& \sum_{j = s+2}^{k}\left( \Gamma _{s+1}^{\Big[\substack{s  \\ s }\Big] \times (j+1)} - \Gamma _{s+1}^{\Big[\substack{s  \\ s }\Big] \times j}   \right)  = \sum_{j = s+2}^{k}  21 \cdot 2^{j+s-1} \nonumber \\
& \Leftrightarrow  \sum_{j = s+3}^{k+1}\Gamma _{s+1}^{\Big[\substack{s  \\ s }\Big] \times j} - \sum_{j = s+2}^{k} \Gamma _{s+1}^{\Big[\substack{s  \\ s }\Big] \times j} = 21\cdot2^{k+s}- 21\cdot2^{2s+1} \nonumber \\
& \Leftrightarrow \Gamma _{s+1}^{\Big[\substack{s  \\ s }\Big] \times (k+1)} - \Gamma _{s+1}^{\Big[\substack{s  \\ s }\Big] \times (s+2)} = 21\cdot2^{k+s}- 21\cdot2^{2s+1}  \nonumber \\
& \Leftrightarrow \Gamma _{s+1}^{\Big[\substack{s  \\ s }\Big] \times (k+1)}= 21\cdot2^{k+s}- 21\cdot2^{2s+1} + 21\cdot(2^{3s-1}- 2^{2s-1}) \nonumber \\
& \Leftrightarrow \Gamma _{s+1}^{\Big[\substack{s  \\ s }\Big] \times (k+1)}=  21[2^{k+s} +2^{3s-1}-5\cdot2^{2s-1}] \quad \text{ if $ k > s +1 $}. \label{eq 13.22}\\
& \nonumber
\end{align}
By \eqref{eq 13.19} we see that \eqref{eq 13.22} holds for k = s+1, then  \eqref{eq 13.22} holds for $ k \geq s+1.$\vspace{0.01 cm}\\

 \end{proof}
  
  \begin{lem}
\label{lem 13.2}In the case m = 1,  we have
\begin{align}
 \Gamma _{s}^{\Big[\substack{s \\ s +1 }\Big] \times s} & = 2^{4s-1} - 3\cdot2^{3s-4} + 2^{2s-3}, \label{eq 13.23} \\
 & \nonumber \\
  \Gamma _{s}^{\Big[\substack{s \\ s +1 }\Big] \times k} & =  2^{k+s-1}+ 21\cdot2^{3s-4} -11 \cdot 2^{2s-3}\quad \text{if $ k > s $},\label{eq 13.24} \\
  & \nonumber \\
  \Gamma _{s+1}^{\Big[\substack{s \\ s +1 }\Big] \times (s+1)} & = 2^{4s+1} - 3\cdot2^{3s- 1} + 2^{2s- 1}, \label{eq 13.25} \\
  & \nonumber \\
  \Gamma _{s+1}^{\Big[\substack{s \\ s +1 }\Big] \times k} & = 11\cdot2^{k+s-1} +21\cdot 2^{3s-1} - 53\cdot2^{2s-1}\quad \text{if $ k> s+1 $}.   \label{eq 13.26}
\end{align}
\end{lem}
\begin{proof}

We proceed as in the proof of Lemma \ref{lem 13.1}.\vspace{0.1 cm}\\

\underline{Proof of \eqref{eq 13.23}}\vspace{0.1 cm}\\
Follows from \eqref{eq 11.26} with m = 1, k = s.\vspace{0.1 cm}\\

\underline{Proof of \eqref{eq 13.24},\eqref{eq 13.25}} \vspace{0.1 cm}\\
Consider the matrix

  $$   \left ( \begin{array} {cccccc}
\alpha _{1} & \alpha _{2} & \alpha _{3} &  \ldots & \alpha _{s}  &  \alpha _{s+1} \\
\alpha _{2 } & \alpha _{3} & \alpha _{4}&  \ldots  &  \alpha _{s+1} &  \alpha _{s+2} \\
\vdots & \vdots & \vdots    &  \vdots & \vdots  &  \vdots \\
\alpha _{s-1} & \alpha _{s} & \alpha _{s +1} & \ldots  &  \alpha _{2s-2} &  \alpha _{2s-1}  \\
\alpha _{s} & \alpha _{s+1} & \alpha _{s +2} & \ldots  &  \alpha _{2s-1} &  \alpha _{2s} \\
 \beta  _{1} & \beta  _{2} & \beta  _{3} & \ldots  &  \beta_{s} &  \beta _{s+1}  \\
\beta  _{2} & \beta  _{3} & \beta  _{4} & \ldots  &  \beta_{s+1} &  \beta _{s+2}  \\
\vdots & \vdots & \vdots    &  \vdots & \vdots  &  \vdots \\
\beta  _{s-1} & \beta  _{s} & \beta  _{s+1} & \ldots  &  \beta_{2s-2} &  \beta _{2s-1}  \\
\beta  _{s} & \beta  _{s+1} & \beta  _{s+2} & \ldots  &  \beta_{2s-1} &  \beta _{2s}\\
\beta  _{s+1} & \beta  _{s+2} & \beta  _{s+3} & \ldots  &  \beta_{2s} &  \beta _{2s +1}
\end{array}  \right). $$ \vspace{0.01 cm}\\

We have respectively  by \eqref{eq 11.1}, \eqref{eq 11.2}, \eqref{eq 11.27} and \eqref{eq 11.28}  with m = 1, k = s+1  \vspace{0.01 cm}\\
  \begin{align}
  \sum_{ i = 0}^{s+1} \Gamma _{i}^{\Big[\substack{s  \\ s +1}\Big] \times (s+1)} & = 2^{4s+1},  \label{eq 13.27}\\
  & \nonumber \\
   \sum_{i = 0}^{s+1} \Gamma _{i}^{\Big[\substack{s  \\ s +1}\Big] \times (s+1)}\cdot2^{-i} & =  2^{3s } + 2^{2s} - 2^{s-1},  \label{eq 13.28}\\
   & \nonumber \\
  \sum_{i = 0}^{s -1} \Gamma _{i}^{\Big[\substack{s  \\ s +1}\Big] \times (s+1)} & = 3\cdot2^{3s-4} - 2^{2s-3} \label{eq 13.29}\\
    &  \text{and}\nonumber \\
   \sum_{i = 0}^{s -1} \Gamma _{i}^{\Big[\substack{s  \\ s +1}\Big] \times (s+1)}\cdot2^{-i} & =  7\cdot2^{2s-4} - 3\cdot2^{s-3}.  \label{eq 13.30}
  \end{align} \vspace{0.01 cm}\\  
  From \eqref{eq 13.27}, \eqref{eq 13.28},  \eqref{eq 13.29} and \eqref{eq 13.30} we deduce after some calculations  \vspace{0.01 cm}\\
  \begin{align}
\Gamma _{s}^{\Big[\substack{s  \\ s+1 }\Big] \times (s+1)} &  =  21\cdot2^{3s-4} - 3\cdot2^{2s-3},  \label{eq 13.31} \\
& \nonumber \\
\Gamma _{s+1}^{\Big[\substack{s  \\ s+1 }\Big] \times (s+1)} & = 2^{4s+1}-3\cdot2^{3s-1} +2^{2s-1}. \label{eq 13.32}
\end{align}
By \eqref{eq 12.33} with j = 0 we get \vspace{0.1 cm}\\
\begin{equation}
\label{eq 13.33}
\Gamma _{s}^{\Big[\substack{s  \\ s +1}\Big] \times (k+1)} - \Gamma _{s}^{\Big[\substack{s  \\ s +1}\Big] \times k }= 2^{k+s-1}\quad \text{if}\quad k > s.
\end{equation}
From \eqref{eq 13.33}, \eqref{eq 13.31} we deduce \vspace{0.01 cm}\\
\begin{align}
& \sum_{j = s+1}^{k}\left( \Gamma _{s}^{\Big[\substack{s  \\ s +1}\Big] \times (j+1)} - \Gamma _{s}^{\Big[\substack{s  \\ s +1}\Big] \times j}   \right)  = \sum_{j = s+1}^{k}   2^{j+s-1} \nonumber \\
& \Leftrightarrow  \sum_{j = s+2}^{k+1}\Gamma _{s}^{\Big[\substack{s  \\ s +1}\Big] \times j} - \sum_{j = s+1}^{k} \Gamma _{s}^{\Big[\substack{s  \\ s +1}\Big] \times j} = 2^{k+s}- 2^{2s} \nonumber \\
& \Leftrightarrow \Gamma _{s}^{\Big[\substack{s  \\ s+1 }\Big] \times (k+1)} - \Gamma _{s}^{\Big[\substack{s  \\ s +1}\Big] \times (s+1)} = 2^{k+s}- 2^{2s} \nonumber \\
& \Leftrightarrow \Gamma _{s}^{\Big[\substack{s  \\ s +1}\Big] \times (k+1)}= 2^{k+s}- 2^{2s} +  21\cdot2^{3s-4} - 3\cdot2^{2s-3} \nonumber \\
& \Leftrightarrow \Gamma _{s}^{\Big[\substack{s  \\ s +1}\Big] \times (k+1)}= 2^{k+s} +  21\cdot2^{3s-4} - 11\cdot2^{2s-3}  \quad \text{ if $ k > s $}. \label{eq 13.34}\\
& \nonumber 
\end{align}
By \eqref{eq 13.31} we see that \eqref{eq 13.34} holds for k = s, then  \eqref{eq 13.34} holds for $ k \geq s.$\vspace{0.01 cm}\\

\underline{Proof of \eqref{eq 13.26}}  \vspace{0.1 cm}\\
Consider the matrix

  $$   \left ( \begin{array} {cccccc}
\alpha _{1} & \alpha _{2} & \alpha _{3} &  \ldots & \alpha _{s+1}  &  \alpha _{s+2} \\
\alpha _{2 } & \alpha _{3} & \alpha _{4}&  \ldots  &  \alpha _{s+2} &  \alpha _{s+3} \\
\vdots & \vdots & \vdots    &  \vdots & \vdots  &  \vdots \\
\alpha _{s-1} & \alpha _{s} & \alpha _{s +1} & \ldots  &  \alpha _{2s-1} &  \alpha _{2s}  \\
\alpha _{s} & \alpha _{s+1} & \alpha _{s +2} & \ldots  &  \alpha _{2s} &  \alpha _{2s+1} \\
 \beta  _{1} & \beta  _{2} & \beta  _{3} & \ldots  &  \beta_{s+1} &  \beta _{s+2}  \\
\beta  _{2} & \beta  _{3} & \beta  _{4} & \ldots  &  \beta_{s+2} &  \beta _{s+3}  \\
\vdots & \vdots & \vdots    &  \vdots & \vdots  &  \vdots \\
\beta  _{s-1} & \beta  _{s} & \beta  _{s+1} & \ldots  &  \beta_{2s-1} &  \beta _{2s}  \\
\beta  _{s} & \beta  _{s+1} & \beta  _{s+2} & \ldots  &  \beta_{2s} &  \beta _{2s +1}\\
\beta  _{s+1} & \beta  _{s+2} & \beta  _{s+3} & \ldots  &  \beta_{2s+1} &  \beta _{2s +2}
\end{array}  \right). $$ \vspace{0.01 cm}\\

We have respectively  by \eqref{eq 11.1}, \eqref{eq 11.2}, \eqref{eq 11.27} and \eqref{eq 11.28} with m = 1 and  k = s+2  \vspace{0.01 cm}\\
  \begin{align}
  \sum_{ i = 0}^{s+2} \Gamma _{i}^{\Big[\substack{s  \\ s +1}\Big] \times (s+2)} & = 2^{4s+3},  \label{eq 13.35}\\
  & \nonumber \\
   \sum_{i = 0}^{s+2} \Gamma _{i}^{\Big[\substack{s  \\ s +1}\Big] \times (s+2)}\cdot2^{-i} & =  2^{3s +1} + 2^{2s+2} - 2^{s},  \label{eq 13.36} \\
   & \nonumber \\
    \sum_{i = 0}^{s -1} \Gamma _{i}^{\Big[\substack{s  \\ s +1}\Big] \times (s+2)} & = 3\cdot2^{3s-4} - 2^{2s-3} \label{eq 13.37} \\
      &  \text{and}\nonumber \\
   \sum_{i = 0}^{s -1} \Gamma _{i}^{\Big[\substack{s  \\ s +1}\Big] \times (s+2)}\cdot2^{-i} & =  7\cdot2^{2s-4} - 3\cdot2^{s-3}.  \label{eq 13.38}
  \end{align} \vspace{0.01 cm}\\
  From \eqref{eq 13.35}, \eqref{eq 13.36}, \eqref{eq 13.37}, \eqref{eq 13.38} and \eqref{eq 13.24} with k = s+2  we obtain  \vspace{0.01 cm}\\
  \begin{align}
& \Gamma _{s+1}^{\Big[\substack{s  \\ s +1}\Big] \times (s+2)} + \Gamma _{s+2}^{\Big[\substack{s  \\ s +1 }\Big] \times (s+2)} \nonumber \\
& = 2^{4s+3} - (3\cdot2^{3s-4} - 2^{2s-3}) - (2^{s+2 +s-1} + 21\cdot2^{3s-4} - 11\cdot2^{2s-3}), \label{eq 13.39}\\
& \nonumber \\
& \Gamma _{s+1}^{\Big[\substack{s  \\ s +1}\Big] \times (s+2)}\cdot2^{-(s+1)} + \Gamma _{s+2}^{\Big[\substack{s  \\ s +1}\Big] \times (s+2)}\cdot 2^{-(s+2)} \nonumber \\
& =   2^{3s +1} + 2^{2s+2} - 2^{s} -(7\cdot2^{2s-4} - 3\cdot2^{s-3} ) -  (2^{s+2 +s-1} + 21\cdot2^{3s-4} -11\cdot2^{2s-3})\cdot2^{-s}. \label{eq 13.40} 
\end{align}\vspace{0.01 cm}\\
Hence by \eqref{eq 13.39},  \eqref{eq 13.40} we deduce after some calculations \vspace{0.01 cm}\\
\begin{align}
\Gamma _{s+1}^{\Big[\substack{s  \\ s +1}\Big] \times (s+2)}& = 21\cdot2^{3s-1}- 9\cdot 2^{2s-1},  \label{eq 13.41}\\
& \nonumber \\
\Gamma _{s+2}^{\Big[\substack{s  \\ s +1}\Big] \times (s+2)}& = 2^{4s+3} -3\cdot2^{3s+2}+ 2^{2s+2}. \label{eq 13.42}\\
& \nonumber
\end{align}
By \eqref{eq 12.33} with j = 1 we get  \vspace{0.1 cm}\\
\begin{equation}
\label{eq 13.43}
\Gamma _{s+1}^{\Big[\substack{s  \\ s +1}\Big] \times (k+1)} - \Gamma _{s+1}^{\Big[\substack{s  \\ s +1}\Big] \times k }= 11\cdot 2^{k+s-1}\quad \text{if}\quad k > s +1.
\end{equation}

From \eqref{eq 13.43}, \eqref{eq 13.41} we deduce  \vspace{0.01 cm}\\
\begin{align}
& \sum_{j = s+2}^{k}\left( \Gamma _{s+1}^{\Big[\substack{s  \\ s +1}\Big] \times (j+1)} - \Gamma _{s+1}^{\Big[\substack{s  \\ s +1}\Big] \times j}   \right)  = \sum_{j = s+2}^{k}11\cdot2^{j+s-1} \nonumber \\
& \Leftrightarrow  \sum_{j = s+3}^{k+1}\Gamma _{s+1}^{\Big[\substack{s  \\ s +1}\Big] \times j} - \sum_{j = s+2}^{k} \Gamma _{s+1}^{\Big[\substack{s  \\ s +1}\Big] \times j} = 11\cdot2^{k+s}- 11\cdot2^{2s+1} \nonumber \\
& \Leftrightarrow \Gamma _{s+1}^{\Big[\substack{s  \\ s+1 }\Big] \times (k+1)} - \Gamma _{s+1}^{\Big[\substack{s  \\ s +1}\Big] \times (s+2)} = 11\cdot2^{k+s}- 11\cdot2^{2s+1} \nonumber \\
& \Leftrightarrow \Gamma _{s+1}^{\Big[\substack{s  \\ s +1}\Big] \times (k+1)}=11\cdot2^{k+s}- 11\cdot2^{2s+1} +  21\cdot2^{3s-1} - 9\cdot2^{2s-1} \nonumber \\
& \Leftrightarrow \Gamma _{s+1}^{\Big[\substack{s  \\ s +1}\Big] \times (k+1)}=  11\cdot2^{k+s} +  21\cdot2^{3s-1} - 53\cdot2^{2s-1} \quad \text{ if $ k > s+1 $}. \label{eq 13.44}\\
& \nonumber
\end{align}

By \eqref{eq 13.41} we see that \eqref{eq 13.44} holds for k = s+1, then  \eqref{eq 13.44} holds for $ k \geq s+1.$\vspace{0.01 cm}\\
\end{proof}
   \begin{lem}
\label{lem 13.3}In the case $ m \geq 2, $ we have 
\begin{align}
 \Gamma _{s}^{\Big[\substack{s \\ s +m }\Big] \times s} & = 2^{4s +m-2} - 3\cdot2^{3s-4} + 2^{2s-3}, \label{eq 13.45} \\
 & \nonumber \\
 \Gamma _{s}^{\Big[\substack{s \\ s +m }\Big] \times k} & =  2^{k+s-1}+ 21\cdot2^{3s-4} -11 \cdot 2^{2s-3}\quad \text{if $ k > s $}, \label{eq 13.46} \\
 & \nonumber \\
  \Gamma _{s+1}^{\Big[\substack{s \\ s +m }\Big] \times (s+1)} & = 2^{4s+m} - 3\cdot2^{3s- 1} + 2^{2s- 1}, \label{eq 13.47} \\
  & \nonumber \\
   \Gamma _{s+1}^{\Big[\substack{s \\ s + m }\Big] \times k} & = 3\cdot2^{k+s-1} +21\cdot[ 2^{3s-1} - 2^{2s-1}]\quad \text{if $ k> s+1 $},  \label{eq 13.48}\\
   & \nonumber \\
 \Gamma _{s+2}^{\Big[\substack{s \\ s +m }\Big] \times (s+2)} & = 2^{4s+m+2} - 3\cdot2^{3s +2} + 2^{2s +2}. \label{eq 13.49} 
\end{align}
\end{lem}

\begin{proof}
We proceed as in the proof of Lemma \ref{lem 13.2}.\vspace{0.1 cm}\\

\underline{Proof of \eqref{eq 13.45}}\vspace{0.1 cm}\\
Follows from \eqref{eq 11.26} with $ m \geq 2,$ k = s.\vspace{0.1 cm}\\

\underline{Proof of \eqref{eq 13.46},  \eqref{eq 13.47}} \vspace{0.1 cm}\\
Consider the matrix

  $$   \left ( \begin{array} {cccccc}
\alpha _{1} & \alpha _{2} & \alpha _{3} &  \ldots & \alpha _{s}  &  \alpha _{s+1} \\
\alpha _{2 } & \alpha _{3} & \alpha _{4}&  \ldots  &  \alpha _{s+1} &  \alpha _{s+2} \\
\vdots & \vdots & \vdots    &  \vdots & \vdots  &  \vdots \\
\alpha _{s-1} & \alpha _{s} & \alpha _{s +1} & \ldots  &  \alpha _{2s-2} &  \alpha _{2s-1}  \\
\alpha _{s} & \alpha _{s+1} & \alpha _{s +2} & \ldots  &  \alpha _{2s-1} &  \alpha _{2s} \\
 \beta  _{1} & \beta  _{2} & \beta  _{3} & \ldots  &  \beta_{s} &  \beta _{s+1}  \\
\beta  _{2} & \beta  _{3} & \beta  _{4} & \ldots  &  \beta_{s+1} &  \beta _{s+2}  \\
\vdots & \vdots & \vdots    &  \vdots & \vdots  &  \vdots \\
\beta  _{s-1} & \beta  _{s} & \beta  _{s+1} & \ldots  &  \beta_{2s-2} &  \beta _{2s-1}  \\
\beta  _{s} & \beta  _{s+1} & \beta  _{s+2} & \ldots  &  \beta_{2s-1} &  \beta _{2s}\\
\vdots & \vdots & \vdots    &  \vdots & \vdots  &  \vdots \\
\beta  _{s+m} & \beta  _{s+m+1} & \beta  _{s+m+2} & \ldots  &  \beta_{2s+m-1} &  \beta _{2s +m}
\end{array}  \right). $$ \vspace{0.01 cm}\\
We have respectively  by \eqref{eq 11.1}, \eqref{eq 11.2}, \eqref{eq 11.27} and \eqref{eq 11.28}  with $ m \geq 2,$ k = s+1  \vspace{0.01 cm}\\
  \begin{align}
  \sum_{ i = 0}^{s+1} \Gamma _{i}^{\Big[\substack{s  \\ s +m}\Big] \times (s+1)} & = 2^{4s+m},  \label{eq 13.50}\\
  & \nonumber \\
   \sum_{i = 0}^{s+1} \Gamma _{i}^{\Big[\substack{s  \\ s +m}\Big] \times (s+1)}\cdot2^{-i} & =  2^{3s +m-1} + 2^{2s} - 2^{s-1},  \label{eq 13.51}\\
   & \nonumber \\
    \sum_{i = 0}^{s -1} \Gamma _{i}^{\Big[\substack{s  \\ s +m}\Big] \times (s+1)} & = 3\cdot2^{3s-4} - 2^{2s-3} \label{eq 13.52}\\
    &  \text{and}\nonumber \\
   \sum_{i = 0}^{s -1} \Gamma _{i}^{\Big[\substack{s  \\ s +m}\Big] \times (s+1)}\cdot2^{-i} & =  7\cdot2^{2s-4} - 3\cdot2^{s-3}.  \label{eq 13.53}
  \end{align} \vspace{0.01 cm}\\  
  From \eqref{eq 13.50}, \eqref{eq 13.51},  \eqref{eq 13.52} and \eqref{eq 13.53} we deduce after some calculations  \vspace{0.01 cm}\\
  \begin{align}
\Gamma _{s}^{\Big[\substack{s  \\ s+m }\Big] \times (s+1)} &  =  21\cdot2^{3s-4} - 3\cdot2^{2s-3},  \label{eq 13.54} \\
& \nonumber \\
\Gamma _{s+1}^{\Big[\substack{s  \\ s+m }\Big] \times (s+1)} & = 2^{4s+ m}-3\cdot2^{3s-1} +2^{2s-1}. \label{eq 13.55}
\end{align}
By \eqref{eq 12.34} with j = 0, $ m\geq 2 $ we get \vspace{0.1 cm}\\
\begin{equation}
\label{eq 13.56}
\Gamma _{s}^{\Big[\substack{s  \\ s +m}\Big] \times (k+1)} - \Gamma _{s}^{\Big[\substack{s  \\ s +m}\Big] \times k }= 2^{k+s-1}\quad \text{if}\quad k > s.
\end{equation}

From \eqref{eq 13.56}, \eqref{eq 13.54} we deduce \vspace{0.01 cm}\\
\begin{align}
& \sum_{j = s+1}^{k}\left( \Gamma _{s}^{\Big[\substack{s  \\ s +m}\Big] \times (j+1)} - \Gamma _{s}^{\Big[\substack{s  \\ s +m}\Big] \times j}   \right)  = \sum_{j = s+1}^{k}   2^{j+s-1} \nonumber \\
& \Leftrightarrow  \sum_{j = s+2}^{k+1}\Gamma _{s}^{\Big[\substack{s  \\ s +m}\Big] \times j} - \sum_{j = s+1}^{k} \Gamma _{s}^{\Big[\substack{s  \\ s +m}\Big] \times j} = 2^{k+s}- 2^{2s} \nonumber \\
& \Leftrightarrow \Gamma _{s}^{\Big[\substack{s  \\ s+m }\Big] \times (k+1)} - \Gamma _{s}^{\Big[\substack{s  \\ s +m}\Big] \times (s+1)} = 2^{k+s}- 2^{2s} \nonumber \\
& \Leftrightarrow \Gamma _{s}^{\Big[\substack{s  \\ s +m}\Big] \times (k+1)}= 2^{k+s}- 2^{2s} +  21\cdot2^{3s-4} - 3\cdot2^{2s-3} \nonumber \\
& \Leftrightarrow \Gamma _{s}^{\Big[\substack{s  \\ s +m}\Big] \times (k+1)}= 2^{k+s} +  21\cdot2^{3s-4} - 11\cdot2^{2s-3}  \quad \text{ if $ k > s $}. \label{eq 13.57}\\
& \nonumber
\end{align}

By \eqref{eq 13.54} we see that \eqref{eq 13.57} holds for k = s, then  \eqref{eq 13.57} holds for $ k \geq s.$\vspace{0.01 cm}\\

\underline{Proof of \eqref{eq 13.48},\eqref{eq 13.49}}  \vspace{0.1 cm}\\
Consider the matrix

  $$   \left ( \begin{array} {cccccc}
\alpha _{1} & \alpha _{2} & \alpha _{3} &  \ldots & \alpha _{s+1}  &  \alpha _{s+2} \\
\alpha _{2 } & \alpha _{3} & \alpha _{4}&  \ldots  &  \alpha _{s+2} &  \alpha _{s+3} \\
\vdots & \vdots & \vdots    &  \vdots & \vdots  &  \vdots \\
\alpha _{s-1} & \alpha _{s} & \alpha _{s +1} & \ldots  &  \alpha _{2s-1} &  \alpha _{2s}  \\
\alpha _{s} & \alpha _{s+1} & \alpha _{s +2} & \ldots  &  \alpha _{2s} &  \alpha _{2s+1} \\
 \beta  _{1} & \beta  _{2} & \beta  _{3} & \ldots  &  \beta_{s+1} &  \beta _{s+2}  \\
\beta  _{2} & \beta  _{3} & \beta  _{4} & \ldots  &  \beta_{s+2} &  \beta _{s+3}  \\
\vdots & \vdots & \vdots    &  \vdots & \vdots  &  \vdots \\
\beta  _{s-1} & \beta  _{s} & \beta  _{s+1} & \ldots  &  \beta_{2s-1} &  \beta _{2s}  \\
\beta  _{s} & \beta  _{s+1} & \beta  _{s+2} & \ldots  &  \beta_{2s} &  \beta _{2s +1}\\
\beta  _{s+1} & \beta  _{s+2} & \beta  _{s+3} & \ldots  &  \beta_{2s+1} &  \beta _{2s +2}\\
\vdots & \vdots & \vdots    &  \vdots & \vdots  &  \vdots \\
\beta  _{s+m} & \beta  _{s+m+1} & \beta  _{s+m+2} & \ldots  &  \beta_{2s+m} &  \beta _{2s +m+1}
\end{array}  \right). $$ \vspace{0.01 cm}\\

We have respectively  by \eqref{eq 11.1}, \eqref{eq 11.2}, \eqref{eq 11.27} and \eqref{eq 11.28} with $ m \geq 2,$ k = s+2  \vspace{0.01 cm}\\
  \begin{align}
  \sum_{ i = 0}^{s+2} \Gamma _{i}^{\Big[\substack{s  \\ s +m}\Big] \times (s+2)} & = 2^{4s+m+2},  \label{eq 13.58}\\
 &  \nonumber \\
   \sum_{i = 0}^{s+2} \Gamma _{i}^{\Big[\substack{s  \\ s +m}\Big] \times (s+2)}\cdot2^{-i} & =  2^{3s +m} + 2^{2s+2} - 2^{s},  \label{eq 13.59} \\
 & \nonumber \\
  \sum_{i = 0}^{s -1} \Gamma _{i}^{\Big[\substack{s  \\ s + m}\Big] \times (s+2)} & = 3\cdot2^{3s-4} - 2^{2s-3} \label{eq 13.60} \\
      &  \text{and}\nonumber \\
   \sum_{i = 0}^{s -1} \Gamma _{i}^{\Big[\substack{s  \\ s +m}\Big] \times (s+2)}\cdot2^{-i} & =  7\cdot2^{2s-4} - 3\cdot2^{s-3}.  \label{eq 13.61}
  \end{align} \vspace{0.01 cm}\\
  From \eqref{eq 13.58}, \eqref{eq 13.59}, \eqref{eq 13.60}, \eqref{eq 13.61} and \eqref{eq 13.46} with k = s+2  we obtain  \vspace{0.01 cm}\\
  \begin{align}
& \Gamma _{s+1}^{\Big[\substack{s  \\ s +m}\Big] \times (s+2)} + \Gamma _{s+2}^{\Big[\substack{s  \\ s +m }\Big] \times (s+2)} \nonumber \\
& = 2^{4s+m+2} - (3\cdot2^{3s-4} - 2^{2s-3}) - (2^{s+2 +s-1} + 21\cdot2^{3s-4} - 11\cdot2^{2s-3}), \label{eq 13.62}\\
& \nonumber \\
& \Gamma _{s+1}^{\Big[\substack{s  \\ s +m}\Big] \times (s+2)}\cdot2^{-(s+1)} + \Gamma _{s+2}^{\Big[\substack{s  \\ s +m}\Big] \times (s+2)}\cdot 2^{-(s+2)} \nonumber \\
& =   2^{3s +m} + 2^{2s+2} - 2^{s} -(7\cdot2^{2s-4} - 3\cdot2^{s-3} ) -  (2^{s+2 +s-1} + 21\cdot2^{3s-4} -11\cdot2^{2s-3})\cdot2^{-s}. \label{eq 13.63} 
\end{align}\vspace{0.01 cm}\\
Hence by \eqref{eq 13.62},\eqref{eq 13.63} we deduce after some calculations \vspace{0.01 cm}\\
\begin{align}
\Gamma _{s+1}^{\Big[\substack{s  \\ s +m}\Big] \times (s+2)}& = 21\cdot2^{3s-1}- 9\cdot 2^{2s-1},  \label{eq 13.64}\\
& \nonumber \\
\Gamma _{s+2}^{\Big[\substack{s  \\ s +m}\Big] \times (s+2)}& = 2^{4s+m+2} -3\cdot2^{3s+2}+ 2^{2s+2}. \label{eq 13.65}\\
& \nonumber
\end{align}
By \eqref{eq 12.34} with $ j = 1$, $ m \geq 2$  we get \vspace{0.1 cm}\\
\begin{equation}
\label{eq 13.66}
\Gamma _{s+1}^{\Big[\substack{s  \\ s +m}\Big] \times (k+1)} - \Gamma _{s+1}^{\Big[\substack{s  \\ s +m}\Big] \times k }= 3 \cdot 2^{k+s-1}\quad \text{if}\quad k > s +1.
\end{equation}
From \eqref{eq 13.66}, \eqref{eq 13.64} we deduce \vspace{0.01 cm}\\

\begin{align}
& \sum_{j = s+2}^{k}\left( \Gamma _{s+1}^{\Big[\substack{s  \\ s +m}\Big] \times (j+1)} - \Gamma _{s+1}^{\Big[\substack{s  \\ s + m}\Big] \times j}   \right)  = \sum_{j = s+2}^{k}3 \cdot2^{j+s-1} \nonumber \\
& \Leftrightarrow  \sum_{j = s+3}^{k+1}\Gamma _{s+1}^{\Big[\substack{s  \\ s +m}\Big] \times j} - \sum_{j = s+2}^{k} \Gamma _{s+1}^{\Big[\substack{s  \\ s +m}\Big] \times j} = 3\cdot2^{k+s}- 3\cdot2^{2s+1} \nonumber \\
& \Leftrightarrow \Gamma _{s+1}^{\Big[\substack{s  \\ s+m }\Big] \times (k+1)} - \Gamma _{s+1}^{\Big[\substack{s  \\ s +m}\Big] \times (s+2)} = 3\cdot2^{k+s}- 3\cdot2^{2s+1} \nonumber \\
& \Leftrightarrow \Gamma _{s+1}^{\Big[\substack{s  \\ s +m}\Big] \times (k+1)}=3\cdot2^{k+s}- 3\cdot2^{2s+1} +  21\cdot2^{3s-1} - 9\cdot2^{2s-1} \nonumber \\
& \Leftrightarrow \Gamma _{s+1}^{\Big[\substack{s  \\ s +m}\Big] \times (k+1)}=  3\cdot2^{k+s} +  21\cdot2^{3s-1} - 21\cdot2^{2s-1} \quad \text{ if $ k > s+1 $}. \label{eq 13.67}\\
& \nonumber
\end{align}
By \eqref{eq 13.64} we see that \eqref{eq 13.67} holds for k = s+1, then  \eqref{eq 13.67} holds for $ k \geq s+1.$\\

\end{proof}

 \section{\textbf{ COMPUTATION OF  $\Gamma_{s+j}^{\left[s\atop s+m\right]\times k}
\;for\; j\in\left\{2,3\right\},\; k \geq s+j ,\;j\leq m $}}
\label{sec 14}

In this section we show the following reduction formulas to be needed in the induction proof in the next section \vspace{0.1 cm}\\
\begin{align*}
  \Gamma _{s +2}^{\Big[\substack{s \\ s +m }\Big] \times k } &  =  8\cdot \Gamma _{s+1}^{\Big[\substack{s \\ s + (m-1) }\Big] \times (k-1)}   &&\text{if  } k\geq s+2,\; m\geq 2, \\
  & \\
   \Gamma _{s + 3}^{\Big[\substack{s \\ s +m }\Big] \times k } &  =  8\cdot \Gamma _{s+1}^{\Big[\substack{s \\ s + (m- 2) }\Big] \times (k- 2)}   &&\text{if  } k\geq s+ 3,\; m\geq 3. \\
   & \\
    \end{align*}
  We recall that the right hand sides in the above equations have already been computed in section \ref{sec 13}.\\
    In fact we have \vspace{0.05 cm}\\
  \begin{equation*}
 \Gamma _{s+1}^{\Big[\substack{s \\ s +(m- (j-1)) }\Big] \times (k- (j-1))} =
  \begin{cases}
   2^{4s+1} - 3\cdot2^{3s -1} + 2^{2s -1}   & \text{if   } \quad j = 2, \; k = s+2, \; m = 2,  \\
   11\cdot2^{k+s-2} + 21\cdot2^{3s- 1} - 53 \cdot2^{2s- 1} & \text{if   } \quad    j = 2, \; k > s+2, \; m = 2,  \\
 2^{4s + m-1} - 3\cdot2^{3s -1} + 2^{2s -1} & \text{if   }\quad    j = 2, \; k = s+2, \; m > 2,  \\
 3\cdot2^{k + s -2} +21\cdot[2^{3s-1}-2^{2s-1}]  & \text{if  }\quad    j = 2, \; k > s+2, \; m > 2, \\
   2^{4s+1} - 3\cdot2^{3s -1} + 2^{2s -1}   & \text{if   } \quad j = 3, \; k = s+3, \; m = 3,  \\
   11\cdot2^{k+s-3} + 21\cdot2^{3s- 1} - 53 \cdot2^{2s- 1} & \text{if   } \quad    j = 3, \; k > s+3, \; m = 3,  \\
 2^{4s + m-2} - 3\cdot2^{3s -1} + 2^{2s -1} & \text{if   }\quad    j = 3, \; k = s+3, \; m > 3,  \\
 3\cdot2^{k + s -3} +21\cdot[2^{3s-1}-2^{2s-1}]  & \text{if  }\quad    j = 3, \; k > s+3, \; m > 3. 
 \end{cases}
\end{equation*}

  \begin{lem}
\label{lem 14.1}We have 
\begin{align}
 \Gamma _{s +2}^{\Big[\substack{s \\ s +m }\Big] \times (s+2)} &  =  8\cdot \Gamma _{s+1}^{\Big[\substack{s \\ s +(m-1) }\Big] \times (s+1)} =  2^{4s+m+2} - 3\cdot2^{3s +2} + 2^{2s +2} &&\text{if\; $ m \geq 2 $}, \label{eq 14.1} \\
 & \nonumber \\
  \Gamma _{s +2}^{\Big[\substack{s \\ s +m }\Big] \times k } &  =  8\cdot \Gamma _{s+1}^{\Big[\substack{s \\ s + (m-1) }\Big] \times (k-1)} = 3 \cdot 2^{k+s+1}  + 21\cdot2^{3s +2} -21 \cdot2^{2s +2} &&\text{if\; $ k > s + 2, \quad  m\geq 3 $}, \label{eq 14.2}\\
  & \nonumber \\
 \Gamma _{s +2}^{\Big[\substack{s \\ s +2 }\Big] \times k } &  =  8\cdot \Gamma _{s+1}^{\Big[\substack{s \\ s + 1 }\Big] \times (k-1)} = 11\cdot 2^{k+s+1}  + 21\cdot2^{3s +2} -53\cdot2^{2s +2} &&\text{if\; $ k > s + 2 $}. \label{eq 14.3}. \\
 & \nonumber 
\end{align}
\end{lem}

\begin{proof}
\underline{ Proof of \eqref{eq 14.1}}\vspace{0.1 cm}\\

By  \eqref{eq 13.65},\eqref{eq 13.47} and \eqref{eq 13.25} we get \vspace{0.1 cm}\\
\begin{align*}
 \Gamma _{s +2}^{\Big[\substack{s \\ s +m }\Big] \times (s+2)} & = 2^{4s+m+2} - 3\cdot2^{3s +2} + 2^{2s +2} \quad \text{if $ m\geq 2 $}, \\
 & \\
  8\cdot \Gamma _{s+1}^{\Big[\substack{s \\ s +(m-1) }\Big] \times (s+1)} & = 8\cdot[2^{4s+(m-1)} - 3\cdot2^{3s- 1} + 2^{2s- 1}] \\
 & =  2^{4s+m+2} - 3\cdot2^{3s +2} + 2^{2s +2} \quad \text{if $ m\geq 2 $}.
\end{align*}

\underline{ Proof of \eqref{eq 14.2}} \vspace{0.1 cm}\\
Consider the matrix

  $$   \left ( \begin{array} {cccccc}
\alpha _{1} & \alpha _{2} & \alpha _{3} &  \ldots & \alpha _{s+2}  &  \alpha _{s+3} \\
\alpha _{2 } & \alpha _{3} & \alpha _{4}&  \ldots  &  \alpha _{s+3} &  \alpha _{s+4} \\
\vdots & \vdots & \vdots    &  \vdots & \vdots  &  \vdots \\
\alpha _{s-1} & \alpha _{s} & \alpha _{s +1} & \ldots  &  \alpha _{2s} &  \alpha _{2s+1}  \\
\alpha _{s} & \alpha _{s+1} & \alpha _{s +2} & \ldots  &  \alpha _{2s+1} &  \alpha _{2s+2} \\
 \beta  _{1} & \beta  _{2} & \beta  _{3} & \ldots  &  \beta_{s+2} &  \beta _{s+3}  \\
\beta  _{2} & \beta  _{3} & \beta  _{4} & \ldots  &  \beta_{s+3} &  \beta _{s+4}  \\
\vdots & \vdots & \vdots    &  \vdots & \vdots  &  \vdots \\
\beta  _{s-1} & \beta  _{s} & \beta  _{s+1} & \ldots  &  \beta_{2s} &  \beta _{2s+1}  \\
\beta  _{s} & \beta  _{s+1} & \beta  _{s+2} & \ldots  &  \beta_{2s+1} &  \beta _{2s +2}\\
\beta  _{s+1} & \beta  _{s+2} & \beta  _{s+3} & \ldots  &  \beta_{2s+2} &  \beta _{2s +3}\\
\vdots & \vdots & \vdots    &  \vdots & \vdots  &  \vdots \\
\beta  _{s+m} & \beta  _{s+m+1} & \beta  _{s+m+2} & \ldots  &  \beta_{2s+m+1} &  \beta _{2s +m+2}
\end{array}  \right). $$ \vspace{0.01 cm}\\
We have respectively  by \eqref{eq 11.1}, \eqref{eq 11.2}, \eqref{eq 11.27} and \eqref{eq 11.28} with $ m \geq 2,$ k = s+3  \vspace{0.01 cm}\\
  \begin{align}
  \sum_{ i = 0}^{s+ 3} \Gamma _{i}^{\Big[\substack{s  \\ s +m}\Big] \times (s+3)} & = 2^{4s+m+4},  \label{eq 14.4}\\
  & \nonumber \\
   \sum_{i = 0}^{s+3} \Gamma _{i}^{\Big[\substack{s  \\ s +m}\Big] \times (s+3)}\cdot2^{-i} & =  2^{3s +m+1} + 2^{2s+4} - 2^{s+1},  \label{eq 14.5} \\
   & \nonumber \\
  \sum_{i = 0}^{s -1} \Gamma _{i}^{\Big[\substack{s  \\ s + m}\Big] \times (s+3)} & = 3\cdot2^{3s-4} - 2^{2s-3} \label{eq 14.6} \\
   &  \text{and}\nonumber \\
   \sum_{i = 0}^{s -1} \Gamma _{i}^{\Big[\substack{s  \\ s +m}\Big] \times (s+3)}\cdot2^{-i} & =  7\cdot2^{2s-4} - 3\cdot2^{s-3}.  \label{eq 14.7}
  \end{align} \vspace{0.01 cm}\\
  From \eqref{eq 14.4}, \eqref{eq 14.5}, \eqref{eq 14.6}, \eqref{eq 14.7},  \eqref{eq 13.46} and  \eqref{eq 13.48} with k = s+3  we obtain \vspace{0.01 cm}\\
  \begin{align}
& \Gamma _{s+2}^{\Big[\substack{s  \\ s +m}\Big] \times (s+3)} + \Gamma _{s+3}^{\Big[\substack{s  \\ s +m }\Big] \times (s+3)} \label{eq 14.8}  \\
& = 2^{4s+m+4} - (3\cdot2^{3s-4} - 2^{2s-3})-(2^{s+3 +s-1} + 21\cdot2^{3s-4} - 11\cdot2^{2s-3}) \nonumber \\
&  - (3\cdot2^{s+3 +s-1} + 21\cdot2^{3s-1} - 21\cdot2^{2s-1}), \nonumber \\
& \nonumber \\
& \Gamma _{s+2}^{\Big[\substack{s  \\ s +m}\Big] \times (s+3)}\cdot2^{-(s+2)} + \Gamma _{s+3}^{\Big[\substack{s  \\ s +m}\Big] \times (s+3)}\cdot 2^{-(s+3)} \label{eq 14.9} \\
& =   2^{3s +m+1} + 2^{2s+4} - 2^{s+1} -(7\cdot2^{2s-4} - 3\cdot2^{s-3} )  \nonumber \\
& -(2^{s+3 +s-1} + 21\cdot2^{3s-4} - 11\cdot2^{2s-3})\cdot2^{-(s+1)}- (3\cdot2^{s+3 +s-1} + 21\cdot2^{3s-1} - 21\cdot2^{2s-1})\cdot2^{-(s+2)}. \nonumber 
\end{align}\vspace{0.01 cm}\\
Hence by \eqref{eq 14.8}, \eqref{eq 14.9} we deduce after some calculations \vspace{0.01 cm}\\
\begin{align}
\Gamma _{s+2}^{\Big[\substack{s  \\ s +m}\Big] \times (s+3)}& = 21\cdot2^{3s +2}- 9\cdot 2^{2s +2},  \label{eq 14.10}\\
& \nonumber \\
\Gamma _{s+3}^{\Big[\substack{s  \\ s +m}\Big] \times (s+3)}& = 2^{4s+m+4} -3\cdot2^{3s+5}+ 2^{2s+5}. \label{eq 14.11}
\end{align}
By \eqref{eq 12.34} with j = 2,\quad $ m \geq 3 $  we get \vspace{0.1 cm}\\
\begin{equation}
\label{eq 14.12}
\Gamma _{s+2}^{\Big[\substack{s  \\ s +m}\Big] \times (k+1)} - \Gamma _{s+2}^{\Big[\substack{s  \\ s +m}\Big] \times k }= 3 \cdot 2^{k+s +1}\quad \text{if}\quad k > s + 2.
\end{equation}

From \eqref{eq 14.12}, \eqref{eq 14.10} we deduce \vspace{0.01 cm}\\
\begin{align}
& \sum_{j = s+3}^{k}\left( \Gamma _{s+2}^{\Big[\substack{s  \\ s +m}\Big] \times (j+1)} - \Gamma _{s+2}^{\Big[\substack{s  \\ s + m}\Big] \times j}   \right)  = \sum_{j = s+3}^{k}3 \cdot2^{j+s +1} \nonumber \\
& \Leftrightarrow  \sum_{j = s+4}^{k+1}\Gamma _{s+2}^{\Big[\substack{s  \\ s +m}\Big] \times j} - \sum_{j = s+3}^{k} \Gamma _{s+2}^{\Big[\substack{s  \\ s +m}\Big] \times j} = 3\cdot2^{k+s +2}- 3\cdot2^{2s+4} \nonumber \\
& \Leftrightarrow \Gamma _{s+2}^{\Big[\substack{s  \\ s+m }\Big] \times (k+1)} - \Gamma _{s+2}^{\Big[\substack{s  \\ s +m}\Big] \times (s+3)} = 3\cdot2^{k+s+2}- 3\cdot2^{2s+4} \nonumber \\
& \Leftrightarrow \Gamma _{s+2}^{\Big[\substack{s  \\ s +m}\Big] \times (k+1)}=3\cdot2^{k+s +2}- 3\cdot2^{2s+4} + 21\cdot2^{3s +2}- 9\cdot 2^{2s +2}    \nonumber \\
& \Leftrightarrow \Gamma _{s+2}^{\Big[\substack{s  \\ s +m}\Big] \times (k+1)}=  3\cdot2^{k+s +2} +  21\cdot2^{3s +2} - 21\cdot2^{2s + 2} \quad \text{ if $ k > s+ 3 $}. \label{eq 14.13}\\
& \nonumber
\end{align}

By \eqref{eq 14.10} we see that \eqref{eq 14.13} holds for k = s+3, then  \eqref{eq 14.13} holds for $ k \geq s+3. $\vspace{0.01 cm}\\

From \eqref{eq 13.48} with $ m \rightarrow m-1, \quad k \rightarrow k-1,\quad  m \geq 3,\quad  k-1> s+1 $ and \eqref{eq 14.13} we get \vspace{0.01 cm}\\
\begin{align}
8\cdot \Gamma _{s+1}^{\Big[\substack{s  \\ s +(m-1)}\Big] \times (k-1)} & = 8\cdot [3\cdot2^{k+s-2} +21\cdot(2^{3s-1} - 2^{2s-1})], \label{eq 14.14}\\
& \nonumber \\
\Gamma _{s+2}^{\Big[\substack{s  \\ s +m}\Big] \times  k} & =  3\cdot2^{k+s +1} +  21\cdot2^{3s +2} - 21\cdot2^{2s + 2}.  \label{eq 14.15} 
\end{align}\vspace{0.01 cm}\\
By \eqref{eq 14.14}, \eqref{eq 14.15} we obtain for $ m\geq 3,\quad k\geq s+3  $\vspace{0.01 cm}\\
$$\Gamma _{s+2}^{\Big[\substack{s  \\ s +m}\Big] \times  k} = 8\cdot \Gamma _{s+1}^{\Big[\substack{s  \\ s +(m-1)}\Big] \times (k-1)}. $$

\underline{ Proof of \eqref{eq 14.3}}\vspace{0.1 cm}\\
By \eqref{eq 12.34} with $ j = 2,\quad m = 2  $ we get \vspace{0.1 cm}\\

\begin{equation}
\label{eq 14.16}
\Gamma _{s+2}^{\Big[\substack{s  \\ s +2}\Big] \times (k+1)} - \Gamma _{s+2}^{\Big[\substack{s  \\ s +2}\Big] \times k }= 11 \cdot 2^{k+s +1}\quad \text{if}\quad k > s + 2.
\end{equation}

From \eqref{eq 14.16}, \eqref{eq 14.10} we deduce \vspace{0.01 cm}\\
\begin{align}
& \sum_{j = s+3}^{k}\left( \Gamma _{s+2}^{\Big[\substack{s  \\ s +2}\Big] \times (j+1)} - \Gamma _{s+2}^{\Big[\substack{s  \\ s + 2}\Big] \times j}   \right)  = \sum_{j = s+3}^{k}11 \cdot2^{j+s +1} \nonumber \\
& \Leftrightarrow  \sum_{j = s+4}^{k+1}\Gamma _{s+2}^{\Big[\substack{s  \\ s +2}\Big] \times j} - \sum_{j = s+3}^{k} \Gamma _{s+2}^{\Big[\substack{s  \\ s +2}\Big] \times j} = 11\cdot2^{k+s +2}- 11\cdot2^{2s+4} \nonumber \\
& \Leftrightarrow \Gamma _{s+2}^{\Big[\substack{s  \\ s+2 }\Big] \times (k+1)} - \Gamma _{s+2}^{\Big[\substack{s  \\ s +2}\Big] \times (s+3)} = 11\cdot2^{k+s+2}- 11\cdot2^{2s+4} \nonumber \\
& \Leftrightarrow \Gamma _{s+2}^{\Big[\substack{s  \\ s +2}\Big] \times (k+1)}=11\cdot2^{k+s +2}- 11\cdot2^{2s+4} + 21\cdot2^{3s +2}- 9\cdot 2^{2s +2}    \nonumber \\
& \Leftrightarrow \Gamma _{s+2}^{\Big[\substack{s  \\ s +2}\Big] \times (k+1)}=  11\cdot2^{k+s +2} +  21\cdot2^{3s +2} - 53 \cdot2^{2s + 2} \quad \text{ if $ k > s+ 3 $}. \label{eq 14.17}
\end{align}\vspace{0.1 cm}\\

By \eqref{eq 14.10} we see that \eqref{eq 14.17} holds for k = s+3, then  \eqref{eq 14.17} holds for $ k \geq s+3. $\vspace{0.05 cm}\\

From \eqref{eq 13.26} ( with  $ k\rightarrow k-1,\quad k-1 > s +1 $ ) and \eqref{eq 14.17} we have \\
\begin{align*}
\Gamma _{s+2}^{\Big[\substack{s  \\ s +2}\Big] \times k}=  11\cdot2^{k+s +1} +  21\cdot2^{3s +2} - 53 \cdot2^{2s + 2} \quad \text{ if $ k \geq  s+ 3 $}, \\
8\cdot\Gamma _{s+1}^{\Big[\substack{s  \\ s+1 }\Big] \times (k-1)} = 8\cdot[11\cdot2^{k+s -2} +  21\cdot2^{3s -1} - 53 \cdot2^{2s -1}] \quad \text{ if $ k -1 > s+ 1 $}.
\end{align*}
\end{proof}
 \begin{lem}
\label{lem 14.2}We have 
\begin{align}
 \Gamma _{s +3}^{\Big[\substack{s \\ s +m }\Big] \times (s+3)} &  =  8^{2}\cdot \Gamma _{s+1}^{\Big[\substack{s \\ s +(m-2) }\Big] \times (s+1)} = 8^{2}[ 2^{4s+m - 2} - 3\cdot2^{3s -1} + 2^{2s -1}] &&\text{if\; $ m \geq 3 $}, \label{eq 14.18} \\
 & \nonumber \\
  \Gamma _{s +3}^{\Big[\substack{s \\ s +m }\Big] \times k } &  =  8^{2}\cdot \Gamma _{s+1}^{\Big[\substack{s \\ s + (m-2) }\Big] \times (k-2)} =8^{2}[ 3 \cdot 2^{k+s -3}  + 21\cdot2^{3s -1} -21 \cdot2^{2s  -1}] &&\text{if\; $ k > s + 3, \quad  m\geq 4 $}, \label{eq 14.19} \\
  & \nonumber \\
   \Gamma _{s +3}^{\Big[\substack{s \\ s +3 }\Big] \times k } &  =  8^{2}\cdot \Gamma _{s+1}^{\Big[\substack{s \\ s + 1 }\Big] \times (k-2)} =8^{2}[ 11\cdot 2^{k+s -3}  + 21\cdot2^{3s  -1} -53\cdot2^{2s  -1}] &&\text{if\; $ k > s + 3 $}. \label{eq 14.20}\\
     & \nonumber 
\end{align}
\end{lem}

\begin{proof}
\underline{ Proof of \eqref{eq 14.18}}\vspace{0.1 cm}\\

By \eqref{eq 14.11}, \eqref{eq 13.47}  and   \eqref{eq 13.25} we get \vspace{0.1 cm}\\
\begin{align*}
\Gamma _{s+3}^{\Big[\substack{s  \\ s +m}\Big] \times (s+3)}& = 2^{4s+m+4} -3\cdot2^{3s+5}+ 2^{2s+5} \quad \text{if \quad $ m\geq 3 $},\\
& \\
\Gamma _{s+1}^{\Big[\substack{s \\ s +(m-2) }\Big] \times (s+1)} &  =  2^{4s+m - 2} - 3\cdot2^{3s -1} + 2^{2s -1}  \quad \text{if \quad $ m -2 \geq 1 $}.\\
& 
\end{align*}
\underline{ Proof of \eqref{eq 14.19}}\vspace{0.1 cm}\\
 Consider the matrix
 
  $$   \left ( \begin{array} {cccccc}
\alpha _{1} & \alpha _{2} & \alpha _{3} &  \ldots & \alpha _{s+3}  &  \alpha _{s+4} \\
\alpha _{2 } & \alpha _{3} & \alpha _{4}&  \ldots  &  \alpha _{s+4} &  \alpha _{s+5} \\
\vdots & \vdots & \vdots    &  \vdots & \vdots  &  \vdots \\
\alpha _{s-1} & \alpha _{s} & \alpha _{s +1} & \ldots  &  \alpha _{2s +1} &  \alpha _{2s+2}  \\
\alpha _{s} & \alpha _{s+1} & \alpha _{s +2} & \ldots  &  \alpha _{2s+2} &  \alpha _{2s+3} \\
 \beta  _{1} & \beta  _{2} & \beta  _{3} & \ldots  &  \beta_{s+ 3} &  \beta _{s+4}  \\
\beta  _{2} & \beta  _{3} & \beta  _{4} & \ldots  &  \beta_{s+4} &  \beta _{s+5}  \\
\vdots & \vdots & \vdots    &  \vdots & \vdots  &  \vdots \\
\beta  _{s-1} & \beta  _{s} & \beta  _{s+1} & \ldots  &  \beta_{2s+1} &  \beta _{2s+2}  \\
\beta  _{s} & \beta  _{s+1} & \beta  _{s+2} & \ldots  &  \beta_{2s+2} &  \beta _{2s +3}\\
\beta  _{s+1} & \beta  _{s+2} & \beta  _{s+3} & \ldots  &  \beta_{2s+3} &  \beta _{2s +4}\\
\vdots & \vdots & \vdots    &  \vdots & \vdots  &  \vdots \\
\beta  _{s+m} & \beta  _{s+m+1} & \beta  _{s+m+2} & \ldots  &  \beta_{2s+m+2} &  \beta _{2s +m+3}
\end{array}  \right). $$ \vspace{0.01 cm}\\
We have  for  $ m \geq 3,$ k = s+4   \vspace{0.01 cm}\\
  \begin{align}
  \sum_{ i = 0}^{s+ 4} \Gamma _{i}^{\Big[\substack{s  \\ s +m}\Big] \times (s+4)} & = 2^{4s+m+ 6}, && \text{by \; \eqref{eq 11.1}}, \label{eq 14.21}\\
  & \nonumber \\
   \sum_{i = 0}^{s + 4} \Gamma _{i}^{\Big[\substack{s  \\ s +m}\Big] \times (s+4)}\cdot2^{-i} & =  2^{3s +m+2} + 2^{2s+6} - 2^{s+2}  && \text{by \; \eqref{eq 11.2}}, \label{eq 14.22} \\
  &  \nonumber \\
    \sum_{i = 0}^{s -1} \Gamma _{i}^{\Big[\substack{s  \\ s + m}\Big] \times (s+4)} & = 3\cdot2^{3s-4} - 2^{2s-3}  && \text{by \; \eqref{eq 11.27}}, \label{eq 14.23} \\
  & \nonumber \\
     \sum_{i = 0}^{s -1} \Gamma _{i}^{\Big[\substack{s  \\ s +m}\Big] \times (s+4)}\cdot2^{-i} & =  7\cdot2^{2s-4} - 3\cdot2^{s-3}  && \text{by \; \eqref{eq 11.28}}, \label{eq 14.24}\\
  & \nonumber \\
    \Gamma _{s}^{\Big[\substack{s  \\ s +m}\Big] \times (s+4)}&  = 21\cdot2^{3s-4} + 53\cdot2^{2s-3}  && \text{by \; \eqref{eq 13.46} with k = s+4}, \label{eq 14.25} \\
   & \nonumber \\
    \Gamma _{s +1}^{\Big[\substack{s  \\ s +m}\Big] \times (s+4)} &  = 21\cdot2^{3s-1} + 27\cdot2^{2s- 1}  && \text{by \; \eqref{eq 13.48} with k = s+4}, \label{eq 14.26} \\
  & \nonumber \\
   \Gamma _{s +2}^{\Big[\substack{s  \\ s +m}\Big] \times (s+4)} &  =  21\cdot2^{3s +2} + 3 \cdot2^{2s +2}  && \text{by \; \eqref{eq 14.2} with k = s+4}. \label{eq 14.27} 
  \end{align} \vspace{0.01 cm}\\ 
 From  $ \eqref{eq 14.21}, \ldots,\eqref{eq 14.27} $ we get   \vspace{0.01 cm}\\ 
 \begin{align}
&  \Gamma _{s+3}^{\Big[\substack{s  \\ s +m}\Big] \times (s+4)}  +  \Gamma _{s+4}^{\Big[\substack{s  \\ s +m}\Big] \times (s+4)} \label{eq 14.28}\\
& =   \sum_{ i = 0}^{s+ 4} \Gamma _{i}^{\Big[\substack{s  \\ s +m}\Big] \times (s+4)} -   \sum_{i = 0}^{s -1} \Gamma _{i}^{\Big[\substack{s  \\ s + m}\Big] \times (s+4)}
 -  \Gamma _{s}^{\Big[\substack{s  \\ s +m}\Big] \times (s+4)}  -  \Gamma _{s +1}^{\Big[\substack{s  \\ s +m}\Big] \times (s+4)} -   \Gamma _{s+2}^{\Big[\substack{s  \\ s +m}\Big] \times (s+4)} \nonumber\\
 & =  2^{4s+m+ 6} -\big( 3\cdot2^{3s-4} - 2^{2s-3} +  21\cdot2^{3s-4} + 53\cdot2^{2s-3} +  21\cdot2^{3s-1} + 27\cdot2^{2s- 1} +  21\cdot2^{3s +2} + 3 \cdot2^{2s +2}  \big), \nonumber  \\
 &          \nonumber    \\
  &  \Gamma _{s+3}^{\Big[\substack{s  \\ s +m}\Big] \times (s+4)}\cdot2^{-(s+3)}  +  \Gamma _{s+4}^{\Big[\substack{s  \\ s +m}\Big] \times (s+4)}\cdot2^{-(s+4)}\label{eq 14.29} \\
 & =  \sum_{i = 0}^{s + 4} \Gamma _{i}^{\Big[\substack{s  \\ s +m}\Big] \times (s+4)}\cdot2^{-i} -  \sum_{i = 0}^{s -1} \Gamma _{i}^{\Big[\substack{s  \\ s +m}\Big] \times (s+4)}\cdot2^{-i} 
 - \Gamma _{s}^{\Big[\substack{s  \\ s +m}\Big] \times (s+4)}\cdot2^{-s}  - \Gamma _{s +1}^{\Big[\substack{s  \\ s +m}\Big] \times (s+4)}\cdot2^{-(s+1)}\nonumber \\
 &  - \Gamma _{s +2}^{\Big[\substack{s  \\ s +m}\Big] \times (s+4)}\cdot2^{-(s+2)}\nonumber \\
 & =  2^{3s +m+2} + 2^{2s+6} - 2^{s+2} - \big( 7\cdot2^{2s-4} - 3\cdot2^{s-3} + (21\cdot2^{3s-4} + 53\cdot2^{2s-3} )\cdot2^{-s} \nonumber  \\
 &  + ( 21\cdot2^{3s-1} + 27\cdot2^{2s- 1} )\cdot2^{-(s+1)}  + ( 21\cdot2^{3s +2} + 3 \cdot2^{2s +2} )\cdot2^{-(s+2)} \big) \nonumber \\
  & =  2^{3s+m+2} + 9\cdot2^{2s+2} - 5\cdot2^{s+2}. \nonumber
 \end{align}\vspace{0.01 cm}\\ 
Hence by \eqref{eq 14.28}, \eqref{eq 14.29} we deduce after some calculations \vspace{0.01 cm}\\ 
\begin{align}
 \Gamma _{s+3}^{\Big[\substack{s  \\ s +m}\Big] \times (s+4)}&  = 21\cdot2^{3s+5} - 9\cdot2^{2s+5} && \text{if $ m\geq 3  $}, \label{eq 14.30}\\
 & \nonumber \\
  \Gamma _{s+4}^{\Big[\substack{s  \\ s +m}\Big] \times (s+4)}&  = 2^{4s+m+6} - 3\cdot2^{3s+8} + 2^{2s+8} && \text{if $ m\geq 3  $}. \label{eq 14.31}
\end{align}\vspace{0.01 cm}\\ 
By \eqref{eq 12.34} with j = 3, $ m \geq 4 $ we get \vspace{0.01 cm}\\ 
\begin{equation}
\label{eq 14.32}
\Gamma _{s+3}^{\Big[\substack{s  \\ s +m}\Big] \times (k+1)} - \Gamma _{s+3}^{\Big[\substack{s  \\ s +m}\Big] \times k }= 3 \cdot 2^{k+s +3}\quad \text{if}\quad k > s + 3.
\end{equation}

From \eqref{eq 14.32}, \eqref{eq 14.30} we deduce \vspace{0.01 cm}\\
\begin{align}
& \sum_{j = s+4}^{k}\left( \Gamma _{s+3}^{\Big[\substack{s  \\ s +m}\Big] \times (j+1)} - \Gamma _{s+3}^{\Big[\substack{s  \\ s + m}\Big] \times j}   \right)  = \sum_{j = s+4}^{k}3 \cdot2^{j+s +3} \nonumber \\
& \Leftrightarrow  \sum_{j = s+5}^{k+1}\Gamma _{s+3}^{\Big[\substack{s  \\ s +m}\Big] \times j} - \sum_{j = s+4}^{k} \Gamma _{s+3}^{\Big[\substack{s  \\ s +m}\Big] \times j} = 3\cdot2^{k+s + 4}- 3\cdot2^{2s+7} \nonumber \\
& \Leftrightarrow \Gamma _{s+3}^{\Big[\substack{s  \\ s+m }\Big] \times (k+1)} - \Gamma _{s+3}^{\Big[\substack{s  \\ s +m}\Big] \times (s+4)} = 3\cdot2^{k+s+4}- 3\cdot2^{2s+7} \nonumber \\
& \Leftrightarrow \Gamma _{s+3}^{\Big[\substack{s  \\ s +m}\Big] \times (k+1)}=3\cdot2^{k+s +4}- 3\cdot2^{2s+7} + 21\cdot2^{3s +5}- 9\cdot 2^{2s +5}    \nonumber \\
& \Leftrightarrow \Gamma _{s+3}^{\Big[\substack{s  \\ s +m}\Big] \times (k+1)}=  3\cdot2^{k+s +4} +  21\cdot2^{3s +5} - 21\cdot2^{2s + 5} \quad \text{ if $ k > s+ 4 $}. \label{eq 14.33}\\
& \nonumber
\end{align}

By \eqref{eq 14.30} we see that \eqref{eq 14.33} holds for k = s+3, then  \eqref{eq 14.33} holds for $ k \geq s+3. $\vspace{0.01 cm}\\

From \eqref{eq 13.48} with $ m\rightarrow m-2,\quad k\rightarrow k-2,\quad m\geq 4,\quad k-2 > s+1 $  and \eqref{eq 14.33} we get \vspace{0.01 cm}\\
  \begin{align}
8^{2}\cdot \Gamma _{s+1}^{\Big[\substack{s  \\ s +(m-2)}\Big] \times (k-2)} & = 8^{2}\cdot [3\cdot2^{k+s-3} +21\cdot(2^{3s-1} - 2^{2s-1})], \label{eq 14.34}\\
 & \nonumber \\
 \Gamma _{s+3}^{\Big[\substack{s  \\ s +m}\Big] \times  k} & =  3\cdot2^{k+s +3} +  21\cdot2^{3s +5} - 21\cdot2^{2s + 5}.  \label{eq 14.35} 
\end{align}\vspace{0.01 cm}\\
   By \eqref{eq 14.34}, \eqref{eq 14.35} we obtain for $ m\geq 4,\quad k\geq s+4  $\vspace{0.01 cm}\\
$$\Gamma _{s+3}^{\Big[\substack{s  \\ s +m}\Big] \times  k} = 8^{2}\cdot \Gamma _{s+1}^{\Big[\substack{s  \\ s +(m-2)}\Big] \times (k- 2)}. $$

\underline{ Proof of \eqref{eq 14.20}}\vspace{0.1 cm}\\
We have \\
\begin{align}
 \Gamma _{s+3}^{\Big[\substack{s  \\ s + 3}\Big] \times (s+4)} &  = 21 \cdot2^{3s+5} - 9 \cdot2^{2s+5} && \text{by \eqref{eq 14.30} with  m = 3 }, \label{eq 14.36} \\
 & \nonumber \\
 \Gamma _{s+3}^{\Big[\substack{s  \\ s + 3}\Big] \times (k+1)} - \Gamma _{s+3}^{\Big[\substack{s  \\ s + 3}\Big] \times k } & = 11 \cdot 2^{k+s +3} && \text{by \eqref{eq 12.34} with  j = m = 3,\; $ k > s + 3 $}. \label{eq 14.37}
\end{align}\vspace{0.01 cm}\\ 
From \eqref{eq 14.36}, \eqref{eq 14.37} we deduce \vspace{0.01 cm} \\ 
\begin{align}
& \sum_{j = s+4}^{k}\left( \Gamma _{s+3}^{\Big[\substack{s  \\ s +3}\Big] \times (j+1)} - \Gamma _{s+4}^{\Big[\substack{s  \\ s + 3}\Big] \times j}   \right)  = \sum_{j = s+4}^{k}11 \cdot2^{j+s +3} \nonumber \\
& \Leftrightarrow  \sum_{j = s+5}^{k+1}\Gamma _{s+3}^{\Big[\substack{s  \\ s +3}\Big] \times j} - \sum_{j = s+4}^{k} \Gamma _{s+3}^{\Big[\substack{s  \\ s +3}\Big] \times j} = 11\cdot2^{k+s +4}- 11\cdot2^{2s+7} \nonumber \\
& \Leftrightarrow \Gamma _{s+3}^{\Big[\substack{s  \\ s+3 }\Big] \times (k+1)} - \Gamma _{s+3}^{\Big[\substack{s  \\ s +3}\Big] \times (s+4)} = 11\cdot2^{k+s+4}- 11\cdot2^{2s+7} \nonumber \\
& \Leftrightarrow \Gamma _{s+3}^{\Big[\substack{s  \\ s +3}\Big] \times (k+1)}=11\cdot2^{k+s +4}- 11\cdot2^{2s+7} + 21\cdot2^{3s +5}- 9\cdot 2^{2s +5}    \nonumber \\
& \Leftrightarrow \Gamma _{s+3}^{\Big[\substack{s  \\ s +3}\Big] \times (k+1)}=  11\cdot2^{k+s +4} +  21\cdot2^{3s +5} - 53 \cdot2^{2s + 5} \quad \text{ if $ k > s+ 4 $}. \label{eq 14.38}
\end{align}\vspace{0.1 cm}\\
By \eqref{eq 14.36} we see that \eqref{eq 14.38} holds for k = s+3,  then  \eqref{eq 14.38} holds for $ k \geq s+3. $\vspace{0.05 cm}\\

From \eqref{eq 13.26} ( with  $ k\rightarrow k-2,\quad k-2 > s +1 $ ) and \eqref{eq 14.38} we have  \vspace{0.05 cm}\\
\begin{align*}
\Gamma _{s+3}^{\Big[\substack{s  \\ s +3}\Big] \times k}=  11\cdot2^{k+s + 3} +  21\cdot2^{3s +5} - 53 \cdot2^{2s + 5} \quad \text{ if $ k \geq  s+ 4 $}, \\
8^{2}\cdot\Gamma _{s+1}^{\Big[\substack{s  \\ s+1 }\Big] \times (k-2)} = 8^{2}\cdot[11\cdot2^{k+s -3} +  21\cdot2^{3s -1} - 53 \cdot2^{2s -1}] \quad \text{ if $ k -2 > s+ 1 $}.\\
& 
\end{align*}
\end{proof}

 \section{\textbf{A REDUCTION FORMULA FOR   $\Gamma_{s+j}^{\left[s\atop s+m\right]\times k}$ in the case $ 1\leq j\leq m,\;k\geq s+j  $}}
 \label{sec 15}
 
 In this section we prove by induction on j the following reduction formula \\
 
\begin{align*}
 \Gamma _{s +j}^{\Big[\substack{s \\ s +m }\Big] \times k} &  =  8^{j-1}\cdot \Gamma _{s+1}^{\Big[\substack{s \\ s +(m-(j-1)) }\Big] \times (k-(j-1))}  &&\text{if   } 1\leq j \leq m, \quad  k\geq s+j. \\
 & 
\end{align*}
Recall that the right hand side in the above equation has been computed in section \ref{sec 13}.\\
In fact we have \vspace{0.05 cm}\\
\begin{equation*}
 \Gamma _{s+1}^{\Big[\substack{s \\ s +(m-(j-1)) }\Big] \times (k-(j-1))} =
  \begin{cases}
  3\cdot2^{k - j +s} +21\cdot[2^{3s-1}-2^{2s-1}]  & \text{if   }\;1\leq j\leq m-1,\;k>s+j,  \\
 2^{4s+(m-(j-1))} - 3\cdot2^{3s -1} + 2^{2s -1} & \text{if   }\;  1\leq j \leq m,    \; k = s+j, \\
11\cdot2^{k- m +s} + 21\cdot2^{3s- 1} - 53 \cdot2^{2s- 1} & \text{if   } \;  j = m,\; k>s+m.  
\end{cases}
\end{equation*}

\begin{lem}
\label{lem 15.1}We have \\

\begin{align}
 \Gamma _{s +l}^{\Big[\substack{s \\ s +m }\Big] \times (s+l)} &  =  8^{l-1}\cdot \Gamma _{s+1}^{\Big[\substack{s \\ s +(m-(l-1)) }\Big] \times (s+1)}  &&\text{if\; $ 1\leq l\leq m $}, \label{eq 15.1} \\
 & \nonumber \\
  \Gamma _{s+1}^{\Big[\substack{s \\ s +(m-(l-1)) }\Big] \times (s+1)} & =  2^{4s+(m-(l-1))} - 3\cdot2^{3s -1} + 2^{2s -1} &&\text{if\; $ 1\leq l\leq m $}, \label{eq 15.2}\\
  & \nonumber \\
  \Gamma _{s +l}^{\Big[\substack{s \\ s +m }\Big] \times k} &  =  8^{l-1}\cdot \Gamma _{s+1}^{\Big[\substack{s \\ s +(m-(l-1)) }\Big] \times (k-(l-1))}  && \text{if\;$1\leq l\leq m,\;k>s+l $}, \label{eq 15.3}\\
  & \nonumber \\ 
  \Gamma _{s+1}^{\Big[\substack{s \\ s +(m-(l-1)) }\Big] \times (k-(l-1))}& = 3\cdot2^{k - l+s} +21\cdot[2^{3s-1}-2^{2s-1}] && \text{if\;$1\leq l\leq m-1,\;k>s+l $}, \label{eq 15.4} \\
   & \nonumber \\
  \Gamma _{s+1}^{\Big[\substack{s \\ s + 1 }\Big] \times (k-(m-1))}& =11\cdot2^{k-m+s} + 21\cdot2^{3s-1} - 53\cdot2^{2s- 1} && \text{if\;$ l = m,\;k>s+m $}. \label{eq 15.5} \\
  & \nonumber
\end{align}
\end{lem}

\begin{proof}
Let $ l $ be a rational integer such that $ 2\leq l\leq m-1; $ we shall prove lemma \ref{lem 15.1} by strong  induction on  $l. $\vspace{0.1 cm}\\
Assume  
\begin{equation}
\label{eq 15.6}
(H_{l-1}) \quad  \Gamma _{s + j}^{\Big[\substack{s \\ s +m }\Big] \times k}   =  8^{j-1}\cdot \Gamma _{s+1}^{\Big[\substack{s \\ s +(m-(j-1)) }\Big] \times (k-(j-1))}\quad for  \quad   1\leq j\leq l-1 \quad k\geq s+j. 
\end{equation}\vspace{0.1 cm}\\
We are going to show that for $ 2\leq l\leq m-1,$ \\
$$ (H_{l-1})\quad \text{implies}\quad  \Gamma _{s +l}^{\Big[\substack{s \\ s +m }\Big] \times k}   =  8^{l-1}\cdot \Gamma _{s+1}^{\Big[\substack{s \\ s +(m-(l-1)) }\Big] \times (k-(l-1))}   \text{for $ \;k\geq s+l $}. $$ \vspace{0.1 cm}\\

By  Lemmas $ \ref{lem 14.1},\ref{lem 14.2} \quad (H_{l-1}) $ holds for l = 3,4. \vspace{0.1 cm}\\

\underline{The case $ k = s+l $} \vspace{0.1 cm}\\
We have \vspace{0.1 cm}\\
\begin{align}
& \Gamma _{ j}^{\Big[\substack{s \\ s +m }\Big] \times (s+l)} = 21\cdot2^{3j-4}-3\cdot2^{2j-3}  \quad   \text{if \quad $1\leq j\leq s-1$\quad by \eqref{eq 11.21}},\label{eq 15.7}\\
 & \nonumber \\
& \Gamma _{ s}^{\Big[\substack{s \\ s +m }\Big] \times (s+l)}  = 2^{2s+l-1} + 21\cdot2^{3s-4} - 11\cdot2^{2s-3} \quad  \text{ by  \eqref{eq 13.46}}, \label{eq 15.8}\\
 & \nonumber \\
& \Gamma _{ s+1}^{\Big[\substack{s \\ s +m }\Big] \times k}  = 3\cdot2^{k+s -1} +21[2^{3s-1} - 2^{2s-1}] \quad   \text{if \quad  $ k > s+1 $ by \eqref{eq 13.48}}, \label{eq 15.9} \\
 & \nonumber \\
&  \Gamma _{s + j}^{\Big[\substack{s \\ s +m }\Big] \times (s+l)}   =  8^{j-1}\cdot \Gamma _{s+1}^{\Big[\substack{s \\ s +(m-(j-1)) }\Big] \times (s+l-(j-1))} \label{eq 15.10}\\
 & = 3\cdot2^{2s+l+2j-3} + 21\cdot[2^{3j+3s-4} - 2^{2s+3j-4}] \quad    \text{ if $ 1 \leq j\leq l-1 $} \nonumber \\
 & \text{by $ (H_{l-1}) $ and ( \eqref{eq 13.48} with $ m \rightarrow m-(j-1) \geq 2,\; k\rightarrow s+l-(j-1),\; s+l-(j-1) > s+1 $)}. \nonumber \\
 & \nonumber
\end{align}

Consider the matrix \\

  $$   \left ( \begin{array} {cccccc}
\alpha _{1} & \alpha _{2} & \alpha _{3} &  \ldots & \alpha _{s+l-1}  &  \alpha _{s+l} \\
\alpha _{2 } & \alpha _{3} & \alpha _{4}&  \ldots  &  \alpha _{s+l} &  \alpha _{s+l+1} \\
\vdots & \vdots & \vdots    &  \vdots & \vdots  &  \vdots \\
\alpha _{s-1} & \alpha _{s} & \alpha _{s +1} & \ldots  &  \alpha _{2s +l-3} &  \alpha _{2s + l-2}  \\
\alpha _{s} & \alpha _{s+1} & \alpha _{s +2} & \ldots  &  \alpha _{2s + l-2} &  \alpha _{2s+l-1} \\
 \beta  _{1} & \beta  _{2} & \beta  _{3} & \ldots  &  \beta_{s+l-1} &  \beta _{s+l}  \\
\beta  _{2} & \beta  _{3} & \beta  _{4} & \ldots  &  \beta_{s+l} &  \beta _{s+ l+ 2}  \\
\vdots & \vdots & \vdots    &  \vdots & \vdots  &  \vdots \\
\beta  _{s-1} & \beta  _{s} & \beta  _{s+1} & \ldots  &  \beta_{2s + l-3} &  \beta _{2s +l-2}  \\
\beta  _{s} & \beta  _{s+1} & \beta  _{s+2} & \ldots  &  \beta_{2s + l -2} &  \beta _{2s +l-1}\\
\beta  _{s+1} & \beta  _{s+2} & \beta  _{s+3} & \ldots  &  \beta_{2s + l-1} &  \beta _{2s +l}\\
\vdots & \vdots & \vdots    &  \vdots & \vdots  &  \vdots \\
\beta  _{s+m} & \beta  _{s+m+1} & \beta  _{s+m+2} & \ldots  &  \beta_{2s+m +l-2} &  \beta _{2s +m+l-1}
\end{array}  \right). $$ \vspace{0.01 cm}\\

We have respectively  by \eqref{eq 11.1} and  \eqref{eq 11.27}  with  k = s+l   \vspace{0.01 cm}\\
  \begin{align}
  \sum_{ i = 0}^{s+l} \Gamma _{i}^{\Big[\substack{s  \\ s +m}\Big] \times (s+l)} & = 2^{4s+m+2l- 2},  \label{eq 15.11}\\
   \sum_{i = 0}^{s -1} \Gamma _{i}^{\Big[\substack{s  \\ s + m}\Big] \times (s+ l)} & = 3\cdot2^{3s-4} - 2^{2s-3}. \label{eq 15.12} 
      \end{align} \vspace{0.01 cm}\\
  From  $ \eqref{eq 15.8},\eqref{eq 15.9},\ldots,\eqref{eq 15.12} $ we get \vspace{0.01 cm}\\
  
     \begin{align}
 \Gamma _{s +l}^{\Big[\substack{s  \\ s +m}\Big] \times (s+l)}& = 2^{4s+m+2l- 2} -  \sum_{ j = 0}^{s -1} \Gamma _{j}^{\Big[\substack{s  \\ s +m}\Big] \times (s+l)}  -
    \Gamma _{s}^{\Big[\substack{s  \\ s +m}\Big] \times (s+l)} - \sum_{ j = 1}^{l-1} \Gamma _{s+j}^{\Big[\substack{s  \\ s +m}\Big] \times (s+l)} \label{eq 15.13} \\
& =    2^{4s+m+2l- 2} -  (3\cdot2^{3s-4} - 2^{2s-3}) - (2^{2s+l-1} + 21\cdot2^{3s-4} - 11\cdot2^{2s-3}) \nonumber \\
& -\sum_{ j = 1}^{l-1}( 3\cdot2^{2s+l+2j-3} + 21\cdot[2^{3j+3s-4} - 2^{2s+3j-4}] ) \nonumber \\
& =  2^{4s+m+2l- 2}-3\cdot2^{3s+3l-4} + 2^{2s+3l-4}, \nonumber \\
 \Gamma _{s + 1}^{\Big[\substack{s  \\ s +(m-(l-1))}\Big] \times (s+1)} & = 2^{4s+(m-(l-1))} - 3\cdot2^{3s-1} + 2^{2s-1} \quad  \text{by \eqref{eq 13.47} since $ m-(l-1)\geq 2 $}. \label{eq 15.14}
\end{align}
From \eqref{eq 15.13}, \eqref{eq 15.14} we get \vspace{0.01 cm}\\
\begin{equation}
\label{eq 15.15}
\Gamma _{s +l}^{\Big[\substack{s \\ s +m }\Big] \times (s+l)}   =  8^{l-1}\cdot \Gamma _{s+1}^{\Big[\substack{s \\ s +(m-(l-1)) }\Big] \times (s+1)}\quad  \text{if\; $ 1\leq l\leq m -1,\; k=s+l $}. 
\end{equation}
\begin{align*}
&
\end{align*}
\underline{The case $ k > s+l $}\vspace{0.01 cm}\\
Consider the matrix \\

  $$   \left ( \begin{array} {cccccc}
\alpha _{1} & \alpha _{2} & \alpha _{3} &  \ldots & \alpha _{s+l}  &  \alpha _{s+l+1} \\
\alpha _{2 } & \alpha _{3} & \alpha _{4}&  \ldots  &  \alpha _{s+l+1} &  \alpha _{s+l+2} \\
\vdots & \vdots & \vdots    &  \vdots & \vdots  &  \vdots \\
\alpha _{s-1} & \alpha _{s} & \alpha _{s +1} & \ldots  &  \alpha _{2s +l-2} &  \alpha _{2s + l-1}  \\
\alpha _{s} & \alpha _{s+1} & \alpha _{s +2} & \ldots  &  \alpha _{2s + l-1} &  \alpha _{2s+l} \\
 \beta  _{1} & \beta  _{2} & \beta  _{3} & \ldots  &  \beta_{s+l} &  \beta _{s+l +1}  \\
\beta  _{2} & \beta  _{3} & \beta  _{4} & \ldots  &  \beta_{s+l+1} &  \beta _{s+ l+ 2}  \\
\vdots & \vdots & \vdots    &  \vdots & \vdots  &  \vdots \\
\beta  _{s-1} & \beta  _{s} & \beta  _{s+1} & \ldots  &  \beta_{2s + l-2} &  \beta _{2s +l-1}  \\
\beta  _{s} & \beta  _{s+1} & \beta  _{s+2} & \ldots  &  \beta_{2s + l-1} &  \beta _{2s +l }\\
\beta  _{s+1} & \beta  _{s+2} & \beta  _{s+3} & \ldots  &  \beta_{2s+l} &  \beta _{2s +l+1}\\
\vdots & \vdots & \vdots    &  \vdots & \vdots  &  \vdots \\
\beta  _{s+m} & \beta  _{s+m+1} & \beta  _{s+m+2} & \ldots  &  \beta_{2s+m +l-1} &  \beta _{2s +m+l}
\end{array}  \right). $$ \vspace{0.01 cm}\\
We get \vspace{0.01 cm}\\
\begin{align}
& \Gamma _{ s}^{\Big[\substack{s \\ s +m }\Big] \times (s+l+1)}  = 2^{2s+l} + 21\cdot2^{3s-4} - 11\cdot2^{2s-3} \quad  \text{ by  \eqref{eq 13.46} },\label{eq 15.16}\\
 & \nonumber \\
&  \Gamma _{s + j}^{\Big[\substack{s \\ s +m }\Big] \times (s+l+1)}   =  8^{j-1}\cdot \Gamma _{s+1}^{\Big[\substack{s \\ s +(m-(j-1)) }\Big] \times (s+l+1-(j-1))} \label{eq 15.17}\\
 & = 3\cdot2^{2s+l+2j-2} + 21\cdot[2^{3j+3s-4} - 2^{2s+3j-4}]  \quad  \text{ if $ 1 \leq j\leq l-1 $} \nonumber \\
 & \text{by (H) and ( \eqref{eq 13.48} with $ m \rightarrow m-(j-1)\geq 2,\; k\rightarrow s+l+1-(j-1),\;s+l+1-(j-1)> s+1  $)}. \nonumber \\
 & \nonumber
\end{align}
We have respectively  by \eqref{eq 11.1}, \eqref{eq 11.2}, \eqref{eq 11.27} and \eqref{eq 11.28} with $ m \geq 2,$ k = s + l +1  \vspace{0.01 cm}\\
  \begin{align}
  \sum_{ i = 0}^{s+ l+1} \Gamma _{i}^{\Big[\substack{s  \\ s +m}\Big] \times (s+ l+1)} & = 2^{4s+m+2l},  \label{eq 15.18}\\
   & \nonumber \\
  \sum_{i = 0}^{s+ l+1} \Gamma _{i}^{\Big[\substack{s  \\ s +m}\Big] \times (s+ l+1)}\cdot2^{-i} & =  2^{3s +m+ l-1} + 2^{2s+2l} - 2^{s+ l-1},  \label{eq 15.19} \\
   & \nonumber \\
    \sum_{i = 0}^{s  - 1} \Gamma _{i}^{\Big[\substack{s  \\ s + m}\Big] \times (s +l+1)} & = 3\cdot2^{3s-4} - 2^{2s-3} \label{eq 15.20} \\
      &  \text{and}\nonumber \\
   \sum_{i = 0}^{s -1} \Gamma _{i}^{\Big[\substack{s  \\ s +m}\Big] \times (s+ l+1)}\cdot2^{-i} & =  7\cdot2^{2s-4} - 3\cdot2^{s-3}.  \label{eq 15.21}
  \end{align} \vspace{0.01 cm}\\
From  $\eqref{eq 15.16}, \ldots,\eqref{eq 15.21} $ we obtain \vspace{0.01 cm}\\ 
  \begin{align}
& \Gamma _{s+ l}^{\Big[\substack{s  \\ s +m}\Big] \times (s+ l+1)} + \Gamma _{s+l+1}^{\Big[\substack{s  \\ s +m }\Big] \times (s+ l+1)} \label{eq 15.22}  \\
& =2^{4s+m+2l} -  \sum_{ j = 0}^{s -1} \Gamma _{j}^{\Big[\substack{s  \\ s +m}\Big] \times (s+l+1)}  -
    \Gamma _{s}^{\Big[\substack{s  \\ s +m}\Big] \times (s+l+1)} - \sum_{ j = 1}^{l-1} \Gamma _{s+j}^{\Big[\substack{s  \\ s +m}\Big] \times (s+l+1)}\nonumber \\
 & =   2^{4s+m+2l} -  (3\cdot2^{3s-4} - 2^{2s-3}) - (2^{2s+l} + 21\cdot2^{3s-4} - 11\cdot2^{2s-3}) \nonumber \\
& -\sum_{ j = 1}^{l-1}( 3\cdot2^{2s+l+2j-2} + 21\cdot[2^{3j+3s-4} - 2^{2s+3j-4}] ) \nonumber \\  
& = 2^{4s +m+2l} - 3\cdot2^{3s +3l - 4} - 2^{2s +3l - 4},\nonumber \\
& \nonumber \\
 & \Gamma _{s+ l}^{\Big[\substack{s  \\ s +m}\Big] \times (s+ l+1)}\cdot2^{-(s+l)} + \Gamma _{s+l+1}^{\Big[\substack{s  \\ s +m }\Big] \times (s+ l+1)}\cdot2^{-(s+l+1)} \label{eq 15.23}  \\
& =   \sum_{ j = 0}^{s  +l+1} \Gamma _{j}^{\Big[\substack{s  \\ s +m}\Big] \times (s+l+1)}\cdot2^{-j}  -  \sum_{ j = 0}^{s -1} \Gamma _{j}^{\Big[\substack{s  \\ s +m}\Big] \times (s+l+1)}\cdot2^{-j}  -
    \Gamma _{s}^{\Big[\substack{s  \\ s +m}\Big] \times (s+l+1)}\cdot2^{-s} - \sum_{ j = 1}^{l-1} \Gamma _{s+j}^{\Big[\substack{s  \\ s +m}\Big] \times (s+l+1)}\cdot2^{-(s+j)}\nonumber \\
 & = 2^{3s+m+ l -1}+ 2^{2s+2l} - 2^{s+l-1}   -  ( 7\cdot2^{2s-4} - 3\cdot2^{s-3}  )   - (2^{s+l} + 21\cdot2^{2s-4} - 11\cdot2^{s-3}) \nonumber \\
& -\sum_{ j = 1}^{l-1}( 3\cdot2^{2s+l+2j-2} + 21\cdot[2^{3j+3s-4} - 2^{2s+3j-4}] )\cdot2^{-(s+j)} \nonumber \\  
& = 2^{4s +m+2l} +  9\cdot2^{3s +3l - 3} - 5\cdot2^{2s +3l - 3}.\nonumber 
  \end{align}\vspace{0.01 cm}\\   
Hence by \eqref{eq 15.22}, \eqref{eq 15.23} we deduce after some calculations \vspace{0.01 cm}\\ 
  \begin{align}
  \Gamma _{s+ l}^{\Big[\substack{s  \\ s +m}\Big] \times (s+ l+1)} & = 21\cdot2^{3s+3l-4} -9\cdot2^{2s+3l-4}, \label{eq 15.24}\\
   & \nonumber \\
  \Gamma _{s+l+1}^{\Big[\substack{s  \\ s +m }\Big] \times (s+ l+1)} & =  2^{4s +m+2l} - 3 \cdot2^{3s +3l - 1} + 2^{2s +3l - 1}.\label{eq 15.25}\\
 & \nonumber
  \end{align}
  
By \eqref{eq 12.34} with $ j = l $, $ m \geq l+1 $ we get \vspace{0.01 cm}\\ 
\begin{equation}
\label{eq 15.26}
\Gamma _{s+ l}^{\Big[\substack{s  \\ s +m}\Big] \times (k+1)} - \Gamma _{s+l}^{\Big[\substack{s  \\ s +m}\Big] \times k }= 3 \cdot 2^{k+s +2l -3}\quad \text{if}\quad k > s + l.
\end{equation}

From \eqref{eq 15.26}, \eqref{eq 15.24} we deduce \vspace{0.01 cm}\\
\begin{align}
& \sum_{j = s+l+1}^{k}\left( \Gamma _{s+l}^{\Big[\substack{s  \\ s +m}\Big] \times (j+1)} - \Gamma _{s+l}^{\Big[\substack{s  \\ s + m}\Big] \times j}   \right)  = \sum_{j = s+l+1}^{k}3 \cdot2^{j+s +2l-3} \nonumber \\
& \Leftrightarrow  \sum_{j = s+l+2}^{k+1}\Gamma _{s+l}^{\Big[\substack{s  \\ s +m}\Big] \times j} - \sum_{j = s+l+1}^{k} \Gamma _{s+l}^{\Big[\substack{s  \\ s +m}\Big] \times j} = 3\cdot2^{k+s +2l+2}- 3\cdot2^{2s+3l-2} \nonumber \\
& \Leftrightarrow \Gamma _{s+l}^{\Big[\substack{s  \\ s+m }\Big] \times (k+1)} - \Gamma _{s+l}^{\Big[\substack{s  \\ s +m}\Big] \times (s+l+1)} = 3\cdot2^{k+s+ 2l+2}- 3\cdot2^{2s+3l-2} \nonumber \\
& \Leftrightarrow \Gamma _{s+ l}^{\Big[\substack{s  \\ s +m}\Big] \times (k+1)}=3\cdot2^{k+s + 2l+2}- 3\cdot2^{2s+ 3l-2} + 21\cdot2^{3s + 3l -4}- 9\cdot 2^{2s +3l -4}    \nonumber \\
& \Leftrightarrow \Gamma _{s+ l}^{\Big[\substack{s  \\ s +m}\Big] \times (k+1)}=  3\cdot2^{k+s +2l -2} +  21\cdot2^{3s +3l -4} - 21\cdot2^{2s + 3l-4} \quad \text{ if $ k > s+ l $}. \label{eq 15.27}\\
& \nonumber
\end{align}
By \eqref{eq 15.24} we see that \eqref{eq 15.27} holds for $ k = s+ l,$ then  \eqref{eq 15.27} holds for $ k \geq s+l. $\vspace{0.01 cm}\\

From \eqref{eq 13.48} with $ m\rightarrow m- (l-1),\quad k\rightarrow k- (l-1),\quad m\geq l+1,\quad k- (l-1) > s+1 $  and \eqref{eq 15.27} we get \vspace{0.01 cm}\\
  \begin{align}
8^{l-1}\cdot \Gamma _{s+1}^{\Big[\substack{s  \\ s +(m-(l-1))}\Big] \times (k-(l-1))} & = 8^{l-1}\cdot [3\cdot2^{k-(l-1)+s-1} +21\cdot(2^{3s-1} - 2^{2s-1})], \label{eq 15.28}\\
 & \nonumber \\
\Gamma _{s+l}^{\Big[\substack{s  \\ s +m}\Big] \times  k} & =  3\cdot2^{k+s +2l-3} +  21\cdot2^{3s + 3l - 4} - 21\cdot2^{2s + 3l-4}.  \label{eq 15.29} 
\end{align}\vspace{0.01 cm}\\
   By \eqref{eq 15.28}, \eqref{eq 15.29} we obtain  \vspace{0.01 cm}\\
    \begin{equation}
\label{eq 15.30}
\Gamma _{s+ l}^{\Big[\substack{s  \\ s +m}\Big] \times  k} = 8^{l-1}\cdot \Gamma _{s+1}^{\Big[\substack{s  \\ s +(m-(l-1))}\Big] \times (k- (l-1))}
\quad  \text{for} \quad  l\leq m-1,\quad k > s+l.
\end{equation}\vspace{0.01 cm}\\
Observe that in view of \eqref{eq 15.30} the formulas \eqref{eq 15.24} and \eqref{eq 15.25} holds for l = m, it suffices to reproduce carefully  
the proof of the above formulas with $ k\rightarrow s+m+1$ instead of  $ k\rightarrow  s+l+1. $\vspace{0.01 cm}\\
In fact we have \vspace{0.01 cm}\\
\begin{align}
  \Gamma _{s+ m}^{\Big[\substack{s  \\ s +m}\Big] \times (s+ m+1)} & = 21\cdot2^{3s+3m-4} -9\cdot2^{2s+3m-4}, \label{eq 15.31}\\
   & \nonumber \\
  \Gamma _{s+m+1}^{\Big[\substack{s  \\ s +m }\Big] \times (s+ m+1)} & =  2^{4s +3m} - 3 \cdot2^{3s +3m - 1} + 2^{2s +3m - 1}.\label{eq 15.32}\\
  & \nonumber
  \end{align}
  \underline{The case $ k = s+m,\;  l = m $}\vspace{0.1 cm}\\
We have  \vspace{0.01 cm}\\
\begin{align}
 \Gamma _{s+ m}^{\Big[\substack{s  \\ s +m }\Big] \times (s+ m)} &  =  2^{4s +3m -2} - 3 \cdot2^{3s +3m - 4} + 2^{2s +3m - 4} \quad \text{by \eqref{eq 15.25} with $ l \rightarrow m-1 $ }, \label{eq 15.33}\\
  & \nonumber \\
  \Gamma _{s + 1}^{\Big[\substack{s  \\ s + 1}\Big] \times (s+1)} & = 2^{4s+ 1} - 3\cdot2^{3s-1} + 2^{2s-1} \quad  \text{by \eqref{eq 13.25}}.  \label{eq 15.34}
\end{align}
From \eqref{eq 15.33}, \eqref{eq 15.34} we obtain \vspace{0.01 cm}\\
$$ \Gamma _{s+ m}^{\Big[\substack{s  \\ s +m }\Big] \times (s+ m)} = 8^{m-1}\cdot  \Gamma _{s + 1}^{\Big[\substack{s  \\ s + 1}\Big] \times (s+1)}. $$\\
\underline{The case $ k  > s+m,\;  l = m $}\vspace{0.01 cm}\\
We have \\
\begin{align}
 \Gamma _{s+m}^{\Big[\substack{s  \\ s + m}\Big] \times (s+m+1)} &  = 21 \cdot2^{3s+3m-4} - 9 \cdot2^{2s+ 3m-4} && \text{by \eqref{eq 15.31}}, \label{eq 15.35} \\
  & \nonumber \\
\Gamma _{s+m}^{\Big[\substack{s  \\ s + m}\Big] \times (k+1)} - \Gamma _{s+m}^{\Big[\substack{s  \\ s + m}\Big] \times k } & = 11 \cdot 2^{k+s +2m-3} && \text{by \eqref{eq 12.34} with  j = m ,\; $ k > s + m $}. \label{eq 15.36}
\end{align}\vspace{0.01 cm}\\ 
From \eqref{eq 15.35}, \eqref{eq 15.36} we deduce \vspace{0.01 cm} \\ 
\begin{align}
& \sum_{j = s+ m +1}^{k}\left( \Gamma _{s+m}^{\Big[\substack{s  \\ s +m}\Big] \times (j+1)} - \Gamma _{s+m}^{\Big[\substack{s  \\ s + m}\Big] \times j}   \right)  = \sum_{j = s+ m+1}^{k}11 \cdot2^{j+s +2m -3} \nonumber \\
& \Leftrightarrow  \sum_{j = s+ m+2}^{k+1}\Gamma _{s+ m}^{\Big[\substack{s  \\ s +m}\Big] \times j} - \sum_{j = s+m+1}^{k} \Gamma _{s+m}^{\Big[\substack{s  \\ s +m}\Big] \times j} = 11\cdot2^{k+s +2m-2}- 11\cdot2^{2s+3m-2} \nonumber \\
& \Leftrightarrow \Gamma _{s+m}^{\Big[\substack{s  \\ s+m }\Big] \times (k+1)} - \Gamma _{s+m}^{\Big[\substack{s  \\ s +m}\Big] \times (s+ m+1)} = 11\cdot2^{k+s+2m-2} - 11\cdot2^{2s+3m-2} \nonumber \\
& \Leftrightarrow \Gamma _{s+m}^{\Big[\substack{s  \\ s +m}\Big] \times (k+1)}=11\cdot2^{k+s + 2m-2}- 11\cdot2^{2s+ 3m-2} + 21\cdot2^{3s + 3m-4} - 9\cdot 2^{2s +3m-4}    \nonumber \\
& \Leftrightarrow \Gamma _{s+m}^{\Big[\substack{s  \\ s + m}\Big] \times (k+1)}=  11\cdot2^{k+s +2m -2} +  21\cdot2^{3s + 3m-4} - 53 \cdot2^{2s + 3m-4} \quad \text{ if $ k > s+ m $}. \label{eq 15.37}
\end{align}\vspace{0.1 cm}\\
By \eqref{eq 15.31} we see that \eqref{eq 15.37} holds for k = s + m,  then  \eqref{eq 15.37} holds for $ k \geq s+ m. $\vspace{0.05 cm}\\
From \eqref{eq 13.26} ( with  $ k\rightarrow k- (m-1),\quad k -(m-1) > s +1 $ ) and \eqref{eq 15.37} with $ k+1 \rightarrow k $ we have  \vspace{0.05 cm}\\
\begin{align*}
\Gamma _{s+m}^{\Big[\substack{s  \\ s +m}\Big] \times k}=  11\cdot2^{k+s + 2m- 3} +  21\cdot2^{3s +3m - 4} - 53 \cdot2^{2s + 3m - 4} \quad \text{ if $ k >  s+ m $}, \\
8^{m-1}\cdot\Gamma _{s+1}^{\Big[\substack{s  \\ s+1 }\Big] \times (k-(m-1))} = 8^{m-1}\cdot[11\cdot2^{k -(m-1)+s-1} +  21\cdot2^{3s -1} - 53 \cdot2^{2s -1}] \quad \text{ if $ k -(m-1) > s+ 1 $}.
\end{align*}
\end{proof}
\begin{lem}
\label{lem 15.2}We have 
\begin{align}
 \Gamma _{s +l}^{\Big[\substack{s \\ s +m }\Big] \times k} &  =  8^{l-1}\cdot \Gamma _{s+1}^{\Big[\substack{s \\ s +(m-(l-1)) }\Big] \times (k-(l-1))}  &&\text{if   \; $ 1\leq l\leq m, \; k\geq s+l $}, \label{eq 15.38} \\
  & \nonumber \\
 \Gamma _{s+l}^{\Big[\substack{s \\ s + m }\Big] \times (s+l)} & =  2^{4s+ 2l+m-2} - 3\cdot2^{3s  +3l -4} + 2^{2s +3l -4} &&\text{if\; $ 1\leq l\leq m $}, \label{eq 15.39}\\
  & \nonumber \\
  \Gamma _{s +l}^{\Big[\substack{s \\ s +m }\Big] \times k} &  = 3\cdot2^{k+2l+s-3} +21\cdot2^{3s+3l-4} -21\cdot2^{2s+3l-4}     && \text{if \;$1\leq l \leq m -1,\;k > s+l $}, \label{eq 15.40}\\
   & \nonumber \\
  \Gamma _{s+ m}^{\Big[\substack{s \\ s + m }\Big] \times  k} & = 11\cdot2^{k + 2m+s-3} +21\cdot2^{3s +3m -4} - 53\cdot2^{2s +3m-4} && \text{if  \; $ k > s+ m $}. \label{eq 15.41} 
  \end{align}
\end{lem}
\begin{proof}
Lemma \ref{lem 15.2} follows from Lemma \ref{lem 15.1}.
\end{proof}

 \section{\textbf{A REDUCTION  FORMULA  FOR  $\Gamma_{s+j+1}^{\left[s\atop s\right]\times k}$ IN THE CASE  $ 1\leq j\leq s-1,\;k\geq s+j +1 $}}
\label{sec 16} 
In this section we prove the following reduction formula by induction on j  \vspace{0.05 cm}\\
\begin{align*}
 \Gamma _{s + j+1}^{\Big[\substack{s \\ s  }\Big] \times k} &  =  8^{2j}\cdot \Gamma _{s -j+1}^{\Big[\substack{s -j\\ s  -j }\Big] \times (k-2j)}  && \text{if\; $ 0\leq j\leq s-1,\;  k\geq s+ j+1 $}. 
\end{align*} \vspace{0.05 cm}\\
We recall once more that the right hand side in the above equation has been computed in section \ref{sec 13}.\\

In fact we have  \vspace{0.05 cm}\\
\begin{equation*}
 \Gamma _{s -j+1}^{\Big[\substack{s -j\\ s  -j }\Big] \times (k-2j)} =                                                                  
  \begin{cases}
 2^{4s -4j} - 3\cdot2^{3s -3j-1} + 2^{2s-2j -1} & \text{if   }\quad 0 \leq j \leq s-1, \; k = s+j+1, \\
  21\cdot[2^{k+s -3j -1} + 2^{3s -3j-1} -5\cdot2^{2s -2j-1}] & \text{if   }\quad 0 \leq j \leq s-2, \; k > s+j+1, \\
  2^{2k-4s+4} - 3\cdot2^{k-2s +2} + 2 & \text{if   }\quad  j = s-1,\;k >2s. 
\end{cases}
\end{equation*}
\begin{align*}
&
\end{align*}

\begin{lem}
\label{lem 16.1}We have \vspace{0.1 cm}\\
\begin{align}
 \Gamma _{s + 2}^{\Big[\substack{s \\ s  }\Big] \times k} &  =  8^{2}\cdot \Gamma _{s}^{\Big[\substack{s -1\\ s  -1 }\Big] \times (k-2)}  && \text{if \; $ k\geq s+2 $}, \label{eq 16.1} \\
 & \nonumber \\
  \Gamma _{s+ 2}^{\Big[\substack{s \\ s  }\Big] \times (s+2)} & =  2^{4s+ 2} - 3\cdot2^{3s +2} + 2^{2s +3} && \text{if \;  k = s+2 }, \label{eq 16.2} \\
  & \nonumber \\
  \Gamma _{s +2}^{\Big[\substack{s \\ s  }\Big] \times k} &  =  21\cdot[2^{k+s+2} + 2^{3s+2} -5\cdot2^{2s+3}] && \text{if \;$ k > s+2 $}. \label{eq 16.3}\\
  & \nonumber
 \end{align}
 \end{lem}
\begin{proof}
\underline{Proof of \eqref{eq 16.1} with k = s+2}\vspace{0.1 cm}\\

From \eqref{eq 13.20} and  \eqref{eq 13.3} (with $ s\rightarrow s-1 $) we get \vspace{0.1 cm}\\
\begin{align*}
\Gamma _{s+ 2}^{\Big[\substack{s \\ s  }\Big] \times (s+2)} &  =  2^{4s+ 2} - 3\cdot2^{3s +2} + 2^{2s +3}, \\
& \\
 8^{2}\cdot \Gamma _{s}^{\Big[\substack{s -1\\ s  -1 }\Big] \times s} &= 8^{2}\cdot[2^{4s-4} -3\cdot2^{3s-4} + 2^{2s-3}].\\
 &
\end{align*}
\underline{Proof of \eqref{eq 16.1} with $ k >  s+2 $}\vspace{0.1 cm}\\
Consider the matrix \\

  $$   \left ( \begin{array} {cccccc}
\alpha _{1} & \alpha _{2} & \alpha _{3} &  \ldots & \alpha _{s +2}  &  \alpha _{s+ 3} \\
\alpha _{2 } & \alpha _{3} & \alpha _{4}&  \ldots  &  \alpha _{s+3} &  \alpha _{s+4} \\
\vdots & \vdots & \vdots    &  \vdots & \vdots  &  \vdots \\
\alpha _{s-1} & \alpha _{s} & \alpha _{s +1} & \ldots  &  \alpha _{2s} &  \alpha _{2s +1}  \\
\alpha _{s} & \alpha _{s+1} & \alpha _{s +2} & \ldots  &  \alpha _{2s +1} &  \alpha _{2s +2} \\
 \beta  _{1} & \beta  _{2} & \beta  _{3} & \ldots  &  \beta_{s +2} &  \beta _{s+3}  \\
\beta  _{2} & \beta  _{3} & \beta  _{4} & \ldots  &  \beta_{s+3} &  \beta _{s+4}  \\
\vdots & \vdots & \vdots    &  \vdots & \vdots  &  \vdots \\
\beta  _{s-1} & \beta  _{s} & \beta  _{s+1} & \ldots  &  \beta_{2s } &  \beta _{2s +1}  \\
\beta  _{s} & \beta  _{s+1} & \beta  _{s+2} & \ldots  &  \beta_{2s +1} &  \beta _{2s +2}
\end{array}  \right). $$ \vspace{0.01 cm}\\

We have \vspace{0.01 cm}\\
  \begin{align}
  \sum_{ i = 0}^{s+3} \Gamma _{i}^{\Big[\substack{s  \\ s }\Big] \times (s+3)} & = 2^{4s +4}  && \text{by \eqref{eq 11.1} with m = 0,\; k = s+3}, \label{eq 16.4}\\
  & \nonumber \\
   \sum_{i = 0}^{s+ 3} \Gamma _{i}^{\Big[\substack{s  \\ s }\Big] \times (s+3)}\cdot2^{-i} & =  2^{3s  +1} + 2^{2s +4} - 2^{s +1}  && \text{by \eqref{eq 11.2} with m = 0,\; k = s+3},  \label{eq 16.5}\\
   & \nonumber \\
    \sum_{i = 0}^{s -1} \Gamma _{i}^{\Big[\substack{s  \\ s }\Big] \times (s+3)} & = 3\cdot2^{3s-4} - 2^{2s-3}&& \text{by \eqref{eq 11.27} with m = 0,\; k = s+3},  \label{eq 16.6}\\
    & \nonumber \\
     \sum_{i = 0}^{s -1} \Gamma _{i}^{\Big[\substack{s  \\ s }\Big] \times (s+ 3)}\cdot2^{-i} & =  7\cdot2^{2s-4} - 3\cdot2^{s-3} && \text{by \eqref{eq 11.28} with m = 0,\; k = s+3},  \label{eq 16.7}\\
     & \nonumber \\
  \Gamma _{s}^{\Big[\substack{s  \\ s }\Big] \times (s+ 3)} &  =  21\cdot2^{3s-4}  +69 \cdot2^{2s-3} && \text{by \eqref{eq 13.2} with k = s+3},  \label{eq 16.8} \\
  & \nonumber \\
    \Gamma _{s+1}^{\Big[\substack{s  \\ s }\Big] \times (s+ 3)} & = 21\cdot[2^{3s-1} +3\cdot2^{2s-1}] && \text{by \eqref{eq 13.4} with  k = s+3}. \label{eq 16.9} \\
    & \nonumber 
  \end{align} \vspace{0.01 cm}\\  
   From $\eqref{eq 16.4},\ldots, \eqref{eq 16.9} $ with k = s+3  we obtain  \vspace{0.01 cm}\\
  \begin{align}
& \Gamma _{s+2}^{\Big[\substack{s  \\ s }\Big] \times (s+3)} + \Gamma _{s+3}^{\Big[\substack{s  \\ s }\Big] \times (s+3)} \label{eq 16.10}  \\
& = 2^{4s+4} - (3\cdot2^{3s-4} - 2^{2s-3}) - ( 21\cdot2^{3s-4}  +69 \cdot2^{2s-3}) -( 21\cdot[2^{3s-1} +3\cdot2^{2s-1}]), \nonumber \\
& \nonumber \\
& \Gamma _{s+2}^{\Big[\substack{s  \\ s }\Big] \times (s+ 3)}\cdot2^{-(s+2)} + \Gamma _{s+3}^{\Big[\substack{s  \\ s }\Big] \times (s+3)}\cdot 2^{-(s+3)}  \label{eq 16.11}  \\
& =   2^{3s +1} + 2^{2s+4} - 2^{s+1} -(7\cdot2^{2s-4} - 3\cdot2^{s-3} ) -  ( 21\cdot2^{3s-4}  +69 \cdot2^{2s-3})\cdot2^{-s} \nonumber \\
&  - (21\cdot[2^{3s-1} +3\cdot2^{2s-1}] )\cdot2^{-(s+1)}. \nonumber 
\end{align}\vspace{0.01 cm}\\
Hence by \eqref{eq 16.10}, \eqref{eq 16.11} we deduce after some calculations  \vspace{0.01 cm}\\
\begin{align}
\Gamma _{s+2}^{\Big[\substack{s  \\ s }\Big] \times (s+3)} &  = 21\cdot2^{3s +2}- 21\cdot 2^{2s +3}, \label{eq 16.12}\\
& \nonumber \\
\Gamma _{s+3}^{\Big[\substack{s  \\ s }\Big] \times (s+3)} &  = 2^{4s+4} -3\cdot2^{3s+5}+ 2^{2s+7}.\label{eq 16.13}\\
 & \nonumber 
\end{align}
By \eqref{eq 12.32} with j = 2,\; $ s\geq 3 $ we get \vspace{0.01 cm}\\
\begin{equation}
\label{eq 16.14}
\Gamma _{s+2}^{\Big[\substack{s  \\ s }\Big] \times (k+1)} - \Gamma _{s+2}^{\Big[\substack{s  \\ s }\Big] \times k }= 21\cdot2^{k+s +2}\quad \text{if}\quad k > s + 2.
\end{equation}

From \eqref{eq 16.14}, \eqref{eq 16.12} we deduce  \vspace{0.01 cm}\\
\begin{align}
& \sum_{j = s+3}^{k}\left( \Gamma _{s+2}^{\Big[\substack{s  \\ s }\Big] \times (j+1)} - \Gamma _{s+2}^{\Big[\substack{s  \\ s }\Big] \times j}   \right)  = \sum_{j = s+3}^{k}  21 \cdot 2^{j+s +2} \nonumber \\
& \Leftrightarrow  \sum_{j = s+4}^{k+1}\Gamma _{s+2}^{\Big[\substack{s  \\ s }\Big] \times j} - \sum_{j = s+2}^{k} \Gamma _{s+2}^{\Big[\substack{s  \\ s }\Big] \times j} = 21\cdot2^{k+s +3}- 21\cdot2^{2s+5} \nonumber \\
& \Leftrightarrow \Gamma _{s+2}^{\Big[\substack{s  \\ s }\Big] \times (k+1)} - \Gamma _{s+2}^{\Big[\substack{s  \\ s }\Big] \times (s+3)} = 21\cdot2^{k+s +3}- 21\cdot2^{2s+5}  \nonumber \\
& \Leftrightarrow \Gamma _{s+2}^{\Big[\substack{s  \\ s }\Big] \times (k+1)}= 21\cdot2^{k+s +3}- 21\cdot2^{2s+5} + 21\cdot(2^{3s +2}-21\cdot2^{2s +3}) \nonumber \\
& \Leftrightarrow \Gamma _{s+2}^{\Big[\substack{s  \\ s }\Big] \times (k+1)}=  21\cdot[2^{k+s+3} +2^{3s +2}-5\cdot2^{2s +3}] \quad \text{ if $ k > s + 2 $}. \label{eq 16.15}\\
& \nonumber
\end{align}
By \eqref{eq 16.12} we see that \eqref{eq 16.15} holds for k = s+ 2, then  \eqref{eq 16.15} holds for $ k \geq s+2.$\vspace{0.01 cm}\\

From \eqref{eq 13.4} (with $ k\rightarrow k-2,\quad s \rightarrow s-1,\quad k-2 > (s-1) +1 $) and \eqref{eq 16.15} we have  \vspace{0.1 cm}\\
\begin{align*}
\Gamma _{s+2}^{\Big[\substack{s  \\ s }\Big] \times  k} & =  21\cdot[2^{k+s+2} +2^{3s +2}-5\cdot2^{2s +3}], \\
8^{2}\cdot\Gamma _{s}^{\Big[\substack{s -1 \\ s -1}\Big] \times  (k-2)} & = 8^{2}\cdot21\cdot[2^{k+s - 4} +2^{3s -4}-5\cdot2^{2s -3}]. 
\end{align*}
\end{proof}

\begin{lem}
\label{lem 16.2}We have for $ s\geq 2 $ \vspace{0.1 cm}\\
\begin{align}
 \Gamma _{s + j+1}^{\Big[\substack{s \\ s  }\Big] \times k} &  =  8^{2j}\cdot \Gamma _{s -j+1}^{\Big[\substack{s -j\\ s  -j }\Big] \times (k-2j)}  && \text{if\; $ 0\leq j\leq s-1,\;  k\geq s+ j+1 $}, \label{eq 16.16} \\
 & \nonumber \\
  \Gamma _{s -j +1}^{\Big[\substack{s -j \\ s -j }\Big] \times (s -j+1)} & =  2^{4s -4j} - 3\cdot2^{3s -3j-1} + 2^{2s -2j-1} && \text{if \;  $ 0\leq j\leq s-1,\; k = s+j+1 $ }, \label{eq 16.17} \\
  & \nonumber \\
   \Gamma _{s -j +1}^{\Big[\substack{s - j \\ s -j }\Big] \times (k-2j)} &  =  21\cdot[2^{k -3j+s-1} + 2^{3s -3j-1} -5\cdot2^{2s -2j-1}] && \text{if \;$ 0\leq j\leq s - 2,\;  k > s + j+1 $}, \label{eq 16.18}\\
   & \nonumber \\
    \Gamma _{s -j +1}^{\Big[\substack{s - j \\ s -j }\Big] \times (k-2j)} &=  \Gamma _{2}^{\Big[\substack{ 1 \\  1}\Big] \times (k-2(s-1))}  = 2^{2k- 4s+4} - 3\cdot2^{k -2s + 2} + 2   && \text{if $ j = s -1,  k >  2s  $}. \label{eq 16.19}\\
    & \nonumber
 \end{align}
 \end{lem}
 
\begin{proof}We proceed as in section \ref{sec 15} by induction on j. \vspace{0.1 cm} \\
Let $l $  be a rational integer such that  $ 2\leq l\leq s-2. $\vspace{0.1 cm}\\
Assume   \vspace{0.1 cm}\\
\begin{equation}
\label{eq 16.20}
(H_{l-1}) \quad  \Gamma _{s + j +1}^{\Big[\substack{s \\ s  }\Big] \times k}   =  8^{2j}\cdot \Gamma _{s -j+1}^{\Big[\substack{s -j \\ s -j }\Big] \times (k- 2j)}\quad for  \quad   0 \leq j\leq l -1,\quad k\geq s+j +1.
\end{equation}\vspace{0.1 cm}\\
We are going to show that $ (H_{l}) $ holds, that is  \vspace{0.1 cm}\\
\begin{equation}
\label{eq 16.21}
 \Gamma _{s + l +1}^{\Big[\substack{s \\ s  }\Big] \times k}   =  8^{2l}\cdot \Gamma _{s - l+1}^{\Big[\substack{s -l \\ s -l }\Big] \times (k- 2l)}\quad for \; k\geq s+ l +1.
\end{equation}\vspace{0.1 cm}\\
By Lemma \ref{lem 16.1}  $ (H_{l-1}) $ holds for $ l = 2  $ (see  \eqref{eq 16.20}).    \vspace{0.1 cm}\\

\underline{Proof of \eqref{eq 16.21} with $ k = s+ l+1 $} \vspace{0.1 cm}\\ 
Consider the matrix \\

  $$   \left ( \begin{array} {cccccc}
\alpha _{1} & \alpha _{2} & \alpha _{3} &  \ldots & \alpha _{s +l}  &  \alpha _{s+ l+1} \\
\alpha _{2 } & \alpha _{3} & \alpha _{4}&  \ldots  &  \alpha _{s+l+1} &  \alpha _{s+ l+2} \\
\vdots & \vdots & \vdots    &  \vdots & \vdots  &  \vdots \\
\alpha _{s-1} & \alpha _{s} & \alpha _{s +1} & \ldots  &  \alpha _{2s +2l-2} &  \alpha _{2s +2l-1}  \\
\alpha _{s} & \alpha _{s+1} & \alpha _{s +2} & \ldots  &  \alpha _{2s +2l-1} &  \alpha _{2s +2l} \\
 \beta  _{1} & \beta  _{2} & \beta  _{3} & \ldots  &  \beta_{s +l} &  \beta _{s+l+1}  \\
\beta  _{2} & \beta  _{3} & \beta  _{4} & \ldots  &  \beta_{s+ l+1} &  \beta _{s+ l+2}  \\
\vdots & \vdots & \vdots    &  \vdots & \vdots  &  \vdots \\
\beta  _{s-1} & \beta  _{s} & \beta  _{s+1} & \ldots  &  \beta_{2s +2l -2} &  \beta _{2s +2l-1}  \\
\beta  _{s} & \beta  _{s+1} & \beta  _{s+2} & \ldots  &  \beta_{2s + 2l-1} &  \beta _{2s + 2l}
\end{array}  \right). $$ \vspace{0.01 cm}\\
We have \vspace{0.01 cm}\\
  \begin{align}
&  \sum_{ i = 0}^{s+ l+1} \Gamma _{i}^{\Big[\substack{s  \\ s }\Big] \times (s+l+1)}  = 2^{4s + 2l} && \text{by \eqref{eq 11.1} with m = 0,\;$ k = s+ l+1 $}, \label{eq 16.22}\\
& \nonumber \\
& \sum_{i = 0}^{s -1} \Gamma _{i}^{\Big[\substack{s  \\ s }\Big] \times (s+ l+1)}  = 3\cdot2^{3s-4} - 2^{2s-3}&& \text{by \eqref{eq 11.27} with m = 0,\;$ k = s+ l+1 $},  \label{eq 16.23}\\
& \nonumber \\
&  \Gamma _{s}^{\Big[\substack{s  \\ s }\Big] \times (s+ l+1)}   = 3\cdot2^{2s+l} + 21\cdot2^{3s-4}  -  27\cdot2^{2s-3} && \text{by \eqref{eq 13.2} with $ k = s+ l+1 $},  \label{eq 16.24} \\
& \nonumber \\
&   \Gamma _{s+ j +1}^{\Big[\substack{s  \\ s }\Big] \times (s+ l+1)}  =  8^{2j}\cdot \Gamma _{s -j+1}^{\Big[\substack{s -j \\ s -j }\Big] \times (s+l+1 - 2j)}  && \text{by \eqref{eq 16.20} $\;\text{for}   \;0\leq j\leq l-1,\; k = s+ l+1 $}  \label{eq 16.25} \\
 & =  2^{6j}\cdot[ 21\cdot(2^{(s+l+1-2j)+(s-j) -1}+ 2^{3(s-j)-1} -5\cdot2^{2(s-j)-1})] && \text{by \eqref{eq 13.4} with $ s\rightarrow s-j,\; k\rightarrow s+l+1-2j > s-j+1 $} \nonumber \\
 & = 21\cdot[2^{2s+3j+l} + 2^{3s+3j-1} - 5\cdot2^{2s+4j-1}]. \nonumber
  \end{align} \vspace{0.01 cm}\\  
From \eqref{eq 16.22}, \eqref{eq 16.23}, \eqref{eq 16.24} and  \eqref{eq 16.25} we get after some calculations  \vspace{0.01 cm}\\  
\begin{align}
 \Gamma _{s+l+1}^{\Big[\substack{s  \\ s }\Big] \times (s+l+1)}& = 2^{4s + 2l}-(3\cdot2^{3s-4} - 2^{2s-3})
 -( 3\cdot2^{2s+l} + 21\cdot2^{3s-4}  -  27\cdot2^{2s-3})\label{eq 16.26} \\
 &  -\left(\sum_{j=0}^{l-1}21\cdot[2^{2s+3j+l} + 2^{3s+3j-1} - 5\cdot2^{2s+4j-1}]\right)\nonumber \\
 & =  2^{4s+2l}  -3\cdot2^{3s +3l -1} + 2^{2s+4l-1}.\nonumber
\end{align}\vspace{0.01 cm}\\  
By \eqref{eq 13.3} with $ s\rightarrow s-l $ and \eqref{eq 16.26} we obtain \vspace{0.01 cm}\\ 

\begin{align*}
\Gamma _{s+l+1}^{\Big[\substack{s  \\ s }\Big] \times (s+l+1)} &= 8^{2l}\cdot\Gamma _{s - l+1}^{\Big[\substack{s -l \\ s -l }\Big] \times (s-l+1)}=2^{6l}\cdot[2^{4(s-l)} - 3\cdot2^{3(s-l)-1} + 2^{2(s-l)-1}] \\
& =  2^{4s+2l} -3\cdot2^{3s +3l -1} + 2^{2s+4l-1}. \quad   \qed
\end{align*}
\underline{Proof of \eqref{eq 16.21} with $ k > s+ l+1 $}\vspace{0.1 cm}\\ 
Consider the matrix \\

  $$   \left ( \begin{array} {cccccc}
\alpha _{1} & \alpha _{2} & \alpha _{3} &  \ldots & \alpha _{s +l+1}  &  \alpha _{s+ l +2} \\
\alpha _{2 } & \alpha _{3} & \alpha _{4}&  \ldots  &  \alpha _{s+l+2} &  \alpha _{s+ l+3} \\
\vdots & \vdots & \vdots    &  \vdots & \vdots  &  \vdots \\
\alpha _{s-1} & \alpha _{s} & \alpha _{s +1} & \ldots  &  \alpha _{2s+l-1} &  \alpha _{2s + l}  \\
\alpha _{s} & \alpha _{s+1} & \alpha _{s +2} & \ldots  &  \alpha _{2s +l} &  \alpha _{2s +l+1} \\
 \beta  _{1} & \beta  _{2} & \beta  _{3} & \ldots  &  \beta_{s + l+1} &  \beta _{s+l+2}  \\
\beta  _{2} & \beta  _{3} & \beta  _{4} & \ldots  &  \beta_{s+ l+2} &  \beta _{s+ l+3}  \\
\vdots & \vdots & \vdots    &  \vdots & \vdots  &  \vdots \\
\beta  _{s-1} & \beta  _{s} & \beta  _{s+1} & \ldots  &  \beta_{2s +l-1} &  \beta _{2s + l}  \\
\beta  _{s} & \beta  _{s+1} & \beta  _{s+2} & \ldots  &  \beta_{2s +l} &  \beta _{2s + l+1}
\end{array}  \right). $$ \vspace{0.01 cm}\\
We have \vspace{0.01 cm}\\
  \begin{align}
&  \sum_{ i = 0}^{s+l+2} \Gamma _{i}^{\Big[\substack{s  \\ s }\Big] \times (s+l+2)} = 2^{4s + 2l+2} && \text{by \eqref{eq 11.1} with m = 0,\;$ k = s+ l+2 $}, \label{eq 16.27}\\
&  \nonumber \\
&  \sum_{i = 0}^{s+ l+2} \Gamma _{i}^{\Big[\substack{s  \\ s }\Big] \times (s+ l+2)}\cdot2^{-i}  =  2^{3s  +l} + 2^{2s +2l+2} - 2^{s + l}  && \text{by \eqref{eq 11.2} with m = 0,\; $ k = s+ l+2 $},  \label{eq 16.28}\\
& \nonumber \\
 &  \sum_{i = 0}^{s -1} \Gamma _{i}^{\Big[\substack{s  \\ s }\Big] \times (s+ l+2)}  = 3\cdot2^{3s-4} - 2^{2s-3}&& \text{by \eqref{eq 11.27} with m = 0,\; $ k = s+ l+2 $},  \label{eq 16.29}\\
 & \nonumber \\
 &   \sum_{i = 0}^{s -1} \Gamma _{i}^{\Big[\substack{s  \\ s }\Big] \times (s+ l+2)}\cdot2^{-i}  =  7\cdot2^{2s-4} - 3\cdot2^{s-3} && \text{by \eqref{eq 11.28} with m = 0,\;$ k = s+l+2 $},  \label{eq 16.30}\\
 & \nonumber \\
 &   \Gamma _{s}^{\Big[\substack{s  \\ s }\Big] \times (s+ l+2)}   = 3\cdot2^{2s+l+1} + 21\cdot2^{3s-4}  -  27\cdot2^{2s-3} && \text{by \eqref{eq 13.2} with $ k = s+ l+2 $},  \label{eq 16.31} \\
 & \nonumber \\
  &   \Gamma _{s+ j +1}^{\Big[\substack{s  \\ s }\Big] \times (s+ l+ 2)}  =  8^{2j}\cdot \Gamma _{s -j+1}^{\Big[\substack{s -j \\ s -j }\Big] \times (s+l+2 - 2j)}  && \text{by \eqref{eq 16.20} $\;\text{for}   \;0\leq j\leq l-1,\; k = s+ l+2 $}  \label{eq 16.32} \\
 & =  2^{6j}\cdot[ 21\cdot(2^{(s+l+2-2j)+(s-j) -1}+ 2^{3(s-j)-1} -5\cdot2^{2(s-j)-1})] && \text{by \eqref{eq 13.4} with $ s\rightarrow s-j,\; k\rightarrow s+l+2-2j > s-j+1 $} \nonumber \\
 & = 21\cdot[2^{2s+3j+l +1} + 2^{3s+3j-1} - 5\cdot2^{2s+4j-1}]. \nonumber
  \end{align} \vspace{0.01 cm}\\  
  From $\eqref{eq 16.27},\ldots, \eqref{eq 16.32} $ with $ k = s+l+2 $ we obtain \vspace{0.01 cm}\\
  \begin{align}
& \Gamma _{s+ l+1}^{\Big[\substack{s  \\ s }\Big] \times (s+l+2)} + \Gamma _{s+l+2}^{\Big[\substack{s  \\ s }\Big] \times (s+ l+2)} \label{eq 16.33}  \\
& = 2^{4s+2l +2} - (3\cdot2^{3s-4} - 2^{2s-3}) - (3\cdot2^{2s+l+1} + 21\cdot2^{3s-4}  -  27\cdot2^{2s-3}) \nonumber \\
& - \sum_{j=0}^{l-1}21\cdot[2^{2s+3j+l+1} + 2^{3s+3j-1} - 5\cdot2^{2s+4j-1}],   \nonumber \\
& \nonumber \\
& \Gamma _{s+ l+1}^{\Big[\substack{s  \\ s }\Big] \times (s+ l+2)}\cdot2^{-(s+l+1)} + \Gamma _{s+l+2}^{\Big[\substack{s  \\ s }\Big] \times (s+l+2)}\cdot 2^{-(s+l+2)}  \label{eq 16.34}  \\
& =   2^{3s +l} + 2^{2s+2l+2} - 2^{s+l} -(7\cdot2^{2s-4} - 3\cdot2^{s-3} ) -  (3\cdot2^{2s+l+1} + 21\cdot2^{3s-4}  -  27\cdot2^{2s-3} )\cdot2^{-s} \nonumber \\
&  -   \sum_{j=0}^{l-1}21\cdot[2^{2s+3j+l+1} + 2^{3s+3j-1} - 5\cdot2^{2s+4j-1}]\cdot2^{-(s+j+1)}.   \nonumber                                
\end{align}\vspace{0.01 cm}\\
Hence by \eqref{eq 16.33}, \eqref{eq 16.34} we deduce after some calculations \vspace{0.01 cm}\\
\begin{align}
\Gamma _{s+l+1}^{\Big[\substack{s  \\ s }\Big] \times (s+l+2)} & = 21\cdot2^{3s +3l-1}- 21\cdot 2^{2s +4l-1}, \label{eq 16.35}\\
& \nonumber \\
\Gamma _{s+ l+2}^{\Big[\substack{s  \\ s }\Big] \times (s+ l+2)} &  = 2^{4s+2l+2} -3\cdot2^{3s+ 3l+2}+ 2^{2s+ 4l+3}.\label{eq 16.36}\\
& \nonumber
\end{align}
By \eqref{eq 12.32} with $ j = l+1,\;   l+1\leq s-1 $ we get \vspace{0.01 cm}\\
\begin{equation}
\label{eq 16.37}
\Gamma _{s+ l+1}^{\Big[\substack{s  \\ s }\Big] \times (k+1)} - \Gamma _{s+ l+1}^{\Big[\substack{s  \\ s }\Big] \times k }= 21\cdot2^{k+s +3(l+1)-4}\quad \text{if}\quad k > s + l+1.
\end{equation}

From \eqref{eq 16.37}, \eqref{eq 16.35} we deduce \vspace{0.01 cm}\\
\begin{align}
& \sum_{j = s+l+2}^{k}\left( \Gamma _{s+l+1}^{\Big[\substack{s  \\ s }\Big] \times (j+1)} - \Gamma _{s+l+1}^{\Big[\substack{s  \\ s }\Big] \times j}   \right)  = \sum_{j = s+l+2}^{k}  21 \cdot 2^{j+s +3l-1} \nonumber \\
& \Leftrightarrow  \sum_{j = s+ l+3}^{k+1}\Gamma _{s+ l+1}^{\Big[\substack{s  \\ s }\Big] \times j} - \sum_{j = s+ l+1}^{k} \Gamma _{s+l+1}^{\Big[\substack{s  \\ s }\Big] \times j} = 21\cdot2^{k+s +3l}- 21\cdot2^{2s+ 4l+1} \nonumber \\
& \Leftrightarrow \Gamma _{s+l+1}^{\Big[\substack{s  \\ s }\Big] \times (k+1)} - \Gamma _{s+l+1}^{\Big[\substack{s  \\ s }\Big] \times (s+ l+2)} = 21\cdot2^{k+s +3l}- 21\cdot2^{2s+4l+1}  \nonumber \\
& \Leftrightarrow \Gamma _{s+ l+1}^{\Big[\substack{s  \\ s }\Big] \times (k+1)}= 21\cdot2^{k+s +3l}- 21\cdot2^{2s+ 4l+1} + 21\cdot(2^{3s +3l-1}-21\cdot2^{2s +4l-1}) \nonumber \\
& \Leftrightarrow \Gamma _{s+l+1 }^{\Big[\substack{s  \\ s }\Big] \times (k+1)}=  21[2^{k+s+3l} + 2^{3s +3l-1} - 5\cdot2^{2s+4l-1}]  \quad \text{ if $ k > s + l+1 $}. \label{eq 16.38}\\
&  \nonumber
\end{align}
By \eqref{eq 16.35} we see that \eqref{eq 16.38} holds for $ k = s+ l+1 $, then  \eqref{eq 16.38} holds for $ k \geq s+ l+1.$\vspace{0.01 cm}\\

From \eqref{eq 16.38} and \eqref{eq 13.4} with $ s\rightarrow s-l, k\rightarrow k-2l > s-l+1 $ we get \vspace{0.01 cm}\\
\begin{align*}
\Gamma _{s+l+1}^{\Big[\substack{s  \\ s }\Big] \times  k} &= 8^{2l}\cdot\Gamma _{s - l+1}^{\Big[\substack{s -l \\ s -l }\Big] \times (k-2l)}=2^{6l}\cdot[21\cdot(2^{k-2l+(s-l)-1} +2^{3(s-l)-1} - 5\cdot2^{2(s-l)-1})] \\
& =  21[2^{k+s+3l-1} + 2^{3s +3l-1} - 5\cdot2^{2s+4l-1}].  \\
 &       \quad   \qed
\end{align*}
We have now established that \vspace{0.01 cm}\\
\begin{equation}
 \Gamma _{s + j+1}^{\Big[\substack{s \\ s  }\Big] \times k}   =  8^{2j}\cdot \Gamma _{s -j+1}^{\Big[\substack{s -j\\ s  -j }\Big] \times (k-2j)}  \quad  \text{if\; $ 0\leq j\leq s-2,\;  k\geq s+ j+1 $}. \label{eq 16.39} \\
\end{equation}\vspace{0.01 cm}\\
\underline{Proof of \eqref{eq 16.16} with $j = s-1,\; k \geq  2s $}\vspace{0.1 cm}\\ 
We shall show that \vspace{0.01 cm}\\
\begin{equation}
\label{eq 16.40}
\Gamma _{2s}^{\Big[\substack{s \\ s  }\Big] \times k}   =  8^{2s-2}\cdot \Gamma _{2}^{\Big[\substack{1 \\ 1 }\Big] \times (k-2(s-1))}  \quad  \text{if $ \; k\geq  2s $}. 
\end{equation}\vspace{0.01 cm}\\
If k =2s,  \eqref{eq 16.40} follows from \eqref{eq 16.21} with $ l = s -1 $ and k = s +(s -1)+1 = 2s.  \vspace{0.01 cm}\\
Let  $ k > 2s.  $  \vspace{0.01 cm}\\
Consider the matrix \\

  $$   \left ( \begin{array} {cccccc}
\alpha _{1} & \alpha _{2} & \alpha _{3} &  \ldots & \alpha _{ 2s}  &  \alpha _{2s +1} \\
\alpha _{2 } & \alpha _{3} & \alpha _{4}&  \ldots  &  \alpha _{2s +1} &  \alpha _{ 2s +2} \\
\vdots & \vdots & \vdots    &  \vdots & \vdots  &  \vdots \\
\alpha _{s-1} & \alpha _{s} & \alpha _{s +1} & \ldots  &  \alpha _{3s -2} &  \alpha _{3s -1}  \\
\alpha _{s} & \alpha _{s+1} & \alpha _{s +2} & \ldots  &  \alpha _{3s -1} &  \alpha _{3s} \\
 \beta  _{1} & \beta  _{2} & \beta  _{3} & \ldots  &  \beta_{ 2s} &  \beta _{ 2s +1}  \\
\beta  _{2} & \beta  _{3} & \beta  _{4} & \ldots  &  \beta_{2s +1} &  \beta _{ 2s +2}  \\
\vdots & \vdots & \vdots    &  \vdots & \vdots  &  \vdots \\
\beta  _{s-1} & \beta  _{s} & \beta  _{s+1} & \ldots  &  \beta_{3s -2}  &  \beta _{3s -1}  \\
\beta  _{s} & \beta  _{s+1} & \beta  _{s+2} & \ldots  &  \beta_{3s -1} &  \beta _{ 3s}
\end{array}  \right). $$ \vspace{0.01 cm}\\
We have \vspace{0.01 cm}\\
  \begin{align}
&  \sum_{ i = 0}^{2s} \Gamma _{i}^{\Big[\substack{s  \\ s }\Big] \times (2s+1)}  = 2^{6s}  && \text{by \eqref{eq 11.1} with m = 0,\; k = 2s +1}, \label{eq 16.41}\\
& \nonumber \\
& \sum_{i = 0}^{s -1} \Gamma _{i}^{\Big[\substack{s  \\ s }\Big] \times (2s +1)}  = 3\cdot2^{3s-4} - 2^{2s-3}&& \text{by \eqref{eq 11.27} with m = 0,\; k =  2s +1},  \label{eq 16.42}\\
& \nonumber \\
&  \Gamma _{s}^{\Big[\substack{s  \\ s }\Big] \times (2s +1)}   = 3\cdot2^{3s} + 21\cdot2^{3s-4}  -  27\cdot2^{2s-3} && \text{by \eqref{eq 13.2} with k = 2s +1},  \label{eq 16.43} \\
& \nonumber \\
&   \Gamma _{s+ j +1}^{\Big[\substack{s  \\ s }\Big] \times (2s +1)}  =  8^{2j}\cdot \Gamma _{s -j+1}^{\Big[\substack{s -j \\ s -j }\Big] \times (2s +1 - 2j)}  && \text{by \eqref{eq 16.20} $\;\text{for}  \;0\leq j\leq s-2 ,\; k >  2s +1 $}  \label{eq 16.44} \\
 & =  2^{6j}\cdot[ 21\cdot(2^{(2s +1-2j)+(s-j) -1}+ 2^{3(s-j)-1} -5\cdot2^{2(s-j)-1})] && \text{by \eqref{eq 13.4} with $ s\rightarrow s-j,\; k\rightarrow 2s+1-2j > s-j+1 $} \nonumber \\
 & = 21\cdot[2^{3s+3j} + 2^{3s+3j-1} - 5\cdot2^{2s+4j-1}]. \nonumber
  \end{align} \vspace{0.01 cm}\\  
From \eqref{eq 16.41}, \eqref{eq 16.42}, \eqref{eq 16.43} and  \eqref{eq 16.44} we get after some calculations  \vspace{0.01 cm}\\  
\begin{align}
 \Gamma _{2s}^{\Big[\substack{s  \\ s }\Big] \times (2s+1)}& = 2^{6s}-(3\cdot2^{3s-4} - 2^{2s-3})
 -( 3\cdot2^{3s} + 21\cdot2^{3s-4}  -  27\cdot2^{2s-3} )\label{eq 16.45} \\
 &  - \sum_{j=0}^{s-2}21\cdot[2^{3s+3j} + 2^{3s+3j-1} - 5\cdot2^{2s+4j-1}]  \nonumber \\
 & =  21\cdot2^{6s-5}.  \nonumber
\end{align}\vspace{0.01 cm}\\  
By \eqref{eq 12.32} with j = s, \; $  k > 2s $ we get \vspace{0.01 cm}\\
\begin{equation}
\label{eq 16.46}
\Gamma _{2s}^{\Big[\substack{s  \\ s }\Big] \times (k+1)} - \Gamma _{2s}^{\Big[\substack{s  \\ s }\Big] \times k }= 3\cdot2^{2k +2s -2} -3\cdot 2^{k+4s-4} \quad \text{if}\quad k > 2s.
\end{equation}
From \eqref{eq 16.46}, \eqref{eq 16.45} we deduce \vspace{0.01 cm}\\
\begin{align}
& \sum_{j = 2s+1}^{k}\left( \Gamma _{2s}^{\Big[\substack{s  \\ s }\Big] \times (j+1)} - \Gamma _{2s}^{\Big[\substack{s  \\ s }\Big] \times j}   \right)  = \sum_{j = 2s+1}^{k} ( 3\cdot2^{2j +2s -2} -3\cdot 2^{j+4s-4} )  \nonumber \\
& \Leftrightarrow  \sum_{j = 2s+2}^{k+1}\Gamma _{2s}^{\Big[\substack{s  \\ s }\Big] \times j} - \sum_{j = 2s +1}^{k} \Gamma _{2s}^{\Big[\substack{s  \\ s }\Big] \times j} = 2^{2k+2s} -3\cdot2^{k+4s-3}  -5\cdot2^{6s-3} \nonumber \\
& \Leftrightarrow \Gamma _{2s}^{\Big[\substack{s  \\ s }\Big] \times (k+1)} - \Gamma _{2s}^{\Big[\substack{s  \\ s }\Big] \times (2s+1)} =  2^{2k+2s} -3\cdot2^{k+4s-3}  -5\cdot2^{6s-3}  \nonumber \\
& \Leftrightarrow \Gamma _{2s}^{\Big[\substack{s  \\ s }\Big] \times (k+1)}=  2^{2k+2s} -3\cdot2^{k+4s-3}  -5\cdot2^{6s-3} +   21\cdot2^{6s-5} \nonumber \\
& \Leftrightarrow \Gamma _{2s }^{\Big[\substack{s  \\ s }\Big] \times (k+1)}=   2^{2k+2s} -3\cdot2^{k+4s-3}+ 2^{6s-5}  \quad \text{ if $ k >  2s $}. \label{eq 16.47}
\end{align}\vspace{0.01 cm}\\
By \eqref{eq 16.45} we see that \eqref{eq 16.47} holds for k =  2s,    then  \eqref{eq 16.47} holds for $  k \geq  2s. $\vspace{0.01 cm}\\
Consider the matrix\\

  $$   \left ( \begin{array} {cccccc}
\alpha _{1} & \alpha _{2} & \alpha _{3} &  \ldots & \alpha _{ k-2s +1}  &  \alpha _{k-2s+2} \\
 \beta  _{1} & \beta  _{2} & \beta  _{3} & \ldots  &  \beta_{ k-2s +1} &  \beta _{ k-2s +2}  
\end{array}  \right). $$ \vspace{0.01 cm}\\
We see easily that  $$ \Gamma _{1}^{\Big[\substack{1 \\ 1 }\Big] \times (k-2(s-1))} = 3\cdot(2^{k-2s+2} -1).  $$ \vspace{0.01 cm}\\
Hence
\begin{equation}
\label{eq 16.48}
 \Gamma _{2}^{\Big[\substack{1 \\ 1 }\Big] \times (k-2(s-1))} = 2^{2k -4s +4} - 3\cdot(2^{k-2s+2} - 1) -1 =  2^{2k -4s +4} - 3\cdot2^{k-2s+2} +2.  
\end{equation}\vspace{0.1 cm}\\
From \eqref{eq 16.48} and \eqref{eq 16.47} with $ k\geq 2s  $ we get \vspace{0.1 cm}\\
$$ \Gamma _{2s }^{\Big[\substack{s  \\ s }\Big] \times k} =  2^{2k+2s-2} -3\cdot2^{k+4s-4}+ 2^{6s-5} = 8^{2(s-1)}\cdot \Gamma _{2}^{\Big[\substack{1 \\ 1 }\Big] \times (k-2(s-1))}. $$
\end{proof}

 \section{\textbf{ A REDUCTION FORMULA FOR  $\Gamma_{s+m+1+j}^{\left[s\atop s+m\right]\times k}$ IN THE CASE  $ 1\leq j\leq s-1,\;k\geq s+m+1+j $}}
\label{sec 17}
In this section we prove the following reduction formula by induction on j \vspace{0.05 cm}\\
\begin{align*}
 \Gamma _{s +m +1 +j}^{\Big[\substack{s \\ s +m }\Big] \times k} &  =  8^{2j + m}\cdot \Gamma _{s-j+1}^{\Big[\substack{s -j\\ s -j }\Big] \times (k- m-2j)}  &&\text{if   \; $ 0\leq j\leq s-1,\;  k\geq s + m+1+j $}. 
\end{align*}\vspace{0.05 cm}\\
We recall again that the right hand side in the above equation has been computed in section \ref{sec 13}.\\
In fact we have \vspace{0.05 cm}\\
\begin{equation*}
 \Gamma _{s -j+1}^{\Big[\substack{s -j\\ s  -j }\Big] \times (k- m-2j)} =                                                                  
  \begin{cases}
 2^{4s -4j} - 3\cdot2^{3s -3j-1} + 2^{2s-2j -1} & \text{if   }\quad 0 \leq j \leq s-1, \; k = s+m+j+1, \\
  21\cdot[2^{k-m +s -3j -1} + 2^{3s -3j-1} -5\cdot2^{2s -2j-1}] & \text{if   }\quad 0 \leq j \leq s-2, \; k > s+m+j+1, \\
  2^{2k-2m-4s+4} - 3\cdot2^{k-m-2s +2} + 2 & \text{if   }\quad  j = s-1,\;k >2s +m.
\end{cases}
\end{equation*}

\begin{lem}
\label{lem 17.1}We have 
\begin{align}
 \Gamma _{s +m +1}^{\Big[\substack{s \\ s +m }\Big] \times k} &  =  8^{m}\cdot \Gamma _{s+1}^{\Big[\substack{s \\ s  }\Big] \times (k- m)}  &&\text{if   \; $  k\geq s + m+1 $}, \label{eq 17.1} \\
 & \nonumber \\
  \Gamma _{s+ m +1}^{\Big[\substack{s \\ s + m }\Big] \times (s+m+1)} & =  2^{4s+3m} - 3\cdot2^{3s  +3m - 1} + 2^{2s +3m - 1}=  8^{m}\cdot \Gamma _{s+1}^{\Big[\substack{s \\ s  }\Big] \times ( s+1)}  &&\text{if   \;  $ k = s + m+1 $}, \label{eq 17.2}\\
  & \nonumber \\
  \Gamma _{s + 1}^{\Big[\substack{s \\ s  }\Big] \times (s+1)} &  = 2^{4s} -3\cdot2^{3s-1} + 2^{2s-1},    \nonumber  \\
  \Gamma _{s+ m +1}^{\Big[\substack{s \\ s + m }\Big] \times  k} & = 21\cdot(2^{k+s+2m-1} + 2^{3s+3m-1} - 5\cdot2^{2s+3m-1})  = 8^{m}\cdot \Gamma _{s+1}^{\Big[\substack{s \\ s  }\Big] \times ( k - m)}  && \text{if  \; $k >s+ m +1 $}, \label{eq 17.3} \\
  \Gamma _{s + 1}^{\Big[\substack{s \\ s  }\Big] \times (k-m)} & =  21\cdot(2^{k -m +s-1} + 2^{3s -1} - 5\cdot2^{2s -1}). \nonumber \\
  & \nonumber 
  \end{align}
\end{lem}
\begin{proof}
\underline{Proof of \eqref{eq 17.2} }\vspace{0.1 cm}\\
By \eqref{eq 15.32} and \eqref{eq 13.3} we get  \vspace{0.1 cm}\\
$$ \Gamma _{s+ m +1}^{\Big[\substack{s \\ s + m }\Big] \times (s+m+1)}  =  2^{4s+3m} - 3\cdot2^{3s  +3m - 1} + 2^{2s +3m - 1}=  8^{m}\cdot \Gamma _{s+1}^{\Big[\substack{s \\ s  }\Big] \times ( s+1)}.  $$

\underline{Proof of \eqref{eq 17.3} }\vspace{0.1 cm}\\
Consider the matrix \\

  $$   \left ( \begin{array} {cccccc}
\alpha _{1} & \alpha _{2} & \alpha _{3} &  \ldots & \alpha _{s+ m +1}  &  \alpha _{s+m+2} \\
\alpha _{2 } & \alpha _{3} & \alpha _{4}&  \ldots  &  \alpha _{s+ m+2} &  \alpha _{s+ m+3} \\
\vdots & \vdots & \vdots    &  \vdots & \vdots  &  \vdots \\
\alpha _{s-1} & \alpha _{s} & \alpha _{s +1} & \ldots  &  \alpha _{2s + m -1} &  \alpha _{2s + m }  \\
\alpha _{s} & \alpha _{s+1} & \alpha _{s +2} & \ldots  &  \alpha _{2s +m } &  \alpha _{2s+ m+ 1} \\
 \beta  _{1} & \beta  _{2} & \beta  _{3} & \ldots  &  \beta_{s+ m +1} &  \beta _{s+ m+2}  \\
\beta  _{2} & \beta  _{3} & \beta  _{4} & \ldots  &  \beta_{s+ m+2} &  \beta _{s+ m +3}  \\
\vdots & \vdots & \vdots    &  \vdots & \vdots  &  \vdots \\
\beta  _{s-1} & \beta  _{s} & \beta  _{s+1} & \ldots  &  \beta_{2s + m-1} &  \beta _{2s +m}  \\
\beta  _{s} & \beta  _{s+1} & \beta  _{s+2} & \ldots  &  \beta_{2s + m} &  \beta _{2s + m +1}\\
\beta  _{s+1} & \beta  _{s+2} & \beta  _{s+3} & \ldots  &  \beta_{2s + m +1} &  \beta _{2s  +m+2}\\
\vdots & \vdots & \vdots    &  \vdots & \vdots  &  \vdots \\
\beta  _{s+m} & \beta  _{s+m+1} & \beta  _{s+m+2} & \ldots  &  \beta_{2s+ 2m} &  \beta _{2s + 2m +1}
\end{array}  \right). $$ \vspace{0.01 cm}\\
We have \vspace{0.01 cm}\\
  \begin{align}
&  \sum_{ i = 0}^{s+ m+2} \Gamma _{i}^{\Big[\substack{s  \\ s +m}\Big] \times (s+ m+2)} = 2^{4s + 3m +2} &&   \text{by \eqref{eq 11.1} with \; k = s+ m +2}, \label{eq 17.4}\\
& \nonumber \\
&  \sum_{i = 0}^{s+ m +2} \Gamma _{i}^{\Big[\substack{s  \\ s +m}\Big] \times (s+ m+2)}\cdot2^{-i}  =  2^{3s  + m} + 2^{2s +2m+2} - 2^{s + m} &&  \text{by \eqref{eq 11.2} with \; k = s+ m+2},  \label{eq 17.5}\\
& \nonumber \\
 &  \sum_{i = 0}^{s -1} \Gamma _{i}^{\Big[\substack{s  \\ s +m}\Big] \times (s+ m+2)}  = 3\cdot2^{3s-4} - 2^{2s-3} &&  \text{by \eqref{eq 11.27} with \; k = s+ m+2},  \label{eq 17.6}\\
 & \nonumber \\
 &   \sum_{i = 0}^{s -1} \Gamma _{i}^{\Big[\substack{s  \\ s +m }\Big] \times (s+ m+2)}\cdot2^{-i}  =  7\cdot2^{2s-4} - 3\cdot2^{s-3} &&    \text{by \eqref{eq 11.28} with \; k = s+m+2},  \label{eq 17.7}\\
 & \nonumber \\
 &   \Gamma _{s}^{\Big[\substack{s  \\ s +m }\Big] \times (s+ m+2)}   = 2^{2s + m+1} + 21\cdot2^{3s-4}  -  11\cdot2^{2s-3}  &&   \text{by \eqref{eq 13.46} with k = s+ m+2},  \label{eq 17.8} \\
 & \nonumber \\
  &   \Gamma _{s+ j }^{\Big[\substack{s  \\ s +m}\Big] \times (s+ m+ 2)}  =   8^{j-1}\cdot \Gamma _{s+1}^{\Big[\substack{s \\ s +(m-(j-1)) }\Big] \times (s+m+2 -(j-1))}  && \text{by \eqref{eq 15.40} $\;\text{for}  \;0\leq j\leq m-1 $}  \label{eq 17.9} \\
  & = 3\cdot2^{2s+m+2j-1} +21\cdot2^{3s+3j-4} - 21\cdot2^{2s+3j-4}  && \text{and \; k = s+ m+2 },  \nonumber \\
  & \nonumber \\
  &   \Gamma _{s+ m }^{\Big[\substack{s  \\ s +m}\Big] \times (s+ m+ 2)} = 11\cdot2^{(s+m+2) +2m +s -3} +21\cdot2^{3s+3m-4} - 53\cdot2^{2s+3m-4}  && \text{by \eqref{eq 15.41} $\;\text{for}  \; k = s+ m+2 $}  \label{eq 17.10} \\
 & = 21\cdot2^{3s+3m-4} + 35\cdot2^{2s+3m-4}.\nonumber
  \end{align} \vspace{0.1 cm}\\  
   From $\eqref{eq 17.4},\ldots, \eqref{eq 17.10} $ with k = s+m+2  we obtain  \vspace{0.01 cm}\\
  \begin{align}
& \Gamma _{s+ m+1}^{\Big[\substack{s  \\ s +m}\Big] \times (s+m+2)} + \Gamma _{s+m+2}^{\Big[\substack{s  \\ s +m}\Big] \times (s+ m+2)} \label{eq 17.11}  \\
& = 2^{4s+3m +2} - (3\cdot2^{3s-4} - 2^{2s-3}) - (2^{2s + m+1} + 21\cdot2^{3s-4}  -  11\cdot2^{2s-3}) \nonumber \\
& - \sum_{j=1}^{m -1} [ 3\cdot2^{2s+m+2j-1} +21\cdot2^{3s+3j-4} - 21\cdot2^{2s+3j-4} ]- ( 21\cdot2^{3s+3m-4} + 35\cdot2^{2s+3m-4}) \nonumber \\
& = 2^{4s+3m+2}  -3\cdot2^{3s+3m-1} -5\cdot2^{2s+3m-1},\nonumber \\
& \nonumber \\
& \Gamma _{s+ m+1}^{\Big[\substack{s  \\ s +m}\Big] \times (s+ m+2)}\cdot2^{-(s+m+1)} + \Gamma _{s+m+2}^{\Big[\substack{s  \\ s +m }\Big] \times (s+m+2)}\cdot 2^{-(s+m+2)}  \label{eq 17.12}  \\
& =   2^{3s +m} + 2^{2s+2m+2} - 2^{s+m} -(7\cdot2^{2s-4} - 3\cdot2^{s-3} ) -  ( 2^{2s + m+1} + 21\cdot2^{3s-4}  -  11\cdot2^{2s-3}  )\cdot2^{-s} \nonumber \\
&  -   \sum_{j= 1}^{m-1}[3\cdot2^{2s+m+2j-1} +21\cdot2^{3s+3j-4} - 21\cdot2^{2s+3j-4}]\cdot2^{-(s+j)} \nonumber  \\
&  - ( 21\cdot2^{3s+3m-4} + 35\cdot2^{2s+3m-4})\cdot2^{-(s+m)}  \nonumber  \\
& =     2^{3s+2m }   + 9\cdot2^{2s+2m - 2} - 13\cdot2^{s+2m - 2}.\nonumber                           
\end{align}\vspace{0.01 cm}\\ 
Hence by \eqref{eq 17.11}, \eqref{eq 17.12} we deduce after some calculations \vspace{0.01 cm}\\
\begin{align}
\Gamma _{s+m+1}^{\Big[\substack{s  \\ s +m}\Big] \times (s+m+2)} &  = 21\cdot2^{3s +3m -1}- 21\cdot 2^{2s + 3m -1}, \label{eq 17.13}\\
& \nonumber \\
\Gamma _{s+ m+2}^{\Big[\substack{s  \\ s +m}\Big] \times (s+ m+2)} &  = 2^{4s+3m +2} -3\cdot2^{3s+ 3m+2}+ 2^{2s+ 3m +3}.\label{eq 17.14}\\
& \nonumber 
\end{align}
By \eqref{eq 12.35} with j = 1,   $ k > s+m+1 $ we get \vspace{0.01 cm}\\
\begin{equation}
\label{eq 17.15}
\Gamma _{s+ m+1}^{\Big[\substack{s  \\ s +m}\Big] \times (k+1)} - \Gamma _{s+ m+1}^{\Big[\substack{s  \\ s +m}\Big] \times k }= 21\cdot2^{k+s + 2m -1}\quad \text{if}\quad k > s + m+1.
\end{equation}

From \eqref{eq 17.15}, \eqref{eq 17.13} we deduce  \vspace{0.01 cm}\\
\begin{align}
& \sum_{j = s+m+2}^{k}\left( \Gamma _{s+m+1}^{\Big[\substack{s  \\ s +m}\Big] \times (j+1)} - \Gamma _{s+m+1}^{\Big[\substack{s  \\ s +m }\Big] \times j}   \right)  = \sum_{j = s+m+2}^{k}  21 \cdot 2^{j+s +2m -1} \nonumber \\
& \Leftrightarrow  \sum_{j = s+ m+3}^{k+1}\Gamma _{s+ m+1}^{\Big[\substack{s  \\ s +m}\Big] \times j} - \sum_{j = s+ m+1}^{k} \Gamma _{s+m+1}^{\Big[\substack{s  \\ s +m}\Big] \times j} = 21\cdot2^{k+s + 2m}- 21\cdot2^{2s+ 3m +1} \nonumber \\
& \Leftrightarrow \Gamma _{s+m+1}^{\Big[\substack{s  \\ s +m}\Big] \times (k+1)} - \Gamma _{s+m+1}^{\Big[\substack{s  \\ s +m}\Big] \times (s+ m+2)} = 21\cdot2^{k+s + 2m}- 21\cdot2^{2s+3m+1}  \nonumber \\
& \Leftrightarrow \Gamma _{s+ m +1}^{\Big[\substack{s  \\ s +m}\Big] \times (k+1)}= 21\cdot2^{k+s + 2m}- 21\cdot2^{2s+ 3m +1} + 21\cdot(2^{3s +3m-1}-21\cdot2^{2s + 3m -1}) \nonumber \\
& \Leftrightarrow \Gamma _{s+m +1 }^{\Big[\substack{s  \\ s +m}\Big] \times (k+1)}=  21[2^{k+s+ 2m} + 2^{3s +3m-1} - 5\cdot2^{2s+ 3m -1}]  \quad \text{ if $ k > s + m+1 $}. \label{eq 17.16}
\end{align}\vspace{0.01 cm}\\
By \eqref{eq 17.13} we see that \eqref{eq 17.16} holds for k = s+ m+1, then  \eqref{eq 17.16} holds for $ k \geq s+ m+1.$\vspace{0.01 cm}\\

From \eqref{eq 17.16} and \eqref{eq 13.4} with $  k\rightarrow k - m > s +1  $ we get \vspace{0.01 cm}\\
\begin{align*}
\Gamma _{s+m+1}^{\Big[\substack{s  \\ s +m}\Big] \times  k} &= 8^{m}\cdot\Gamma _{s +1}^{\Big[\substack{s  \\ s  }\Big] \times (k- m)}=2^{3m}\cdot[21\cdot(2^{k- m + s -1} +2^{3s -1} - 5\cdot2^{2s -1})] \\
& =  21[2^{k+s+2m -1} + 2^{3s +3m-1} - 5\cdot2^{2s+ 3m -1}].                           
\end{align*}
\end{proof}

\begin{lem}
\label{lem 17.2}We have 
\begin{align}
 \Gamma _{s +m +1 +j}^{\Big[\substack{s \\ s +m }\Big] \times k} &  =  8^{2j + m}\cdot \Gamma _{s-j+1}^{\Big[\substack{s -j\\ s -j }\Big] \times (k- m-2j)}  &&\text{if   \; $ 0\leq j\leq s-1,\;  k\geq s + m+1+j $}, \label{eq 17.17} \\
 & \nonumber \\
  \Gamma _{s+ m +1+j}^{\Big[\substack{s \\ s + m }\Big] \times (s+m+1+j)} & =   8^{2j + m}\cdot \Gamma _{s - j+1}^{\Big[\substack{s -j\\ s -j }\Big] \times ( s - j+1)} &&\text{if   \; $ 0\leq j\leq s-1,\;  k =  s + m+1+j $} \label{eq 17.18}  \\
 & = 8^{2j + m}\cdot(2^{4s-4j} -3\cdot2^{3s-3j-1} + 2^{2s-2j-1} ) \nonumber \\
&   = 2^{4s +3m  +2j} -3\cdot2^{3s+3m+3j-1} + 2^{2s +3m+ 4j-1},\nonumber \\
& \nonumber \\
 \Gamma _{s+ m +1+j }^{\Big[\substack{s \\ s + m }\Big] \times  k} & =   8^{2j + m}\cdot \Gamma _{s-j+1}^{\Big[\substack{s -j\\ s -j }\Big] \times (k- m-2j)}  \label{eq 17.19}  \\
   & =8^{2j +m}\cdot 21\cdot[2^{k-m -3j+s-1} + 2^{3s-3j-1} -5\cdot2^{2s-2j-1}] &&\text{ if $ \; 0\leq j\leq s-2,\; k > s+m+ 1+j $} \nonumber\\
   & = 21\cdot[2^{k +2m  + 3j+s-1} + 2^{3s +3m +3j -1} -5\cdot2^{2s +3m +4j -1}], \nonumber \\
   & \nonumber \\
 \Gamma _{2s +m }^{\Big[\substack{s \\ s +m }\Big] \times k} &  = 8^{2s+m-2}\cdot \Gamma _{2}^{\Big[\substack{1 \\ 1  }\Big] \times (k-m -2s +2)} &&\text{if   \; $ j =  s-1,\;  k > 2s + m $} \label{eq 17.20} \\
  & =  8^{2s+m-2}\cdot[2^{2(k-m) -4s +4} -3\cdot2^{k-m -2s +2} +2] \nonumber  \\
  & = 2^{2k+2s+m-2} -3\cdot2^{k+2m+4s -4} +2^{6s+3m-5}.\nonumber \\
  & \nonumber 
   \end{align}
\end{lem}
\begin{proof}We will do the proof by induction on j.\vspace{0.1 cm}\\
Let $l$ be a rational integer such that $ 1\leq l\leq s-2. $\vspace{0.1 cm}\\
Let $ (H_{l-1}) $ denote the following statement \vspace{0.1 cm}\\
\begin{align}
\Gamma _{s +m +1 +j}^{\Big[\substack{s \\ s +m }\Big] \times k} &  =  8^{2j + m}\cdot \Gamma _{s-j+1}^{\Big[\substack{s -j\\ s -j }\Big] \times (k- m-2j)}  &&\text{if   \; $ 0\leq j\leq l -1,\;  k\geq s + m+1+j $}. \label{eq 17.21}\\
 & \nonumber
\end{align}
By Lemma \ref{lem 17.1}  $ (H_{l-1}) $ holds for $ l = 1. $\vspace{0.1 cm}\\
Assume that  $ (H_{l-1}) $ holds. \vspace{0.1 cm}\\
We are going to show that $(H_{l})$ holds, that is \vspace{0.1 cm}\\
\begin{align}
\Gamma _{s +m +1 +l}^{\Big[\substack{s \\ s +m }\Big] \times k} &  =  8^{2l + m}\cdot \Gamma _{s-l+1}^{\Big[\substack{s -l\\ s -l }\Big] \times (k- m-2l)}  &&\text{if   \; $  k\geq s + m+l+1 $}.\label{eq 17.22}\\
& \nonumber
\end{align}
\underline{The case $ k = s+m+l+1 $}\vspace{0.1 cm}\\
Consider the matrix\\

  $$   \left ( \begin{array} {cccccc}
\alpha _{1} & \alpha _{2} & \alpha _{3} &  \ldots & \alpha _{s+ m +l}  &  \alpha _{s+m+ l+1} \\
\alpha _{2 } & \alpha _{3} & \alpha _{4}&  \ldots  &  \alpha _{s+ m+l+1} &  \alpha _{s+ m+ l+2} \\
\vdots & \vdots & \vdots    &  \vdots & \vdots  &  \vdots \\
\alpha _{s-1} & \alpha _{s} & \alpha _{s +1} & \ldots  &  \alpha _{2s + m  +l-2} &  \alpha _{2s + m +l-1}  \\
\alpha _{s} & \alpha _{s+1} & \alpha _{s +2} & \ldots  &  \alpha _{2s +m +l-1} &  \alpha _{2s+ m+ l} \\
 \beta  _{1} & \beta  _{2} & \beta  _{3} & \ldots  &  \beta_{s+ m + l} &  \beta _{s+ m+ l+1}  \\
\beta  _{2} & \beta  _{3} & \beta  _{4} & \ldots  &  \beta_{s+ m+l+1} &  \beta _{s+ m +l+2}  \\
\vdots & \vdots & \vdots    &  \vdots & \vdots  &  \vdots \\
\beta  _{s-1} & \beta  _{s} & \beta  _{s+1} & \ldots  &  \beta_{2s + m+ l-2} &  \beta _{2s +m +l-1}  \\
\beta  _{s} & \beta  _{s+1} & \beta  _{s+2} & \ldots  &  \beta_{2s + m +l-1} &  \beta _{2s + m + l}\\
\beta  _{s+1} & \beta  _{s+2} & \beta  _{s+3} & \ldots  &  \beta_{2s + m + l} &  \beta _{2s  +m+ l+1}\\
\vdots & \vdots & \vdots    &  \vdots & \vdots  &  \vdots \\
\beta  _{s+m} & \beta  _{s+m+1} & \beta  _{s+m+2} & \ldots  &  \beta_{2s+ 2m +l-1} &  \beta _{2s + 2m +l}
\end{array}  \right). $$ \vspace{0.01 cm}\\
We have \vspace{0.01 cm}\\
  \begin{align}
&  \sum_{ i = 0}^{s+ m+l+1} \Gamma _{i}^{\Big[\substack{s  \\ s +m}\Big] \times (s+ m+ l+1)} = 2^{4s + 3m +2l} &&   \text{by \eqref{eq 11.1} with \; $ k = s+ m +l+1 $}, \label{eq 17.23}\\
& \nonumber \\
 &  \sum_{i = 0}^{s -1} \Gamma _{i}^{\Big[\substack{s  \\ s +m}\Big] \times (s+ m+ l+1)}  = 3\cdot2^{3s-4} - 2^{2s-3} &&  \text{by \eqref{eq 11.27} with \;  $ k = s+ m+l+1 $},  \label{eq 17.24}\\
 & \nonumber \\
  &   \Gamma _{s}^{\Big[\substack{s  \\ s +m }\Big] \times (s+ m+ l+1)}   = 2^{2s + m+l} + 21\cdot2^{3s-4}  -  11\cdot2^{2s-3}  &&   \text{by \eqref{eq 13.46} with $ k = s+ m+l+1 $},  \label{eq 17.25} \\
  & \nonumber \\
   &   \Gamma _{s+ j }^{\Big[\substack{s  \\ s +m}\Big] \times (s+ m+ l+1)}  =   8^{j-1}\cdot \Gamma _{s+1}^{\Big[\substack{s \\ s +(m-(j-1)) }\Big] \times (s+m+ l+1 -(j-1))}  && \text{by \eqref{eq 15.40} $\;\text{for}   \;1\leq j\leq m-1 $}  \label{eq 17.26} \\
  & = 3\cdot2^{2s+m+l +2j- 2} +21\cdot2^{3s+3j-4} - 21\cdot2^{2s+3j-4}  && \text{and \;$ k = s+ m+ l+1$ },  \nonumber \\
  & \nonumber \\
   &   \Gamma _{s+ m }^{\Big[\substack{s  \\ s +m}\Big] \times (s+ m+ l+1)} = 11\cdot2^{(s+m+l+1) +2m +s -3} +21\cdot2^{3s+3m-4} - 53\cdot2^{2s+3m-4}  && \text{by \eqref{eq 15.41} $\;\text{for}  \; k = s+ m+2 $}  \label{eq 17.27} \\
 & =   11\cdot2^{2s+3m+l-2} +21\cdot2^{3s+3m-4} - 53\cdot2^{2s+3m-4},    \nonumber\\
 & \nonumber \\
  &  \Gamma _{s+ m +1+j }^{\Big[\substack{s \\ s + m }\Big] \times s+m+l+1  }  =   8^{2j + m}\cdot \Gamma _{s-j+1}^{\Big[\substack{s -j\\ s -j }\Big] \times ((s+m+l+1)- m-2j)}    && \text{by $ (H_{l-1}) $ ( see   \eqref{eq 17.21}) }\label{eq 17.28} \\
   & =8^{2j +m}\cdot 21\cdot[2^{(s+m+l+1) -m -3j+s-1} + 2^{3s-3j-1} -5\cdot2^{2s-2j-1}] &&\text{ if $ \; 0\leq j\leq l-1,\; k = s+m+l+1 $} \nonumber\\
   & = 21\cdot[2^{2s +3m +l+3j} + 2^{3s +3m +3j -1}  -5\cdot2^{2s +3m +4j -1}]. \nonumber \\
   & \nonumber 
  \end{align} \vspace{0.1 cm} 
   From $\eqref{eq 17.23}, \ldots,\eqref{eq 17.28}\; with\; k = s+m+l+1 $ we obtain after some calculations  \vspace{0.1 cm} \\
   \begin{align}
& \Gamma _{s+m+l+1}^{\Big[\substack{s  \\ s +m}\Big] \times (s+ m+ l+1)}=  2^{4s + 3m +2l} - (3\cdot2^{3s-4} - 2^{2s-3} )
 - ( 2^{2s + m+l} + 21\cdot2^{3s-4}  -  11\cdot2^{2s-3} ) \label{eq 17.29}\\
 & - \sum_{j = 1}^{m - 1} [3\cdot2^{2s+m+l +2j- 2} +21\cdot2^{3s+3j-4} - 21\cdot2^{2s+3j-4}] -( 11\cdot2^{2s+3m+l-2} +21\cdot2^{3s+3m-4} - 53\cdot2^{2s+3m-4}  ) \nonumber\\
 & - \sum_{j = 0}^{l-1} 21\cdot[2^{2s +3m +l+3j} + 2^{3s +3m +3j -1}  -5\cdot2^{2s +3m +4j -1}] \nonumber \\
 & = 2^{4s+3m+2l} -3\cdot2^{3s+3m+3l-1} +2^{2s+3m+4l-1}.\nonumber 
\end{align}
By \eqref{eq 13.3} with $ s\rightarrow s-l $ and \eqref{eq 17.29} we get \vspace{0.1 cm} \\
\begin{align*}
& \Gamma _{s+m+l+1}^{\Big[\substack{s  \\ s +m}\Big] \times (s+ m+ l+1)}=8^{2l+m}\cdot[2^{4(s-l)} - 3\cdot2^{3(s-l)-1} + 2^{2(s-l)-1}]
= 2^{6l+3m}\cdot\Gamma _{s - l +1}^{\Big[\substack{s -l  \\ s  - l}\Big] \times (s - l+1)}. 
\end{align*}
 \underline{The case $ k >  s+m+l+1 $}\vspace{0.1 cm}\\
Consider the matrix
 $$   \left ( \begin{array} {cccccc}
\alpha _{1} & \alpha _{2} & \alpha _{3} &  \ldots & \alpha _{s+ m +l +1}  &  \alpha _{s+m+ l+2} \\
\alpha _{2 } & \alpha _{3} & \alpha _{4}&  \ldots  &  \alpha _{s+ m+l+2} &  \alpha _{s+ m+ l+3} \\
\vdots & \vdots & \vdots    &  \vdots & \vdots  &  \vdots \\
\alpha _{s-1} & \alpha _{s} & \alpha _{s +1} & \ldots  &  \alpha _{2s + m  +l- 1} &  \alpha _{2s + m +l}  \\
\alpha _{s} & \alpha _{s+1} & \alpha _{s +2} & \ldots  &  \alpha _{2s +m +l} &  \alpha _{2s+ m+ l+1} \\
 \beta  _{1} & \beta  _{2} & \beta  _{3} & \ldots  &  \beta_{s+ m + l +1} &  \beta _{s+ m+ l+2}  \\
\beta  _{2} & \beta  _{3} & \beta  _{4} & \ldots  &  \beta_{s+ m+l+2} &  \beta _{s+ m +l+ 3}  \\
\vdots & \vdots & \vdots    &  \vdots & \vdots  &  \vdots \\
\beta  _{s-1} & \beta  _{s} & \beta  _{s+1} & \ldots  &  \beta_{2s + m+ l- 1} &  \beta _{2s +m +l}  \\
\beta  _{s} & \beta  _{s+1} & \beta  _{s+2} & \ldots  &  \beta_{2s + m +l} &  \beta _{2s + m + l +1}\\
\beta  _{s+1} & \beta  _{s+2} & \beta  _{s+3} & \ldots  &  \beta_{2s + m + l+1} &  \beta _{2s  +m+ l+2}\\
\vdots & \vdots & \vdots    &  \vdots & \vdots  &  \vdots \\
\beta  _{s+m} & \beta  _{s+m+1} & \beta  _{s+m+2} & \ldots  &  \beta_{2s+ 2m +l} &  \beta _{2s + 2m +l +1}
\end{array}  \right). $$ \vspace{0.01 cm}\\
We have \\
 \begin{align}
&  \sum_{ i = 0}^{s+ m+l+2} \Gamma _{i}^{\Big[\substack{s  \\ s +m}\Big] \times (s+ m+ l+2)} = 2^{4s + 3m +2l +2} &&   \text{by \eqref{eq 11.1}\; with \;$ k = s+m+l+2 $}, \label{eq 17.30}\\
& \nonumber \\
&  \sum_{i = 0}^{s+ m + l +2} \Gamma _{i}^{\Big[\substack{s  \\ s +m}\Big] \times (s+ m+ l +2)}\cdot2^{-i}  =  2^{3s  + 2m + l} + 2^{2s +2m+2l +2} - 2^{s + m + l} &&  \text{by \eqref{eq 11.2}\; with \;$ k = s+m+l+2 $},  \label{eq 17.31}\\
& \nonumber \\
&   \sum_{i = 0}^{s -1} \Gamma _{i}^{\Big[\substack{s  \\ s +m }\Big] \times (s+ m+ l+2)}\cdot2^{-i}  =  7\cdot2^{2s-4} - 3\cdot2^{s-3} &&    \text{by \eqref{eq 11.28} with \; $ k = s+m+l+2 $},  \label{eq 17.32}\\
& \nonumber \\
 &  \sum_{i = 0}^{s -1} \Gamma _{i}^{\Big[\substack{s  \\ s +m}\Big] \times (s+ m+ l+2)}  = 3\cdot2^{3s-4} - 2^{2s-3} &&  \text{by \eqref{eq 11.27} with \;$ k = s+m+l+2 $},  \label{eq 17.33}\\
 & \nonumber \\
  &   \Gamma _{s}^{\Big[\substack{s  \\ s +m }\Big] \times (s+ m+ l+2)}   = 2^{2s + m+l +1} + 21\cdot2^{3s-4}  -  11\cdot2^{2s-3}  &&   \text{by \eqref{eq 13.46} with $ k = s+m+l+2 $},  \label{eq 17.34} \\
  & \nonumber \\
  &   \Gamma _{s+ j }^{\Big[\substack{s  \\ s +m}\Big] \times (s+ m+ l+2)}  =   8^{j-1}\cdot \Gamma _{s+1}^{\Big[\substack{s \\ s +(m-(j-1)) }\Big] \times (s+m+ l+2 -(j-1))}  && \text{by \eqref{eq 15.40} $\;\text{for}   \;1\leq j\leq m-1 $}  \label{eq 17.35} \\
  & = 3\cdot2^{2s+m+l +2j- 1} +21\cdot2^{3s+3j-4} - 21\cdot2^{2s+3j-4}  && \text{and \;$ k = s+m+l+2 $},  \nonumber \\
  & \nonumber \\
   &   \Gamma _{s+ m }^{\Big[\substack{s  \\ s +m}\Big] \times (s+ m+ l+2)} = 11\cdot2^{(s+m+l+2) +2m +s -3} +21\cdot2^{3s+3m-4} - 53\cdot2^{2s+3m-4}  && \text{by \eqref{eq 15.41} $\;\text{for}  \; k = s+ m+ l+2 $}  \label{eq 17.36} \\
 & =   11\cdot2^{2s+3m + l - 1} +21\cdot2^{3s+3m-4} - 53\cdot2^{2s+3m-4},    \nonumber\\
 & \nonumber \\
  &  \Gamma _{s+ m +1+j }^{\Big[\substack{s \\ s + m }\Big] \times s+m+l+2  }  =   8^{2j + m}\cdot \Gamma _{s-j+1}^{\Big[\substack{s -j\\ s -j }\Big] \times ((s+m+l+2)- m-2j)}    && \text{by $ (H_{l-1}) $ ( see   \eqref{eq 17.21}) }\label{eq 17.37} \\
   & =8^{2j +m}\cdot 21\cdot[2^{(s+m+l+2) -m -3j+s-1} + 2^{3s-3j-1} -5\cdot2^{2s-2j-1}] &&\text{ if $ \; 0\leq j\leq l-1,\; k = s+m+l+2 $} \nonumber\\
   & = 21\cdot[2^{2s +3m +l+3j +1} + 2^{3s +3m +3j -1}  -5\cdot2^{2s +3m +4j -1}]. \nonumber \\
   & \nonumber
  \end{align} \vspace{0.1 cm}
    From $\eqref{eq 17.30},\ldots, \eqref{eq 17.37} $ with $ k = s+m+l+2 $ we obtain after some calculations   \vspace{0.01 cm}\\
  \begin{align}
& \Gamma _{s+ m+ l+1}^{\Big[\substack{s  \\ s +m}\Big] \times (s+m+ l+2)} + \Gamma _{s+m+ l+2}^{\Big[\substack{s  \\ s +m}\Big] \times (s+ m+ l+2)} \label{eq 17.38}  \\
& = 2^{4s + 3m +2l +2} - (3\cdot2^{3s-4} - 2^{2s-3}) - ( 2^{2s + m + l +1} + 21\cdot2^{3s-4}  -  11\cdot2^{2s-3}) \nonumber \\
& - \sum_{j=1}^{m -1} [3\cdot2^{2s+m+l +2j- 1} +21\cdot2^{3s+3j-4} - 21\cdot2^{2s+3j-4}  ]- ( 11\cdot2^{2s+3m + l - 1} +21\cdot2^{3s+3m-4} - 53\cdot2^{2s+3m-4}  ) \nonumber \\
& -\sum_{j = 0}^{l-1}21\cdot[2^{2s +3m +l+3j +1} + 2^{3s +3m +3j -1}  -5\cdot2^{2s +3m +4j -1}]  \nonumber \\
& = 2^{4s+3m+2l+2} - 3\cdot2^{3s+3m+3l-1} - 5\cdot2^{2s+3m+4l-1},\nonumber \\
& \nonumber \\
& \Gamma _{s+ m+ l+1}^{\Big[\substack{s  \\ s +m}\Big] \times (s+ m+ l+2)}\cdot2^{-(s+m+ l+1)} + \Gamma _{s+m+ l+2}^{\Big[\substack{s  \\ s +m }\Big] \times (s+m+ l+2)}\cdot 2^{-(s+m +l+2)}  \label{eq 17.39}  \\
& =   2^{3s  + 2m + l} + 2^{2s +2m+2l +2} - 2^{s + m + l} -(7\cdot2^{2s-4} - 3\cdot2^{s-3} ) -  (  2^{2s + m + l +1} + 21\cdot2^{3s-4}  -  11\cdot2^{2s-3}  )\cdot2^{-s} \nonumber \\
&  -   \sum_{j= 1}^{m-1}[3\cdot2^{2s+m+l+2j-1} +21\cdot2^{3s+3j-4} - 21\cdot2^{2s+3j-4}]\cdot2^{-(s+j)} \nonumber  \\
&  - ( 11\cdot2^{2s+3m + l - 1} +21\cdot2^{3s+3m-4} - 53\cdot2^{2s+3m-4})\cdot2^{-(s+m)}   \nonumber \\
& -   \sum_{j = 0}^{l-1}21\cdot[2^{2s +3m +l+3j +1} + 2^{3s +3m +3j -1}  -5\cdot2^{2s +3m +4j -1}]\cdot2^{-(s+m+1+j)} \nonumber   \\
& = 2^{3s+2m+l} + 9\cdot2^{2s+2m+2l-2} - 13\cdot2^{s+2m+3l-2}. \nonumber                          
\end{align}\vspace{0.01 cm}\\ 
 Hence by \eqref{eq 17.38}, \eqref{eq 17.39} we deduce  \vspace{0.01 cm}\\
\begin{align}
\Gamma _{s+m+ l+1}^{\Big[\substack{s  \\ s +m}\Big] \times (s+m+ l+2)} &  = 21\cdot2^{3s +3m +3l -1}- 21\cdot 2^{2s + 3m  +4l-1}, \label{eq 17.40}\\
& \nonumber \\
\Gamma _{s+ m+ l+2}^{\Big[\substack{s  \\ s +m}\Big] \times (s+ m+ l+2)} & = 2^{4s+3m +2l +2} -3\cdot2^{3s+ 3m+ 3l +2}+ 2^{2s+ 3m + 4l +3}.\label{eq 17.41}\\
& \nonumber
\end{align}
By \eqref{eq 12.35} with $ j =  l+1 (  \leq s-1),    k > s+m+l+1 $ we get \vspace{0.01 cm}\\
\begin{equation}
\label{eq 17.42}
\Gamma _{s+ m+ l+1}^{\Big[\substack{s  \\ s +m}\Big] \times (k+1)} - \Gamma _{s+ m+ l+1}^{\Big[\substack{s  \\ s +m}\Big] \times k }= 21\cdot2^{k+s + 2m + 3l-1 }\quad \text{if}\quad k > s + m+ l+1.
\end{equation}

From \eqref{eq 17.42}, \eqref{eq 17.40} we deduce \vspace{0.01 cm}\\
\begin{align}
& \sum_{j = s+m+ l+ 2}^{k}\left( \Gamma _{s+m+l+1}^{\Big[\substack{s  \\ s +m}\Big] \times (j+1)} - \Gamma _{s+m+ l+1}^{\Big[\substack{s  \\ s +m }\Big] \times j}   \right)  = \sum_{j = s+m+ l+2}^{k}  21 \cdot 2^{j+s +2m  +3l-1} \nonumber \\
& \Leftrightarrow  \sum_{j = s+ m+ l+3}^{k+1}\Gamma _{s+ m+1}^{\Big[\substack{s  \\ s +m}\Big] \times j} - \sum_{j = s+ m+ l+2}^{k} \Gamma _{s+m+1}^{\Big[\substack{s  \\ s +m}\Big] \times j} = 21\cdot2^{k+s + 2m +3l}- 21\cdot2^{2s+ 3m + 4l+1} \nonumber \\
& \Leftrightarrow \Gamma _{s+m+ l+1}^{\Big[\substack{s  \\ s +m}\Big] \times (k+1)} - \Gamma _{s+m+l+1}^{\Big[\substack{s  \\ s +m}\Big] \times (s+ m+ l+2)} = 21\cdot2^{k+s + 2m +3l} - 21\cdot2^{2s+3m+4l +1}  \nonumber \\
& \Leftrightarrow \Gamma _{s+ m +l +1}^{\Big[\substack{s  \\ s +m}\Big] \times (k+1)}= 21\cdot2^{k+s + 2m +3l}- 21\cdot2^{2s+ 3m +4l +1} + 21\cdot(2^{3s +3m +3l-1} - \cdot2^{2s + 3m  +4l-1}) \nonumber \\
& \Leftrightarrow \Gamma _{s+m + l+1 }^{\Big[\substack{s  \\ s +m}\Big] \times (k+1)}=  21[2^{k+s+ 2m +3l} + 2^{3s +3m +3l-1} - 5\cdot2^{2s+ 3m +4l-1}]  \quad \text{ if $ k > s + m+ l+1 $}. \label{eq 17.43}\\
& \nonumber
\end{align}\vspace{0.01 cm}\\

By \eqref{eq 17.40} we see that \eqref{eq 17.43} holds for $ k = s+ m+ l+1 $, then  \eqref{eq 17.43} holds for $ k \geq s+ m+ l+1.$\vspace{0.01 cm}\\
From \eqref{eq 17.43} and \eqref{eq 13.4} with $ s\rightarrow s-l, \; k\rightarrow k - m -2l \; ( > s - l +1 )  $ we get \vspace{0.01 cm}\\
\begin{align*}
\Gamma _{s+m+ l+1}^{\Big[\substack{s  \\ s +m}\Big] \times  k} &= 8^{2l + m}\cdot\Gamma _{s - l +1}^{\Big[\substack{s -l \\ s -l }\Big] \times (k- m -2l) }=2^{3m +6l}\cdot[21\cdot(2^{(k- m -2l) + (s-l) -1} +2^{3(s-l) -1} - 5\cdot2^{2(s-l) -1})] \\
& =  21[2^{k+s+2m  +3l-1} + 2^{3s +3m +3l-1} - 5\cdot2^{2s+ 3m +4l-1}]  \quad \text{ if \;$ k > s+m+l+1 $}. \\
&        \qed                 
\end{align*}
We have now established that \vspace{0.01 cm}\\
$$ \Gamma _{s +m +1 +j}^{\Big[\substack{s \\ s +m }\Big] \times k}   =  8^{2j + m}\cdot \Gamma _{s-j+1}^{\Big[\substack{s -j\\ s -j }\Big] \times (k- m-2j)} \quad  \text{if   \; $ 0\leq j\leq s-2,\;  k\geq s + m+1+j $}.  $$\vspace{0.01 cm}\\
It remains to prove \vspace{0.01 cm}\\
\begin{equation}
\label{eq 17.44}
 \Gamma _{2s +m }^{\Big[\substack{s \\ s +m }\Big] \times k}   =  8^{2s + m -2}\cdot \Gamma _{ 2}^{\Big[\substack{1 \\  1}\Big] \times (k- m-2(s-1))} \quad  \text{if  $ \; k\geq 2s+m $}.
\end{equation}\vspace{0.01 cm}\\

\underline{The case k = 2s + m}\vspace{0.1 cm}\\
\eqref{eq 17.41} holds for $ l = s -2,  $ we then obtain by \eqref{eq 16.48} with $ k \rightarrow   2s  $ \vspace{0.1 cm}\\  
$$  8^{2(s-1) + m}\cdot \Gamma _{2}^{\Big[\substack{ 1 \\ 1 }\Big] \times 2}=
   8^{2(s-1) + m}\cdot[2^{ 4} -3\cdot2^{2} +2] = 2^{6s+3m-6}\cdot 6 = 3\cdot2^{6s+3m-5}= \Gamma _{2s +m }^{\Big[\substack{s \\ s +m }\Big] \times (2s+m)}.$$ \qed  \vspace{0.1 cm}\\ 
    
 \underline{The case $ k > 2s + m $}\vspace{0.1 cm}\\  
  We proceed as in the case $ k > s+m+l+1 $ with $ l = s-1.  $\vspace{0.1 cm}\\
 From \eqref{eq 17.38}, with $ l = s-1 $ we get \vspace{0.1 cm}\\
 \begin{align}
& \Gamma _{2s+m}^{\Big[\substack{s  \\ s +m}\Big] \times (2s+m+1)}
 =  2^{4s+3m+2(s-1)+2} - 3\cdot2^{3s+3m+3(s-1)-1} - 5\cdot2^{2s+3m+4(s-1)-1}\label{eq 17.45} \\
 & =  2^{6s+3m} - 3\cdot2^{6s+3m - 4} - 5\cdot2^{6s+3m -5}= 21\cdot2^{6s +3m -5}.\nonumber
\end{align}
By \eqref{eq 12.35} with $ j = s,    k > 2s+m $ we get \vspace{0.01 cm}\\
\begin{equation}
\label{eq 17.46}
\Gamma _{2s+m}^{\Big[\substack{s  \\ s +m}\Big] \times (k+1)} - \Gamma _{2s+m}^{\Big[\substack{s  \\ s +m}\Big] \times k }=3\cdot2^{2k+2s+m-2} - 3\cdot2^{k+4s+2m-4}    \quad \text{if}\quad k >  2s+m
\end{equation}

From \eqref{eq 17.46}, \eqref{eq 17.45} we deduce \vspace{0.01 cm}\\
\begin{align}
& \sum_{j = 2s+m+1}^{k}\left( \Gamma _{2s+m}^{\Big[\substack{s  \\ s +m}\Big] \times (j+1)} - \Gamma _{2s+m}^{\Big[\substack{s  \\ s +m }\Big] \times j}   \right)  = \sum_{j = 2s+m+1}^{k} 3\cdot2^{2j+2s+m-2} - 3\cdot2^{j+4s+2m-4}   \nonumber \\
& \Leftrightarrow  \sum_{j = 2s+m+2}^{k+1}\Gamma _{2s+m}^{\Big[\substack{s  \\ s +m}\Big] \times j} - \sum_{j = 2s+m+1}^{k} \Gamma _{2s+m}^{\Big[\substack{s  \\ s +m}\Big] \times j} = 2^{2k+2s+m} - 3\cdot2^{k+4s+2m-3} - 5\cdot2^{6s+3m-3} \nonumber \\
& \Leftrightarrow \Gamma _{2s+m}^{\Big[\substack{s  \\ s +m}\Big] \times (k+1)} - \Gamma _{2s+m}^{\Big[\substack{s  \\ s +m}\Big] \times (2s+m+1)} = 2^{2k+2s+m} - 3\cdot2^{k+4s+2m-3} - 5\cdot2^{6s+3m-3}   \nonumber \\
& \Leftrightarrow \Gamma _{2s+m}^{\Big[\substack{s  \\ s +m}\Big] \times (k+1)}=  2^{2k+2s+m} - 3\cdot2^{k+4s+2m-3} - 5\cdot2^{6s+3m-3} +    21\cdot2^{6s +3m -5}     \nonumber \\
& \Leftrightarrow \Gamma _{2s+m }^{\Big[\substack{s  \\ s +m}\Big] \times (k+1)}=    2^{2k+2s+m} - 3\cdot2^{k+4s+2m-3} + 2^{6s +3m -5}    \quad \text{ if $ k > 2s+m  $}. \label{eq 17.47}
\end{align}\vspace{0.01 cm}\\
By \eqref{eq 17.45} we see that \eqref{eq 17.47} holds for k = 2s+m, then  \eqref{eq 17.47} holds for $ k \geq 2s+m. $\vspace{0.01 cm}\\
From \eqref{eq 16.48}with $ k\rightarrow k-m,$ and \eqref{eq 17.47} we obtain \vspace{0.01 cm}\\
\begin{align*}
8^{2s+m-2}\cdot \Gamma _{2}^{\Big[\substack{ 1 \\ 1 }\Big] \times ( k-m-2(s-1))}& = 2^{6s+3m-6}\cdot[2^{2(k-m) - 4s +4} -3\cdot2^{k-m-2s+2} +2] \\
&=  2^{2k +m +2s-2} - 3\cdot2^{k+2m+4s - 4} +2^{6s+3m-5} \nonumber \\
& = \Gamma _{2s+m}^{\Big[\substack{s  \\ s +m}\Big] \times k }.\nonumber
\end{align*}
  \end{proof}
  
\section{\textbf{PROOF  OF THEOREMS  $3.1,  3.2,\ldots, 3.14, 3.15 $}}
\label{sec 18}
\subsection{\textbf{PROOF  OF  THEOREM  3.1}}
\label{subsec 18.1}
Follows from \eqref{eq 4.1} and \eqref{eq 4.19}.
 \subsection{\textbf{PROOF  OF  THEOREM  3.2}}
\label{subsec 18.2}  
  Follows from Lemma \ref{lem 5.14}.
 \subsection{\textbf{PROOF  OF  THEOREM  3.3}}
\label{subsec 18.3}  
  Follows from Lemma \ref{lem 8.3}.
 \subsection{\textbf{PROOF  OF  THEOREM  3.4}}
\label{subsec 18.4}
Follows from Lemma \ref{lem 8.4}.  
 \subsection{\textbf{PROOF  OF  THEOREM  3.5}}
\label{subsec 18.5}
 Follows from Lemma \ref{lem 8.5}.   
 \subsection{\textbf{PROOF  OF  THEOREM  3.6}}
\label{subsec 18.6} We have \\
\eqref{eq 3.10} follows from \eqref{eq 11.21},\\
\eqref{eq 3.11} follows from \eqref{eq 12.1},\\
\eqref{eq 3.12} follows from \eqref{eq 12.15},\\
\eqref{eq 3.13} follows from \eqref{eq 12.20},\\
\eqref{eq 3.14} follows from \eqref{eq 12.21}.
 \subsection{\textbf{PROOF  OF  THEOREM  3.7}}
\label{subsec 18.7}Follows from Lemma \ref{lem 15.2}.
  
 \subsection{\textbf{PROOF  OF  THEOREM  3.8}}
\label{subsec 18.8} We have \\
 \eqref{eq 3.19} follows from \eqref{eq 17.17},\\
\eqref{eq 3.20} follows from \eqref{eq 16.17},\\
\eqref{eq 3.21} follows from \eqref{eq 16.18} with $ k\rightarrow k-m, $ \\
\eqref{eq 3.22} follows from \eqref{eq 17.20}.
 \subsection{\textbf{PROOF  OF  THEOREM  3.9}}
\label{subsec 18.9} We have in the following cases : \\
 From \eqref{eq 11.21} we get   $ \Gamma _{i}^{\Big[\substack{s  \\ s }\Big] \times k } $   in the case  $1\leq i\leq s-1,\; k > i.$\\
From \eqref{eq 13.2} we get   $ \Gamma _{i}^{\Big[\substack{s  \\ s }\Big] \times k } $ in the case    $ i = s, \; k > s.$\\
From \eqref{eq 16.16}, \eqref{eq 16.18} and \eqref{eq 16.19} we deduce    $ \Gamma _{i}^{\Big[\substack{s  \\ s }\Big] \times k } $  in the case  $ s+1\leq i\leq2s,\; k> i. $

 \subsection{\textbf{PROOF  OF  THEOREM  3.10}}
\label{subsec 18.10} 
We have in the following two  cases  :\\
 From \eqref{eq 11.26} we get   $ \Gamma _{i}^{\Big[\substack{s  \\ s }\Big] \times i } $   in the case  $1\leq i\leq s.$ \\
From \eqref{eq 16.16}, \eqref{eq 16.17}  we deduce    $ \Gamma _{i}^{\Big[\substack{s  \\ s }\Big] \times  i } $  in the case  $ s+1\leq i\leq2s. $ 
 \subsection{\textbf{PROOF  OF  THEOREM  3.11}}
\label{subsec 18.11}  We have in the following cases :\\
 From \eqref{eq 13.21} we get   $ \Gamma _{i}^{\Big[\substack{s  \\ s +1}\Big] \times k } $   in the case  $1\leq i\leq s-1,\; k > i.$\\
From \eqref{eq 13.24} we get   $ \Gamma _{i}^{\Big[\substack{s  \\ s +1}\Big] \times k } $ in the case    $ i = s, \; k > s.$\\
From \eqref{eq 13.26} we get   $ \Gamma _{i}^{\Big[\substack{s  \\ s +1}\Big] \times k } $ in the case    $ i = s +1, \; k > s +1.$\\
From \eqref{eq 17.19} with j = i-s-2 we get   $ \Gamma _{i}^{\Big[\substack{s  \\ s +1}\Big] \times k } $ in the case    $s+2\leq i\leq 2s,\; k> i. $ \\
From \eqref{eq 17.20} with m = 1 we get   $ \Gamma _{i}^{\Big[\substack{s  \\ s +1}\Big] \times k } $  in the case  $i = 2s+1,\;k >2s+1.   $

 \subsection{\textbf{PROOF  OF  THEOREM  3.12}}
\label{subsec 18.12}  We have in the following two  cases :\\
 From \eqref{eq 13.25} and \eqref{eq 11.26} with m = 1 and $ k\rightarrow i  $ we get   $ \Gamma _{i}^{\Big[\substack{s  \\ s +1}\Big] \times i } $   in the case  $1\leq i\leq s+1.$ \\
From \eqref{eq 17.18}  with m = 1 and j = i - s - 2  we deduce    $ \Gamma _{i}^{\Big[\substack{s  \\ s +1}\Big] \times  i } $  in the case  $ s+2 \leq i\leq 2s +1. $ 

 \subsection{\textbf{PROOF  OF  THEOREM  3.13}}
\label{subsec 18.13}  
We have in the following cases :\\
 From \eqref{eq 13.21} we get   $ \Gamma _{i}^{\Big[\substack{s  \\ s +m}\Big] \times k } $   in the case  $1\leq i\leq s-1,\; k > i.$\\
From \eqref{eq 13.46} we get   $ \Gamma _{i}^{\Big[\substack{s  \\ s +m}\Big] \times k } $ in the case    $ i = s, \; k > s.$\\
From \eqref{eq 15.40} with $ l = i -s $ we get   $ \Gamma _{i}^{\Big[\substack{s  \\ s +m}\Big] \times k } $ in the case    $ s+1\leq i\leq s+m-1,\;k>i. $\\
From \eqref{eq 15.41} we get   $ \Gamma _{i}^{\Big[\substack{s  \\ s +m}\Big] \times k } $ in the case    $i = s+m, \; k> i. $ \\
From \eqref{eq 17.19} with j = i-s-m-1 we get   $ \Gamma _{i}^{\Big[\substack{s  \\ s + m}\Big] \times k } $  in the case  $s+m+1\leq i\leq 2s+m-1,\;k > i . $ \\
From \eqref{eq 17.20} we get   $ \Gamma _{i}^{\Big[\substack{s  \\ s + m}\Big] \times k } $  in the case  $i = 2s+m,\;k > i. $

 \subsection{\textbf{PROOF  OF  THEOREM  3.14}}
\label{subsec 18.14}  
 We have in the following three  cases :\\
 From \eqref{eq 11.26}  with $k\rightarrow i $  and \eqref{eq 13.47}   we get   $ \Gamma _{i}^{\Big[\substack{s  \\ s +m}\Big] \times i } $   in the case  $1\leq i\leq s+1.$ \\
From \eqref{eq 15.39}  with $l\rightarrow i $  and \eqref{eq 17.2}  we get   $ \Gamma _{i}^{\Big[\substack{s  \\ s +m}\Big] \times i } $   in the case  $s+2\leq i\leq s+m+1.$ \\
From \eqref{eq 17.18}  with  j = i-s-m-1   we  get   $ \Gamma _{i}^{\Big[\substack{s  \\ s +m}\Big] \times i } $   in the case  $s+m+2\leq i\leq 2s+m.$ 
 \subsection{\textbf{PROOF  OF  THEOREM  3.15}}
\label{subsec 18.15}  
Follows from Lemma \ref{lem 4.13}.

 \end{document}